\documentclass[11pt]{amsart}

\usepackage[margin=1in]{geometry}
\usepackage{amsmath,amssymb,amsthm}
\usepackage{mathtools}
\usepackage{enumitem}
\usepackage{microtype}
\usepackage{capt-of}
\usepackage[hidelinks]{hyperref}
\usepackage[nameinlink,noabbrev]{cleveref}
\newenvironment{datadisplay}[1]
  {\par\smallskip\noindent\begin{minipage}{\linewidth}
   \captionof{table}{#1}}
  {\end{minipage}\par\smallskip}

\setlist[enumerate]{leftmargin=2.2em,itemsep=2pt,topsep=4pt}
\setlist[itemize]{leftmargin=2em,itemsep=2pt,topsep=4pt}

\newtheorem{theorem}{Theorem}[section]
\newtheorem{proposition}[theorem]{Proposition}
\newtheorem{lemma}[theorem]{Lemma}
\newtheorem{corollary}[theorem]{Corollary}
\theoremstyle{definition}
\newtheorem{definition}[theorem]{Definition}
\theoremstyle{remark}
\newtheorem{remark}[theorem]{Remark}
\crefname{lemma}{lemma}{lemmas}
\Crefname{lemma}{Lemma}{Lemmas}
\crefname{proposition}{proposition}{propositions}
\Crefname{proposition}{Proposition}{Propositions}

\newcommand{\V}{\mathcal V}
\newcommand{\F}{\mathcal F}
\newcommand{\N}{\mathcal N}
\newcommand{\conv}{\operatorname{conv}}
\newcommand{\aff}{\operatorname{aff}}
\newcommand{\R}{\mathbb R}
\newcommand{\eps}{\varepsilon}

\title[Assignment obstructions beyond free joins]
{Sharp Vertex--Facet Assignment Obstructions Beyond Free Joins}

\author{Qihang Wang}
\address{School of Mathematical Sciences, Peking University,
Beijing 100871, China}
\email{2501110049@stu.pku.edu.cn}
\thanks{\raggedright School of Mathematical Sciences, Peking University, Beijing 100871,
China. Correspondence: \texttt{2501110049@stu.pku.edu.cn}.
ORCID: \href{https://orcid.org/0009-0003-2504-8756}{0009-0003-2504-8756}.}

\hypersetup{
  pdftitle={Sharp Vertex--Facet Assignment Obstructions Beyond Free Joins},
  pdfauthor={Qihang Wang},
  pdfsubject={Vertex--facet assignments for convex polytopes},
  pdfkeywords={convex polytope, vertex--facet assignment, matching deficiency,
    nonincidence graph, extremal geometry, free join}
}


\begin{document}

\begin{abstract}
A vertex--facet assignment of a convex polytope is a matching between
vertices and nonincident facets that covers the smaller of the two
classes. Jahn and Winter proved that assignments exist in dimension at
most six, constructed free-join counterexamples in every higher
dimension, and asked whether there are seven-dimensional
counterexamples other than free joins. Since free joins have
disconnected nonincidence graphs, we study counterexamples whose
nonincidence graph is connected. Every such seven-polytope has at least
seventeen vertices and seventeen facets, and an integral polytope
attains both bounds. At least sixteen rational combinatorial types
attain this minimum, and two equality branches are classified
completely. The lower-bound proof analyzes how the ambient incidences
of a face carrying an optimal Hall obstruction can coexist with
connectedness; the extremal examples are convex hulls of exposed faces
sharing one direction and coplanar one-point extensions of joins. The same
constructions give connected integral examples of every positive
deficiency in every dimension at least seven. For deficiency at least
\(h\) in dimension seven, the minimum larger class size lies between
\(3h+14\) and \(5h+12\). Finally, for strict join--stack constructions
with an inherited Hall certificate, we determine the sharp size floor
\(2d+2h+6\), its first-attainment dimension, the complete
deficiency-one sequence, and further exact optima in fixed dimension.

This work used the generative artificial intelligence models OpenAI
GPT-5.6 Sol, Anthropic Claude Fable 5, xAI Grok 4.6, OpenAI GPT-6.1 Sol,
OpenAI GPT-6 Astra, and Anthropic Claude Opus 5.5. Further details
appear in the disclosure at the end of the paper.
\end{abstract}

\maketitle

\section*{Acknowledgements}

The author acknowledges GeWu SciLab for providing the
collaborative research environment in which this project was developed. Bohrium's Large Knowledge Model (LKM) was used at the
problem-selection stage to locate the question in Section~6 of Jahn and
Winter's paper~\cite{jahnWinter2021} about seven-dimensional vertex--facet
assignment obstructions beyond free joins. This source-linked question
provided the starting point for the investigation.

\section{Introduction}
\label{sec:introduction}

Which face lattices allow vertices to be paired with nonincident facets?
For a convex \(d\)-polytope \(P\), a \emph{vertex--facet assignment}
is a matching that covers the smaller of the vertex and facet classes.
Equivalently, its bipartite nonincidence graph \(\N(P)\) has a matching
of size \(\min\{f_0(P),f_{d-1}(P)\}\). This formulation connects
face incidence with matching theory; its geometric motivation comes
from reduced polytopes and antipodality
\cite{averkovMartini2008,gonzalezMerinoJahnPolyanskiiWachsmuth2018}.

Jahn and Winter proved that assignments always exist through dimension
six and constructed counterexamples by free joins in every higher
dimension \cite{jahnWinter2021}. The nonincidence graph of a free join
is the disjoint union of the factor graphs. Opposite vertex--facet
imbalances in different components can therefore prevent an assignment,
even if each factor has one. Requiring the graph to be connected removes
this source of failure. In the first possible dimension, seven, Jahn and
Winter asked whether counterexamples exist beyond free joins, and
suggested taking the convex hull of a free join and a point
\cite[Section~6]{jahnWinter2021}.

We show that connected counterexamples exist and determine their
minimum size exactly: seventeen vertices and seventeen facets are
necessary, and one integral polytope attains both bounds. At least
sixteen rational combinatorial types attain the minimum. The geometric issue is how a face carrying an optimal
Hall obstruction sits inside the polytope: its ambient incidences must
retain the deficiency without disconnecting the graph. We use these
incidences both to exclude smaller examples and to describe equality.
Two further questions then arise naturally: how to prescribe the
deficiency, and how small a construction can be when a specified
free-join factor must remain its Hall witness.

\subsection{The unrestricted connected problem}

Write
\[
 \delta(P)=\min\{f_0(P),f_{d-1}(P)\}-\nu(\N(P))
\]
for the matching deficiency, where \(\nu(\N(P))\) is the maximum size
of a matching in \(\N(P)\). For integers \(d\ge7\), \(h\ge1\), let
\(\omega_{d,h}\) be the minimum of
\(\max\{f_0(P),f_{d-1}(P)\}\) over all \(d\)-polytopes with
connected nonincidence graph and \(\delta(P)\ge h\). We call \(P\)
\emph{balanced} if \(f_0(P)=f_{d-1}(P)\).
For a face \(\sigma\), let \(\kappa_P(\sigma)\) count the facets
containing it; \(\sigma\) is a \emph{maximizing Hall face} if it
maximizes \(f_0(\sigma)+\kappa_P(\sigma)\), and \(\sigma^\Delta\) denotes
the corresponding face of the polar \(P^*\). If \(\sigma\) and
\(\sigma^\Delta\) have vertex--facet counts \((a,\alpha)\) and
\((c,\gamma)\), their \emph{count slacks} are
\((f_0(P)-a-\gamma,\,f_{d-1}(P)-c-\alpha)\). A \emph{minimal coface} of
\(\sigma\) is a face of dimension \(\dim\sigma+1\) containing it, and its
vertices outside \(\sigma\) are \emph{outside vertices}.

\begin{theorem}[Sharp connected minimum and prescribed deficiency]
\label{thm:intro-connected}
Every seven-polytope with connected nonincidence graph and positive
matching deficiency has at least seventeen vertices and seventeen facets.
An integral seven-polytope attains both bounds and has matching number
sixteen. Thus
\[
 \omega_{7,1}=17.
\]
More generally, for every \(h\ge1\),
\[
 3h+14\le\omega_{7,h}\le5h+12.
\]
The upper bound has a balanced integral attainer of deficiency exactly
\(h\). For every \(d\ge7\) and \(h\ge1\), an integral connected
example of deficiency exactly \(h\) exists. At deficiency one,
\(\omega_{8,1}\le22\) and \(\omega_{9,1}\le27\).
\end{theorem}

\begin{proof}
Theorem~\ref{thm:seven-second-size-exclusion} supplies the lower bound;
Corollaries~\ref{cor:shared-prescribed} and~\ref{cor:sharp-seventeen}
supply the integral attainers and the remaining upper bounds.
\end{proof}

\begin{theorem}[Sixteen sharp equality types]
\label{thm:intro-sixteen}
At least sixteen pairwise combinatorially nonisomorphic rational
seven-polytopes have connected nonincidence graph and
\[
 (f_0,f_6,\nu)=(17,17,16).
\]
\end{theorem}

\begin{proof}[Proof of Theorem~\ref{thm:intro-sixteen}]
Theorem~\ref{thm:sixteen-extremizers} separates seven coplanar
join extensions, their seven polars, and the shared-coordinate pair by
facet-size and incidence invariants.
\end{proof}

The sharp example is the convex hull of an exposed hexagonal prism and
an exposed four-simplex sharing one direction. The prism is the
maximizing Hall face, and the shared direction prevents a free-join
decomposition. In
Section~\ref{sec:shared-coordinate} this becomes one family
\(W_{n,m,\ell}\), where \(\ell\ge1\), \(d=\ell+6\), \(b=m+\ell+1\),
\(m\ge3\) and \(n\ge b+1\), with exact counts
\[
 (f_0,f_{d-1},\nu)=
 \bigl(2n+b,\;n+m(\ell+1)+b,\;n+2b\bigr)
\]
and deficiency \(\min\{n-b,\ell(m-1)-1\}\).
The proof gives the complete inequality system, a maximum matching
and connectedness. For \(\ell=1\), the prism is the unique maximizing
Hall face (Proposition~\ref{prop:shared-unique-hall}). At the sharp minimum,
seventeen vertices and seventeen facets, any zero-slack maximizing pair forces the
entire face lattice to be that of \(W_{6,3,1}\), up to polarity
(Theorem~\ref{thm:zero-seventeen-rigidity}). Thus the construction
also identifies one complete equality class.

Coplanar extensions provide equality examples with positive count
slacks.
Theorem~\ref{thm:coplanar-extremizer} constructs them from a marked
six-facet three-polytope, with an exposed Hall face of counts \((10,7)\)
and dual counts \((8,6)\). Their count slacks are \((1,2)\), or
\((2,1)\) after polarity. Proposition~\ref{prop:coplanar-cube} gives
seventeen explicit rational inequalities and a sixteen-edge matching.
The associated triple-coface branch has a complete converse.

\begin{theorem}[The triple-coface equality branch]
\label{thm:intro-triple-branch}
Let \(P\) be a seven-polytope with connected nonincidence graph and
\((f_0,f_6,\nu)=(17,17,16)\). Suppose that a maximizing Hall face \(A\)
and its dual face \(D=A^\Delta\) are three-dimensional with counts
\((8,6)\) and \((10,7)\), and that the minimal cofaces of \(A\) are six
pyramids and one four-dimensional coface with three outside vertices. Then \(P^*\) is
combinatorially equivalent to one of the seven coplanar types of
Theorem~\ref{thm:sixteen-extremizers}. Every one of these
seven types occurs in this way.
\end{theorem}

\begin{proof}[Proof of Theorem~\ref{thm:intro-triple-branch}]
The stated counts give count slacks \((2,1)\), so
Corollary~\ref{cor:final-seven-boundary} makes both Hall faces simple.
Theorem~\ref{thm:triple-global-classification} gives the necessary
marked seed and all seven rational converses.
\end{proof}

At the boundary size \(3h+14\), the triple allocation is impossible
for every \(h\ge2\) (Theorem~\ref{thm:triple-coface-exclusion}). The other allocations
at the positive-slack boundary remain to be characterized.

The minimum \(\omega_{7,1}\) is taken over all seven-polytopes with
connected nonincidence graph. For \(h\ge2\) the bounds leave a gap; by
Corollary~\ref{cor:final-seven-boundary}, a polytope attaining
\(3h+14\) has a maximizing pair of simple Hall faces with counts
\((2h+6,h+5)\) and \((2h+8,h+6)\).
The eight- and
nine-dimensional examples give upper bounds in the unrestricted problem.

Theorem~\ref{thm:incidence-support-injection} gives an incidence-only
form of the local support obstruction. In three dimensions it yields
\(\alpha-3\le s\) for a simple base facet with \(\alpha\) facets, one
exceptional facet and \(s\) outside vertices, and the almost-simple
obstruction holds without sign assumptions
(Corollary~\ref{cor:sign-free-almost-simple}).

\subsection{The factor-certified extremal problem}

The unrestricted extremal problem allows every polytope. To isolate the
cost of a particular construction, we now retain a factor as a Hall
witness while connecting the graph by strict stacking. This additional
requirement defines a different minimum, which can be determined exactly
in a substantial range.
Let \(A,B\) have dimensions \(r,m\), vertex--facet counts
\((a,\alpha),(b,\beta)\), and let a facet \(F\) of \(A\) have
\(\rho\) ridges. Add one point strictly beyond \(F*B\) in \(A*B\)
and strictly beneath every other facet. The resulting
\(d=r+m+1\) polytope \(Q\) has
\[
 f_0(Q)=a+b+1,\qquad f_{d-1}(Q)=\alpha+2\beta+\rho-1.
\]
The surviving copy of \(A\) certifies deficiency at least \(h\) exactly
when
\begin{equation}\label{eq:intro-certificate}
 a-\alpha\ge\beta+\rho+h-1,\qquad \beta-b\ge h+1.
\end{equation}
Let \(\mathcal C_{d,h}\) be the class with this certificate, and let
\(\eta_{d,h}\) be the minimum of the larger vertex--facet class size,
with value \(+\infty\) when \(\mathcal C_{d,h}\) is empty.
The certificate is part of the definition of this class. Thus
\(\omega_{7,1}=17\) and \(\eta_{7,1}=33\) measure different
extremal problems, not competing bounds for the same one.

\begin{theorem}[Sharp size and exact first dimension]\label{thm:intro-sharp}
For an integer \(h\ge1\), define
\[
 \sigma_h=\min_{\substack{u,v\ge1\\uv=h+2}}(u+v),
 \qquad
 D_h=\sigma_h+1+\left\lceil2\sqrt{2h+5+\sigma_h}\right\rceil,
\]
where the minimum is over integers. Every \(Q\in\mathcal C_{d,h}\)
has at least \(2d+2h+6\) vertices and at least \(2d+2h+6\) facets.
Moreover,
\begin{equation}\label{eq:intro-sharp}
 \eta_{d,h}=2d+2h+6
 \quad\Longleftrightarrow\quad d\ge D_h.
\end{equation}
For every \(d\ge D_h\), a balanced rational polytope attains both bounds,
has connected nonincidence graph, and has matching number \(2d+h+6\).
Its deficiency is exactly \(h\), and its face lattice has no nontrivial
free-join decomposition.
\end{theorem}

\begin{proof}[Proof of Theorem~\ref{thm:intro-sharp}]
Combine Proposition~\ref{prop:sharp-floor},
Theorem~\ref{cor:stable-h}, and
Theorem~\ref{thm:sharp-construction}.
\end{proof}

The threshold has two sources. Equality in the complementary factor
forces
\[
 B^*\cong\operatorname{Pyr}^{m-u-v}(\Delta_u\times\Delta_v),
 \qquad uv=h+2,\quad u+v\le m
\]
(Proposition~\ref{prop:sharp-equality}); this determines
\(\sigma_h\). The remaining dimension must support a witness of
sufficient capacity with a simplex stacking facet. Its sharp capacity
bound yields the square-root term. A controlled vertex cut then attains
the required integer count, and a three-block matching attains the Hall
bound. The threshold therefore includes realizability and exact
deficiency, rather than count feasibility alone.

At fixed dimension, the optimum reflects the capacities of both
factors. Theorems~\ref{thm:seven-deficiency},
\ref{thm:eight-deficiency} and~\ref{thm:nine-deficiency}
give closed integer formulas and equality factor counts in dimensions
seven, eight and nine, respectively. Their leading terms are
\(7h\), \(4h\) and \(2h\). Euler's inequalities control dimension
seven; quadratic witness capacity enters in dimension eight, and
quadratic complementary capacity enters in dimension nine.
The face-closure formula in dimension seven, and matching saturation
in dimensions eight and nine, determine the exact deficiency of every
extremizer.

In higher fixed dimension, the two factors contribute inverse
face-number capacities of different orders.
For fixed \(d\ge9\), put \(t=\lfloor(d-1)/4\rfloor\).
Theorem~\ref{thm:second-order} gives
\begin{equation}\label{eq:intro-second-order}
 \eta_{d,h}=2h+K_d h^{1/t}+L_d h^{1/(t+1)}+O_d(1)
 \qquad(h\longrightarrow\infty),
\end{equation}
where \(K_d\) is the optimal second-order coefficient;
\(L_d=2((t+1)!)^{1/(t+1)}\) for \(d=4t+3,4t+4\), and is zero
otherwise. Balanced rational connected examples attain this profile
with exact deficiency \(h\). Corollary~\ref{cor:capacity-stability}
identifies the oriented factor dimensions and bounds the remaining
count, ridge and certificate-slack budgets near the optimum.

The integer refinements separate witness parity from complementary
face-pair feasibility. Theorems~\ref{thm:even-integer}
and~\ref{thm:odd-witness-integer} determine the exact eventual optimum
when the optimal complement is even-dimensional, including the equality
factor data and exact deficiency. Classical coplanar stacking removes
the odd-witness parity loss. With an odd-dimensional complement,
Theorem~\ref{thm:odd-integer} reduces the eventual optimum to one
specified vertex--facet pair. The classical facet gap produces the
ten-dimensional plateaus of Corollary~\ref{cor:ten-plateau} and the
higher-dimensional density-one families of
Theorem~\ref{thm:density-plateaus}. The latter density is over capacity
indices, rather than over all deficiencies.

For deficiency one, the low-dimensional constructions complete the
entire sequence. Write \(\eta_d=\eta_{d,1}\).

\begin{theorem}[Complete deficiency-one optimum]\label{thm:intro-complete}
The class \(\mathcal C_{d,1}\) is empty for \(d<7\), and
\[
 \eta_7=33,\qquad \eta_8=\eta_9=28,\qquad
 \eta_{10}=29,\qquad \eta_{11}=31,\qquad
 \eta_d=2d+8\quad(d\ge12).
\]
Each finite value is attained by a balanced rational polytope with
connected nonincidence graph and matching number \(\eta_d-1\).
\end{theorem}

\begin{proof}
This is Theorem~\ref{thm:complete-one}.
\end{proof}

\subsection{The common mechanism: Hall faces and extensions}

The matching problem becomes geometric after closing a Hall witness
under face incidence. For \(n=f_0(P)\) and \(m=f_{d-1}(P)\), put
\[
 H(P)=\max_{\sigma\le P}\{f_0(\sigma)+\kappa_P(\sigma)\},
\]
where \(\kappa_P(\sigma)\) counts the facets containing \(\sigma\);
the maximum includes the empty face and \(P\) itself.
Jahn and Winter's criterion is
\(H(P)\le\max\{n,m\}\) \cite[Theorem~6]{jahnWinter2021}.
The same closure argument, followed by K\H{o}nig's theorem, gives the
elementary quantitative form
\begin{equation}\label{eq:intro-exact}
 \nu(\N(P))=n+m-H(P),\qquad
 \delta(P)=H(P)-\max\{n,m\}.
\end{equation}
Section~\ref{sec:defect} also separates within-component deficiency
from cancellation between oppositely unbalanced components. This
explains how a free join can fail an assignment even when its factors
admit assignments.

Preserving such a face and making the graph connected are separate
requirements. Theorem~\ref{thm:horizon-transfer} treats both for
nondegenerate vertex-preserving one-point extensions in the strict
visible/retained setting: a surviving face transfers the Hall bound,
and a rooted horizon graph characterizes connectedness. In particular,
strict one-facet stacking always connects the graph
(Corollary~\ref{cor:stacking-transfer}).

An explicit coordinate family records this mechanism in every
dimension at least seven.

\begin{theorem}\label{thm:main}
For integers \(d\ge7\) and \(n\ge2^{d-4}+6\), there is a rational
\(d\)-polytope \(Q_{d,n}\) with connected nonincidence graph such that
\[
 f_0(Q_{d,n})=2n+2d-7,\qquad
 f_{d-1}(Q_{d,n})=n+5+2^{d-3},
\]
and
\[
 \nu(\N(Q_{d,n}))=n+2^{d-4}+2d-2
 <\min\{f_0(Q_{d,n}),f_{d-1}(Q_{d,n})\}.
\]
Its face lattice has no nontrivial free-join decomposition.
\end{theorem}

For fixed \(d\), the vertex count strictly increases with \(n\),
so the family consists of pairwise combinatorially nonisomorphic
polytopes. The proof of Theorem~\ref{thm:main} is completed in
Section~\ref{sec:matching}.

\begin{corollary}\label{cor:seven}
For every \(n\ge14\), the seven-polytope \(Q_{7,n}\) has \(2n+7\)
vertices, \(n+21\) facets and matching number \(n+20\).
In particular, seven-dimensional counterexamples exist beyond free joins.
\end{corollary}

\subsection{Relation to earlier work}

The matching ingredients are Hall's and K\H{o}nig's theorems
\cite{hall1935,koenig1931,lovaszPlummer2009}.
The geometric constructions use classical operations.
Free joins also occur in self-duality constructions
\cite[Corollary~2]{tiwaryElbassioni2014}, and subdirect sums are
projections of joins \cite[Section~15.1.3]{henkRichterGebertZiegler2017}.
Our shared-coordinate hull belongs to this setting. Its contribution
to the assignment problem is the complete facet budget and matching
structure. Throughout, indecomposability concerns the face lattice,
rather than Minkowski decomposability.

The capacity inputs are McMullen's Upper Bound Theorem
\cite{mcmullen1970}, his few-facet classification \cite{mcmullen1971},
and the facet types recorded in
\cite[Lemma~1 and Remark~6]{pinedaTritamaWangYost2026}.
Few-facet constructions provide wide ranges of vertex counts
\cite[Section~6]{pinedaUgonYost2019}. Billera--Lee sufficiency
\cite{billeraLee1980,billeraLee1981}, with necessity due to Stanley
\cite{stanley1980}, supplies the simplicial realizations.
Sj\"oberg and Ziegler's interpolation of the top \(g\)-entry and
coplanar stacking supply the count and parity adjustments
\cite[Theorems~3.2--3.3 and their proofs]{sjobergZiegler2020}.
Our optimization couples these classical capacities to a retained
Hall certificate and a matching that saturates its bound.

The distinction is also visible in the lower bounds.
Recent work minimizes face numbers at fixed vertex count
\cite{pinedaUgonYost2019,pinedaTritamaWangYost2026}.
Guillermo Pineda-Villavicencio and Jie Wang continue this direction in
their preprint \emph{A lower bound theorem for \(d\)-polytopes with at
most \(3d-1\) vertices} (2025), arXiv:2512.07456v1,
\url{https://doi.org/10.48550/arXiv.2512.07456}.
Its Section~6(2) already describes simple-vertex truncations of pyramids
over \(\Delta_2\times\Delta_{s-2}\), a constituent family used here.
Our extremal objective instead couples the larger vertex--facet class
to a specified Hall certificate, its realization and a saturating matching.
Almost-simplicial bounds with one distinguished
nonsimplex facet assume that every other facet is a simplex
\cite[Theorem~1.2]{nevoPinedaUgonYost2018}.
Our deleted-polar estimate applies the Upper Bound Theorem to a hull
and a vertex figure without that global hypothesis.
For odd-dimensional complements, the nonsimplicial facet gap is
a classical input, as derived by pulling in
\cite[proof of Theorem~3.3]{sjobergZiegler2020}.
The extremal plateaus arise when this gap interacts with the Hall
certificate and its equality slacks; the density argument uses
Weyl's equidistribution theorem \cite{weyl1916}.

For the unrestricted lower bound, face-number inequalities alone leave
small candidate budgets. The additional constraints come from the
incidence relation between paired Hall faces and their ambient cofaces.
The slack-matrix rank and row-space language is classical
\cite[Theorem~14]{gouveiaEtAl2013} and
\cite[Definition~2.4 and Theorem~2.5]{gouveiaEtAl2019}.
The local exclusions decide which candidate budgets can actually be
realized. The support injection uses incidence zeros rather than
positivity; the separate ambient argument establishes the nonvanishing
needed to apply it. This distinction is essential to the sharp bound.

\subsection{Proof architecture}

The proof of the unrestricted minimum separates necessity from
realization. After the extension principle and quantitative Hall identity
in Sections~\ref{sec:preliminaries} and~\ref{sec:defect},
Section~\ref{sec:connected-extremal} eliminates the small dual Hall-face
budgets. Ambient nonvanishing and support injection provide the local
incidence obstruction. Section~\ref{sec:shared-coordinate} develops
the opposite direction: shared-coordinate hulls attain the bound and
prescribe deficiency, while their equality analysis and the coplanar
extension analysis give the zero-slack and triple-coface converses.

The certified theory follows the same necessity--attainment order.
The balanced construction in Section~\ref{sec:balanced} leads to the
size floor and first-attainment dimension in Section~\ref{sec:sharp},
followed by fixed-dimensional capacity laws and their integer
refinements. Section~\ref{sec:lowdim} closes the remaining dimensions
at deficiency one. The coordinate family and explicit maximum matching
in Sections~\ref{sec:construction} and~\ref{sec:matching} make the
extension mechanism concrete. The appendices give alternative
constructions and describe verification programs for several of the
explicit families.

\section{Assignments, joins, and strict stacking}
\label{sec:preliminaries}

For a polytope \(P\), write \(\V(P)\) and \(\F(P)\) for its vertex and facet
sets.  Its \emph{nonincidence graph} \(\N(P)\) is the bipartite graph on
\(\V(P)\sqcup\F(P)\) in which \(vF\) is an edge precisely when \(v\notin F\).
A vertex--facet assignment is a matching in \(\N(P)\) of size
\(\min\{|\V(P)|,|\F(P)|\}\).
Write \(\nu(X)\) for the matching number of a graph \(X\).

We shall use the following elementary counting bound, which underlies the
necessity part of Hall's theorem \cite{hall1935}; compare
\cite[Proposition~4 and Theorem~6]{jahnWinter2021}.

\begin{lemma}[Face-neighborhood bound]\label{lem:face-bound}
Let \(G\) be a face of \(P\), with \(r\) vertices, and suppose that exactly
\(c\) facets of \(P\) contain \(G\).  Every matching in \(\N(P)\) has size at
most
\[
  |\F(P)|-c+|\V(P)|-r.
\]
\end{lemma}

\begin{proof}
Put \(S=\V(G)\).  A facet has no neighbor in \(S\) exactly when it contains
every vertex of \(G\), hence exactly when it contains \(G\).  Thus
\(|N(S)|=|\F(P)|-c\).  At most this many matching edges can meet \(S\), and
at most \(|\V(P)|-r\) further edges can meet vertices outside \(S\).
\end{proof}

We realize the free join of polytopes \(P_1\subset\R^{d_1}\) and
\(P_2\subset\R^{d_2}\) as
\[
 P_1*P_2=\conv\bigl((P_1,0,0)\cup(0,P_2,1)\bigr)
 \subset\R^{d_1+d_2+1}.
\]
It has dimension \(d_1+d_2+1\), and its facets are \(F_1*P_2\) and
\(P_1*F_2\), where \(F_i\) ranges over the facets of \(P_i\); see
\cite[Section~15.1.3]{henkRichterGebertZiegler2017}. Here a point has the
empty face as its only facet, and a free join is \emph{nontrivial} if both
factors are nonempty.

\begin{lemma}[Free joins disconnect nonincidence]\label{lem:join-disconnects}
If the face lattice of a positive-dimensional polytope \(P\) is a nontrivial
free join, then \(\N(P)\) is disconnected.
\end{lemma}

\begin{proof}
For \(P=P_1*P_2\), the vertices split into the two factor classes and the
facets split into the two induced classes.  Every vertex from \(P_1\) belongs
to every facet \(P_1*F_2\), and every vertex from \(P_2\) belongs to every
facet \(F_1*P_2\).  Therefore every nonincidence edge stays within one of the
two blocks
\[
 \V(P_1)\sqcup\{F_1*P_2:F_1\in\F(P_1)\},
 \qquad
 \V(P_2)\sqcup\{P_1*F_2:F_2\in\F(P_2)\}.
\]
Each block contains the vertices of one factor and is therefore nonempty,
so \(\N(P)\) is disconnected. (If \(P_1\) is a point, its block consists of
that vertex and the facet \(P_2\).)
\end{proof}

Finally, we record the standard strict-stacking rule
\cite[Section~5.2]{gruenbaum2003}.  For a facet with defining inequality
\(g\le0\), where \(g\) is affine, a point \(x\) is \emph{beyond} the facet if
\(g(x)>0\) and \emph{beneath} it if \(g(x)<0\).

\begin{lemma}[Strict one-facet stacking]\label{lem:stacking}
Let \(P\) be a \(d\)-polytope with \(d\ge2\), and let \(p\) lie strictly beyond one facet
\(S\) and strictly beneath every other facet.  Then
\(\conv(P\cup\{p\})\) retains every old facet except \(S\) and has one new
facet \(\conv(p,R)\) for each ridge \(R\) of \(S\).  It has no other facets,
and every old vertex remains a vertex.
\end{lemma}

If \(g_0\le0\) supports \(S\), \(g_1\le0\) supports an old facet meeting
\(S\) in a ridge, and \(\lambda=g_0(p)>0\), \(\mu=g_1(p)<0\), then the new
cone facet over that ridge has supporting form
\begin{equation}\label{eq:cone-form}
  \ell_{g_1}=\lambda g_1-\mu g_0\le0.
\end{equation}
Indeed, \(\ell_{g_1}\le0\) on \(P\), it vanishes at \(p\), and among the
old vertices it vanishes precisely at the vertices of the ridge
\(S\cap\{g_1=0\}\).  Together with Lemma~\ref{lem:stacking}, this gives an
explicit complete facet description of a stacked polytope.

We next treat a general nondegenerate one-point extension.  Let
\(p\notin P\) lie on no facet hyperplane of \(P\); call the facets beyond
which \(p\) lies \emph{visible} and the others \emph{retained}; since
\(p\notin P\), at least one facet is visible.  Assume that
\(Q=\conv(P\cup\{p\})\) is vertex-preserving, that is, every vertex of \(P\)
is a vertex of \(Q\).  Its new facets are the cones \(C_R=\conv(p,R)\) over the horizon ridges \(R\),
namely the ridges incident to one visible and one retained facet
\cite[Section~5.2]{gruenbaum2003}.  Define the \emph{rooted horizon graph} \(\Gamma(P,p)\) on
a root \(*\), the old vertices of \(P\), and the horizon ridges by putting
\[
 *\sim v
 \quad\Longleftrightarrow\quad
 v\text{ misses a retained facet},
 \qquad
 v\sim R
 \quad\Longleftrightarrow\quad v\notin R.
\]

\begin{theorem}[One-point transfer through the horizon]
\label{thm:horizon-transfer}
In the setting above, \(\N(Q)\) is connected if and only if
\(\Gamma(P,p)\) is connected.  Suppose in addition that \(G\) is a nonempty
proper face of \(P\), with \(r\) vertices, equal to the intersection of
a nonempty set of retained facets, and suppose that exactly \(c\)
retained facets in total contain it. Assume also that no horizon ridge
contains \(G\). If there are \(v\) old vertices, \(f_R\) retained
facets, and \(q\) horizon ridges, then \(G\) survives in \(Q\), is
contained in exactly \(c\) facets, and
\[
  \nu(\N(Q))\le f_R+q-c+v+1-r.
\]
\end{theorem}

\begin{proof}
The apex is nonincident to every retained facet and incident to every new
facet.  Hence the apex and all retained facet nodes span a connected star;
contracting it to one root preserves connectedness in both directions.  An old vertex is adjacent to this root precisely when it misses a
retained facet.  The node of \(C_R\) is adjacent to an old vertex \(v\)
precisely when \(v\notin R\).  The contracted graph is therefore exactly
\(\Gamma(P,p)\), proving the connectivity equivalence.

The retained supporting hyperplanes in the hypothesis still cut out \(G\).
Indeed, every point of \(Q\) has the form \(tp+(1-t)x\) with \(x\in P\)
and \(t\in[0,1]\). If \(g_F(tp+(1-t)x)=0\) for a retained facet \(F\),
then \(g_F(p)<0\) and \(g_F(x)\le0\) force \(t=0\) and \(x\in F\).
A new facet \(C_R\) does not
contain \(G\), since \(C_R\cap P=R\) and no horizon ridge contains \(G\).
Hence \(G\) is contained in exactly the stated \(c\) retained facets.
The extension has \(v+1\) vertices and \(f_R+q\) facets, so
Lemma~\ref{lem:face-bound} gives the bound.
\end{proof}

For strict one-facet stacking the horizon graph is always connected, so
the nonincidence graph of the stacked polytope is automatically connected.

\begin{corollary}[Connected Hall transfer under strict stacking]
\label{cor:stacking-transfer}
Let \(P\) be a \(d\)-polytope with \(d\ge2\), with \(v\) vertices and
\(f\) facets, let \(S\) be a facet of \(P\) with \(\rho\) ridges, and
let \(Q\) be obtained by strict stacking beyond \(S\).
Then \(\N(Q)\) is connected.  If a nonempty proper face \(G\) of \(P\) has
\(r\) vertices, is contained in exactly \(c\) facets of \(P\), and is not
contained in \(S\), then
\(G\) survives in \(Q\), is contained in exactly \(c\) facets of \(Q\), and
\[
  \nu(\N(Q))\le v+f+\rho-r-c.
\]
Consequently, \(Q\) admits no nontrivial combinatorial free-join
decomposition, and if the right-hand side is smaller than
\(\min\{v+1,f-1+\rho\}\), then \(Q\) has no vertex--facet assignment.
\end{corollary}

\begin{proof}
For connectivity, first consider a ridge \(R\) of \(S\).  Choose a vertex
\(u\in\V(S)\setminus\V(R)\). Some facet \(K\) of \(P\) misses \(u\);
since \(u\in S\), this facet is retained. Thus the horizon node \(R\) has
the path \(R\sim u\sim *\) in \(\Gamma(P,p)\). Every vertex of \(S\)
misses a retained facet and is adjacent to \(*\). A vertex outside \(S\)
lies in no ridge of \(S\), so it is adjacent to every horizon ridge node,
and these were already shown to lie in the root component.  Hence \(\Gamma(P,p)\), and therefore
\(\N(Q)\), is connected by Theorem~\ref{thm:horizon-transfer}.

Every facet containing \(G\) is retained because \(G\not\subset S\), and
these facets cut out \(G\).  No ridge of \(S\) contains \(G\).  Apply
Theorem~\ref{thm:horizon-transfer} with \(f_R=f-1\) and \(q=\rho\) to obtain
the displayed bound.  The final assertions follow from the definition of an
assignment and Lemma~\ref{lem:join-disconnects}.
\end{proof}

\section{Exact defect and its component decomposition}
\label{sec:defect}

The face-neighborhood bound is sufficient to certify a failed assignment,
but it does not by itself say whether the chosen face is optimal or how many
matching edges are lost.  Let \(P\) be a \(d\)-polytope with
\(n=f_0(P)\) vertices and \(m=f_{d-1}(P)\) facets.  For a face
\(\sigma\) of \(P\), including \(\varnothing\) and \(P\), write
\[
 v(\sigma)=|\V(\sigma)|,\qquad
 \kappa_P(\sigma)=|\{F\in\F(P):\sigma\subseteq F\}|,
\]
so that \(v(\varnothing)=0\),
\(\kappa_P(\varnothing)=m\), and
\(\kappa_P(P)=0\).  Define the \emph{face-closure number}
\begin{equation}\label{eq:face-closure-number}
H(P)=\max_{\sigma\le P}
 \{v(\sigma)+\kappa_P(\sigma)\}.
\end{equation}
Jahn and Winter already proved that an assignment exists if and only if
\(H(P)\le\max\{n,m\}\)~\cite[Theorem~6]{jahnWinter2021}.  We now show that
this same number is the independence number of \(\N(P)\), so that
K\H{o}nig's theorem makes the matching number exact.

\begin{theorem}[Exact face-closure formula]
\label{thm:exact-defect}
For every convex polytope \(P\),
\begin{equation}\label{eq:exact-matching-formula}
 \nu(\N(P))=n+m-H(P).
\end{equation}
Consequently, the matching deficiency
\begin{equation}\label{eq:deficiency-definition}
 \delta(P):=\min\{n,m\}-\nu(\N(P))
\end{equation}
is exactly
\begin{equation}\label{eq:exact-deficiency-formula}
 \delta(P)=H(P)-\max\{n,m\}.
\end{equation}
In particular, \(P\) has a vertex--facet assignment if and only if
\begin{equation}\label{eq:exact-assignment-criterion}
 v(\sigma)+\kappa_P(\sigma)\le \max\{n,m\}
 \qquad\text{for every face }\sigma,
\end{equation}
which is the criterion of Jahn and Winter~\cite[Theorem~6]{jahnWinter2021}.
\end{theorem}

\begin{proof}
Let \(I\) be an independent set of the bipartite nonincidence graph, and
write \(I=A\sqcup B\), where \(A\subseteq\V(P)\) and \(B\subseteq\F(P)\).
Independence says that every vertex in \(A\) is incident to every facet in
\(B\).  If \(A,B\ne\varnothing\), put
\(\sigma=\bigcap_{F\in B}F\).  Then \(\sigma\) is a face,
\(A\subseteq\V(\sigma)\), and every member of \(B\) contains \(\sigma\), so
\[
 |I|=|A|+|B|\le v(\sigma)+\kappa_P(\sigma)\le H(P).
\]
If \(A=\varnothing\), use \(\sigma=\varnothing\); if \(B=\varnothing\), use
\(\sigma=P\).  Thus the independence number satisfies
\(\operatorname{ind}(\N(P))\le H(P)\).

Conversely, for every face \(\sigma\), the set
\[
 \V(\sigma)\sqcup\{F\in\F(P):\sigma\subseteq F\}
\]
is independent: each displayed vertex belongs to each displayed facet.
Therefore \(\operatorname{ind}(\N(P))=H(P)\).  The complement of a maximum
independent set is a minimum vertex cover, and by K\H{o}nig's
theorem~\cite{koenig1931} a minimum vertex cover of a bipartite graph has
the same size as a maximum matching.  Hence
\(\operatorname{ind}(\N(P))+\nu(\N(P))=n+m\), proving
\eqref{eq:exact-matching-formula}.  Substituting it into
\eqref{eq:deficiency-definition} and using
\(\min\{n,m\}+\max\{n,m\}=n+m\) gives
\eqref{eq:exact-deficiency-formula}, and
\eqref{eq:exact-assignment-criterion} follows.
\end{proof}

\begin{corollary}[Closed maximum witnesses]
\label{cor:closed-witnesses}
Every maximum independent set of \(\N(P)\) is of the form
\[
 \V(\sigma)\sqcup\{F\in\F(P):\sigma\subseteq F\}
\]
for some face \(\sigma\).  In particular, an \emph{optimal Hall
obstruction}, that is, a maximum independent set of \(\N(P)\), may
always be chosen to consist of all vertices of one face together with all
facets containing that face.
\end{corollary}

\begin{proof}
The proof of Theorem~\ref{thm:exact-defect} shows that every independent
set is contained in one of the displayed sets, which are themselves
independent. For a maximum independent set the containment is an
equality.
\end{proof}

A face \(\sigma\) attaining the maximum in \eqref{eq:face-closure-number}
is called a \emph{maximizing Hall face}, and \(\sigma^\Delta\) its
\emph{dual Hall face}. When \(\delta(P)>0\), both are nonempty proper
faces, and by polarity their roles may be exchanged. We call
\((\sigma,\sigma^\Delta)\) a \emph{maximizing Hall pair}, and its
\emph{count slacks} are
\[
 u=f_0(P)-f_0(\sigma)-f_{\dim\sigma^\Delta-1}(\sigma^\Delta),\qquad
 v=f_{d-1}(P)-f_0(\sigma^\Delta)-f_{\dim\sigma-1}(\sigma).
\]

The exact formula also separates two mechanisms that can cause a matching
defect.  Let the connected components of \(\N(P)\) be
\(C_1,\ldots,C_t\).  Write \(n_i\) and \(m_i\) for the sizes of the vertex
and facet sides of \(C_i\), respectively, and put
\[
 \delta_i=\min\{n_i,m_i\}-\nu(C_i),\qquad
 \Delta_i=n_i-m_i.
\]

\begin{theorem}[Intrinsic--cancellation decomposition]
\label{thm:component-defect}
The total deficiency decomposes as
\begin{equation}\label{eq:component-defect}
 \delta(P)=\sum_{i=1}^t\delta_i+
 \frac{\sum_{i=1}^t|\Delta_i|-\left|\sum_{i=1}^t\Delta_i\right|}{2}.
\end{equation}
\end{theorem}

\begin{proof}
A matching is the disjoint union of matchings in the components, so
\[
 \nu(\N(P))=\sum_i
 \bigl(\min\{n_i,m_i\}-\delta_i\bigr).
\]
Using \(\min\{x,y\}=(x+y-|x-y|)/2\) first for each component and then for
the two total class sizes gives
\[
 \min\Bigl\{\sum_i n_i,\sum_i m_i\Bigr\}-\sum_i\min\{n_i,m_i\}
 =\frac{\sum_i|\Delta_i|-|\sum_i\Delta_i|}{2},
\]
which proves the formula.
\end{proof}

The sum \(\sum_i\delta_i\) measures a failure already present inside a
connected component; we call it the \emph{intrinsic defect}.  The second
term measures cancellation between components whose two bipartition classes
are unbalanced in opposite directions.  It vanishes for a connected graph,
so every defect of a connected counterexample is intrinsic.

\begin{corollary}\label{cor:component-assignment}
The polytope \(P\) has a vertex--facet assignment if and only if
\(\delta_i=0\) for every component and all nonzero \(\Delta_i\) have the
same sign.  If \(P\) is combinatorially equivalent to \(A*B\), both
factors admit assignments, and their imbalances have opposite signs, then
\begin{equation}\label{eq:join-cancellation-defect}
 \delta(P)=\min\{|f_0(A)-f_{\dim A-1}(A)|,
                    |f_0(B)-f_{\dim B-1}(B)|\}.
\end{equation}
\end{corollary}

\begin{proof}
The first assertion follows from \eqref{eq:component-defect}, since both
terms are nonnegative. For the second, the nonincidence graph of a free
join is the disjoint union of the two factor nonincidence graphs, with
the factor facet nodes relabeled by the corresponding join facets.
Since both factors admit assignments,
\(\nu(\N(P))=\min\{n_A,m_A\}+\min\{n_B,m_B\}\), where
\((n_A,m_A)\) and \((n_B,m_B)\) are the factor class sizes.
Interchanging \(A\) and \(B\) if necessary, write
\(n_A-m_A=p>0\) and \(m_B-n_B=q>0\). Then
\(\min\{n_A+n_B,m_A+m_B\}=m_A+n_B+\min\{p,q\}\), so
\(\delta(P)=\min\{p,q\}\).
\end{proof}

Thus free-join examples can be obstructed entirely by cancellation between
otherwise assignment-admitting blocks.  By contrast, the defect of every
connected example constructed below is intrinsic; for instance, the
polytopes \(R_d\) of Theorem~\ref{thm:balanced} have \(\delta_1=1\).

\section{The connected seven-dimensional lower bound}
\label{sec:connected-extremal}

We now pass from the matching identity to the extremal geometry.
The bound \(3h+8\) of Proposition~\ref{prop:global-seven} is sharp
among all seven-polytopes, but its equality cases are free joins
(Theorem~\ref{thm:seven-rigidity}). The subsequent exclusions measure the additional
cost of connectedness. They culminate in the bound \(3h+14\) of
Theorem~\ref{thm:seven-second-size-exclusion} and the boundary data of
Corollary~\ref{cor:final-seven-boundary}.

\subsection{The free-join benchmark}

Jahn and Winter's criterion, combined with Euler's formula, already
gives a sharp size bound in dimension seven before connectedness is
imposed.

\begin{proposition}[Seven-dimensional size bound]\label{prop:global-seven}
Let \(h\ge1\) be an integer. Every seven-polytope with matching deficiency at least
\(h\) has at least \(3h+8\) vertices and at least \(3h+8\) facets.
Both bounds are attained simultaneously by a rational free join with
deficiency exactly \(h\).
\end{proposition}

\begin{proof}
Let \(P\) be a seven-polytope with \(\delta(P)\ge h\). Choose a face
\(\sigma\) realizing \(H(P)\), and let \(\sigma^\Delta\) be its dual
face in a polar polytope. Since \(\delta(P)\ge1\), the face \(\sigma\)
is nonempty and proper. Write
\[
 a=f_0(\sigma),\quad \alpha=f_{\dim\sigma-1}(\sigma),\quad
 c=f_0(\sigma^\Delta)=\kappa_P(\sigma),\quad
 \gamma=f_{\dim\sigma^\Delta-1}(\sigma^\Delta),
\]
so that \(\alpha\) and \(\gamma\) are the facet counts of the two faces.
Distinct facets of \(\sigma\) are cut out by distinct facets of \(P\)
not containing \(\sigma\), and the dual assertion holds for
\(\sigma^\Delta\). Thus, for \(n=f_0(P)\) and \(m=f_6(P)\),
\[
 m\ge c+\alpha,\qquad n\ge a+\gamma.
\]
The identity \(a+c=H(P)\ge\max\{n,m\}+h\) gives
\(a-\alpha\ge h\) and \(c-\gamma\ge h\).
Polytopes of dimension at most two, including a point, have equally
many vertices and facets, so \(\dim\sigma,\dim\sigma^\Delta\ge3\); since
\(\dim\sigma+\dim\sigma^\Delta=6\), both faces are three-dimensional
(compare \cite[proof of Corollary~8]{jahnWinter2021}). In particular
\(\alpha=f_2(\sigma)\) and \(\gamma=f_2(\sigma^\Delta)\).
Euler's formula and minimum vertex degree three give
\(\alpha\ge a/2+2\) and \(\gamma\ge c/2+2\).
Consequently \(a,c\ge2h+4\) and \(\alpha,\gamma\ge h+4\),
which prove the asserted two bounds.

For attainment, take a rational prism \(A\) over an \((h+2)\)-gon.
It has \(2h+4\) vertices and \(h+4\) facets. Its polar \(B\) has
the exchanged counts, so \(A*B\) has \(3h+8\) vertices and facets.
Both factors admit assignments because they are simple and simplicial,
respectively \cite[Corollary~9]{jahnWinter2021}.
Equation~\eqref{eq:join-cancellation-defect} gives deficiency \(h\).
This is an instance of the free-join construction of
\cite[Theorem~11]{jahnWinter2021}.
\end{proof}

Equality in this benchmark is rigid. We first isolate the face-lattice
argument, which is independent of dimension seven. For a face \(\sigma\)
of \(P\), a facet \(F\) of \(P\) \emph{induces} the face \(F\cap\sigma\)
of \(\sigma\). A \emph{minimal coface} of \(\sigma\) is a face of
dimension \(\dim\sigma+1\) containing \(\sigma\); these correspond to the
facets of \(\sigma^\Delta\). An \emph{outside vertex} is a vertex of
\(P\) not in \(\sigma\).

\begin{lemma}[Equality forces a join]\label{lem:face-equality-join}
Let \(\sigma\) be a proper face of a \(d\)-polytope \(P\), and suppose
that both \(\sigma\) and its polar dual face \(\sigma^\Delta\) have
positive dimension. Write \(a,\alpha\) for the numbers of vertices and
facets of \(\sigma\), and \(c,\gamma\) for those of
\(\sigma^\Delta\). If
\[
 f_0(P)=a+\gamma,\qquad f_{d-1}(P)=c+\alpha,
\]
then \(P\) is combinatorially a free join \(\sigma*B\), where \(B\)
is a face with vertex set \(\V(P)\setminus\V(\sigma)\).
\end{lemma}

\begin{proof}
Each facet \(\tau\) of \(\sigma\) is induced by a facet of \(P\)
not containing \(\sigma\). Such a facet can induce at most one
\(\tau\). There are exactly \(\alpha\) facets not containing
\(\sigma\), so each \(\tau\) has a unique inducing facet
\(G_\tau\), and these exhaust that group.

By duality, the \(\gamma\) facets of \(\sigma^\Delta\) correspond
to the minimal cofaces \(\Pi\) of \(\sigma\). Distinct such cofaces
intersect in \(\sigma\), and each contains an outside vertex.
There are exactly \(\gamma\) outside vertices. Consequently every
minimal coface has exactly one, say \(v\), and every outside vertex
occurs in this way. Thus \(\Pi=\sigma*v\) is a pyramid.
For each facet \(\tau\) of \(\sigma\), its join with \(v\) is a
facet of \(\Pi\). A facet of \(P\) cutting \(\tau*v\) from \(\Pi\)
cannot contain \(\sigma\), and it induces \(\tau\) on \(\sigma\);
by uniqueness it is \(G_\tau\). Hence every outside vertex belongs to
every \(G_\tau\).

It follows that \(B=\bigcap_\tau G_\tau\) has exactly the outside
vertices: no vertex of \(\sigma\) lies in every facet of \(\sigma\).
All other facets of \(P\) contain \(\sigma\). The first group
therefore induces the facets of \(\sigma\) and contains \(B\);
the second contains \(\sigma\). Each facet of \(B\) is \(B\cap F\) for a
facet \(F\) not containing \(B\); as every \(G_\tau\) contains \(B\),
such an \(F\) belongs to the second group. Hence, for a set \(X\) of
facets in the first group and a set \(Y\) in the second, the vertex set
of \(\bigcap X\cap\bigcap Y\) is
\(\V(\bigcap_{G_\tau\in X}\tau)\cup\V(B\cap\bigcap Y)\), and the two
parts range independently over the vertex sets
\(\V(\sigma')\cup\V(B')\), for faces \(\sigma'\le\sigma\),
\(B'\le B\), including empty faces.
Conversely, every face of \(P\) is an intersection of facets and has
this form. Its face lattice is therefore the product of the two
factor face lattices, which is the face lattice of their free join.
\end{proof}

The following refinement allows \(s\) extra facets while the vertex
count stays extremal, provided that few facets induce exceptionally and
the dual face has few nonsimple vertices.

\begin{lemma}[Near-equality rigidity]\label{lem:near-equality-join}
Let \(\sigma\) be a proper face of a \(d\)-polytope \(P\), with
\(r=\dim\sigma\ge1\) and \(t=\dim\sigma^\Delta\ge1\).
Write \((a,\alpha)\) and \((c,\gamma)\) for the vertex--facet counts
of \(\sigma\) and \(D=\sigma^\Delta\), respectively.
Let \(\iota(\sigma)\) be the largest number of facets of \(\sigma\)
containing one vertex, and let \(e(D)\) count the nonsimple vertices of
\(D\). Suppose that, for an integer \(s\ge0\),
\[
 f_0(P)=a+\gamma,\qquad f_{d-1}(P)=c+\alpha+s,
 \qquad \alpha-s>\iota(\sigma),\qquad 2s+e(D)\le t+1.
\]
Then \(s=0\), and \(P\) is a geometric free join \(\sigma*B\),
where \(B\) is a face with vertex set
\(\V(P)\setminus\V(\sigma)\).
\end{lemma}

\begin{proof}
The \(\gamma\) minimal cofaces of \(\sigma\) correspond to the
facets of \(D\). They have dimension \(r+1\), intersect pairwise in
\(\sigma\), and each has an outside vertex. Since there are exactly
\(\gamma\) outside vertices, each coface has exactly one, every
outside vertex occurs in a coface, and these cofaces are pyramids
\(\sigma*v\).

For each facet \(\tau\) of \(\sigma\), let \(q_\tau\ge1\) count
the facets of \(P\) inducing exactly \(\tau\) on \(\sigma\).
Let \(b_0\) count the remaining facets not containing \(\sigma\),
whose intersections with \(\sigma\) are smaller or empty. Then
\[
 s=b_0+\sum_\tau(q_\tau-1).
\]
At least \(\alpha-s\) facets \(\tau\) therefore have a unique
inducing facet \(G_\tau\). For every outside vertex \(v\),
the facet \(\tau*v\) of the pyramid \(\sigma*v\) is cut out by a
facet of \(P\); that facet induces \(\tau\) on \(\sigma\), so it is
necessarily \(G_\tau\). Thus all outside
vertices belong to each uniquely inducing \(G_\tau\).
Their intersection is a nonempty face \(B\). It has no vertex of
\(\sigma\), since more than \(\iota(\sigma)\) distinct facets
of \(\sigma\) occur in the intersection. Its vertices are precisely
the outside vertices, so \(P=\conv(\sigma\cup B)\).
Writing \(k=\dim B\), we have
\(d=r+t+1\le r+k+1\), hence \(k\ge t\).

There are at most \(2s\) facets not containing \(\sigma\) other
than the unique \(G_\tau\): a group of size \(q_\tau\ge2\)
contributes at most \(2(q_\tau-1)\), and each of the \(b_0\) remaining
facets contributes one. The unique facets contain \(B\) and cannot induce
its facets. A facet \(K\) containing \(\sigma\) corresponds to a
vertex \(w\) of \(D\). Its outside vertices correspond to the
facets of \(D\) containing \(w\). If \(w\) is simple, there are
exactly \(t\) such vertices.

If \(k>t\), every facet of \(B\) has at least \(t+1\) vertices,
so it cannot be induced by such a simple-vertex facet \(K\).
Only the \(e(D)\) nonsimple-vertex facets and the at most \(2s\)
exceptional facets remain. Every facet of \(B\) is induced by an
ambient facet, and each ambient facet induces at most one. Consequently
\[
 f_{k-1}(B)\le 2s+e(D)\le t+1<k+1,
\]
contrary to the minimum facet count of a \(k\)-polytope.
Thus \(k=t\). The affine hulls of \(\sigma\) and \(B\) span
the maximum dimension \(r+t+1\), so they are skew and
\(P=\conv(\sigma\cup B)\) is a geometric free join.
Its face figure at \(\sigma\) is \(B\), which is dual to \(D\).
Therefore \(B\) has \(c\) facets, the join has \(\alpha+c\)
facets, and comparison with the assumed count gives \(s=0\).
\end{proof}

\begin{theorem}[Seven-dimensional equality cases]\label{thm:seven-rigidity}
For an integer \(h\ge1\), a seven-polytope of deficiency at least \(h\) has
larger class size \(3h+8\) if and only if it is combinatorially
\(A*B\), where \(A\) is a simple three-polytope with
\(2h+4\) vertices and \(h+4\) facets, and \(B\) is a simplicial
three-polytope with \(h+4\) vertices and \(2h+4\) facets.
Every such join has both class sizes \(3h+8\) and deficiency exactly
\(h\).
\end{theorem}

\begin{proof}
Use the maximizing face and notation of
Proposition~\ref{prop:global-seven}. Since both class sizes are at
least \(3h+8\) and the larger one equals \(3h+8\), both equal \(3h+8\).
Equality in both size bounds forces
\[
 a=c=2h+4,\qquad \alpha=\gamma=h+4,\qquad
 f_0(P)=a+\gamma,\quad f_6(P)=c+\alpha.
\]
Equality in the minimum-degree form of Euler's bound makes both
three-faces simple. Lemma~\ref{lem:face-equality-join} gives
\(P\cong A*B\), with \(A=\sigma\). The face figure of \(\sigma\) in
\(\sigma*B\) is \(B\), so \(\sigma^\Delta\cong B^*\); as
\(\sigma^\Delta\) is simple, \(B\) is simplicial with the claimed
exchanged counts. Conversely, the simple and simplicial factors
admit assignments \cite[Corollary~9]{jahnWinter2021}. Their opposite
imbalances both have magnitude \(h\), and
\eqref{eq:join-cancellation-defect} proves the assertion.
\end{proof}

\begin{corollary}[Connected seven-dimensional lower bound]
\label{cor:connected-seven-floor}
A seven-polytope with connected nonincidence graph and deficiency at
least \(h\), where \(h\ge1\) is an integer, has at least \(3h+10\)
vertices and at least \(3h+10\) facets.
\end{corollary}

\begin{proof}
By polarity it suffices to exclude vertex counts \(3h+8\) and
\(3h+9\). Continue with the maximizing face notation above, and put
\(n=f_0(P)\), \(M=f_6(P)\). We have
\[
 a,c\ge2h+4,\quad \alpha,\gamma\ge h+4,\quad
 n\ge a+\gamma,\quad M\ge c+\alpha,\quad
 a+c\ge\max\{n,M\}+h.
\]
If \(n=3h+8\), then \(a=2h+4\), \(\gamma=h+4\).
Euler's bound \(c\le2\gamma-4\) forces \(c=2h+4\), and the
last inequality forces \(M=3h+8\). Theorem~\ref{thm:seven-rigidity}
then gives a free join.

If \(n=3h+9\), the pair \((a,\gamma)=(2h+4,h+4)\) is impossible,
since it would give \(a+c\le4h+8<n+h\). There are only two other
pairs. If \((a,\gamma)=(2h+4,h+5)\), the inequalities
\(\alpha\ge h+4\), \(a-\alpha\ge h\) give \(\alpha=h+4\).
Then \(c+\alpha\le M\le a+c-h=c+\alpha\), so both face-count
bounds are equalities. If \((a,\gamma)=(2h+5,h+4)\), Euler's bound
gives \(c=2h+4\) and \(\alpha\ge h+5\). Hence
\(3h+9\le c+\alpha\le M\le a+c-h=3h+9\), again giving both
equalities. Lemma~\ref{lem:face-equality-join} applies in either case.
Every nontrivial free join has disconnected nonincidence graph by
Lemma~\ref{lem:join-disconnects}, a contradiction.
\end{proof}

Near-equality rigidity improves Corollary~\ref{cor:connected-seven-floor}
by excluding the next two sizes, \(3h+10\) and \(3h+11\).

\begin{theorem}[Small seven-dimensional obstructions are joins]
\label{thm:small-seven-join}
Let \(h\ge1\) be an integer, and let \(P\) be a seven-polytope with
\(\delta(P)\ge h\). If either \(f_0(P)\le3h+11\) or
\(f_6(P)\le3h+11\), then \(P\) is combinatorially a nontrivial free
join.
In particular, if \(\N(P)\) is connected, then
\[
 f_0(P)\ge3h+12,\qquad f_6(P)\ge3h+12.
\]
\end{theorem}

\begin{proof}
By polarity it suffices to treat \(n=f_0(P)\le3h+11\).
Choose a face \(\sigma\) attaining \(H(P)\), and put
\(D=\sigma^\Delta\), \(M=f_6(P)\). As in the proof of
Proposition~\ref{prop:global-seven}, both \(\sigma\) and \(D\) are
three-dimensional. With their vertex--facet counts
\((a,\alpha),(c,\gamma)\), introduce the nonnegative integer slacks
\[
 u=n-a-\gamma,\qquad v=M-c-\alpha.
\]
The witness inequality \(a+c\ge\max\{n,M\}+h\) gives
\(a-\alpha\ge h+v\) and \(c-\gamma\ge h+u\).
Euler's bound gives
\(a,c\ge2h+4\), \(\alpha,\gamma\ge h+4\). Write
\[
 a=2h+4+x,\quad \alpha=h+4+y,\quad
 c=2h+4+z,\quad \gamma=h+4+w.
\]
The four increments are nonnegative integers satisfying
\[
 x-y\ge v,\qquad z-w\ge u,\qquad x\le2y,\qquad z\le2w.
\]
In particular, \(x\ge2v\), \(w\ge u\), and
\[
 2(u+v)\le x+w+u=n-(3h+8)\le3.
\]
Thus \(u+v\le1\). If \(u=v=0\),
Lemma~\ref{lem:face-equality-join} gives the claimed join.

For a three-polytope \(S\), Euler's formula expresses the total
vertex-degree excess as
\[
 E(S):=\sum_{v_0\in\V(S)}(\deg_S v_0-3)=2f_2(S)-f_0(S)-4.
\]
The number of facets at a vertex equals its degree. Thus \(E(S)\le1\)
implies \(\iota(S)\le4\) and at most one nonsimple vertex.
If \(u=0,v=1\), then \(x\ge2\), \(w\le1\), and
\(1\le y\le x-1\), so
\[
 0\le E(\sigma)=2y-x\le x-2\le1,
 \qquad 0\le E(D)=2w-z\le w\le1.
\]
Moreover \(\alpha\ge h+5\ge6\).
Lemma~\ref{lem:near-equality-join} applies with \(r=t=3\) and \(s=1\):
\(\alpha-1\ge5>\iota(\sigma)\) and \(2+e(D)\le3\le4\).
It forces \(s=0\), a contradiction.

If \(u=1,v=0\), then \(w\ge1\), \(x+w\le2\), \(y\le x\),
and \(z\ge w+1\). Hence
\[
 0\le E(\sigma)=2y-x\le x\le1,
 \qquad 0\le E(D)=2w-z\le w-1\le1,
\]
and \(\gamma\ge h+5\ge6\). Apply the same lemma in the polar
polytope, with face \(D\), whose dual face is \(\sigma\): the vertex
count \(M=c+\alpha\) is an equality, the facet count
\(n=a+\gamma+1\) has slack \(s=1\), \(\gamma-1\ge5>4\ge\iota(D)\), and
\(2+e(\sigma)\le3\le4\). This again gives a contradiction.
Only \(u=v=0\) remains. Every nontrivial free join has disconnected
nonincidence graph by Lemma~\ref{lem:join-disconnects}, proving the
connected bounds.
\end{proof}

Thus, among seven-polytopes without an assignment, the minimum larger
class size is \(11\); by Theorem~\ref{thm:seven-rigidity} the attainers
are exactly the polytopes combinatorially equivalent to the free join of
a triangular prism and a triangular bipyramid. The following subsections strengthen the
connected bound to \(17\) on both sides, and
Section~\ref{sec:shared-coordinate} attains it.

\subsection{Structure at the first possible connected size}
\label{subsec:first-size}

At the first size not excluded by Theorem~\ref{thm:small-seven-join},
namely \(3h+12\), the face data are already far more rigid than for an
arbitrary connected obstruction.
The following two results describe that data and its geometric allocation.
Theorem~\ref{thm:seven-first-size-exclusion} will show that these
necessary conditions cannot hold simultaneously.

\begin{theorem}[First-size Hall faces]\label{thm:first-seven-size}
Let \(h\ge1\) be an integer, and let \(P\) be a seven-polytope with
connected nonincidence graph and \(\delta(P)\ge h\).
If either \(f_0(P)\le3h+12\) or \(f_6(P)\le3h+12\), then
\[
 f_0(P)=f_6(P)=3h+12,\qquad \delta(P)=h.
\]
Every face \(\sigma\) attaining \(H(P)\), and its dual face
\(D=\sigma^\Delta\), are simple three-polytopes with
\[
 f_0(\sigma)=f_0(D)=2h+6,\qquad
 f_2(\sigma)=f_2(D)=h+5.
\]
Both face-count slacks are one:
\(f_0(P)-f_0(\sigma)-f_2(D)=
f_6(P)-f_0(D)-f_2(\sigma)=1\).
\end{theorem}

\begin{proof}
By polarity suppose \(n=f_0(P)\le3h+12\), and put \(M=f_6(P)\).
Theorem~\ref{thm:small-seven-join} gives \(n=3h+12\).
Choose any maximizing face \(\sigma\); both it and \(D\) are
three-dimensional, as in the proof of
Proposition~\ref{prop:global-seven}.
As in the preceding proof, write their counts as
\[
 a=2h+4+x,\quad \alpha=h+4+y,\quad
 c=2h+4+z,\quad \gamma=h+4+w,
\]
and put \(u=n-a-\gamma\), \(v=M-c-\alpha\).
All increments and slacks are nonnegative integers, with
\[
 x-y\ge v,\quad z-w\ge u,\quad x\le2y,\quad z\le2w,
 \quad x\ge2v,\quad w\ge u,\quad x+w+u=4.
\]
In particular \(u+v\le2\). The case \(u=v=0\) forces a join by
Lemma~\ref{lem:face-equality-join}, contrary to connectedness.

For a three-polytope \(S\), its degree excess
\(E(S)=2f_2(S)-f_0(S)-4\) bounds both the number of nonsimple
vertices and \(\iota(S)-3\). Suppose \(u=0\), \(v>0\).
Then \(E(\sigma)=2y-x\) and \(E(D)=2w-z\le w\), so
\[
 \alpha-v-\iota(\sigma)
 \ge h+1+x-y-v\ge h+1>0,
 \qquad 2v+e(D)\le2v+w\le x+w=4.
\]
Lemma~\ref{lem:near-equality-join}, with \(s=v\) and \(t=3\),
forces \(v=0\), a contradiction.
If \(v=0\), \(u>0\), then \(y\le x\) and
\(E(\sigma)\le x\). The polar budgets are
\[
 \gamma-u-\iota(D)\ge h+1+z-w-u\ge h+1>0,
 \qquad 2u+e(\sigma)\le2u+x\le u+w+x=4.
\]
The same lemma applied at \(D\) in the polar polytope forces
\(u=0\), again a contradiction.

Thus \(u=v=1\). From \(x\ge2\), \(w\ge1\) and
\(x+w+u=4\) we get \(x=2\), \(w=1\).
The other inequalities give \(y=1\), \(z=2\).
Both degree excesses vanish, so both faces are simple, with the stated
counts. Finally \(M=c+\alpha+v=3h+12\) and
\(H(P)=a+c=4h+12\), proving \(\delta(P)=h\).
The maximizing face was arbitrary, and polarity covers the other count
assumption.
\end{proof}

\begin{theorem}[Coface allocation]\label{thm:first-size-cofaces}
Under the hypotheses of Theorem~\ref{thm:first-seven-size}, fix any
face \(\sigma\) attaining \(H(P)\). Every vertex outside \(\sigma\)
belongs to a minimal coface of \(\sigma\). Exactly one minimal coface
contains two outside vertices; every other minimal coface is a pyramid
over \(\sigma\). Dually, every facet not containing \(\sigma\)
induces a facet of \(\sigma\). Exactly one facet of \(\sigma\)
has two inducing facets of \(P\); all others have one.
\end{theorem}

The proof excludes a vertex belonging to none of the minimal cofaces.
Such a vertex first forces a complementary three-face. Viewing the
polytope as a one-point extension of the resulting join then produces
too many inducing facets. Unlike Lemma~\ref{lem:stacking}, the argument
does not assume that the added vertex lies strictly beneath every
facet \(\sigma*J\); it uses only those \(J\) for which this holds.

\begin{proof}
Write \(a=c=2h+6\), \(\alpha=\gamma=h+5\) for the counts of
\(\sigma\) and \(D\). There are \(\gamma\) minimal cofaces,
each of dimension four, and \(\gamma+1\) outside vertices.
Distinct minimal cofaces intersect in \(\sigma\).
Thus either one coface has two outside vertices and all others have
one, or all are pyramids \(\sigma*v_i\) and one outside vertex
\(x_0\) belongs to none. We exclude the latter alternative.

\smallskip\noindent\emph{The quotient face figure.}
Choose a point in the relative interior of \(\sigma\), and project
the cone of feasible directions of \(P\) at this point modulo the directions of
\(\aff\sigma\). This is a pointed four-dimensional cone.
A strictly positive linear functional gives a three-dimensional
section \(T\), the face figure at \(\sigma\)
\cite[Chapter~2]{ziegler1995}. Its face lattice is dual to that of
\(D\), so \(T\) is simplicial. Its vertices correspond to the
\(\gamma\) minimal cofaces. The normalized image of \(x_0\)
cannot be a vertex of \(T\), since an image on an extreme ray places
\(x_0\) in the corresponding minimal coface.
A nonvertex point of a three-polytope lies on at most two facets: its
minimal face is an edge, a two-face, or the whole polytope.
Consequently \(x_0\) belongs to at most two facets of \(P\)
containing \(\sigma\).

\smallskip\noindent\emph{A complementary face.}
For each facet \(\tau\) of \(\sigma\), let \(q_\tau\ge1\)
count its inducing facets, and let \(b_0\) count the facets not
containing \(\sigma\) that induce a smaller or empty face. Since
there are \(\alpha+1\) facets not containing \(\sigma\),
\[
 1=b_0+\sum_\tau(q_\tau-1).
\]
At least \(\alpha-1\) facets have a unique inducing facet
\(G_\tau\). Each such \(G_\tau\) contains every \(v_i\):
the pyramid facet \(\tau*v_i\) must be induced by \(G_\tau\).
Let \(B\) be the intersection of these unique inducing facets.
It contains every \(v_i\) and possibly \(x_0\), but no vertex
of \(\sigma\), because \(\sigma\) is simple and
\(\alpha-1\ge5>3\). Put \(k=\dim B\).
The projected directions of the \(v_i\) span four dimensions, so
their affine span has dimension at least three; hence \(k\ge3\).

If \(k>3\), each facet of \(B\) needs at least four vertices.
A facet of \(P\) containing \(\sigma\) contains precisely three
of the \(v_i\), because the corresponding vertex of \(D\) is
simple. Only the at most two such facets containing \(x_0\)
can therefore induce a facet of \(B\).
Unique inducing facets contain \(B\); the remaining noncontaining
facets number at most two. Every facet of \(B\) is induced by an
ambient facet, so \(f_{k-1}(B)\le4<k+1\), a contradiction.
Thus \(k=3\). If \(x_0\in B\), then
\(P=\conv(\sigma\cup B)\) has the maximal join dimension seven,
and is a geometric free join, contrary to connectedness.
Hence \(B=\conv\{v_1,\ldots,v_\gamma\}\).

\smallskip\noindent\emph{The one-point contradiction.}
The polytope \(P_0=\conv(\sigma\cup B)\) has dimension seven:
the directions of \(\sigma\) and the four-dimensional quotient
cone span the ambient space. Thus \(P_0=\sigma*B\), and
\(P=\conv(P_0\cup\{x_0\})\).
The quotient cones of \(P_0\) and \(P\) at \(\sigma\) agree,
because the extreme rays already come from the \(v_i\).
All \(c\) facets \(\sigma*J\), for facets \(J\) of \(B\),
remain supporting facets of \(P\), and \(x_0\) is on at most
two of their hyperplanes. A unique inducing facet has the supporting
hyperplane of \(\tau*B\), since it contains this six-face of \(P_0\).
If all groups were unique, \(x_0\) would satisfy every facet
inequality of \(P_0\), contradicting that it is a vertex outside
\(P_0\). Therefore exactly one facet \(\tau_0\) has
\(q_{\tau_0}=2\), and \(b_0=0\).
The only facet of \(P_0\) beyond which \(x_0\) can lie is
\(S=\tau_0*B\), and it must lie strictly beyond it.

Choose affine facet inequalities \(g_0\le0\) for \(S\) and
\(g_J\le0\) for \(\sigma*J\) on \(P_0\).
Put \(\lambda=g_0(x_0)>0\). For every \(J\) with
\(\mu_J=g_J(x_0)<0\), the inequality
\[
 \lambda g_J-\mu_J g_0\le0
\]
holds on \(P\) (compare \eqref{eq:cone-form}). Its equality set on \(P_0\) is exactly the ridge
\(\tau_0*J\); adjoining \(x_0\), which is beyond \(S\),
gives a six-dimensional facet of \(P\). It induces \(\tau_0\)
on \(\sigma\). Distinct \(J\) give distinct facets, since their
equality sets on \(P_0\) differ. There are at least
\(c-2\ge6\) such facets, contrary to \(q_{\tau_0}=2\).
This excludes \(x_0\) and proves the vertex allocation.

\smallskip\noindent\emph{The dual allocation.}
Apply the same argument at \(D\) in the polar polytope. A polar
vertex outside \(D\) corresponds to a facet \(G\) not containing
\(\sigma\); it belongs to a minimal coface of \(D\) precisely
when \(G\cap\sigma\) is a facet of \(\sigma\).
Thus no smaller intersection occurs, and the single surplus polar
vertex gives exactly one inducing multiplicity two.
\end{proof}

Thus a hypothetical connected deficiency-one example with fifteen
vertices or fifteen facets would have both counts fifteen. Each
maximizing Hall three-face and its dual would have eight vertices and
six facets, with one double-vertex coface and one double-induced facet.
The next argument excludes this coupled configuration.

\subsection{Excluding the first connected size}
\label{subsec:first-size-exclusion}

The face counts and allocation in the preceding two theorems are not
realizable. The obstruction uses the rank of a slack matrix, rather
than a further face-count inequality. We recall that a slack matrix
of a \(d\)-polytope has rank \(d+1\)
\cite[Theorem~14]{gouveiaEtAl2013}. Positive row and column scalings
preserve its rank and its incidence zeros; a row-scaled matrix need
not itself be a slack matrix of the same affine realization.

\begin{theorem}[Exclusion of the first connected size]
\label{thm:seven-first-size-exclusion}
For every integer \(h\ge1\), a seven-polytope \(P\) with connected
nonincidence graph and \(\delta(P)\ge h\) satisfies
\[
 f_0(P)\ge3h+13,\qquad f_6(P)\ge3h+13.
\]
\end{theorem}

We use two elementary properties of simple three-polytopes.
If \(A\) is simple with at least six facets, deleting any one
column from its facet-slack matrix leaves rank four. Indeed, the
corresponding simplicial polar cannot be a pyramid: the base of a
simplicial three-dimensional pyramid is a triangle, giving only four
vertices. Thus its remaining homogeneous vertices span four dimensions.
Also, any three distinct facet-slack columns are linearly independent,
since three vertices of a convex polytope cannot be collinear.

\begin{proof}
Suppose that either count is at most \(3h+12\).
Theorems~\ref{thm:first-seven-size} and~\ref{thm:first-size-cofaces}
give a maximizing Hall face \(A=\sigma\) and dual face \(D\),
both simple three-polytopes with
\[
 a=c=2h+6,\qquad \alpha=\gamma=h+5.
\]
There are \(\gamma-1\) pyramid cofaces \(A*v_i\) and one
four-dimensional coface \(H=\conv(A\cup\{p,q\})\).
Every facet \(\tau\) of \(A\) has a unique inducing facet
\(G_\tau\), except \(\tau_0\), which has two, \(G_0,G_1\).
Write \(I=\{v_i\}\), and let \(J\) index the unique inducing
facets. Each \(G_\tau\), \(\tau\ne\tau_0\), contains all of
\(I\), since it induces the pyramid facet \(\tau*v_i\).

\smallskip\noindent\emph{Rank and the exceptional corner.}
Partition a slack matrix by rows \(\V(A),I,\{p,q\}\) and
columns \(C,J,\{G_0,G_1\}\), where \(C\) consists of the facets
containing \(A\). Its upper-left block is zero.
The block \(X\) on \(\V(A)\) and the columns outside \(C\)
is, up to positive column scaling, a facet-slack matrix of \(A\) in
which the \(\tau_0\) column occurs twice: each inducing facet restricts
to \(\aff A\) as a positive multiple of the corresponding facet form of
\(A\). Scale the two \(\tau_0\) columns positively so that they agree.
The block \(Y\) on the outside vertices and \(C\) is a slack
matrix of the pointed four-dimensional quotient cone at \(A\).
The rows for \(p,q\) represent the same extreme ray; scale them
positively to make them agree. All remaining entries stay nonnegative.
Denote the outside-vertex/noncontaining-facet block by \(Z\).

The restrictions \(X_J\) of \(X\) to the columns \(J\) and \(Y_I\) of
\(Y\) to the rows \(I\) both have rank four, by the preceding
observation applied to \(A\) and \(D\). For \(Y_I\), this means that
the \(\gamma-1\) rays indexed by \(I\) span the four-dimensional
quotient, because the simplicial figure at \(A\), dual to \(D\) and
with \(\gamma\ge6\) vertices, is not a pyramid. The submatrix on
rows \(\V(A),I\) and columns \(C,J\) is
\[
 \begin{pmatrix}0&X_J\\Y_I&0\end{pmatrix},
\]
of rank eight, the rank of the full matrix.
Consequently, for \(i\in\{p,q\}\) there is a unique vector
\(d_i\in\operatorname{rowspan}X\) whose restriction to \(J\)
is \(Z_{i,J}\). Likewise, for \(j\in\{0,1\}\) there is a
unique vector \(e_j\in\operatorname{colspan}Y\) whose
restriction to \(I\) is \(Z_{I,G_j}\).
Let \(\xi_i\) be the duplicated-column coordinate of \(d_i\),
and \(\zeta_j\) the duplicated-row coordinate of \(e_j\).
Then
\begin{equation}\label{eq:first-size-corner}
 Z_{i,G_j}=\xi_i+\zeta_j
 \qquad(i=p,q;\ j=0,1).
\end{equation}
To see this explicitly, express the common row \(Y_p=Y_q\)
as a (not necessarily unique) linear combination \(\ell Y_I\).
Writing \(e_j=Y\beta_j\), so that \(Z_{I,G_j}=Y_I\beta_j\) and
\(\zeta_j=Y_p\beta_j\), the same combination of the rows indexed by
\(I\) has corner entries \(\ell Y_I\beta_j=Y_p\beta_j=\zeta_j\),
independently of the choice of \(\ell\).
Subtracting it from row \(i\) leaves zero on \(C\). Write the
difference as a combination of the rows indexed by \(\V(A)\cup I\);
its \(I\)-part vanishes on \(C\), and on \(J\) because every
\(G_\tau\) contains \(I\). A row-space vector vanishing on \(C\cup J\)
is zero, by the rank-eight submatrix above; so the difference is a
combination of the rows indexed by \(\V(A)\). Its restriction
to \(J\) determines \(d_i\), and its corner entries are \(\xi_i\).
The rank-eight submatrix makes this extension unique.

\smallskip\noindent\emph{No exceptional row or column vanishes.}
Suppose \(d_p=0\). The face
\(B=\bigcap_{\tau\ne\tau_0}G_\tau\) contains \(I\) and \(p\),
and possibly \(q\), but no vertex of \(A\): its defining
\(\alpha-1\ge5\) facets restrict to distinct facets of the simple
three-polytope \(A\).
Their homogeneous forms have rank at least four, since \(X_J\)
does. Thus \(\dim B\le3\). Its vertices represent all
\(\gamma\) quotient extreme rays, which span four dimensions,
so \(\dim B\ge3\).
If \(q\in B\), then \(P=\conv(A\cup B)\) has the maximal
join dimension seven and is a free join, contrary to connectedness.
Otherwise \(B=\conv(I\cup\{p\})\), and
\(P_0=\conv(A\cup B)=A*B\) is a full-dimensional join.
Its quotient face figure at \(A\) is the original figure, so
\(B\) is simplicial with \(\gamma\) vertices and
\(c=2\gamma-4\) facets.

The remaining vertex \(q\) projects to the same quotient ray as
\(p\). Hence every facet \(A*K\) of \(P_0\) stays supporting
for \(P\), and \(q\) is on its hyperplane exactly when \(p\in K\).
Each unique \(G_\tau\) has the hyperplane of \(\tau*B\), so
\(q\) satisfies its inequality. Since \(q\notin P_0\), it must
lie strictly beyond the sole remaining factor facet \(\tau_0*B\).
The supporting-inequality construction in the proof of
Theorem~\ref{thm:first-size-cofaces} now gives a distinct facet
of \(P\) inducing \(\tau_0\) for every facet \(K\) of \(B\)
not containing \(p\). There are at least
\[
 c-\deg_B(p)\ge(2\gamma-4)-(\gamma-1)
 =\gamma-3\ge3
\]
such facets, contradicting the inducing multiplicity two.
Thus \(d_p,d_q\ne0\). Applying the argument in the polar
polytope, at \(D\), also gives \(e_0,e_1\ne0\).

\smallskip\noindent\emph{The three mixed facets.}
Every facet of \(H\) other than \(A\) meets \(A\) in a facet:
it is induced by a facet of \(P\) not containing \(A\), to which
the allocation theorem applies. Each facet \(\tau\) of \(A\) is a
ridge of \(H\), hence lies in exactly two facets of \(H\): \(A\) and
one further facet \(L_\tau\).
Every vertex of \(A\) is therefore simple in \(H\), with three
edges in \(A\) and exactly one edge to \(p\) or \(q\).
Color it by this unique outside neighbor. Each color has at least
three vertices: each of \(p,q\) has degree at least four in \(H\),
and can have at most one neighbor outside \(A\).

Call a facet of \(A\) mixed if it has both colors. For a mixed
\(\tau\), the facet \(L_\tau\) contains both \(p,q\), since
the outside edge at a vertex of \(\tau\) lies in \(L_\tau\).
For a unique inducing facet this gives
\(d_p(\tau)=d_q(\tau)=0\).
For every unique \(\tau\), mixed or not, at least one of these
coordinates is zero: \(L_\tau\) is three-dimensional and cannot
have all its vertices in the two-dimensional \(\tau\).

There cannot be three unique mixed facets. Their three columns in
\(X\) are independent, so the nonzero vectors \(d_p,d_q\)
would belong to one common one-dimensional subspace. Their
restrictions to \(J\) are nonnegative and nonzero, making the
proportionality factor positive. Since at least one coordinate
vanishes at every column of \(J\), both restrictions would then
vanish identically, a contradiction.

There are at least three mixed facets altogether. Since the graph of
\(A\) is connected and both colors occur, some edge joins the two
colors, and each such edge belongs to two mixed facets. If there were
at most two, every such edge would lie in their intersection, so there
would be exactly one. That edge would be a bridge in the connected graph of
\(A\), impossible because every edge lies on a polygonal facet cycle.
Thus there are exactly three mixed facets: \(\tau_0,F_1,F_2\),
with \(F_1,F_2\) uniquely induced.
Each mixed polygon has at least two edges between the colors.
Every such edge belongs to two mixed facets, while a pair of facets
has at most one edge in common. It follows that there are exactly
three cross-color edges, one for each pair of mixed facets.

These three facets have empty common intersection. Otherwise,
at a common vertex \(v\), simplicity makes their three pairwise
intersection edges precisely the three edges at \(v\), all of
which cross between the colors. The graph of \(A\) with \(v\)
removed is connected: the neighbors of \(v\) are joined pairwise
by the boundary paths avoiding \(v\) in its three incident
polygonal facets, and every component must contain such a neighbor.
All remaining vertices would then have the opposite color from \(v\),
contradicting the minimum color-class size three.
In particular the edge \(E=F_1\cap F_2\) is disjoint from \(\tau_0\).

\smallskip\noindent\emph{Positive exceptional coordinates.}
Let \(u,v\) be the endpoints of \(E\). Since \(d_i\) vanishes
on the independent columns for \(F_1,F_2\), it lies in the
two-dimensional subspace of the row space of \(X\) vanishing on these
two columns, which is spanned by the rows \(X_u,X_v\):
\[
 d_i=\lambda_iX_u+\mu_iX_v.
\]
The third facet at \(u\) does not contain \(v\); the third at
\(v\) does not contain \(u\). Neither is \(\tau_0\), so both
are unique inducing columns. Nonnegativity there gives
\(\lambda_i,\mu_i\ge0\). Since \(d_i\ne0\), at least one is positive.
Both endpoints are outside \(\tau_0\), yielding \(\xi_i>0\).
The same color and edge argument at \(D\) in the polar polytope
gives \(\zeta_0,\zeta_1>0\). Positive row and column scalings
do not change these conclusions.
Equation~\eqref{eq:first-size-corner} now makes all four exceptional
corner entries strictly positive.

Finally, \(L_{\tau_0}\) has dimension three and contains the
two-dimensional \(\tau_0\), so it contains \(p\) or \(q\).
An ambient facet inducing it must be one of \(G_0,G_1\).
At least one exceptional corner entry must therefore be zero,
the required contradiction.
\end{proof}

For deficiency one this theorem gives sixteen vertices and sixteen
facets as necessary bounds. The stronger bound below and the integral
attainer in Section~\ref{sec:shared-coordinate} determine the connected
minimum at seventeen. (The value \(\eta_{7,1}=33\) of
Theorem~\ref{thm:intro-complete} concerns the different,
certificate-restricted class \(\mathcal C_{7,1}\).)

\subsection{Capacity rigidity and simple dual faces}
\label{subsec:capacity-rigidity}

The facet budget in Lemma~\ref{lem:near-equality-join} can be
compared with a vertex capacity, rather than only with the minimum
number of facets. This gives a dimension-independent refinement.
For \(q\ge k+1\) and \(k\ge2\), let \(U_k(q)\) be the number
of facets of the cyclic \(k\)-polytope on \(q\) vertices. Explicitly,
\[
 U_{2p}(q)=\binom{q-p}{p}+\binom{q-p-1}{p-1},\qquad
 U_{2p+1}(q)=2\binom{q-p-1}{p}.
\]
For integers \(t\ge1\), \(q\ge0\), put
\[
 \Theta_t(q)=\max_{t+1\le k\le q-1}U_k(q),
\]
where an empty maximum is zero. McMullen's Upper Bound Theorem
\cite{mcmullen1970}, applied after polarity, bounds the number of
vertices of a \(k\)-polytope with at most \(q\) facets by \(U_k(q)\),
since \(U_k(q)\) is nondecreasing in \(q\).
In particular,
\[
 \Theta_t(q)=0\quad(q\le t+1),\qquad
 \Theta_t(t+2)=t+2,\qquad
 \Theta_t(t+3)=\left\lfloor\frac{(t+3)^2}{4}\right\rfloor.
\]

\begin{lemma}[Capacity-enhanced near-equality rigidity]
\label{lem:capacity-rigidity}
Let \(\sigma\) be a proper face of a \(d\)-polytope \(P\), with
\(r=\dim\sigma\ge1\), \(t=\dim D\ge1\), and
\(D=\sigma^\Delta\). Write \((a,\alpha)\), \((c,\gamma)\)
for the vertex--facet counts of \(\sigma,D\), respectively.
Let \(e(D)\) count the nonsimple vertices of \(D\), and let
\(\iota(\sigma)\) be the maximum number of its facets containing
one vertex. Suppose that, for an integer \(s\ge0\),
\[
 f_0(P)=a+\gamma,\qquad f_{d-1}(P)=c+\alpha+s,
 \qquad \alpha-s>\iota(\sigma).
\]
If \(\gamma>\Theta_t(2s+e(D))\), then \(s=0\) and \(P\)
is a geometric free join \(\sigma*B\), where \(B\) is a face
with vertex set \(\V(P)\setminus\V(\sigma)\).
\end{lemma}

\begin{proof}[Proof of Lemma~\ref{lem:capacity-rigidity}]
There are \(\gamma\) minimal cofaces of \(\sigma\), each with
an outside vertex, and distinct cofaces meet in \(\sigma\).
The outside count is exactly \(\gamma\), so every coface is a
pyramid \(\sigma*v\), and every outside vertex is used.
If \(q_\tau\ge1\) counts ambient facets inducing a facet
\(\tau\) of \(\sigma\), and \(b_0\) counts the remaining
noncontaining facets, then
\[
 s=b_0+\sum_\tau(q_\tau-1).
\]
At least \(\alpha-s\) facets have a unique inducer \(G_\tau\).
Each \(G_\tau\) contains every outside vertex, since it must
induce the pyramid facet \(\tau*v\).
Their intersection \(B\) has no vertex of \(\sigma\), by
\(\alpha-s>\iota(\sigma)\), and has exactly \(\gamma\) vertices.
Thus \(P=\conv(\sigma\cup B)\). With \(k=\dim B\),
the equality \(d=r+t+1\) gives \(k\ge t\).

There are at most \(2s\) noncontaining facets other than the unique
inducers: a group of size \(q_\tau\ge2\) contributes at most
\(2(q_\tau-1)\), and each of the \(b_0\) remaining facets contributes one.
The unique inducers contain \(B\) and do not induce its facets.
A facet containing \(\sigma\) and corresponding to a simple
vertex of \(D\) contains exactly \(t\) outside vertices.
It cannot induce a facet of \(B\) if \(k>t\).
Consequently, for \(k>t\),
\[
 f_{k-1}(B)\le N:=2s+e(D).
\]
Then \(k+1\le N\), and the classical vertex capacity gives
\[
 \gamma=f_0(B)\le U_k(N)\le\Theta_t(N),
\]
a contradiction. Hence \(k=t\). The affine spans have maximal
join dimension \(r+t+1\), so \(P=\sigma*B\).
The face figure at \(\sigma\) is dual to \(D\), hence \(B\)
has \(c\) facets. Comparing the join's \(\alpha+c\) facets
with the assumed count proves \(s=0\).
\end{proof}

For \(2s+e(D)=t+2\), the hypothesis is simply \(\gamma>t+2\).
The only remaining higher-dimensional complementary face would be
a \((t+1)\)-simplex, which has exactly \(t+2\) vertices.
This equality case accounts for the strict inequality.

The slack argument also extends beyond the balanced counts in
Theorem~\ref{thm:seven-first-size-exclusion}.
Positive row and column scalings give generalized slack matrices in
the sense of \cite[Definition~2.4 and Theorem~2.5]{gouveiaEtAl2019}.
Only rank, incidence zeros, and nonnegativity enter the following result.

\begin{proposition}[Simple dual faces exclude two unit slacks]
\label{prop:simple-double-slack}
Let \(A\) be a three-face of a seven-polytope \(P\), and put
\(D=A^\Delta\). Suppose that \(A,D\) are simple and have
vertex--facet counts \((a,\alpha)\), \((c,\gamma)\), with
\(\alpha,\gamma\ge6\). Then the two equalities
\[
 f_0(P)=a+\gamma+1,\qquad f_6(P)=c+\alpha+1
\]
cannot hold simultaneously.
\end{proposition}

The proposition does not require a maximizing Hall face or a connected
nonincidence graph. We first prove the coface allocation without
assuming \(a=c\) or \(\alpha=\gamma\), and then check each bound used
in the proof of Theorem~\ref{thm:seven-first-size-exclusion}.

\begin{proof}[Proof of Proposition~\ref{prop:simple-double-slack}]
Suppose both equalities hold. There are \(\gamma\) minimal cofaces
and \(\gamma+1\) outside vertices. If one vertex \(x_0\) lies
in none of the cofaces, all cofaces are pyramids \(A*v_i\).
At least \(\alpha-1\) facets of \(A\) have unique inducers.
Their intersection \(B\) contains all \(v_i\) and possibly \(x_0\);
since \(\alpha-1\ge5>3=\iota(A)\), it contains no vertex of \(A\).
Deleting one facet-slack column of simple \(A\) leaves rank four,
as proved before Theorem~\ref{thm:seven-first-size-exclusion}.
The homogeneous forms defining \(B\) therefore have rank at least
four, so \(\dim B\le3\). The quotient ray span gives
\(\dim B\ge3\).
If \(x_0\in B\), then \(P=A*B\) is a full-dimensional join.
Its figure at \(A\) would have \(\gamma+1\) vertices, contrary
to the actual \(\gamma\) extreme rays.

Otherwise \(P_0=A*B\) has the same quotient cone at \(A\) as
\(P\), and all its \(c\) containing facets remain supporting.
The nonvertex quotient image of \(x_0\) lies on at most two
facets of the three-dimensional figure.
Every unique inducing plane is the plane of \(\tau*B\).
If all groups were unique, \(x_0\) would satisfy every facet
inequality of \(P_0\). Thus exactly one facet \(\tau_0\)
has two inducers, and \(x_0\) lies strictly beyond \(\tau_0*B\).
The supporting forms used in Theorem~\ref{thm:first-size-cofaces}
give a distinct inducer of \(\tau_0\) for every containing
facet whose plane misses \(x_0\). There are at least
\(c-2=2\gamma-6\ge6\) such inducers, contradicting the fact that
\(\tau_0\) has exactly two.
Thus one minimal coface has two outside vertices and all others are
pyramids. Polarity gives the corresponding dual allocation, with one
double inducing group and all others unique.

Use the blocks \(X,Y,Z\) and extension vectors \(d_i,e_j\)
from the proof of Theorem~\ref{thm:seven-first-size-exclusion}.
The blocks \(X_J\) and \(Y_I\) have rank four each, and hence
exhaust the ambient rank eight. The exceptional corner again satisfies
\(Z_{i,G_j}=\xi_i+\zeta_j\).
The nonzero-vector argument needs only \(\alpha,\gamma\ge6\).
Indeed, a vanishing \(d_p\) gives a complementary three-face
containing all quotient ray representatives. If it also contains
the other exceptional vertex, the resulting join has too many figure
vertices. Otherwise the one-point join argument gives at least
\[
 c-\deg_B(p)\ge(2\gamma-4)-(\gamma-1)=\gamma-3\ge3
\]
inducers of the double-induced facet. Both alternatives are impossible.
The polar bound is \(\alpha-3\ge3\), so all four extension
vectors are nonzero.

The mixed-facet part of that proof now applies with the same hypotheses:
local base vertices have one outside neighbor; each color has at least
three vertices; three distinct facet-slack columns are independent;
and each unique inducing column has a zero in at least one exceptional
row. These facts force exactly three mixed facets
\(\tau_0,F_1,F_2\) with empty triple intersection, and the edge
\(F_1\cap F_2\) misses \(\tau_0\). Nonnegative
endpoint coefficients then give \(\xi_i>0\) on both rows and,
dually, \(\zeta_j>0\) on both columns.
All exceptional entries are positive, whereas the facet
\(L_{\tau_0}\) of the double-vertex coface requires one to vanish. This is a contradiction.
No step used \(\alpha=\gamma\) or connectedness.
\end{proof}

\subsection{The next boundary and its exceptional vertex figure}
\label{subsec:next-seven-boundary}

At size \(3h+13\), Lemma~\ref{lem:capacity-rigidity} and
Proposition~\ref{prop:simple-double-slack} force a single Hall-face
type. We first derive its count data and exceptional vertex figure.
The next subsection excludes the pair using only its simple member.

\begin{theorem}[Hall faces at the next connected size]
\label{thm:next-seven-size}
Let \(h\ge1\) be an integer, and let \(P\) be a seven-polytope
with connected nonincidence graph and \(\delta(P)\ge h\).
If either \(f_0(P)\le3h+13\) or \(f_6(P)\le3h+13\), then
\[
 f_0(P)=f_6(P)=3h+13,\qquad
 \delta(P)=h,\qquad \nu(\N(P))=2h+13.
\]
For every face attaining \(H(P)\), its pair with the dual face,
up to interchanging the two faces, has counts
\[
 (2h+6,h+5),\qquad(2h+7,h+6).
\]
The first face is simple. The second has exactly one vertex of degree
four, and every other vertex has degree three. Both face-count slacks
are one.
\end{theorem}

The proof first excludes a zero count slack using the capacity of a
four-simplex. The remaining integer budget has three alternatives;
Proposition~\ref{prop:simple-double-slack} removes one, and the other two
are exchanged by polarity.

\begin{proof}[Proof of Theorem~\ref{thm:next-seven-size}]
By polarity suppose \(n=f_0(P)\le3h+13\), and put \(M=f_6(P)\).
Theorem~\ref{thm:seven-first-size-exclusion} gives
\(n=3h+13\), \(M\ge3h+13\).
Choose any face \(\sigma\) attaining \(H(P)\).
Both \(\sigma\) and \(D=\sigma^\Delta\) are three-dimensional, as in
the proof of Proposition~\ref{prop:global-seven}.
With their vertex--facet counts \((a,\alpha)\), \((c,\gamma)\),
write
\[
 a=2h+4+x,\quad \alpha=h+4+y,\quad
 c=2h+4+z,\quad \gamma=h+4+w,
\]
and set \(u=n-a-\gamma\), \(v=M-c-\alpha\).
As in Theorem~\ref{thm:first-seven-size}, all these increments
and slacks are nonnegative integers satisfying
\[
 x-y\ge v,\quad z-w\ge u,\quad x\le2y,\quad z\le2w,
 \quad x\ge2v,\quad w\ge u,\quad x+w+u=5.
\]
Thus \(u+v\le2\). If \(u=v=0\),
Lemma~\ref{lem:face-equality-join} contradicts connectedness.
For a three-polytope \(S\), the degree excess
\(E(S)=2f_2(S)-f_0(S)-4\) bounds both its number of nonsimple
vertices and \(\iota(S)-3\). Here
\(E(\sigma)=2y-x\), \(E(D)=2w-z\).

\smallskip\noindent\emph{Excluding a zero slack.}
Suppose \(u=0\), \(v>0\). Then
\[
 \alpha-v-\iota(\sigma)\ge h+1+x-y-v\ge h+1>0,
 \qquad q:=2v+e(D)\le2v+E(D)\le2v+w\le x+w=5.
\]
If \(q\le4\), Lemma~\ref{lem:near-equality-join} gives a join.
For \(q=5\), every inequality in the latter chain is an equality:
\(x=2v\), \(w=5-2v\), \(z=w\), and \(e(D)=E(D)=w\).
Necessarily \(v=1\) or \(2\), so \(\gamma=h+4+w>5\).
Lemma~\ref{lem:capacity-rigidity}, with
\(t=3\) and \(\Theta_3(5)=5\), again gives a join.
Both conclusions contradict connectedness.

Suppose instead \(v=0\), \(u>0\). The polar bounds are
\[
 \gamma-u-\iota(D)\ge h+1+z-w-u\ge h+1>0,
 \qquad q':=2u+e(\sigma)\le2u+E(\sigma)
 \le2u+x\le u+w+x=5.
\]
Again \(q'\le4\) is excluded by Lemma~\ref{lem:near-equality-join}
in the polar polytope.
If \(q'=5\), equality gives \(w=u\), \(x=5-2u\),
\(y=x\), and \(z=2u\). Thus \(u=1\) or \(2\), and
\(\alpha=h+4+x>5\). The capacity lemma in the polar polytope
gives the same contradiction. Consequently \(u=v=1\).

\smallskip\noindent\emph{The remaining counts.}
Now \(x+w=4\), \(x\ge2\), \(w\ge1\).
For \((x,w)=(3,1)\), the inequalities force \((y,z)=(2,2)\).
For \((x,w)=(2,2)\), they give \(y=1\) and \(3\le z\le4\),
by \(z-w\ge u\) and \(z\le2w\).
The alternative \(z=4\) makes both faces simple, with
\(\alpha=h+5\ge6\), \(\gamma=h+6\ge7\).
Proposition~\ref{prop:simple-double-slack} excludes it.
Thus \(z=3\). The two surviving data sets are polar equivalents,
with \(M=3h+13\), \(H(P)=a+c=4h+13\), and degree excesses
zero and one. Excess one means one degree-four vertex and all other
degrees three. The face formula gives \(\delta(P)=h\) and
\(\nu(\N(P))=n-h=2h+13\). The maximizing face was arbitrary.
\end{proof}

\begin{proposition}[The exceptional prism figure]
\label{prop:boundary-prism}
Under the hypotheses of Theorem~\ref{thm:next-seven-size}, replace
\(P\) by \(P^*\) if necessary, so that the maximizing face \(\sigma\) is
simple and \(D=\sigma^\Delta\) has the unique degree-four vertex \(b\).
At \(D\) in the polar polytope, every outside vertex belongs to
a minimal coface. Exactly one such coface \(K\) contains two
outside vertices; all other cofaces are pyramids over \(D\).
The vertex figure of \(b\) in \(K\) is a triangular prism.
In particular, \(b\) is adjacent in \(K\) to both outside vertices.
\end{proposition}

We verify that allocation and the rank-eight corner still hold for
the almost-simple face. If the exceptional vertex had only one
outside neighbor, the sign argument at the end of the proof of
Theorem~\ref{thm:seven-first-size-exclusion} would apply on both sides
and make every exceptional corner entry positive, contrary to the
required incidence zero.

\begin{proof}[Proof of Proposition~\ref{prop:boundary-prism}]
Write \(a=2h+6\), \(\alpha=h+5\), \(c=2h+7\),
\(\gamma=h+6\).

\smallskip\noindent\emph{Allocation and normal ranks.}
The polar \(D^*\) has \(\gamma\ge7\) vertices, one
quadrilateral facet, and all other facets triangular.
It is not a pyramid: its base would have \(\gamma-1\ge6\)
vertices, exceeding its maximum facet size four.
Deleting any vertex of \(D^*\) therefore preserves affine dimension
three, so deleting any facet-slack column of \(D\) leaves rank four.
The same deletion property holds for simple \(\sigma\).
There are at least \(\alpha-1\ge5>3\) unique inducers at
\(\sigma\), and at least \(\gamma-1\ge6>4\) at \(D\).
Their intersections have no base vertex on either side.
As in the allocation part of Proposition~\ref{prop:simple-double-slack},
the restriction rank bounds their dimension by three, while their
quotient rays give the reverse bound. A join containing all outside
vertices would have too many figure vertices.
If a vertex belonged to no minimal coface, its nonvertex quotient
image would lie on at most two facets of the three-dimensional figure.
The support forms would then produce at least \(c-2\) inducers
at \(\sigma\), or at least \(a-2\) in the polar argument.
Both exceed two. Thus the two allocations hold, each with one double
coface and one double inducing group.

The normal slack blocks have ranks four and four. The exceptional
corner has the additive form \eqref{eq:first-size-corner}.
The nonzero-vector argument also survives. On the nonsimplicial
quotient side, the complementary three-face \(B\) has
\(\gamma\) vertices and \(c=2\gamma-5\) facets. At a vertex
of a three-polytope, the numbers of incident facets and edges agree,
by its polygonal link. Hence
\[
 c-\deg_B(p)\ge c-(\gamma-1)=\gamma-4=h+2\ge3.
\]
These would be distinct inducers of the same double-induced facet.
The opposite side retains the bound \(\alpha-3\ge3\).
Consequently all exceptional extension vectors are nonzero.

\smallskip\noindent\emph{The two vertex figures.}
Every facet of \(K\) other than \(D\) meets \(D\) in a facet,
with one local facet per base facet.
At \(b\), exactly five facets of \(K\) are incident: the base
and four side facets. Its three-dimensional vertex figure has five
facets and a quadrilateral base. Euler's formula and minimum degree
three bound its vertex count by six. Four vertices would give a
tetrahedron, so there are five or six.
With five vertices, the unique degree-four vertex is joined to the
other four, which form a cycle; the figure is a quadrilateral pyramid.
With six, the figure is simple. Its simplicial polar has five vertices
and graph \(K_5\) minus one edge, so it is a triangular bipyramid;
the figure is a triangular prism.

\smallskip\noindent\emph{Excluding the pyramid figure.}
Suppose that the figure is a quadrilateral pyramid. Then \(b\)
has exactly one edge to an outside vertex of \(K\), lying in all
four side facets through \(b\). Every other base vertex is simple
in \(K\) and also has one outside edge.
Color the vertices of \(D\) by these unique outside neighbors.
Each color has at least three vertices, as in the earlier proof.
The mixed-facet rank argument gives at most two unique mixed facets;
the edge-cycle argument gives exactly three mixed facets and exactly
three cross-color edges, one for each pair.
They cannot share a simple vertex, by the singleton-color argument.
Nor can they meet at \(b\): its four incident base facets have a
quadrilateral link, and no three are pairwise adjacent along edges.
Thus the edge of the two unique mixed facets misses the exceptional
base facet.

The two endpoint slack rows span the kernel of those two independent
columns. At either endpoint, including \(b\) if it is an endpoint,
an additional base facet misses the other endpoint. It is a unique
inducing column because the edge is disjoint from the exceptional
facet. Its inequality forces the corresponding endpoint coefficient
to be nonnegative. The two coefficients are not both zero, so the
exceptional coordinates are positive on this side.
On the simple side \(\sigma\), the endpoint-coefficient argument of
the proof of Theorem~\ref{thm:seven-first-size-exclusion} applies
verbatim.
Equation~\eqref{eq:first-size-corner} would make all exceptional
entries positive, contradicting the incidence zero required by the
facet \(L_{\tau_0}\) of the double-vertex coface.
The pyramid figure is impossible. The remaining prism figure has two
vertices outside its quadrilateral base, giving the two outside edges
at \(b\) asserted in the proposition.
\end{proof}

These theorems specify the local data forced at size \(3h+13\).
The single-side support obstruction below shows that the pair itself
cannot occur, without using its exceptional prism figure.

\subsection{A single-side support obstruction}
\label{subsec:single-side-obstruction}

The positive endpoint representation from the first exclusion has a
stronger consequence on the simple face itself. The edge used there
has endpoints lying in only four facets of the simple face. A further unique facet missing both gives
a contradiction to the local coface incidence, without requiring
positivity on the opposite Hall face.

\begin{lemma}[One-exception support obstruction]
\label{lem:single-side-support}
Let \(A\) be a simple three-polytope with \(\alpha\ge6\) facets
and vertex--facet slack matrix \(X\). Let \(H\) be a four-polytope
with facet \(A\) and exactly two vertices \(p,q\) outside \(A\).
Suppose that its other facets are \(L_\tau\), one for each facet
\(\tau\) of \(A\), with \(L_\tau\cap A=\tau\).
Fix a facet \(\tau_0\) of \(A\).
There do not exist \(d_p,d_q\in\operatorname{rowspan}X\), both
nonzero, such that, for every \(\tau\ne\tau_0\) and \(i\in\{p,q\}\),
\[
 d_i(\tau)\ge0,\qquad i\in L_\tau\ \Longrightarrow\ d_i(\tau)=0.
\]
\end{lemma}

The proof identifies the three mixed facets, represents both vectors
by the endpoint rows of the edge shared by the two mixed facets other
than \(\tau_0\), and then uses the six-facet threshold.

\begin{proof}[Proof of Lemma~\ref{lem:single-side-support}]
Each vertex of \(A\) lies in four facets of \(H\): \(A\)
and the three \(L_\tau\) through it. It is simple in \(H\),
with three edges in \(A\) and one to \(p\) or \(q\).
Color it by that outside neighbor, and call a facet of \(A\)
\emph{mixed} if its vertices have both colors. Both colors have at
least three vertices, since an outside vertex has degree at least four
and at most one neighbor outside \(A\).
Each \(L_\tau\) is three-dimensional and must contain \(p\)
or \(q\). Thus at least one of \(d_p(\tau),d_q(\tau)\) is
zero when \(\tau\ne\tau_0\).
For a mixed facet \(\tau\), the outside edges at its two colors
lie in \(L_\tau\), which therefore contains both outside vertices.

Put \(J=\{\tau:\tau\ne\tau_0\}\).
Deleting the \(\tau_0\) column of \(X\) leaves rank four,
by the nonpyramid argument for simple three-polytopes preceding
Theorem~\ref{thm:seven-first-size-exclusion}. Hence each nonzero
\(d_i\) has nonzero restriction to \(J\).
If three mixed facets lay in \(J\), their independent columns
would constrain \(d_p,d_q\) to the same line in the four-dimensional
row space. Their nonnegative, nonzero restrictions to \(J\) force
positive proportionality. Since one coordinate is zero at every
column of \(J\), both restrictions would vanish, a contradiction.
There are at most two mixed facets in \(J\).

There are at least three mixed facets altogether. Since the graph of
\(A\) is connected and both colors occur, some edge is cross-color, and
every cross-color edge lies in two mixed facets. With at most two
mixed facets, all cross-color edges would lie in their intersection,
so there would be exactly one. That edge would be a bridge, impossible because it lies on a
polygonal facet cycle. Thus there are exactly three mixed facets,
\(\tau_0,F_1,F_2\). Each mixed polygon has at least two cross-color
edges, each shared by two mixed facets. Since a pair of facets has
at most one common edge, there are exactly three cross-color edges,
one for each pair.

The three facets have empty common intersection. Otherwise a common
vertex \(v\) is simple, and their three pairwise intersection edges
are its three incident edges, all cross-color. The graph with \(v\)
removed is connected: its neighbors are joined by boundary paths
avoiding \(v\) in the three incident polygonal facets, and every
component contains a neighbor. All remaining vertices would have
the other color, contrary to both color classes having at least three
vertices. Therefore \(E=F_1\cap F_2\) is an edge disjoint from
\(\tau_0\). Let its endpoints be \(u_0,v_0\).

The two independent columns for \(F_1,F_2\) have a two-dimensional
common kernel in \(\operatorname{rowspan}X\), spanned by the
independent endpoint rows \(X_{u_0},X_{v_0}\). Write
\[
 d_i=\lambda_iX_{u_0}+\mu_iX_{v_0}\qquad(i=p,q).
\]
The third facet at \(u_0\) misses \(v_0\), and the third at
\(v_0\) misses \(u_0\). Both are in \(J\), so nonnegativity
there gives \(\lambda_i,\mu_i\ge0\). At least one is positive.
The two endpoints belong to three facets each and share precisely
two. Only four facets contain either endpoint. Since \(\alpha\ge6\),
at least two miss both. One is \(\tau_0\); another, \(\kappa\),
lies in \(J\). Both endpoint slacks at \(\kappa\) are positive,
and hence \(d_p(\kappa),d_q(\kappa)>0\). This contradicts the
required zero for the outside vertex of \(L_\kappa\).
\end{proof}

The six-facet threshold is necessary for this lemma. Take
\(H=\Delta_3\times[0,1]\) and its triangular-prism facet
\(A=\Delta_2\times[0,1]\), which has five facets.
Choose \(\tau_0\) as one rectangular side and \(E=\{x\}\times[0,1]\)
as the intersection of the other two. Let \(w\) be the vertex of
\(\Delta_3\) outside \(\Delta_2\), and put \(p=(w,0)\), \(q=(w,1)\),
\(u_0=(x,0)\), \(v_0=(x,1)\). The rows \(d_p=X_{u_0}\),
\(d_q=X_{v_0}\) are nonnegative and satisfy all six required zero
conditions on the four columns other than \(\tau_0\).
There is no second facet missing both endpoints.

\begin{proposition}[One simple face excludes unit count slacks]
\label{prop:one-simple-double-slack}
Let \(A\) be a three-face of a seven-polytope \(P\), and put
\(D=A^\Delta\). Write \((a,\alpha),(c,\gamma)\) for their
vertex--facet counts. If \(A\) is simple and
\(a-\alpha\ge2\), \(c-\gamma\ge2\), then
\[
 f_0(P)=a+\gamma+1,\qquad f_6(P)=c+\alpha+1
\]
cannot both hold.
\end{proposition}

No connectedness, Hall-maximization, balance, or simplicity of \(D\)
is assumed.

\begin{proof}[Proof of Proposition~\ref{prop:one-simple-double-slack}]
Suppose both equalities hold. Simplicity gives \(a=2\alpha-4\),
so \(\alpha\ge6\). Euler's inequality \(c\le2\gamma-4\)
and \(c\ge\gamma+2\) give \(\gamma\ge6\).
The degree excess of \(D\) is \(2\gamma-c-4\), and hence
\[
 \iota(D)\le3+(2\gamma-c-4)\le\gamma-3.
\]
Its polar is not a pyramid: a base facet would give a vertex of \(D\)
in \(\gamma-1>\iota(D)\) facets. Deleting any vertex of that polar preserves dimension three.
Thus deleting one facet-slack column of \(D\) leaves rank four,
as it does for simple \(A\).

\smallskip\noindent\emph{The two allocations.}
There are \(\gamma\) minimal cofaces and \(\gamma+1\) outside
vertices at \(A\). Suppose that one vertex \(x_0\) belongs to
none of these cofaces. All cofaces are pyramids \(A*v_i\).
At least \(\alpha-1\) facets have unique inducers; their
intersection \(B\) contains all \(v_i\), no vertex of \(A\),
and possibly \(x_0\). Restriction rank four gives \(\dim B\le3\),
and the four-dimensional quotient span gives the reverse inequality.
If \(x_0\in B\), the full join \(P=A*B\) has a face figure
with \(\gamma+1\) vertices, contrary to the actual \(\gamma\) rays.
Otherwise \(P_0=A*B\) contains all vertices except \(x_0\).
Its quotient cone at \(A\) is unchanged. All \(c\) containing-facet
planes remain supporting, and the nonvertex quotient image of \(x_0\)
lies on at most two facets of the three-dimensional figure.

Every unique inducer has the plane of \(\tau*B\).
If all groups were unique, \(x_0\) would satisfy every facet inequality
of \(P_0\). Thus one \(\tau_0\) has two inducers and the others
are unique; \(x_0\) lies strictly beyond \(\tau_0*B\).
For supporting forms \(g_0\le0\) there and \(h_K\le0\) at
\(A*K\), set \(\lambda=g_0(x_0)>0\) and, for each facet \(K\) of
\(B\) such that the hyperplane of \(A*K\) misses \(x_0\),
\(\mu_K=h_K(x_0)<0\).
The form
\[
 \lambda h_K-\mu_Kg_0\le0
\]
supports a facet \(\conv(x_0,\tau_0*K)\) inducing \(\tau_0\).
Distinct such \(K\) give distinct facets, and there are at least
\(c-2\ge6\) of them, contrary to multiplicity two.
Hence one minimal coface has two outside vertices \(p,q\), and
the other \(\gamma-1\) cofaces are pyramids.
The polar argument uses \(\gamma-1>\iota(D)\), its deletion
rank four, and \(a-2\ge6\), and gives the dual allocation.
In particular every noncontaining ambient facet induces a facet of
\(A\); exactly one group \(\tau_0\) is double.

\smallskip\noindent\emph{The exceptional rows.}
Let \(I\) index the single-apex vertices, \(J\) the unique
inducers, and \(C\) the facets containing \(A\).
Use the normalized blocks \(X,Y,Z\) from
Theorem~\ref{thm:seven-first-size-exclusion}.
Positive scaling equalizes the duplicated \(\tau_0\) columns
of \(X\) and the two quotient-ray rows of \(Y\).
The submatrix on rows \(\V(A),I\) and columns \(C,J\) is
\[
 \begin{pmatrix}0&X_J\\Y_I&0\end{pmatrix},
\]
of rank eight, the full slack rank. Thus \(Z_{p,J},Z_{q,J}\)
extend uniquely to \(d_p,d_q\in\operatorname{rowspan}X\).
On \(J\) they are nonnegative and retain the ambient incidence zeros.

If \(d_p=0\), the face \(B=\bigcap_{\tau\in J}G_\tau\)
contains \(I,p\), possibly \(q\), and no vertex of \(A\).
Its dimension is three by the same rank and quotient bounds.
If \(q\in B\), the full join again has too many figure vertices.
Otherwise \(P_0=A*B\) contains all vertices except \(q\), and
\(B\) is projectively the original quotient figure, with
\(\gamma\) vertices and \(c\) facets.
The point \(q\) has the same quotient ray as \(p\) and lies beyond
only \(\tau_0*B\). The support-form construction gives a distinct
\(\tau_0\) inducer for every facet of \(B\) missing \(p\).
The number of incident facets at a three-polytope vertex equals its
edge degree, which is at most \(\gamma-1\). At least
\(c-(\gamma-1)\ge3\) inducers are required, a contradiction.
Therefore \(d_p\ne0\); likewise \(d_q\ne0\).

\smallskip\noindent\emph{Applying the support obstruction.}
Let \(H=\conv(A,p,q)\) be the double coface.
Every facet \(G\ne A\) of \(H\) has the form \(G=H\cap F\) for an
ambient facet \(F\not\supseteq H\). Such an \(F\) cannot contain \(A\):
it contains \(p\) or \(q\), and a face containing \(A\) and an outside
vertex of \(H\) contains the minimal coface \(H\). By the dual
allocation, \(G\cap A=F\cap A\) is a facet of \(A\).
Each ridge \(\tau\) lies in exactly two facets of \(H\), one
being \(A\), so there is one \(L_\tau\) per \(\tau\).
For a unique \(\tau\), its ambient inducer contains \(L_\tau\).
A local incidence with \(p\) or \(q\) is therefore the corresponding
zero of \(d_i(\tau)\).
Remove one duplicated column of \(X\) to obtain its ordinary
facet-slack matrix. The two nonzero extensions satisfy
Lemma~\ref{lem:single-side-support}, which gives the contradiction.
\end{proof}

\begin{theorem}[Stronger connected seven-dimensional bound]
\label{thm:seven-second-size-exclusion}
For every integer \(h\ge1\), a seven-polytope \(P\) with connected
nonincidence graph and \(\delta(P)\ge h\) satisfies
\[
 f_0(P)\ge3h+14,\qquad f_6(P)\ge3h+14.
\]
\end{theorem}

\begin{proof}[Proof of Theorem~\ref{thm:seven-second-size-exclusion}]
If either count were at most \(3h+13\),
Theorem~\ref{thm:next-seven-size} would give a maximizing pair,
up to polarity, with a simple face \((2h+6,h+5)\) and a face
\((2h+7,h+6)\), both count slacks one.
Both vertex--facet excesses equal \(h+1\ge2\).
Proposition~\ref{prop:one-simple-double-slack} excludes the pair.
The integer counts must therefore be at least \(3h+14\).
\end{proof}

For deficiency one this gives seventeen vertices and seventeen facets
as necessary bounds. Section~\ref{sec:shared-coordinate} attains both
bounds simultaneously. Before constructing the attainer, we refine
the equality data using the shadows of complementary four-faces.

\subsection{Six-facet shadows and the sharp lower boundary}
\label{subsec:shadow-boundary}

At the next count budget a complementary four-face can survive the
capacity argument. Its three-dimensional shadow supplies an additional
restriction, which uses that the shadow is vertex-bijective: general
polytope extensions may have hidden vertices. Pashkovich and Weltge studied this
phenomenon for minimum extensions of general-position heptagons
\cite[Section~2.2]{pashkovichWeltge2015}. We use the classical few-facet
classification with a different, three-dimensional shadow constraint.

\begin{lemma}[Three-dimensional shadows of six-facet polytopes]
\label{lem:six-facet-shadow}
Let \(B\) have dimension at least four, at most six facets and at least
seven vertices. Let \(L:\aff B\to\R^3\) and
\(\ell:\aff B\to\R\) be affine maps with \(\ell>0\) on \(B\).
Suppose that \(\pi(x)=L(x)/\ell(x)\) maps the vertices of \(B\)
bijectively onto the vertices of a three-polytope \(T\). Then
\[
 f_0(T)\in\{7,8,9\},\qquad f_2(T)\le9.
\]
\end{lemma}

\begin{proof}
Each edge of \(T\) lifts to an edge of \(B\). Indeed, pull back a
supporting affine form whose zero face is that edge, and multiply its
pullback by \(\ell\). The resulting affine supporting
form on \(B\) has exactly the two preimage vertices in its zero face.
That face is their segment. The graph of \(T\) is therefore a spanning
subgraph of the graph of \(B\).

Six facets imply \(\dim B\le5\). Dimension five would give a simplex
with six vertices, and a four-polytope with five facets also has only
five vertices. Thus \(B\) is four-dimensional with six facets.
The classical classification \cite{mcmullen1971}, in the form of
\cite[Lemma~1]{pinedaTritamaWangYost2026}, leaves the following types
with at least seven vertices:
\[
 \operatorname{Pyr}(\Delta_1\times\Delta_2),\qquad
 \Delta_1\times\Delta_3,\qquad \Delta_2\times\Delta_2.
\]
They have respectively \((7,15),(8,16),(9,18)\) vertices and edges.
The other six-facet type, \(\operatorname{Pyr}^2(\Delta_1\times\Delta_1)\),
has six vertices.

The first graph is a cone over the triangular-prism graph. Contracting
one triangular layer leaves a \(K_5\) minor, formed by the other
triangle, the contracted vertex and the apex. In the second graph,
contracting one of its two \(K_4\) layers also leaves a \(K_5\) minor.
These graphs are nonplanar: contractions preserve planarity, whereas
\(K_5\) violates the planar edge bound. A three-polytope shadow must
omit an edge, and therefore has at most fourteen or fifteen edges,
respectively.

The third graph is \(K_3\mathbin{\square}K_3\). It has exactly six
triangles, its three rows and three columns; each edge lies in exactly
one of them. If the graph of \(T\) had eighteen edges, Euler's
formula would give eleven faces. If \(q\) of them were triangular,
the boundary-length sum would give
\(36\ge3q+4(11-q)\), so \(q\ge8\), impossible.
With seventeen edges it would have ten faces and
\(34\ge3q+4(10-q)\), so \(q\ge6\). The omitted edge destroys
one of the six triangles, again impossible. Hence the shadow has at
most sixteen edges. In all three cases Euler's formula now gives
\(f_2(T)=f_1(T)-f_0(T)+2\le9\).
\end{proof}

\begin{theorem}[Preliminary Hall-pair types at the connected lower boundary]
\label{thm:seven-second-boundary}
Let \(h\ge1\) be an integer, and let \(P\) be a seven-polytope
with connected nonincidence graph and \(\delta(P)\ge h\).
If either class size is at most \(3h+14\), then
\[
 f_0(P)=f_6(P)=3h+14,\qquad
 \delta(P)=h,\qquad \nu(\N(P))=2h+14.
\]
Every maximizing Hall face \(A\) and its dual \(D=A^\Delta\)
are three-dimensional. Write \((a,\alpha),(c,\gamma)\) for their
vertex--facet counts and
\[
 u=f_0(P)-a-\gamma,\qquad v=f_6(P)-c-\alpha.
\]
Up to polarity, the pair has one of the following types:
\begin{enumerate}
\item Both faces have counts \((2h+7,h+6)\), each with one
degree-four vertex and all other degrees three; \((u,v)=(1,1)\).
\item Both faces are simple, with counts \((2h+6,h+5)\) and
\((2h+8,h+6)\), respectively; \((u,v)=(2,1)\).
\item One of the three exceptional rows in Table~\ref{tab:exceptional-hall}, with \(u=0\):
\begin{datadisplay}{Exceptional zero-slack Hall-pair data}\label{tab:exceptional-hall}
\[
\begin{array}{c|c|c|c|c}
 h&(a,\alpha)&(c,\gamma)&v&B\\\hline
 1&(10,7)&(8,7)&2&\operatorname{Pyr}(\Delta_1\times\Delta_2)\\
 1&(12,8)&(6,5)&3&\Delta_4\\
 2&(14,9)&(8,6)&3&\operatorname{Pyr}^2(\Delta_1\times\Delta_1)
\end{array}
\]
\end{datadisplay}
Here \(A\) is simple, and \(B\) is the four-face on the vertices
outside \(A\). In the first row \(D\) has two degree-four vertices
and all other degrees three; in the other two rows it is simple.
\end{enumerate}
In particular, for \(h\ge3\) only the first two types are possible.
\end{theorem}

This reduction supplies necessary data for every maximizing pair.
Corollary~\ref{cor:zero-boundary-exclusion} below eliminates two of the
three exceptional rows. Section~\ref{sec:shared-coordinate} realizes
the four-simplex row at \(h=1\) and determines its extremal rigidity.

\begin{proof}[Proof of Theorem~\ref{thm:seven-second-boundary}]
By polarity put \(n=f_0(P)\le3h+14\), \(M=f_6(P)\).
Theorem~\ref{thm:seven-second-size-exclusion} gives
\(n=3h+14\), \(M\ge n\).
Both members of a maximizing pair are three-dimensional by the
violating-face argument of Proposition~\ref{prop:global-seven}.
For any maximizing pair, introduce the count increments
\[
 a=2h+4+x,\quad\alpha=h+4+y,\quad
 c=2h+4+z,\quad\gamma=h+4+w.
\]
The face inequalities and Euler's bound give nonnegative integers with
\[
 x-y\ge v,\quad z-w\ge u,\quad x\le2y,\quad z\le2w,
 \quad x\ge2v,\quad w\ge u,\quad x+w+u=6.
\]
Consequently \(u+v\le3\). If \(u=v=0\),
Lemma~\ref{lem:face-equality-join} makes \(P\) a free join, contrary
to connectedness.

\smallskip\noindent\emph{Zero slack and the shadow.}
Suppose \(u=0,v>0\). The incidence and excess estimates give
\[
 \alpha-v-\iota(A)\ge h+1>0,\qquad
 q:=2v+e(D)\le2v+(2w-z)\le2v+w\le6.
\]
For \(q\le4\), Lemma~\ref{lem:near-equality-join} forces a join.
For \(q=5\), one has \(\gamma>5\): otherwise \(h=1,w=z=0\),
and \(q=2v\) is even. Since \(\Theta_3(5)=5\),
Lemma~\ref{lem:capacity-rigidity} again forces a join.
Connectedness excludes both possibilities. Thus \(q=6\), and all
inequalities in the upper chain are equalities. They give
\[
 x=2v,\quad y=v,\quad w=z=6-2v,\quad
 e(D)=2\gamma-c-4=w,\qquad v\in\{1,2,3\}.
\]
Hence \(A\) is simple, \(c-\gamma=h\), and \(M=n\).
Every nonsimple vertex of \(D\) has degree four.

Take the intersection \(B\) of the unique inducers of facets of \(A\),
as in the capacity-rigidity proof. The zero vertex slack assigns all
\(\gamma\) outside vertices to distinct pyramid cofaces.
Every unique inducer contains them, and at least
\(\alpha-v=h+4>3=\iota(A)\) such inducers exclude all vertices
of \(A\). Thus \(B\) is precisely the face on the outside vertices.
Its quotient-ray span gives \(\dim B\ge3\); equality would make
\(P=A*B\), contradicting connectedness.
The facet budget from that proof gives \(f_3(B)\le q=6\).
Moreover, any three distinct facet-slack columns of \(A\) are
independent: after polar normalization they represent three distinct
vertices, which cannot be collinear. Three unique inducers therefore
have independent ambient homogeneous forms, so \(\dim B\le7-3=4\).
Consequently \(B\) is four-dimensional.

The quotient map normalized by a positive cone section maps its vertices
bijectively onto those of the figure \(T\) at \(A\): each outside
vertex represents its own extreme ray. Its denominator is positive
on their convex hull. The figure has \(\gamma\) vertices and \(c\)
facets. For \(h\ge3\), one has \(\gamma\ge h+4\ge7\) and
\(c=\gamma+h\ge10\), contrary to Lemma~\ref{lem:six-facet-shadow}.
For \(h\in\{1,2\}\) and \(v=1\), one has \(\gamma=h+8\).
At \(h=1\), the shadow has nine vertices and ten facets, impossible;
at \(h=2\), \(\gamma=10>9=\Theta_3(6)\) excludes it.
For \(v=2\), \(\gamma=h+6\), and \(h=2\) gives \(c=10>9\),
excluded by the shadow bound. Only \(h=1\) remains.
For \(v=3\), only \(h=1,2\) remain. Substitution gives exactly
the three exceptional rows. Their four-face types follow from the
same classical classification; five vertices force a four-simplex.
If \(v=0,u>0\), the polar chain reads
\(\gamma-u-\iota(D)\ge h+1\) and
\(q':=2u+e(A)\le2u+(2y-x)\le2u+x\le u+w+x=6\), since \(y\le x\).
The cases \(q'\le5\) are excluded as before: if \(q'=5\) and
\(\alpha\le5\), then \(h=1\), \(x=y=0\) and \(q'=2u\) is even.
Equality at \(q'=6\) gives \(y=x=6-2u\), \(w=u\) and \(z=2u\), so
\(D\) is simple and \(M=n\). The shadow argument at \(D\) in the polar
polytope now applies verbatim.

\smallskip\noindent\emph{Positive slacks.}
Now \((u,v)\) is \((1,1),(1,2)\), or \((2,1)\).
If \(u=v=1\), then \(x+w=5\), \(x\ge2\), \(w\ge1\).
For \((x,w)=(2,3)\), \(y=1\) and \(A\) is simple.
The slack inequalities give both vertex--facet excesses at least
\(h+1\ge2\), so Proposition~\ref{prop:one-simple-double-slack}
excludes this case. For \((x,w)=(4,1)\), \(z=2\) and \(D\)
is simple; polarity excludes it in the same way.
For \((x,w)=(3,2)\), \(y=2\), and \(z-w\ge u\), \(z\le2w\) give
\(3\le z\le4\). The value \(z=4\) makes \(D\) simple and is
again excluded. Hence \(z=3\), both faces have counts
\((2h+7,h+6)\), and both degree excesses are one.

For \((u,v)=(1,2)\), the budget forces
\((x,w,y,z)=(4,1,2,2)\); both faces are simple with counts
\((2h+8,h+6)\), \((2h+6,h+5)\).
For \((u,v)=(2,1)\), it forces \((2,2,1,4)\), the polar data.
Every remaining positive or zero row gives \(M=n=3h+14\) and
\(H(P)=a+c=4h+14\).
The face formula gives \(\delta(P)=h\) and
\(\nu(\N(P))=2h+14\). The maximizing pair was arbitrary.
\end{proof}

\subsection{Two-face facet accounting}
\label{subsec:two-face-budget}

The zero-slack alternatives have an additional geometric constraint:
the outside vertices form a face. We count all ambient facets when two
disjoint faces exhaust the vertex set. This is a projected-join setting,
related to the classical subdirect-sum construction
\cite[Section~15.1.3]{henkRichterGebertZiegler2017}. The count below
will eliminate two of the preliminary alternatives.

\begin{lemma}[Shared-line normal form]
\label{lem:two-face-normal}
Let \(P\) be an \((r+k)\)-polytope, where \(r,k\ge2\), and let
\(A,B\) be disjoint faces of dimensions \(r,k\) whose vertex sets
partition \(\V(P)\). Then \(P\) has a projectively equivalent bounded
realization in which \(\aff A\) and \(\aff B\) are disjoint and their
direction spaces meet in a line \(L\). It admits coordinates \((x,y,z,t)\), with
\(x\in\R^{r-1}\), \(y\in\R^{k-1}\), such that
\[
 A\subset\{y=0,t=0\},\qquad B\subset\{x=0,t=1\},
\]
and the shared direction is the \(z\)-direction.
\end{lemma}

\begin{proof}
Write \(U=\aff(A)\), \(V=\aff(B)\). Their combined affine span
is all of \(\aff(P)\). If they are disjoint, the affine dimension
formula gives \(\dim(\operatorname{dir}U\cap\operatorname{dir}V)=1\).
If they meet, their intersection is a single point \(o\) and their
direction spaces intersect trivially. Since
\(P\cap\aff(A)=A\), \(P\cap\aff(B)=B\), and \(A\cap B=\varnothing\),
the point \(o\) lies outside \(P\). Strict separation gives a
hyperplane through \(o\) disjoint from the compact polytope \(P\).
Send this hyperplane to infinity. The projective denominator has constant
sign and is bounded away from zero on \(P\), so its image is a bounded
polytope with the same face lattice. The transformed affine hulls meet
only at their common infinite point; their direction intersection is
therefore a line \(L\).

In either case the sum of the direction spaces has codimension one.
Choose an affine coordinate \(t\) whose linear part vanishes on
\(\operatorname{dir}U+\operatorname{dir}V\), with \(t=0\) on \(U\) and
\(t=1\) on \(V\). Write \(\operatorname{dir}U=L\oplus X\) and
\(\operatorname{dir}V=L\oplus Y\), take \(z\) along \(L\), and let
\(x,y\) be coordinates on \(X,Y\). Translating the origin inside \(U\)
along \(X\) makes the \(x\)-coordinate vanish on \(V\), which gives
the stated form.
\end{proof}

\begin{proposition}[Complete facet count and allocation budget]
\label{prop:two-face-budget}
In the normal form of Lemma~\ref{lem:two-face-normal}, let \(R,T\)
be the projections of \(A,B\) along \(L\). Put
\[
 p=f_{r-2}(R),\quad q=f_{k-2}(T),\quad
 \alpha=f_{r-1}(A),\quad\beta=f_{k-1}(B).
\]
Partition the facets of each face by the sign of the \(z\)-coefficient
in an outward supporting form, with counts
\(\alpha_+,\alpha_-,\alpha_0\) and
\(\beta_+,\beta_-,\beta_0\). Then
\begin{equation}\label{eq:two-face-count}
 f_{r+k-1}(P)=p+q+\alpha_+\beta_++\alpha_-\beta_-.
\end{equation}
Moreover \(p\ge\alpha_0\), and the allocation excess satisfies
\begin{align}
 v&:=f_{r+k-1}(P)-q-\alpha\notag\\
  &=(p-\alpha_0)+(\beta_+-1)\alpha_+
                       +(\beta_--1)\alpha_-
    \ge\beta-\beta_0-2.\label{eq:two-face-budget}
\end{align}
If \(B\to T\) maps vertices bijectively onto vertices and \(T\)
is simplicial, then \(\beta_0=0\), hence \(v\ge\beta-2\).
\end{proposition}

\begin{proof}
Every facet of \(P\) meets both vertex classes: otherwise its dimension
is at most \(r\) or \(k\), strictly smaller than \(r+k-1\).
A facet containing \(A\) has a supporting form vanishing on
\(\operatorname{dir}A\), including \(L\). Its intersection with \(B\)
projects to a face of \(T\); the dimension count below shows that it
is a facet. Conversely, a facet inequality \(\ell(y)\le\ell_0\) of
\(T\) extends to \(\ell(y)-\ell_0t\le0\) on \(P\).
The lifted face of \(B\) has dimension \(k-2\) without \(L\) in its
direction space, or dimension \(k-1\) with \(L\). Joining it to \(A\)
has dimension \(r+k-1\) in both cases. These are exactly the \(q\)
facets containing \(A\). Interchanging the two faces gives exactly
\(p\) facets containing \(B\). No facet contains both entire faces.

For a facet containing neither face, its intersections with \(A,B\)
must both be facets, and their direction intersection must be zero:
their dimensions are at most \(r-1,k-1\), and their convex hull has
dimension \(r+k-1\). The global supporting form restricts to the same
coefficient on \(L\) in both faces. If it were zero, both facet
direction spaces would contain \(L\), reducing the hull dimension to
at most \(r+k-2\). Thus the two coefficients are nonzero with the same
sign.

Conversely, scale any same-sign pair of facet forms positively so that
both \(z\)-coefficients are \(+1\) or both are \(-1\). Combine their
\(x\)- and \(y\)-coefficients, and choose the constant and
\(t\)-coefficient to reproduce their two offsets. This global form is
valid on \(A\) and \(B\). Its two zero faces have direction intersection
zero, since neither contains \(L\), and their hull has dimension
\(r+k-1\). Distinct pairs give distinct intersections with \(A,B\).
This exhausts the facets and proves \eqref{eq:two-face-count}.

Every zero-type facet of \(A\) projects to a distinct facet of \(R\),
so \(p\ge\alpha_0\). Boundedness implies
\(\alpha_+,\alpha_-,\beta_+,\beta_-\ge1\): inequalities with only
one sign cannot bound a fiber parallel to \(L\) in both directions.
Subtracting \(q+\alpha\) from \eqref{eq:two-face-count} now gives
\eqref{eq:two-face-budget}.

Finally, a zero-type facet of \(B\) projects to a facet of \(T\)
with the same vertex set under a vertex-bijective projection. It has
at least \(k\) vertices, whereas a facet of the simplicial
\((k-1)\)-polytope \(T\) has exactly \(k-1\). Therefore no such
facet exists, proving the last assertion.
\end{proof}

\begin{corollary}[Elimination of the zero-slack boundary alternatives]
\label{cor:zero-boundary-exclusion}
Under the hypotheses of Theorem~\ref{thm:seven-second-boundary}, the
only zero-slack alternative, up to polarity, is
\[
 h=1,\quad (a,\alpha)=(12,8),\quad(c,\gamma)=(6,5),
 \quad(u,v)=(0,3),\quad B\cong\Delta_4.
\]
For every \(h\ge2\), only the two positive-slack types of
Theorem~\ref{thm:seven-second-boundary} remain.
Consequently, if \(P\) is a seven-polytope with connected
nonincidence graph, \(\delta(P)\ge h\ge2\), and some maximizing Hall
pair has a zero count slack, then \(f_0(P),f_6(P)\ge3h+15\).
\end{corollary}

\begin{proof}
In a zero-slack row, \(A\) and the outside four-face \(B\) exhaust
the vertices, and the figure \(T\) at \(A\) is a vertex-bijective
three-dimensional shadow of \(B\). The projective normalization preserves
this quotient description. Apply Proposition~\ref{prop:two-face-budget}
with \(r=3,k=4\).

For the \(h=2\) row, \(B\) has six facets and \(T\) has six vertices
and eight facets. Euler's formula gives twelve edges, and the boundary
length sum is twenty-four. All eight facets are triangles. Thus
\(\beta_0=0\), and \(v\ge6-2=4\), contradicting \(v=3\).

For the remaining \(h=1\) prism-pyramid row, \(B\) has seven vertices,
six facets and fifteen edges. Its facets have vertex counts
\(6,4,4,5,5,5\); the two tetrahedral facets share only the apex and
no edges. The dual counts and degrees of \(D\) show that \(T\) has
two quadrilateral facets and six triangular facets. It has thirteen
edges. Since \(v=2\), \eqref{eq:two-face-budget} forces
\(\beta_0\ge2\). Only the two tetrahedral facets can project to
facets of \(T\) with the same vertex set, so both project to
quadrilaterals. Each loses its two diagonal edges from the shadow graph.
Those four edges are distinct, and every edge of \(T\) lifts to an
edge of \(B\), as in Lemma~\ref{lem:six-facet-shadow}. Hence
\(f_1(T)\le15-4=11\), contrary to thirteen.

The simplex row is the sole zero-slack survivor. The last assertion
follows from Theorem~\ref{thm:seven-second-size-exclusion}: a class
size below \(3h+15\) would equal \(3h+14\), where the preceding
exclusions rule out a zero-slack maximizing pair. Polarity gives the
bound for both classes.
\end{proof}

\subsection{The exceptional quadrilateral and the final boundary data}
\label{subsec:almost-simple}

The remaining unit-slack type has a single degree-four vertex on each
Hall face. Coloring vertices by an outside neighbor is classical in the
facet-plus-edge setting: see the Shephard-property argument in
\cite[Lemma~20]{pinedaUgonYost2022}. Here the additional row-space
conditions give a support obstruction. Two mixed facets may meet only
at the degree-four vertex; in the polar they then span a diagonal of
the quadrilateral facet, and that case needs the full incidence data
of the coface.

\begin{lemma}[One exceptional vertex in the support obstruction]
\label{lem:almost-simple-support}
Let \(A\) be a three-polytope with \(\alpha\ge7\) facets, all
vertices of degree three except at most one of degree four, and let
\(X\) be its vertex--facet slack matrix. Let \(H\) be a four-polytope
with facet \(A\) and exactly two outside vertices \(p,q\). Suppose
its other facets are \(L_\tau\), one per facet \(\tau\) of \(A\),
with \(L_\tau\cap A=\tau\). For a fixed facet \(\tau_0\), there
are no nonzero \(d_p,d_q\in\operatorname{rowspan}X\) such that,
for every \(\tau\ne\tau_0\) and \(i\in\{p,q\}\),
\[
 d_i(\tau)\ge0,\qquad i\in L_\tau\ \Longrightarrow\ d_i(\tau)=0.
\]
\end{lemma}

The proof reduces the mixed facets to a triangle. Genuine edges give
positive endpoint representations; the possible diagonal instead forces
both vectors onto the exceptional vertex row.

\begin{proof}[Proof of Lemma~\ref{lem:almost-simple-support}]
If \(A\) is simple, Lemma~\ref{lem:single-side-support} applies.
Suppose therefore that \(v_*\) is its unique degree-four vertex.

\smallskip\noindent\emph{Ranks and the mixed facets.}
The polar \(S=A^*\) has exactly one quadrilateral facet \(Q\) and
all other facets triangular. Deleting one vertex leaves dimension three:
otherwise \(S\) is a pyramid with a base on \(\alpha-1\ge6\)
vertices, incompatible with this facet list. Thus deleting the
\(\tau_0\) column leaves \(X_J\) of rank four, where
\(J=\F(A)\setminus\{\tau_0\}\). Every nonzero row-space vector
has nonzero restriction to \(J\). Any three different columns are
independent, since their positive normalizations represent distinct
polar vertices, which cannot be collinear.

Here, unlike in Lemma~\ref{lem:single-side-support}, call \(\tau\)
mixed if \(L_\tau\) contains both outside vertices.
Every \(L_\tau\) contains at least one, since its intersection with
\(A\) is only two-dimensional. Hence at least one of
\(d_p(\tau),d_q(\tau)\) vanishes at every \(\tau\in J\).
Three mixed facets in \(J\) would constrain both vectors to one
common line. Their nonnegative, nonzero restrictions would make them
positively proportional, and the zero at each coordinate would then
make both restrictions vanish. Thus there are at most two mixed facets
in \(J\), and at most three altogether.

Every ordinary vertex of \(A\) lies in exactly four facets of \(H\).
It is simple in \(H\), with exactly one outside edge; color it by
its neighbor \(p\) or \(q\). Each outside vertex has degree at least
four, with at most one outside neighbor and at most the neighbor
\(v_*\), so each color occurs at least twice among ordinary vertices.
The graph of \(A\) with \(v_*\) deleted is connected: facet boundary
paths avoiding \(v_*\) join consecutive neighbors around its cyclic
vertex figure, and every component contains a neighbor. It has no
bridges. An edge avoiding \(v_*\) belongs to two facets, at least one
of which misses \(v_*\); otherwise their intersection edge would
contain \(v_*\). That facet cycle bypasses the edge after deletion.

Every cross-color ordinary edge belongs to two mixed facets: the unique
outside edge at a simple vertex of \(H\) lies in all three facets
\(L_\tau\) through it. With at most two mixed facets, every
cross-color edge would be their single common edge, a bridge in the
connected ordinary graph. This is impossible. There are therefore
exactly three mixed facets, one of them \(\tau_0\).

\smallskip\noindent\emph{The polar mixed triangle.}
In \(H^*\), write \(o\) for the vertex dual to \(A\),
\(z_\tau\) for that dual to \(L_\tau\), and \(F_p,F_q\)
for the facets dual to the outside vertices. Their intersection is
the triangle \(C\) on the three mixed vertices. As a ridge, it lies
in exactly \(F_p,F_q\). Each edge \([z_\tau,z_\sigma]\) lies
in at least one further facet \(G_x\), dual to an \(A\)-vertex
\(x\). Thus \(\tau\cap\sigma\ne\varnothing\). If this
intersection is not an edge, it is \(v_*\), and the two facets are
opposite among the four through \(v_*\). A triple of mixed facets
can contain at most one such opposite pair.

Suppose it contains one. The corresponding edge belongs to exactly
\(F_p,F_q,G_{v_*}\), so its two incident faces in \(G_{v_*}\)
are \(F_p\cap G_{v_*}\) and \(F_q\cap G_{v_*}\).
The facet \(G_{v_*}\) has five vertices, \(o\) and the four
neighbors corresponding to the facets through \(v_*\), and \(o\)
has degree four within it. Euler's formula and planarity give eight
or nine edges. With eight, the other four vertices form a cycle and
\(G_{v_*}\) is a quadrilateral pyramid with apex \(o\); with
nine, its graph is \(K_5\) minus an edge and it is a triangular
bipyramid with \(o\) equatorial. Both indicated incident faces miss
\(o\), whereas a quadrilateral pyramid has only one such facet.
Thus the bipyramid case is forced. Write \(z_1,z_2\) for its other
equatorial vertices and \(z_U,z_V\) for its apices. After relabeling,
\[
 F_p\cap G_{v_*}=\conv(z_1,z_2,z_U),\qquad
 F_q\cap G_{v_*}=\conv(z_1,z_2,z_V).
\]
The third mixed vertex \(z_0\) is not a vertex of \(G_{v_*}\):
otherwise \(z_0\in C\subseteq F_p\cap F_q\) would lie in both
displayed triangles, whose only common vertices are \(z_1,z_2\). The pair \(1,2\) corresponds to a diagonal of \(Q\),
not an \(A\)-edge; the other two pairs give genuine edges.

If no opposite pair occurs, all three pairs give \(A\)-edges, and
the three facets have empty common intersection. A common ordinary
vertex would give a tetrahedral \(G_x\) containing \(C\), contrary
to ridge incidence. A common \(v_*\) would put three facets among
the four around it, including an opposite pair.

\smallskip\noindent\emph{Positive endpoint representations.}
In the no-diagonal case, the two mixed facets in \(J\) meet in an
edge \(E\) disjoint from \(\tau_0\). Its independent endpoint
rows \(X_u,X_w\) span the two-dimensional kernel of the two columns.
Write \(d_i=\lambda_iX_u+\mu_iX_w\). Each endpoint has a facet missing the
other; these facets lie in \(J\), since \(\tau_0\) misses \(E\).
Nonnegativity gives \(\lambda_i,\mu_i\ge0\), with one positive.
The endpoints belong together to at most five facets: at most one has
degree four, and they share exactly two facets. Since \(\alpha\ge7\),
at least two facets miss both. Besides \(\tau_0\), one lies in
\(J\), where both \(d_i\) are positive, contradicting the required zero.

If a diagonal occurs and \(\tau_0\) is one of its endpoints,
the other two mixed facets meet in a genuine edge. The third mixed
facet misses \(v_*\), so the edge has ordinary endpoints and is
disjoint from \(\tau_0\). The same argument applies, now with only
four facets containing either endpoint.

\smallskip\noindent\emph{The diagonal with the exceptional column outside \(Q\).}
It remains that \(\tau_0\) is the third mixed facet, outside the
four facets through \(v_*\). All four columns \(1,2,U,V\) lie
in \(J\), and the local incidence triangles give
\[
 d_p(1)=d_p(2)=d_p(U)=0,\qquad
 d_q(1)=d_q(2)=d_q(V)=0.
\]
The four columns have rank three, with every triple independent, since
they represent the quadrilateral \(Q\). Hence the kernel in
\(\operatorname{rowspan}X\) of each triple \(\{1,2,U\}\),
\(\{1,2,V\}\) is the line spanned by \(X_{v_*}\), which vanishes on
all four columns. Both vectors are therefore multiples of that row. Nonnegative, nonzero restrictions to \(J\) make the
multiples positive. The zero at every \(J\)-coordinate again forces
both restrictions to vanish, a contradiction.
\end{proof}

\begin{proposition}[An almost-simple face excludes unit count slacks]
\label{prop:almost-simple-unit-slack}
Let \(A\) be a three-face of a seven-polytope \(P\), and put
\(D=A^\Delta\), with vertex--facet counts \((a,\alpha),(c,\gamma)\).
Suppose \(A\) has exactly one degree-four vertex and all other degrees
three, \(\alpha\ge7\), and \(c-\gamma\ge2\). Then
\[
 f_0(P)=a+\gamma+1,\qquad f_6(P)=c+\alpha+1
\]
cannot both hold.
\end{proposition}

\begin{proof}[Proof of Proposition~\ref{prop:almost-simple-unit-slack}]
Suppose both counts hold. Euler's formula gives \(a=2\alpha-5\ge9\)
and \(\iota(A)=4\). Its one-column deletion has rank four by the
preceding proof. The inequalities \(c\le2\gamma-4\),
\(c\ge\gamma+2\) imply \(\gamma\ge6\) and
\[
 \iota(D)\le3+(2\gamma-c-4)\le\gamma-3.
\]
Thus \(D^*\) cannot be a pyramid with \(\gamma-1\) base vertices;
its one-column deletion also has rank four.

\smallskip\noindent\emph{The two allocations.}
There are \(\gamma\) minimal cofaces and \(\gamma+1\) outside
vertices at \(A\). If an outside vertex \(x_0\) belongs to none,
all cofaces are pyramids \(A*v_i\). At least \(\alpha-1\) unique
inducers intersect in a face \(B\) containing every \(v_i\), perhaps
\(x_0\), and no \(A\)-vertex, since \(\alpha-1\ge6>\iota(A)\).
Restricted normal rank gives \(\dim B\le3\), and the quotient span
gives the reverse inequality. If \(x_0\in B\), the full join has
\(\gamma+1\) figure vertices instead of \(\gamma\). Otherwise
\(P_0=A*B\) contains all vertices except \(x_0\), with the same
quotient rays and all \(c\) containing facets supporting. The
nonvertex normalized quotient image of \(x_0\) lies on at most two
facets of the three-dimensional figure.

Every unique inducer is the plane of \(\tau*B\). If all groups
were unique, \(x_0\) would satisfy every facet inequality of \(P_0\).
Thus one group \(\tau_0\) is double; \(x_0\) lies beyond
\(\tau_0*B\), while the other planes remain supporting. For forms
\(g_0\le0\) there and \(h_K\le0\) at \(A*K\) strictly missing
\(x_0\), put \(\lambda=g_0(x_0)>0\), \(\mu_K=h_K(x_0)<0\).
The form \(\lambda h_K-\mu_Kg_0\le0\) supports a distinct facet
\(\conv(x_0,\tau_0*K)\) inducing \(\tau_0\) for every such
\(K\). There are at least \(c-2\ge6\), contrary to multiplicity
two. Hence one minimal coface contains two outside vertices \(p,q\),
and all others are pyramids.

In polarity the same argument uses \(\gamma-1>\iota(D)\), deletion
rank four and \(a-2\ge7\). It shows that every noncontaining ambient
facet induces a facet of \(A\), with exactly one double group
\(\tau_0\) and all others unique.

\smallskip\noindent\emph{The nonzero row-space extensions.}
Let \(I\) index the single-apex outside vertices, \(C\) the facets
containing \(A\), and \(J\) the unique inducers. Normalize the
quotient-ray rows and inducing columns as in the proof of
Theorem~\ref{thm:seven-first-size-exclusion}. The block on rows
\(\V(A),I\) and columns \(C,J\) is
\[
 \begin{pmatrix}0&X_J\\Y_I&0\end{pmatrix},
\]
of rank eight, the full slack rank. It uniquely extends
\(Z_{p,J},Z_{q,J}\) to \(d_p,d_q\in\operatorname{rowspan}X\),
nonnegative on \(J\) with their ambient incidence zeros.

If \(d_p=0\), the unique inducers intersect in a three-face \(B\)
containing \(I,p\), perhaps \(q\), and no \(A\)-vertex. Including
\(q\) gives the excessive figure count again. Excluding it gives
\(P_0=A*B\), with \(\gamma\) vertices and \(c\) facets in \(B\).
The point \(q\) has the same quotient ray as \(p\), lies beyond
only \(\tau_0*B\), and forces an inducer for every \(B\)-facet
missing \(p\), by the same supporting-form construction. At least
\(c-(\gamma-1)\ge3\) such facets are needed, contrary to the double
group. Hence \(d_p\ne0\), and likewise \(d_q\ne0\).

\smallskip\noindent\emph{The local coface.}
The double minimal coface \(H=\conv(A,p,q)\) is four-dimensional.
Every facet \(G\ne A\) of \(H\) contains an outside vertex \(i\).
Extend it to an ambient facet not containing \(H\). That facet cannot
contain \(A\): a face containing \(A\) and \(i\) contains their
minimal coface \(H\). The dual allocation makes its intersection with
\(A\) an \(A\)-facet lying in \(G\). Thus \(G\cap A\) is a
facet of \(A\), and ridge incidence gives one \(L_\tau\) for each
\(\tau\). For unique \(\tau\), the ambient inducer contains
\(L_\tau\), so local incidence with \(p,q\) gives the required zero
of \(d_i(\tau)\). Lemma~\ref{lem:almost-simple-support} contradicts
these two nonzero extensions.
\end{proof}

\begin{corollary}[Final Hall-pair data at the lower boundary]
\label{cor:final-seven-boundary}
Under the hypotheses of Theorem~\ref{thm:seven-second-boundary}, every
positive-slack maximizing pair consists, up to polarity, of simple faces
with counts
\[
 (a,\alpha)=(2h+6,h+5),\quad(c,\gamma)=(2h+8,h+6),
 \qquad(u,v)=(2,1).
\]
For \(h\ge2\) these are the only possible boundary data. For \(h=1\),
the additional zero-slack alternative has \(P\cong W_{6,3,1}\) or its
polar, as in Theorem~\ref{thm:zero-seventeen-rigidity}.
Every boundary attainer is balanced of deficiency exactly \(h\), with
matching number \(2h+14\).
\end{corollary}

\begin{proof}
In the almost-simple alternative of Theorem~\ref{thm:seven-second-boundary},
\(\alpha=\gamma=h+6\ge7\), \(c-\gamma=h+1\ge2\) and \(u=v=1\).
Proposition~\ref{prop:almost-simple-unit-slack} excludes it.
Corollary~\ref{cor:zero-boundary-exclusion} leaves only the stated simple
type for \(h\ge2\) and the simplex zero-slack row for \(h=1\).
Theorem~\ref{thm:zero-seventeen-rigidity} identifies the latter's
combinatorial type. The counts, deficiency and matching number are those
of Theorem~\ref{thm:seven-second-boundary}.
\end{proof}

\subsection{Incidence-only support injection}
\label{subsec:incidence-injection}

The support obstruction has a dimension-independent form. At a simple
vertex of the base facet, the unique outside edge identifies an outside
vertex. Incidence zeros then determine a line in the slack row space.
Nonvanishing, rather than positivity, makes this correspondence injective.
Coloring base vertices by their outside neighbors is classical
\cite[Lemma~20]{pinedaUgonYost2022}; the row-space condition below gives
the additional restriction.

\begin{theorem}[Incidence-only support injection]
\label{thm:incidence-support-injection}
Let \(A\) be an \(r\)-polytope, \(r\ge2\), with vertex--facet
slack matrix \(X\). Let \(H\) be an \((r+1)\)-polytope having
\(A\) as a facet, and let \(U=\V(H)\setminus\V(A)\).
Suppose that all other facets of \(H\) are indexed by \(\F(A)\),
one facet \(L_\tau\) for each \(\tau\in\F(A)\), with
\(L_\tau\cap A=\tau\). Let \(\mathcal E\subseteq\F(A)\).
For each \(i\in U\), suppose there is a nonzero vector
\(d_i\in\operatorname{rowspan}X\) such that
\[
 i\in L_\tau\ \Longrightarrow\ d_i(\tau)=0
 \qquad(\tau\in\F(A)\setminus\mathcal E).
\]
Then the simple vertices of \(A\) outside
\(\bigcup_{\tau\in\mathcal E}\tau\) inject into \(U\).
In particular, their number is at most \(|U|\).
\end{theorem}

\begin{proof}[Proof of Theorem~\ref{thm:incidence-support-injection}]
Let \(v\) be a simple vertex of \(A\). Exactly \(r\) facets
of \(A\) contain \(v\), so exactly \(r+1\) facets of \(H\)
contain it: \(A\) and the corresponding \(L_\tau\).
Thus \(v\) is simple in \(H\). Of its \(r+1\) edges,
\(r\) lie in \(A\); the remaining edge joins it to an outside
vertex \(i(v)\in U\). This edge is the intersection of the
\(r\) facets \(L_\tau\) through \(v\), and consequently
\(i(v)\) belongs to each of them.

Suppose that \(v\) misses every facet in \(\mathcal E\).
The vector \(d_{i(v)}\) vanishes on the \(r\) facet columns
incident with \(v\). Those columns have rank \(r\): their
homogeneous facet normals annihilate \((1,v)\), while their ordinary
normals span \(\R^r\). The usual full-rank slack factorization
identifies these column ranks with their ranks in \(X\), and
\(\operatorname{rowspan}X\) has dimension \(r+1\).
The common zero subspace in this row space is therefore a line.
It contains the vertex row \(X_v\), so
\[
 d_{i(v)}=c_vX_v,\qquad c_v\ne0.
\]

If another simple vertex \(w\) missing the exceptional union has
\(i(w)=i(v)\), the same vector vanishes on every facet through
\(w\). Since \(c_v\ne0\), all these facets contain \(v\).
Their intersection in \(A\) is \(\{w\}\), and hence \(v=w\).
Thus \(v\mapsto i(v)\) is injective.
\end{proof}

The theorem permits arbitrary signs in every \(d_i\). Its nonzero
hypothesis is essential: zero vectors would satisfy every incidence
condition without restricting the number of base vertices. In an ambient
Hall-pair application, proving that the extension vectors are nonzero is
a separate geometric step, as in
Propositions~\ref{prop:one-simple-double-slack}
and~\ref{prop:almost-simple-unit-slack}.

\begin{corollary}[Simple three-dimensional threshold]
\label{cor:sign-free-simple-support}
Under the hypotheses of Theorem~\ref{thm:incidence-support-injection},
suppose that \(A\) is a simple three-polytope with \(\alpha\)
facets, \(\mathcal E=\{\tau_0\}\), and \(|U|=s\). Then
\[
 \alpha-3\le s.
\]
Consequently, if \(\alpha\ge s+4\), the required nonzero vectors
cannot exist, without any sign assumptions.
\end{corollary}

\begin{proof}[Proof of Corollary~\ref{cor:sign-free-simple-support}]
Euler's formula and simplicity give \(f_0(A)=2\alpha-4\).
Every edge of the polygon \(\tau_0\) belongs to a different
other facet of \(A\), since two facets intersect in a single face.
Thus \(f_0(\tau_0)\le\alpha-1\). At least \(\alpha-3\)
vertices miss \(\tau_0\), and the injection bounds their number
by \(s\).
\end{proof}

For \(s=2\), this removes the nonnegativity hypothesis from
Lemma~\ref{lem:single-side-support}. For \(s=3\), the threshold
becomes seven facets. The same incidence mechanism also handles the
exceptional degree-four vertex.

\begin{corollary}[Sign-free almost-simple obstruction]
\label{cor:sign-free-almost-simple}
Let \(A\) be a three-polytope with \(\alpha\ge7\) facets,
all vertices of degree three except at most one of degree four,
and vertex--facet slack matrix \(X\). Let \(H\) be a four-polytope
with facet \(A\) and exactly two outside vertices \(p,q\), whose
other facets are \(L_\tau\), one for each \(\tau\in\F(A)\), with
\(L_\tau\cap A=\tau\). For a fixed \(\tau_0\in\F(A)\),
there are no nonzero \(d_p,d_q\in\operatorname{rowspan}X\) satisfying
\[
 i\in L_\tau\ \Longrightarrow\ d_i(\tau)=0
 \qquad(i\in\{p,q\},\ \tau\ne\tau_0).
\]
\end{corollary}

\begin{proof}[Proof of Corollary~\ref{cor:sign-free-almost-simple}]
If \(A\) is simple, Corollary~\ref{cor:sign-free-simple-support}
applies. Otherwise let \(v_*\) be its unique degree-four vertex.
Euler's formula gives \(f_0(A)=2\alpha-5\), and the facet-edge
argument above gives \(f_0(\tau_0)\le\alpha-1\).
There are at least \(\alpha-5\) ordinary, degree-three vertices
outside \(\tau_0\). When \(\alpha\ge8\), this exceeds two,
contrary to Theorem~\ref{thm:incidence-support-injection}.

It remains to treat \(\alpha=7\); then \(f_0(A)=9\).
If \(v_*\in\tau_0\), at least three ordinary vertices miss
\(\tau_0\), again a contradiction. Hence \(v_*\notin\tau_0\).
There are exactly two ordinary vertices \(u,w\) outside \(\tau_0\),
and \(\tau_0\) is a hexagon. Their distinct outside neighbors
exhaust \(\{p,q\}\). By the proof of the injection theorem,
after relabeling, \(d_p\) and \(d_q\) are nonzero multiples of
\(X_u\) and \(X_w\), respectively.

Every hexagon vertex \(z\) is ordinary and lies in two facets of
\(A\) other than \(\tau_0\). Its outside \(H\)-neighbor is
\(p\) or \(q\), and the corresponding vector vanishes on both
of these facet columns. Thus \(u\) or \(w\) lies in both facets,
hence on the edge of \(A\) in which they meet; this edge contains
\(z\). That vertex is the unique \(A\)-neighbor of
\(z\) outside the hexagon. On the other hand, \(v_*\) has at
most two neighbors outside the hexagon, namely \(u,w\), so at least
two of its four neighbors lie on the hexagon. At either such neighbor
\(z\), its unique edge leaving the hexagon joins it to \(v_*\),
contradicting the preceding conclusion.
\end{proof}

\subsection{The triple-coface boundary}
\label{subsec:triple-coface}

At the boundary where both Hall faces are simple and the count slacks
are \((u,v)=(2,1)\), the two outside vertices not needed as pyramid
apices can be allocated to one minimal coface. The incidence-only injection excludes
this allocation at every larger deficiency. The nonzero hypotheses
require an ambient argument: adding two points to a join need not
preserve the cone facets formed by either point alone.

\begin{theorem}[Triple-coface exclusion]
\label{thm:triple-coface-exclusion}
Let \(h\ge1\), and let \(A\) be a simple three-face of a
seven-polytope \(P\) whose dual face \(D=A^\Delta\) of \(P^*\) is also
simple, with counts
\[
 (a,\alpha)=(2h+6,h+5),\qquad
 (c,\gamma)=(2h+8,h+6).
\]
Suppose \(f_0(P)=f_6(P)=3h+14\). Suppose that the minimal
cofaces of \(A\) consist of \(\gamma-1\) pyramids with distinct
outside apices and one four-coface
\[
 H=\conv(A\cup\{p_1,p_2,p_3\})
\]
with exactly three outside vertices. Then \(h=1\).
Every facet of \(P\) not containing \(A\) induces a facet of \(A\);
one facet \(\tau_0\) has two inducers and all others have one.
The facets of \(H\) are \(A\) and exactly one \(L_\tau\) per
\(\tau\in\F(A)\), with \(L_\tau\cap A=\tau\).
\end{theorem}

We first exclude a hidden polar vertex, then prove nonvanishing by
counting the second facets at retained ridges. This supplies the full
local incidence data needed for the injection.

\begin{proof}[Proof of Theorem~\ref{thm:triple-coface-exclusion}]
\emph{The polar allocation.}
Since \(f_0(P)=a+\gamma+2\) and the minimal cofaces of \(A\) contain
\(\gamma+2\) distinct outside vertices, every outside vertex lies in a
minimal coface. Dually, in \(P^*\) every facet not containing \(D\)
induces a facet of \(D\):
the facet \(\sigma\) corresponding to \(H\) has three inducers,
and each of the other \(\gamma-1\) facets has one. There are
\(\alpha+1\) vertices outside \(D\), and its figure has
\(\alpha\) extreme rays. If one vertex \(x_0\) belongs to no
minimal coface at \(D\), all those cofaces are pyramids with
distinct apices.

Intersect the \(\gamma-1\) unique inducing facets. The resulting
face \(B_0\) contains all \(\alpha\) apices, possibly \(x_0\),
and no \(D\)-vertex, since \(\gamma-1>3\).
Deleting the \(\sigma\) column from the slack matrix of simple
\(D\) leaves rank four, as before
Theorem~\ref{thm:seven-first-size-exclusion}. Thus
\(\dim B_0\le3\), while its quotient rays span four dimensions,
giving \(\dim B_0\ge3\). Since \(\dim D+\dim B_0+1=7\), the affine
hulls of \(D\) and \(B_0\) are skew. If \(x_0\in B_0\), the full join would have
\(\alpha+1\) figure vertices instead of \(\alpha\).
Otherwise the vertices other than \(x_0\) span the join
\(P_0=D*B_0\), whose face figure at \(D\) is \(B_0\); thus \(B_0\) has
\(\alpha\) vertices and \(a\) facets.

Every containing plane \(D*K\) remains supporting, and every
unique inducing plane is \(\theta*B_0\), \(\theta\ne\sigma\).
Hence \(x_0\), being outside \(P_0\), is strictly beyond only
\(\sigma*B_0\). Its normalized quotient image is a nonvertex
point and lies on at most two facets of the three-dimensional figure.
For each of the at least \(a-2\) strictly missed planes \(D*K\),
the supporting-form combination at the ridge \(\sigma*K\) gives
a different facet inducing \(\sigma\), as in the proof of
Proposition~\ref{prop:one-simple-double-slack}.
But \(a-2\ge6\), exceeding its three inducers. This contradiction
excludes \(x_0\). There is consequently one double coface at \(D\)
in polarity and all other cofaces are pyramids. Dualizing gives
the asserted inducing groups at \(A\), with one double group
\(G_0,G_1\) over \(\tau_0\).

\smallskip\noindent\emph{Rank and nonvanishing.}
Let \(I\) index the \(\gamma-1\) single apices at \(A\),
\(C\) the containing facets, and \(J\) the unique inducers.
Normalize the duplicated \(\tau_0\) columns and the three equal
quotient-ray rows as in
Theorem~\ref{thm:seven-first-size-exclusion}. Write \(X\) for
the ordinary facet-slack matrix of \(A\), with the duplicate merged,
and \(Y\) for the quotient-ray block. Positive column scaling
corresponds to rescaling facet inequalities, so \(X\) is still a slack
matrix of \(A\). The submatrix
\[
 \begin{pmatrix}0&X_J\\Y_I&0\end{pmatrix}
\]
has rank eight, exhausting the ambient slack rank
\cite[Theorem~14]{gouveiaEtAl2013}. Both component ranks are four,
by one-column deletion for the simple faces. Each restricted row
\(Z_{p_i,J}\) therefore extends uniquely to
\(d_i\in\operatorname{rowspan}X\). On \(J\) it is nonnegative
and retains the ambient incidence zeros.

Suppose \(d_1=0\). The face
\[
 B=\bigcap_{\tau\ne\tau_0}G_\tau
\]
contains \(I,p_1\), possibly another \(p_i\), and no \(A\)-vertex.
Restriction rank four and the quotient span give \(\dim B=3\).
Since \(\dim A+\dim B+1=7\), the affine hulls of \(A\) and \(B\) are
skew, so quotient normalization is projectively injective on \(B\).
The three \(p_i\) have the same quotient ray; thus neither
\(p_2\) nor \(p_3\) belongs to \(B\).
Consequently \(B=\conv(I\cup\{p_1\})\), with \(\gamma\)
vertices and \(c\) facets, and \(P_0=A*B\) contains every
vertex except \(p_2,p_3\).

Each of \(p_2,p_3\) satisfies the inequalities of all facets \(A*K\)
and \(\tau*B\) \((\tau\ne\tau_0)\) of \(P_0\).
It lies outside \(P_0\), and hence strictly beyond \(\tau_0*B\).
For every \(B\)-facet \(K\) missing \(p_1\), both new points
strictly miss the plane of \(A*K\); this facet remains unchanged.
Its ridge \(\tau_0*K\) therefore has a second adjacent facet
in the final polytope \(P\), inducing \(\tau_0\).
Different \(K\) require different second facets. Otherwise one
facet would contain \(\tau_0\) and two distinct \(B\)-facets,
which span the old six-plane \(\tau_0*B\). That plane is not
supporting because both new points lie strictly beyond it.
The number of required inducers is at least
\[
 c-\deg_B(p_1)\ge c-(\gamma-1)=h+3\ge4,
\]
contrary to multiplicity two. Here the number of incident facets
at a three-polytope vertex equals its edge degree, at most
\(\gamma-1\). Thus \(d_1\ne0\), and the same proof gives
\(d_2,d_3\ne0\). The retained-ridge argument does not assume
that a cone facet survives the addition of the other new point.

\smallskip\noindent\emph{The local list and injection.}
Every nonbase facet of \(H\) contains some \(p_i\).
Extend it to an ambient facet cutting \(H\). That facet cannot
contain \(A\), since the smallest face containing \(A,p_i\)
is \(H\). By the polar allocation its intersection with \(A\)
is an \(A\)-facet. Each such ridge lies in exactly two \(H\)-facets,
one being \(A\); hence there is exactly one \(L_\tau\) per
\(\tau\). For \(\tau\ne\tau_0\), the unique ambient inducer
contains \(L_\tau\), so
\[
 p_i\in L_\tau\ \Longrightarrow\ d_i(\tau)=0.
\]
Corollary~\ref{cor:sign-free-simple-support} applies with
three outside vertices and gives
\(\alpha-3=h+2\le3\). Thus \(h=1\).
\end{proof}

\begin{proposition}[Local rigidity of the triple coface]
\label{prop:triple-local-rigidity}
Under the hypotheses of
Theorem~\ref{thm:triple-coface-exclusion}, the exceptional facet
\(\tau_0\) is a pentagon. The four-coface \(H\) is simple and
is combinatorially obtained by two successive vertex truncations
of a four-simplex. Its facet \(L_{\tau_0}\) contains all three
outside vertices. Moreover, if \(v_1,v_2,v_3\) are the vertices of
\(A\) outside \(\tau_0\), then after relabeling the extension vectors
\(d_i\) from the proof of Theorem~\ref{thm:triple-coface-exclusion}
satisfy \(d_i=c_iX_{v_i}\) with \(c_i>0\); in particular
\(d_i(\tau_0)>0\).
\end{proposition}

\begin{proof}[Proof of Proposition~\ref{prop:triple-local-rigidity}]
The preceding theorem gives \(h=1\), so \(A\) has eight vertices
and six facets. The injection bounds the number of vertices outside
\(\tau_0\) by three; the facet-size bound gives at least three.
Thus \(\tau_0\) is a pentagon, and its three outside vertices
\(v_1,v_2,v_3\) have distinct neighbors \(p_1,p_2,p_3\) in \(H\).
The proof of Theorem~\ref{thm:incidence-support-injection} gives
\[
 d_i=c_iX_{v_i},\qquad c_i>0.
\]
The sign follows from nonnegativity and nonzero restriction on
the five unique columns.

The simplicial polar \(S=A^*\) has six vertices. Write \(z_0\)
for the vertex dual to \(\tau_0\). It is adjacent to all other
five vertices. Its star is a cone over their five-cycle, while
the complementary disk is a triangulated pentagon with no interior
vertex. Its two diagonals share a vertex. Label the cycle
\(z_1,\ldots,z_5\) so that its three triangles are
\[
 z_1z_2z_3,\qquad z_1z_3z_4,\qquad z_1z_4z_5.
\]
Relabel so that these triangles correspond to \(v_1,v_2,v_3\),
relabeling the \(p_i\) and \(d_i\) accordingly.

By Theorem~\ref{thm:triple-coface-exclusion}, the facets of \(H\) are
\(A\) and the \(L_\tau\), so every vertex of \(A\) is simple in \(H\)
and has a unique neighbor outside \(A\); color it by this neighbor. A pentagon vertex lies in two nonexceptional facets.
Its color \(p_i\) makes both columns vanish on \(d_i\), and
therefore on \(X_{v_i}\). Their intersection edge joins that
pentagon vertex to \(v_i\). It is the unique edge of \(A\)
leaving the pentagon. Thus the color is determined by the fan
triangle incident with the corresponding cycle edge.

The vertices of \(L_\tau\) outside \(A\) are exactly the colors
present among the vertices of \(\tau\). One inclusion follows
from the outside edges. For the reverse inclusion, any outside
vertex of the three-polytope \(L_\tau\) has degree at least
three and at most two outside neighbors, so it has a neighbor
in \(\tau\), of that color. Each fan triangle has a boundary
edge, so all three colors occur on the pentagon.
In particular \(L_{\tau_0}\) contains all three \(p_i\).

In \(H^*\), let \(o\) be dual to \(A\), and use \(z_\tau\)
for the vertex dual to \(L_\tau\). The eight facets through
\(o\) are \(o\) joined to the eight triangular facets of \(S\).
The three other facets, dual to \(p_1,p_2,p_3\), are
\[
 z_0z_1z_2z_3,\qquad z_0z_1z_3z_4,\qquad z_0z_1z_4z_5.
\]
All eleven facets are tetrahedra. Their list is exactly that
obtained from the four-simplex on \(o,z_0,z_1,z_2,z_3\) by
stacking \(z_4\) on \(oz_0z_1z_3\), then \(z_5\) on
\(oz_0z_1z_4\). Thus \(H^*\) is simplicial, \(H\) is
simple, and duality gives the stated two vertex truncations.
\end{proof}

\section{Sharp constructions and equality types}
\label{sec:shared-coordinate}

The preceding lower bound is attained without strict stacking. The two
exposed faces in the following construction share one direction. This
reduces their join dimension by one, while the exposed prism still
yields the Hall bound. A complete facet description makes both the matching
number and connectedness explicit.

\begin{theorem}[Integral shared-coordinate family]
\label{thm:shared-coordinate}
Let \(\ell\ge1\), \(m\ge3\) and \(n\ge m+\ell+2\) be integers,
and put \(d=\ell+6\), \(b=m+\ell+1\). There is an integral
\(d\)-polytope \(W_{n,m,\ell}\) with connected nonincidence graph and
\begin{align*}
 f_0(W_{n,m,\ell})&=2n+b,\\
 f_{d-1}(W_{n,m,\ell})&=n+m(\ell+1)+b,\\
 \nu(\N(W_{n,m,\ell}))&=n+2b,\\
 \delta(W_{n,m,\ell})&=\min\{n-b,\ell(m-1)-1\}.
\end{align*}
In particular its deficiency is positive and its face lattice admits
no nontrivial free-join decomposition.
\end{theorem}

\begin{proof}[Proof of Theorem~\ref{thm:shared-coordinate}]
\emph{Vertices and dimension.}
Let \(P_n=\conv\{(k,k^2):0\le k<n\}\), and let \(K_m\) be the
polygon with vertices
\[
 v_j=\bigl(6j-3(m-1),\;6j^2-(m-1)(2m-1)\bigr),
 \qquad 0\le j<m.
\]
Its vertex centroid is the origin, which is therefore an interior
point. Write \(\lambda_j(u)\le1\)
for its normalized facet inequalities, indexed by the edges
\([v_j,v_{j+1}]\), with subscripts modulo \(m\).
These rational forms are explicit: for \(j<m-1\),
\[
 \lambda_j(u)=\frac{(2j+1)u_1-u_2}{\beta_j},\qquad
 \beta_j=6j(j+1)-3(2j+1)(m-1)+(m-1)(2m-1),
\]
and
\[
 \lambda_{m-1}(u)=\frac{-(m-1)u_1+u_2}{(m-1)(m-2)}.
\]
The denominators are positive: each numerator vanishes at the
interior point \(0\) of \(K_m\), so its maximum over \(K_m\), attained
on the corresponding edge, is positive.

In \(\R^\ell\), take the simplex
\(S_\ell=\conv\{q_0,\ldots,q_\ell\}\), with
\(q_0=(-1,\ldots,-1)\) and \(q_i=e_i\) for \(i\ge1\).
Its normalized facet forms are
\[
 g_0(w)=\sum_{r=1}^\ell w_r,\qquad
 g_i(w)=\sum_{r=1}^\ell w_r-(\ell+1)w_i\quad(1\le i\le\ell).
\]
Thus \(g_i(q_i)=-\ell\), \(g_i(q_k)=1\) for \(k\ne i\), and
\(\sum_i g_i=0\).
In coordinates \((x,y,z,u_1,u_2,w_1,\ldots,w_\ell,t)\), set
\begin{align*}
 a_{k,s}&=(k,k^2,s,0,0,0_\ell,0)
       &&(0\le k<n,\ s\in\{-1,1\}),\\
 b_j&=(0,0,0,v_j,0_\ell,1)&&(0\le j<m),\\
 c_i&=(0,0,1,0,0,q_i,1)&&(0\le i\le\ell),
\end{align*}
and take their convex hull \(W\).
The exposed face \(A=W\cap\{t=0\}\) is \(P_n\times[-1,1]\).
The exposed face \(B=W\cap\{t=1\}\) is the free join
\(K_m*S_\ell\) in coordinates \((z,u,w)\). It has dimension
\(\ell+3\), with \(b\) vertices and \(b\) facets.
All the displayed points remain vertices of these exposed faces,
and there are no other vertices of their convex hull.
The direction spaces of \(A,B\) intersect in the \(z\)-direction;
the additional \(t\)-direction separates them. Consequently
\(\dim W=3+(\ell+3)-1+1=\ell+6\), and \(f_0(W)=2n+b\).

\smallskip\noindent\emph{Complete facet description.}
The following inequalities describe \(W\):
\begin{alignat}{2}
 G_k:\quad &(2k+1)x-y-k(k+1)(1-t)\le0
       &&\quad(0\le k<n-1),\notag\\
 G_{n-1}:\quad &y-(n-1)x\le0,&&\notag\\
 C_{j,i}:\quad &\lambda_j(u)+g_i(w)-t\le0
       &&\quad(0\le j<m,\ 0\le i\le\ell),\label{eq:shared-facets}\\
 U_j:\quad &\lambda_j(u)+z-1\le0
       &&\quad(0\le j<m),\notag\\
 L_i:\quad &g_i(w)-z+t-1\le0
       &&\quad(0\le i\le\ell).\notag
\end{alignat}
Every displayed vertex satisfies them. To prove the reverse inclusion,
let a point satisfy all the inequalities. The three side inequalities
\(G_0,G_{n-2},G_{n-1}\) give
\[
 0\le x\le(n-1)(1-t).
\]
Hence \(t\le1\), and the full side system says
\((x,y)\in(1-t)P_n\), including \((x,y)=0\) when \(t=1\).
Put
\[
 r_- =\max_i g_i(w),\qquad
 r_+ =\min_j\{t-\lambda_j(u)\}.
\]
The \(C\) inequalities give \(r_-\le r_+\).
Since \(\sum_i g_i=0\), one has \(r_-\ge0\).
Boundedness and the interior origin of \(K_m\) imply a positive
linear dependence among its facet normals. Averaging the corresponding
values gives \(\min_j(t-\lambda_j(u))\le t\).
Thus \(0\le r_-\le r_+\le t\), in particular \(t\ge0\).
The \(L,U\) inequalities now give
\[
 r_--(1-t)\le z\le r_++(1-t).
\]
Choose \(z_B\in[r_-,r_+]\) with
\(z_A=z-z_B\in[-(1-t),1-t]\). Then
\[
 w\in z_B S_\ell,\qquad u\in(t-z_B)K_m.
\]
The point \((z_B,u,w,t)\) is a sum of a \(z_B\)-weighted simplex
point and a \((t-z_B)\)-weighted polygon point of \(B\).
The remaining coordinates form a \((1-t)\)-weighted point of \(A\).
Their weights sum to one. At a zero weight, the corresponding coordinates
vanish by boundedness of the factor. This proves that the point lies in
\(\conv(A\cup B)=W\), and so proves completeness of the system.

Each inequality defines a facet. A \(G\) equality contains the
two-dimensional edge-prism of \(A\) and all of \(B\). Their direction
spaces share exactly one dimension, giving
\(2+(\ell+3)-1+1=\ell+5\).
A \(C_{j,i}\) equality contains all of \(A\) and, in \(B\),
the polygon edge joined to the simplex facet opposite \(q_i\).
That \((\ell+1)\)-face projects injectively to \((u,w)\): its
\(z\)-coordinate is determined by \(g_i(w)=z\).
Its directions are therefore disjoint from those of \(A\), giving
\(3+(\ell+1)+1=\ell+5\).
A \(U_j\) equality contains the upper polygon of \(A\) and the
\((\ell+2)\)-face of \(B\) joining the polygon edge to all of
\(S_\ell\). Projection to \((u,w)\) is injective because
\(z=1-\lambda_j(u)\). It gives \(2+(\ell+2)+1=\ell+5\).
Likewise \(L_i\) contains the lower polygon and the face joining
all of \(K_m\) to the simplex facet; its projection is injective
because \(z=g_i(w)\), and has the same dimension.
This argument includes \(\ell=1\), when that simplex facet is a point.
Their vertex sets are pairwise distinct, so these facets are distinct.
Thus the system is irredundant,
and
\[
 f_{d-1}(W)=n+m(\ell+1)+b.
\]
Exactly the \(m(\ell+1)\) facets \(C_{j,i}\) contain \(A\).

\smallskip\noindent\emph{Maximum matching.}
Only the \(b\) vertices outside \(A\) can match its containing facets;
equivalently, Lemma~\ref{lem:face-bound} applied to \(A\) gives
\[
 \nu(\N(W))\le f_{d-1}(W)-m(\ell+1)+b=n+2b.
\]
The following nonincidences attain this bound. First match
\[
 U_j\longleftrightarrow a_{j,-1}\quad(0\le j<m),\qquad
 L_i\longleftrightarrow a_{i,1}\quad(0\le i\le\ell).
\]
Next match the side facets by
\begin{align*}
 G_k&\longleftrightarrow a_{k+\ell+1,1}
       &&(0\le k\le n-\ell-2),\\
 G_{n-\ell-1+r}&\longleftrightarrow a_{m+r,-1}
       &&(0\le r\le\ell).
\end{align*}
The first row uses the remaining upper indices \(\ell+1,\ldots,n-1\).
Since \(\ell\ge1\), the assigned index \(k+\ell+1\) is not an endpoint
\(k,k+1\) of \(G_k\). In the second row with \(r<\ell\), the side
\(G_{n-\ell-1+r}\) has endpoints \(n-\ell-1+r\) and \(n-\ell+r\), both
larger than \(m+r\) because \(n\ge m+\ell+2\).
The closing side misses index \(m+\ell\), which is neither zero nor
\(n-1\). These lower indices are disjoint from those assigned to the
\(U\) facets. We have matched all \(n+b\) noncontaining facets.
Finally use
\begin{align*}
 b_k&\longleftrightarrow C_{k+1,0}&&(0\le k\le m-2),\\
 b_{m-1}&\longleftrightarrow C_{0,1},&\qquad
 c_i&\longleftrightarrow C_{0,i}&&(i\ne1),\\
 c_1&\longleftrightarrow C_{1,1}.&&
\end{align*}
Every polygon vertex misses the indicated next edge, and \(c_i\)
misses the simplex facet opposite it. The facet endpoints are distinct;
index one exists since \(\ell\ge1\). These \(b\) additional pairs
are disjoint from the earlier pairs. Hence \(\nu(\N(W))=n+2b\),
and subtraction from the smaller class gives the asserted deficiency.

\smallskip\noindent\emph{Connectedness.}
Since \(n\ge m+\ell+2\ge6\), any two side facets miss a common vertex
of \(A\):
their edges use at most four polygon indices. Every vertex of \(A\)
misses a side, so these vertices and all side facets form one component.
Each \(U_j\) joins it through lower vertices of \(A\), and each
\(L_i\) through upper vertices. Each \(b_j\) misses \(m-2\ge1\)
of the \(U\) facets, each \(c_i\) misses \(L_i\), and each
\(C_{j,i}\) misses \(c_i\). Thus every graph vertex belongs to this
component. A nontrivial free join has disconnected nonincidence graph
(Lemma~\ref{lem:join-disconnects}), so no such face-lattice
decomposition is possible.
\end{proof}

\begin{corollary}[Prescribed deficiency in every admissible dimension]
\label{cor:shared-prescribed}
For integers \(d\ge7\), \(h\ge1\), put
\[
 \ell=d-6,\qquad
 m=\max\left\{3,1+\left\lceil\frac{h+1}{\ell}\right\rceil\right\},
 \qquad b=m+\ell+1,\qquad n=b+h.
\]
Then \(W_{n,m,\ell}\) has deficiency exactly \(h\), connected
nonincidence graph, and vertex--facet counts
\[
 N=2h+3b,\qquad M=h+(\ell+3)m+2\ell+2\ge N.
\]
For fixed \(d\), its larger class has size
\((2+3/(d-6))h+O_d(1)\).
In dimension seven it is balanced of size \(5h+12\), with matching
number \(4h+12\). At deficiency one the dimensions seven, eight and
nine give respectively
\[
 (f_0,f_{d-1},\nu)=(17,17,16),\quad(20,22,19),\quad(23,27,22).
\]
\end{corollary}

\begin{proof}
The chosen parameters satisfy \(n\ge m+\ell+2\),
\(n-b=h\) and \(\ell(m-1)-1\ge h\).
Theorem~\ref{thm:shared-coordinate} gives exact deficiency \(h\).
The count difference is \(M-N=\ell(m-1)-1-h\ge0\).
For \(h\ge\ell\), this difference is at most \(\ell-1\),
and \(m=h/\ell+O_\ell(1)\), proving the size law.
The stated specializations follow by substituting \(\ell=1,2,3\).
\end{proof}

\begin{corollary}[Sharp connected seven-dimensional minimum]
\label{cor:sharp-seventeen}
Every seven-polytope with connected nonincidence graph and no
vertex--facet assignment has at least seventeen vertices and seventeen
facets. Both bounds are attained simultaneously by the integral polytope
\(W_{6,3,1}\), whose matching number is sixteen.
\end{corollary}

\begin{proof}
The lower bounds are Theorem~\ref{thm:seven-second-size-exclusion}
with \(h=1\), and the attainer is
Corollary~\ref{cor:shared-prescribed} with \((d,h)=(7,1)\).
\end{proof}

The shared direction is what keeps this attainer from being a free
join. Its exposed faces
are a hexagonal prism \(A\), with counts \((12,8)\), and a
four-simplex \(B\). The three-dimensional figure at \(A\) has five
vertices and six facets: its \((u,w)\)-image is the free sum of a
triangle and a segment, with the origin interior in both factors.
By Proposition~\ref{prop:shared-unique-hall}, \(A\) is the unique
maximizing Hall face. Its dual Hall face has counts \((6,5)\), and the
ambient count slacks
are \((0,3)\). This realizes the simplex row of
Theorem~\ref{thm:seven-second-boundary}, not a free-join equality case.
For fixed \(d\), distinct deficiencies give different vertex counts
in Corollary~\ref{cor:shared-prescribed}, so the construction also
provides infinitely many pairwise nonisomorphic polytopes.

\subsection{Unique Hall faces and zero-slack extremizers}
\label{subsec:shared-rigidity}

In dimension seven the shared-coordinate family has a unique maximizing
Hall face. Together with the facet budget, this determines the
combinatorial type of every zero-slack extremizer with seventeen vertices
and seventeen facets.

\begin{proposition}[Unique Hall face of the seven-dimensional family]
\label{prop:shared-unique-hall}
For \(m\ge3\), \(n\ge m+3\), the exposed \(n\)-gonal prism \(A\)
is the unique face attaining \(H(W_{n,m,1})\). Every other three-face
has equal numbers of vertices and facets.
\end{proposition}

\begin{proof}
Use the coordinates of Theorem~\ref{thm:shared-coordinate}. A three-face
contained in \(A\) is \(A\) itself. A three-face contained in
\(B=K_m*S_1\) is a tetrahedron or a pyramid over \(K_m\),
and therefore has equal vertex and facet counts.

Consider a three-face meeting both \(A\) and \(B\), with intersections
of dimensions \(i,j\) in \(A,B\). Choose a supporting form exposing it. If this form has
nonzero \(z\)-coefficient,
neither intersection direction space contains the shared line. The
hull of the two intersections has dimension \(i+j+1=3\), so \((i,j)\) is \((2,0),(1,1)\)
or \((0,2)\). It is a polygonal pyramid or a tetrahedron, with equal
vertex and facet counts.

If the \(z\)-coefficient is zero, the intersection with \(A\) is a
vertical edge or lateral quadrilateral. It cannot be all of \(A\),
since meeting \(B\) would then give dimension at least four. The
projection of \(B\) along the shared line is the simplicial bipyramid
over \(K_m\). A proper face exposed by a zero-\(z\) form lifts a
vertex, edge or triangle of this bipyramid. The projection is
vertex-bijective, so these lifts preserve dimension and have no shared
line in their direction spaces. The dimension sum again makes the
mixed three-face a tetrahedron or a quadrilateral pyramid. This exhausts
all three-faces other than \(A\).

Theorem~\ref{thm:shared-coordinate} gives positive deficiency.
By the violating-face argument of Proposition~\ref{prop:global-seven},
every maximizing Hall face \(S\) and its dual are three-dimensional.
The ambient facet count is at least
\(f_0(S^\Delta)+f_2(S)\). Since
\[
 H(W)=f_0(S)+f_0(S^\Delta)>f_6(W),
\]
one has \(f_0(S)>f_2(S)\). Only \(A\) has this strict excess.
Hence \(A\) is the unique maximizer.
\end{proof}

\begin{theorem}[Zero-slack extremizers of size seventeen]
\label{thm:zero-seventeen-rigidity}
Let \(P\) be a seven-polytope with connected nonincidence graph and
positive deficiency, and suppose either \(f_0(P)\le17\) or
\(f_6(P)\le17\). If some maximizing Hall pair has a zero count slack,
then \(P\) is combinatorially equivalent to \(W_{6,3,1}\) or its
polar. Its maximizing Hall face is unique.
\end{theorem}

\begin{proof}
Theorem~\ref{thm:seven-second-boundary} gives
\(f_0(P)=f_6(P)=17\), \(\delta(P)=1\). Orient by polarity so that
\(u=0\). Corollary~\ref{cor:zero-boundary-exclusion} leaves only the
simplex row: \(A\) is simple with twelve vertices and eight facets,
\(B\) is a four-simplex, and its shadow \(T\) has five vertices and
six facets. Euler's formula makes \(T\) simplicial; it is a triangular
bipyramid. Proposition~\ref{prop:two-face-budget} gives \(\beta_0=0\).

\smallskip\noindent\emph{The simplex sign partition.}
Let \(b_0,\ldots,b_4\) project to the vertices \(t_0,\ldots,t_4\)
of \(T\). Their unique affine dependence is
\[
 \sum_{i=0}^4\theta_i t_i=0,\qquad\sum_{i=0}^4\theta_i=0.
\]
No coefficient vanishes. Otherwise four projected vertices would be
coplanar; all are vertices, so they would form a quadrilateral and the
fifth point its pyramid, contrary to simpliciality. Neither sign occurs
once, since that would express a vertex as a convex combination of the
others. The signs therefore split two and three. The vector
\(l=\sum_i\theta_i b_i\) is a nonzero vector along \(L\): the
\(b_i\) are affinely independent. A simplex facet form opposite
\(b_i\) has value zero at every other vertex and strictly negative
value \(-s_i\) at \(b_i\). Its linear part evaluated at \(l\) is
\(-s_i\theta_i\). Orient \(L\) by \(l\); then the \(z\)-coefficient
of this form has the sign of \(-\theta_i\). Thus
\(\{\beta_+,\beta_-\}=\{2,3\}\).

\smallskip\noindent\emph{The prism and the full facet incidences.}
Here \(v=3\). In \eqref{eq:two-face-budget} both signed coefficients
\(\beta_+-1,\beta_--1\) are positive, with sum three, while
\(\alpha_+,\alpha_-\ge1\) and \(p-\alpha_0\ge0\).
Equality forces
\(\alpha_+=\alpha_-=1\), \(p=\alpha_0=6\).
The projection \(R\) of \(A\) is a hexagon. Its six edge inequalities
are precisely the six vertical facets of \(A\). The two remaining
facet inequalities have the form
\[
 A=\{(x,z):x\in R,\ l_-(x)\le z\le l_+(x)\},
\]
with affine functions \(l_-,l_+\). Their difference is nonnegative
on \(R\), since \(R\) is the projection of \(A\). At a hexagon vertex, equality
\(l_-=l_+\) would give a lifted vertex incident to two vertical and
both end facets, contradicting simplicity. Thus the two end vertices over
each corner are distinct; these twelve points exhaust the \(a=12\)
vertices, so \(A\) is combinatorially a hexagonal prism.

The complete facet count now gives six facets containing \(A\), six
containing \(B\), and five mixed facets. The first six meet \(B\)
in the triangles of its bipyramid shadow, and the next six meet \(A\)
in its lateral quadrilaterals. Among the mixed facets, three each contain one
hexagonal end and one of the three simplex facets of one sign class;
the other two each contain the other end and one of the two remaining
simplex facets.
The sign classes are the two parts of the unique circuit of \(T\).
They identify its two apices and its three equatorial vertices, so this
description fixes all facet vertex sets, up to relabeling and reversal
of the ends. They are exactly those of \(W_{6,3,1}\).
Vertex--facet incidences determine the face lattice by intersections of
facets. Hence \(P\cong W_{6,3,1}\), and
Proposition~\ref{prop:shared-unique-hall} gives uniqueness of its Hall
face. Undoing polarity proves the alternative in the statement.
\end{proof}

The theorem classifies the zero-slack extremizers combinatorially.
The simple positive-slack alternative from
Corollary~\ref{cor:final-seven-boundary} is realized by the coplanar
construction below. Projective uniqueness is a separate question.

\subsection{Coplanar one-point extremizers}
\label{subsec:coplanar-extremizers}

The sharp minimum also admits the one-point mechanism suggested by
Jahn and Winter \cite[Section~6]{jahnWinter2021}. The added point lies
on several retained facet planes, so the strict transfer theorem
(Theorem~\ref{thm:horizon-transfer}) does not apply. We give the complete facet and matching descriptions instead.
Joins, truncation and stacking are classical operations
\cite{ziegler1995}; coplanar one-point modifications also occur in
\cite[proof of Theorem~3.3]{sjobergZiegler2020}.

\begin{lemma}[A fixed six-vertex factor]
\label{lem:coplanar-fixed-factor}
Let \(B\subset\R^3\) be the convex hull of
\[
 p=(0,-2,1),\quad q=(0,-2,-1),\quad
 r_1=(-1,-1,0),\ r_2=(-1,1,0),\ r_3=(1,1,0),\ r_4=(1,-1,0).
\]
Its facets are the eight triangles
\[
\begin{gathered}
 P_i=\conv(p,r_i,r_{i+1}),\quad Q_i=\conv(q,r_i,r_{i+1})
 \quad(1\le i\le3),\\
 E_1=\conv(p,q,r_1),\qquad E_4=\conv(p,q,r_4).
\end{gathered}
\]
They are defined by
\begin{alignat*}{3}
 P_1:&\ -y_1+y_3\le1,&\quad P_2:&\ y_2+3y_3\le1,&\quad P_3:&\ y_1+y_3\le1,\\
 Q_1:&\ -y_1-y_3\le1,& Q_2:&\ y_2-3y_3\le1,& Q_3:&\ y_1-y_3\le1,\\
 E_1:&\ -y_1-y_2\le2,& E_4:&\ y_1-y_2\le2.&&
\end{alignat*}
Its graph is the join of the edge \([p,q]\) and the path
\(r_1r_2r_3r_4\); exactly \(Q_1,Q_2,Q_3\) miss \(p\).
\end{lemma}

\begin{proof}[Proof of Lemma~\ref{lem:coplanar-fixed-factor}]
Substitution gives the stated equality triangles, with every remaining
point strictly beneath each plane. The normals for \(P_1,P_2,P_3\)
are independent and fix \(p\); those for \(Q_1,Q_2,Q_3\) fix \(q\).
The independent triples \((P_1,Q_1,E_1)\), \((P_1,Q_1,P_2)\),
\((P_3,Q_3,P_2)\), \((P_3,Q_3,E_4)\) fix \(r_1,r_2,r_3,r_4\),
respectively. Thus all six points are extreme; their hull is
three-dimensional because \(p,q,r_1,r_2\) are affinely independent.
The Euler inequality \(f_2\le2f_0-4\) bounds its facets by eight;
the eight distinct supporting triangles exhaust them. Their edges give
the asserted graph and the three facets missing \(p\).
\end{proof}

Let \(T_0\) be a rational three-polytope with six facets and a marked
vertex \(v\) of degree \(\rho\in\{3,4,5\}\), every other
vertex being simple. Truncate \(v\) by a rational plane strictly
separating \(v\) from the other vertices of \(T_0\), obtaining \(D\),
and denote the new facet by \(\sigma\).
Euler's formula gives \(f_0(T_0)=11-\rho\). Truncation replaces
\(v\) by \(\rho\) simple vertices, so \(D\) is simple with
\[
 f_0(D)=10,\qquad f_2(D)=7,\qquad f_0(\sigma)=\rho.
\]
The facets adjacent to \(\sigma\) are precisely the \(\rho\)
old facets whose planes contain \(v\).

Embed \(D*B\) in coordinates \((x,y,t)\in\R^3\times\R^3\times\R\)
by \((x,0,0)\) and \((0,y,1)\), and put
\[
 w_*=(v/2,p/2,1/2),\qquad
 Q^{\mathrm{cop}}_{T_0,v}=\conv\bigl((D*B)\cup\{w_*\}\bigr).
\]
The notation suppresses the realization of \(T_0\) and the truncating
plane; the conclusions below hold for every such rational choice.

\begin{theorem}[Sharp coplanar join construction]
\label{thm:coplanar-extremizer}
The polytope \(Q=Q^{\mathrm{cop}}_{T_0,v}\) is rational and
seven-dimensional, with
\[
 f_0(Q)=f_6(Q)=17,\qquad \nu(\N(Q))=16,\qquad\delta(Q)=1.
\]
Its nonincidence graph is connected. Its complete facet list consists of:
\begin{enumerate}
\item the six retained facets \(\tau*B\), \(\tau\ne\sigma\),
with \(w_*\) added to their vertex sets exactly when
\(v\in\aff\tau\);
\item the eight retained facets \(D*K\), \(K\in\F(B)\),
with \(w_*\) added exactly when \(p\in K\);
\item the three facets
\(F_i=\conv\bigl((\sigma*Q_i)\cup\{w_*\}\bigr)\), \(1\le i\le3\).
\end{enumerate}
The exposed face \(D\) and its dual face have counts \((10,7)\)
and \((8,6)\), respectively, and count slacks \((1,2)\).
\end{theorem}

\begin{proof}[Proof of Theorem~\ref{thm:coplanar-extremizer}]
\emph{Facet planes and completeness.}
The point \(w_*\) is beyond only \(\sigma*B\). It lies on
the \(\rho\) adjacent \(D\)-factor planes and the five
\(B\)-factor planes through \(p\), and strictly beneath all other
factor planes. All fourteen retained facets remain six-dimensional.

A facet of \(Q\) through \(w_*\) meets the old join in a face of
dimension at least five. A six-dimensional old face is a retained
coplanar facet. Otherwise the face is a ridge, whose normal cone is
generated by its two incident facet forms. A positive combination can
vanish at \(w_*\) only when one form is \(g\le0\) for
\(\sigma*B\), with \(g(w_*)>0\), and the other is a form
\(\phi\le0\) with \(\phi(w_*)<0\). Two zero values put \(w_*\)
in the ridge's affine hull and cannot give a facet; a positive
and a zero value cannot cancel. The ridges of \(\sigma*B\)
are \(R*B\) for edges \(R\) of \(\sigma\), and
\(\sigma*K\) for facets \(K\) of \(B\).
The other plane at each \(R*B\) is coplanar with \(w_*\).
Among the \(\sigma*K\), exactly the three \(K=Q_i\) have
a negative value. The supporting form
\[
 g(w_*)\phi-\phi(w_*)g\le0
\]
then has old equality face \(\sigma*Q_i\) of dimension five,
and \(w_*\) outside its affine hull. This gives exactly the three
listed new facets, proving the complete seventeen-facet description.

\smallskip\noindent\emph{Vertices.}
The eight \(D*K\) planes intersect in the affine span of \(D\).
For a \(D\)-vertex off \(\sigma\), its three incident retained
\(D\)-planes fix it in that span. For a vertex on \(\sigma\),
its two other incident \(D\)-planes and a new \(F_i\) recover
the missing \(\sigma\)-plane and again fix it. Thus all ten
\(D\)-vertices remain extreme. The six retained \(D\)-factor
planes intersect in the affine span of \(B\): deleting one facet
column from the slack matrix of simple \(D\) leaves rank four.
Otherwise the six corresponding vertices of the simplicial polar
\(D^*\) would be coplanar, and \(D^*\) would be a pyramid over a
hexagonal facet.
Within the affine span of \(B\), its incident \(B\)-planes fix
each old vertex. All six remain extreme. The point \(w_*\) is
outside the old join and is a new vertex. Hence \(f_0(Q)=17\);
the old join already spans dimension seven.

\smallskip\noindent\emph{Connectedness.}
Every vertex node in \(\N(D)\) has four facet neighbors.
Every component therefore has at least four facet nodes; two cannot
fit among seven. Facet nodes are not isolated, so \(\N(D)\)
is connected. For \(B^*\), every vertex node has three neighbors
among six facet nodes. If disconnected, there are two groups of three
facets. All vertices in one group belong to the three facets of the
other group; three distinct facets of a three-polytope have at most
one common vertex. There would be at most two vertices, instead of
eight. Thus \(\N(B)\), the transpose of \(\N(B^*)\), is connected.

The ten \(D\)-vertices and the nine facets not containing \(D\)
induce \(\N(D)\) with its \(\sigma\)-node replaced by three
copies. The six \(B\)-vertices and the eight \(D*K\) facets
induce \(\N(B)\). The vertex \(p\) misses all three \(F_i\)
and joins these blocks. The vertex \(w_*\) misses every \(D*Q_i\)
and belongs to the same component. This accounts for all thirty-four
nodes of \(\N(Q)\).

\smallskip\noindent\emph{Matching saturation.}
The face \(D\) has ten vertices and exactly eight containing
facets, so Lemma~\ref{lem:face-bound} gives \(\nu(\N(Q))\le16\).
We match the nine facets not containing \(D\) to distinct
\(D\)-vertices. Their nonincidences are those of its six facets
other than \(\sigma\) and three copies of \(\sigma\).
Hall's condition follows in three cases. A set using one distinct
facet has at most three nodes and at least four neighbors, because
a facet of simple \(D\) has at most six vertices. A set using two
distinct facets has at most four nodes and at least eight neighbors,
because their intersection has at most two vertices. A set using at
least three distinct facets has at least nine neighbors and at most
nine nodes, because their common intersection has at most one vertex.
For the other seven pairs use
\[
\begin{gathered}
 p\longleftrightarrow D*Q_1,\quad q\longleftrightarrow D*P_1,\quad
 r_1\longleftrightarrow D*P_2,\quad r_2\longleftrightarrow D*E_4,\\
 r_3\longleftrightarrow D*E_1,\quad r_4\longleftrightarrow D*Q_2,
 \quad w_*\longleftrightarrow D*Q_3.
\end{gathered}
\]
These are distinct nonincidences, disjoint from the nine earlier pairs.
Hence \(\nu(\N(Q))=16\), and \(\delta(Q)=1\).

Finally, the quotient cone at \(D\) has the six rays of \(B\).
The new point has quotient direction \((p/2,1/2)\), on the old
\(p\)-ray. Its dual face is therefore combinatorially \(B^*\),
with eight vertices and six facets. The count slacks are
\(17-10-6=1\) and \(17-8-7=2\).
\end{proof}

\subsection{A seventeen-inequality rational example}
\label{subsec:coplanar-cube}

The cube input makes every coordinate and matching edge explicit.
Take \(T_0=[-1,1]^3\), \(v=(1,1,1)\), and
\[
 D=[-1,1]^3\cap\{x_1+x_2+x_3\le2\}.
\]
Its vertices are the seven cube vertices other than \(v\), together
with \(a_1=(0,1,1)\), \(a_2=(1,0,1)\), \(a_3=(1,1,0)\).
In the join embedding, add
\[
 w_*=(1/2,1/2,1/2,0,-1,1/2,1/2).
\]

\begin{proposition}[Complete rational H/V description]
\label{prop:coplanar-cube}
The seventeen points consisting of the ten embedded \(D\)-vertices,
the six embedded \(B\)-vertices, and \(w_*\) are exactly the
vertices of a seven-polytope whose irredundant facet inequalities are
\begin{alignat}{2}
 &\pm x_j+t\le1 &&\quad(1\le j\le3),\notag\\
 &-y_1+y_3\le t,\quad y_2+3y_3\le t,\quad y_1+y_3\le t,&&\notag\\
 &-y_1-y_3\le t,\quad y_2-3y_3\le t,\quad y_1-y_3\le t,&&\notag\\
 &-y_1-y_2\le2t,\quad y_1-y_2\le2t,&&\label{eq:coplanar-cube-hrep}\\
 &2(x_1+x_2+x_3)-y_1-y_3+3t\le4,&&\notag\\
 &6(x_1+x_2+x_3)+y_2-3y_3+11t\le12,&&\notag\\
 &2(x_1+x_2+x_3)+y_1-y_3+3t\le4.&&\notag
\end{alignat}
Its nonincidence graph is connected and has a matching of size sixteen.
\end{proposition}

\begin{proof}[Proof of Proposition~\ref{prop:coplanar-cube}]
Theorem~\ref{thm:coplanar-extremizer} identifies all vertices and
facets. We also verify directly that every solution of the
displayed system lies in \(\conv((D*B)\cup\{w_*\})\). Put
\[
 g=x_1+x_2+x_3-2+2t,
\]
the deleted join-facet form. Write
\(h_1=-y_1-y_3-t\), \(h_2=y_2-3y_3-t\),
\(h_3=y_1-y_3-t\). The last three inequalities are
\(h_i+2b_i g\le0\), where \((b_1,b_2,b_3)=(1,3,1)\).
They are precisely the three new supporting forms from the theorem.
Every displayed vertex satisfies all inequalities.

If a feasible point has \(g\le0\), it satisfies the entire
old join system and belongs to \(D*B\). Suppose \(g>0\).
The positive coordinate bounds imply
\[
 g\le1-t,\qquad g\le x_j\quad(1\le j\le3).
\]
Adding \(-y_1+y_3\le t\) to \(h_3+2g\le0\) gives
\(g\le t\). Hence \(\theta=2g\) satisfies \(0<\theta\le1\).
Subtract \(\theta(1,w_*)\) from the homogeneous point
\((1,x,y,t)\). Its homogeneous constant is \(1-\theta\),
and its value of the homogeneous form \(g\) is zero.
The coplanar old forms remain nonpositive. The three \(h_i\)
remain nonpositive because \(h_i(w_*)=-b_i\) and
\(h_i(x,y,t)\le-2b_i g\).
For an old negative coordinate form \(-x_j+t-1\), its value
increases by \(\theta\), but stays nonpositive since
\[
 2g\le2x_j\le1-t+x_j.
\]
Thus every old join inequality holds for the homogeneous residual.
If \(\theta<1\), dividing by \(1-\theta\) gives a point
of \(D*B\). If \(\theta=1\), the residual lies in the
recession cone of the bounded old join, and is zero. In both cases
the original point belongs to \(\conv((D*B)\cup\{w_*\})\).
The theorem's facet dimensions prove irredundancy.

An explicit matching of the nine non-\(D\)-containing facets is
\[
\begin{array}{c|c@{\qquad}c|c@{\qquad}c|c}
 +x_1&a_1&+x_2&a_2&+x_3&a_3\\
 -x_1&(1,-1,-1)&-x_2&(1,1,-1)&-x_3&(1,-1,1)\\
 F_1&(-1,-1,-1)&F_2&(-1,-1,1)&F_3&(-1,1,-1).
\end{array}
\]
Here \(\pm x_j\) names the corresponding facet
\(\pm x_j+t=1\), and each \(D\)-point is embedded at
\((y,t)=(0,0)\). Together with the seven pairs in the proof of
Theorem~\ref{thm:coplanar-extremizer}, these are sixteen distinct
nonincidences. They leave \((-1,1,1)\) and \(D*P_3\)
unmatched. Connectedness and the matching upper bound follow from
that theorem.
\end{proof}

\subsection{At least sixteen sharp combinatorial types}
\label{subsec:coplanar-types}

Seven rational marked inputs distinguish seven members of the coplanar
family. Their polars supply seven more, and the shared-coordinate example
and its polar are different from all fourteen. The separation uses only
the explicit facet vertex sets.

\begin{theorem}[Richness at the sharp connected minimum]
\label{thm:sixteen-extremizers}
There are at least sixteen pairwise combinatorially nonisomorphic rational
seven-polytopes with connected nonincidence graph, seventeen vertices,
seventeen facets and matching number sixteen.
\end{theorem}

\begin{proof}[Proof of Theorem~\ref{thm:sixteen-extremizers}]
\emph{Seven marked rational inputs.}
Each input is a rational three-polytope \(S_0\) with six vertices and a
distinguished \(\rho\)-gonal facet, all other facets being triangles.
After a rational interior point is translated to the origin,
\(T_0=S_0^*\) has six facets, and the vertex \(v\) dual to the
distinguished facet has degree \(\rho\); all other vertices are simple.
The choices in Table~\ref{tab:marked-seeds} have the indicated vertex degrees. The \emph{common-edge chain} is the combinatorial type of the polytope
\(B\) of Lemma~\ref{lem:coplanar-fixed-factor}, formed by the three
tetrahedra \(pqr_ir_{i+1}\) glued along the edge \(pq\); the facet labels
\(E_1,P_1,P_2\) refer to that lemma.
\begin{datadisplay}{Seven marked seeds and their vertex degrees}\label{tab:marked-seeds}
\[
\begin{array}{c|l|c|c|c}
 &S_0\text{ and distinguished facet}&\rho&
 \text{degrees on facet}&\text{other degrees}\\\hline
 1&\text{octahedron, triangle}&3&4,4,4&4,4,4\\
 2&\text{common-edge chain, }E_1&3&5,5,3&4,4,3\\
 3&\text{common-edge chain, }P_1&3&5,4,3&5,4,3\\
 4&\text{common-edge chain, }P_2&3&5,4,4&5,3,3\\
 5&\text{stacked square pyramid, base}&4&4,4,3,3&5,3\\
 6&\text{square with two upper vertices, base}&4&4,4,3,3&4,4\\
 7&\text{pentagonal pyramid, base}&5&3,3,3,3,3&5
\end{array}
\]
\end{datadisplay}
The octahedron and the fixed chain are rational. Input five is obtained
by strict rational stacking on a triangular side of a square pyramid.
For input six, use the square vertices
\(s_1=(-1,-1,0)\), \(s_2=(1,-1,0)\), \(s_3=(1,1,0)\),
\(s_4=(-1,1,0)\), and upper vertices
\(u_1=(-1,1,1)\), \(u_2=(1,-1,1)\). Its facets are the square
and the six triangles
\[
 s_1s_2u_2,\quad s_2s_3u_2,\quad s_3s_4u_1,\quad s_4s_1u_1,\quad
 s_1u_1u_2,\quad s_3u_1u_2.
\]
They lie in the supporting planes for
\(z\ge0\), \(\pm x\le1\), \(\pm y\le1\), and
\(\pm(x+y)+2z\le2\). All six points are extreme.
A six-vertex three-polytope with a quadrilateral facet has at most
seven facets: Euler gives \(e=f+4\), while \(2e\ge3f+1\).
Thus these seven supporting facets exhaust its faces and give the
listed degrees. Input seven is a pyramid over any rational pentagon.
Each resulting \(T_0\) satisfies the marked-vertex assumptions, and
any rational truncation of \(v\) as above gives a member of
Theorem~\ref{thm:coplanar-extremizer}.

\smallskip\noindent\emph{Facet-size separation.}
A vertex of \(S_0\) of degree \(k\) corresponds to a
\(T_0\)-facet with \(k\) vertices. If it lies on the distinguished
facet, truncation increases that facet's size by one and the coplanar
extension adds \(w_*\). Its \(Q^{\mathrm{cop}}\)-facet therefore
has size \(k+8\). Otherwise its size is \(k+6\).
The five \(D*K\) facets through \(p\) have size fourteen;
the other three have size thirteen. The three new \(F_i\) have
size \(\rho+4\). Applying these formulas to Table~\ref{tab:marked-seeds}
gives the multisets in Table~\ref{tab:coplanar-profiles}:
\begin{datadisplay}{Facet-size multisets of the seven coplanar extensions}\label{tab:coplanar-profiles}
\[
\begin{array}{c|l}
 \text{input}&\text{facet-size multiset}\\\hline\noalign{\vskip 3pt}
 1&7^{[3]},10^{[3]},12^{[3]},13^{[3]},14^{[5]}\\
 2&7^{[3]},9^{[1]},10^{[2]},11^{[1]},13^{[5]},14^{[5]}\\
 3&7^{[3]},9^{[1]},10^{[1]},11^{[2]},12^{[1]},13^{[4]},14^{[5]}\\
 4&7^{[3]},9^{[2]},11^{[1]},12^{[2]},13^{[4]},14^{[5]}\\
 5&8^{[3]},9^{[1]},11^{[3]},12^{[2]},13^{[3]},14^{[5]}\\
 6&8^{[3]},10^{[2]},11^{[2]},12^{[2]},13^{[3]},14^{[5]}\\
 7&9^{[3]},11^{[6]},13^{[3]},14^{[5]}
\end{array}
\]
\end{datadisplay}
Here \(k^{[m]}\) means \(m\) facets with \(k\) vertices,
not exponentiation. The seven multisets are distinct.

\smallskip\noindent\emph{Polarity and the remaining pair.}
The old \(D\)-vertices off \(\sigma\) lie in eleven facets,
and its \(\rho\) vertices on \(\sigma\) lie in thirteen.
The old \(B\)-vertices \(p,q,r_1,r_2,r_3,r_4\) lie in
\(11,14,10,12,12,10\) facets, respectively; the new vertex
lies in \(\rho+8\). Thus the vertex--facet degree multiset is
\[
 10^{[2]},\quad11^{[11-\rho]},\quad12^{[2]},\quad
 13^{[\rho]},\quad14^{[1]},\quad(\rho+8)^{[1]},
\]
with repeated values combined. Every polar has minimum facet size ten,
whereas every original has a facet of size at most nine. No original
is isomorphic to any polar; the seven polars are mutually distinct
because duality preserves and reflects combinatorial isomorphism.
The total vertex--facet incidence count is \(187+3\rho\), namely
\(196,199,202\).

For \(W_{6,3,1}\), Theorem~\ref{thm:shared-coordinate} gives facet
sizes \(9^{[6]},10^{[5]},15^{[6]}\): the \(G\)-facets contain
four prism vertices and five opposite-face vertices; the \(U,L\)
facets contain six prism vertices and four opposite-face vertices;
the \(C\)-facets contain all twelve prism vertices and three
opposite-face vertices. The same facet lists give vertex--facet degrees
\(10^{[6]},11^{[6]},13^{[2]},14^{[3]}\), which are its polar's
facet sizes. These two types differ from each other, and their total
incidence count is \(194\), different from every coplanar type
and its polar. There are therefore at least \(7+7+2=16\) types.
All polars are taken after an interior rational translation and are
rational. Their nonincidence graphs are transposes of the originals,
so connectedness and matching number sixteen are preserved.
\end{proof}

\subsection{Classification of the triple-coface branch}
\label{subsec:triple-classification}

The seven marked inputs above are exhaustive for this allocation.
We first classify their small marked seeds, then recover them from an
arbitrary ambient triple coface. The small-polytope geometry is classical
\cite{gruenbaum2003,ziegler1995}.

\begin{lemma}[The seven marked seeds]
\label{lem:seven-marked-seeds}
Let \(S_0\) be a three-polytope with six vertices and a distinguished
\(\rho\)-gonal facet, \(\rho\in\{3,4,5\}\), all other facets
being triangles. Its marked combinatorial type is one of the seven
inputs in the proof of Theorem~\ref{thm:sixteen-extremizers}.
\end{lemma}

\begin{proof}[Proof of Lemma~\ref{lem:seven-marked-seeds}]
If \(\rho=5\), the sixth vertex is the apex of a pentagonal
pyramid. This gives one marked type.

If \(\rho=4\), there are two vertices \(p,q\) outside the
quadrilateral. Euler's formula and the facet-edge count give
seven facets and eleven edges. A base diagonal cannot be an edge,
since an edge contained in that facet is one of its boundary edges.
The vertices \(p,q\) must be adjacent. Otherwise every other
facet would use one outside vertex and a base edge; there are
only four such boundary edges, fewer than the required six triangles.
After removing the four base edges and \(pq\), six edges join
outside vertices to base vertices. Every base vertex has an
outside neighbor, so exactly two base vertices are common neighbors
of \(p,q\). They are also the third vertices in the two triangles
at \(pq\). The outside degrees sum to eight and are at least
three, giving \((3,5)\), \((4,4)\) or \((5,3)\).
In the first and third cases the two common neighbors must be
adjacent in the base, since the link of the degree-three vertex
is a triangle. This gives a triangular stacking on a square
pyramid. In the \((4,4)\) case the common neighbors must be
opposite: if adjacent, the edge between the other two base
vertices would have neither \(p\) nor \(q\) as a possible third
vertex of its other facet. Opposite common neighbors give
exactly the square with two upper vertices. Their full degree
multisets are \((5,4,4,3,3,3)\) and \((4,4,4,4,3,3)\).
Thus there are exactly two marked quadrilateral types.

If \(\rho=3\), the polytope is simplicial, with twelve edges
and eight facets. If all degrees are at least four, their sum
twenty-four makes every degree four. Its graph is \(K_6\) minus
a perfect matching, the octahedral graph; since the graph of a
three-polytope is planar and \(3\)-connected, it determines the face
lattice. Its triangular facets
form one orbit.
Otherwise a degree-three vertex can be unstacked: its triangular
link is replaced by the triangle on its neighbors. The visible cap
of the hull of the other vertices has this triangular boundary and
no interior vertices, since such a vertex would cease to be extreme
after stacking back. Thus the cap is one triangle. This leaves
the simplicial five-vertex triangular bipyramid. All its facets
are equivalent, so stacking back gives the common-edge chain.
In the notation of Lemma~\ref{lem:coplanar-fixed-factor}, its
eight triangles have three orbits: \(E_1,E_4\), of degree
triple \((5,5,3)\); \(P_1,Q_1,P_3,Q_3\), of triple
\((5,4,3)\); and \(P_2,Q_2\), of triple \((5,4,4)\).
Exchanging the two universal vertices and reversing the path
acts transitively within each group; different degree triples
separate the groups. Thus there are four marked triangular types.
Together the cases give \(1+2+4=7\). Rational realizations
of all seven were supplied in
Theorem~\ref{thm:sixteen-extremizers}.
\end{proof}

\begin{theorem}[The classified coplanar equality branch]
\label{thm:triple-global-classification}
Under the hypotheses of
Theorem~\ref{thm:triple-coface-exclusion}, \(h=1\) and \(P^*\)
is combinatorially equivalent to one of the seven coplanar types
\(Q^{\mathrm{cop}}_{T_0,v}\) in
Theorem~\ref{thm:sixteen-extremizers}. Conversely, the polar of
each of these seven types satisfies those hypotheses, has connected
nonincidence graph and matching number sixteen. Thus the polytopes
satisfying these hypotheses form exactly seven combinatorial types,
the polars of the seven coplanar types.
\end{theorem}

The proof reconstructs a complementary join face from a vanishing
polar extension. The sign pattern of the other extension then recovers
the marked seed and the full ambient facet list.

\begin{proof}[Proof of Theorem~\ref{thm:triple-global-classification}]
\emph{The polar exceptional extensions.}
By Theorem~\ref{thm:triple-coface-exclusion}, \(h=1\). As in its
proof, let \(\sigma\) be the facet of \(D\) corresponding to \(H\).
Let \(G_0,G_1\) be the double inducers of \(\tau_0\), viewed
as outside vertices in \(P^*\). They lie in a double minimal
coface \(K\) of \(D\). All other minimal cofaces there are
pyramids. The noncontaining inducing groups at \(D\) consist of
three facets over \(\sigma\), corresponding to \(p_1,p_2,p_3\),
and one over each other \(D\)-facet. The rank-eight
extension of that proof, applied in \(P^*\) at \(D\), gives ordinary \(D\)-facet row-space vectors \(e_0,e_1\),
nonnegative off \(\sigma\), with the actual incidence zeros.

Every nonbase facet of \(K\) meets \(D\) in a facet: extend it
to a noncontaining ambient facet and use this inducing allocation.
Ridge incidence gives exactly one local facet per \(D\)-facet.
For \(\theta\ne\sigma\), the unique inducer of \(\theta\) contains the
local facet over \(\theta\), so an incidence of \(G_j\) with that local
facet gives \(e_j(\theta)=0\). If both \(e_j\) were nonzero,
Theorem~\ref{thm:incidence-support-injection} would inject all
vertices of \(D\) outside \(\sigma\) into \(\{G_0,G_1\}\).
But \(D\) is simple with ten vertices and seven facets, so any
facet has at most six vertices. At least four miss \(\sigma\).
Hence one extension vanishes; relabel so that \(e_0=0\).

Normalize the exceptional rows and columns as in the proof of
Theorem~\ref{thm:seven-first-size-exclusion}. Its rank-eight
corner identity, now with three rows and two columns, gives
\[
 Z_{p_i,G_j}=\xi_i+\zeta_j,\qquad
 \xi_i=d_i(\tau_0),\quad \zeta_j=e_j(\sigma).
\]
Proposition~\ref{prop:triple-local-rigidity} gives \(\xi_i>0\).
Since \(e_0=0\), the facet \(G_0\) misses every \(p_i\).
The local facet \(L_{\tau_0}\) contains all three, so its
ambient inducer is \(G_1\). Consequently
\[
 \zeta_1=-\xi_1=-\xi_2=-\xi_3<0.
\]
In particular \(e_1\ne0\).

\smallskip\noindent\emph{Recovering the old join.}
In \(P^*\), intersect the unique inducing facets of \(D\).
Their rank-four restrictions and the quotient span give a
three-face \(B\) containing the five single apices and \(G_0\).
The vertex \(G_1\) cannot also belong: its quotient ray agrees
with \(G_0\), whereas quotient normalization is projectively
injective on this complementary three-face.
Thus these are precisely its six vertices, and all old vertices
form the full join \(D*B\).
The figure at \(D\) is combinatorially \(A^*\), so \(B\)
is the common-edge chain identified by the pentagonal fan in
Proposition~\ref{prop:triple-local-rigidity}. Its marked vertex
\(b_0=G_0\), dual to \(\tau_0\), has degree five.

The remaining point \(G_1\) satisfies all containing planes
\(D*K'\) and is on such a plane exactly when \(b_0\in K'\):
it has the same positive quotient ray. Every unique factor
plane \(\theta*B\), \(\theta\ne\sigma\), remains supporting.
Since \(G_1\) is outside the old join, it is strictly beyond
only \(\sigma*B\).

\smallskip\noindent\emph{The marked seed.}
The local facet of \(K\) over \(\sigma\) contains \(G_1\)
but not \(G_0\), since each of its three ambient inducers
contains \(G_1\) and misses \(G_0\).
Every \(D\)-vertex on \(\sigma\) has outside neighbor \(G_1\).
No other vertex can have that neighbor. Otherwise the injection
proof would make \(e_1\) a nonzero multiple of its vertex row.
Nonnegativity off \(\sigma\) makes the multiple positive,
contrary to \(e_1(\sigma)<0\).
Hence the other \(D\)-vertices have neighbor \(G_0\).

A \(D\)-facet adjacent to \(\sigma\) has both colors, giving
the zero \(e_1(\theta)=0\). A disjoint facet has only
\(G_0\)-colored vertices. Its local three-facet cannot also
contain \(G_1\): an outside vertex there would have degree
at least three and at most one outside neighbor, forcing a
neighbor in \(D\) of the other color. Therefore
\(e_1(\theta)>0\) for every disjoint facet.
In a simple three-polytope two facets sharing a vertex share
an edge, so these cases exhaust the facets other than \(\sigma\).

In the simplicial three-polytope \(S=D^*\), interpret \(e_1\)
as a linear form on normalized homogeneous vertex columns.
It is negative at the vertex \(s\) dual to \(\sigma\), zero
on all its \(\rho\) neighbors, and positive on all other vertices.
The latter set is nonempty: deleting \(s\) leaves homogeneous
rank four, so the nonzero form cannot vanish on all six columns.
Thus \(3\le\rho\le5\).
Cutting \(S\) by \(e_1\ge0\) introduces no new vertices,
since every edge leaving \(s\) ends at a zero neighbor.
The result \(S_0\) has exactly the other six vertices and a
marked \(\rho\)-gonal facet. All its other facets are retained
triangles. No old triangular facet can consist only of zero
neighbors, because its plane would have \(s\) on one side
and a positive vertex on the other and would not support \(S\).
The six retained columns span four dimensions, so \(S_0\)
is three-dimensional.

\smallskip\noindent\emph{Ambient reconstruction and converse.}
Lemma~\ref{lem:seven-marked-seeds} classifies \(S_0\).
Its polar \(T_0\) has six facets and a marked vertex \(v\)
of degree \(\rho\); all other vertices are simple.
The original \(D\) is combinatorially its marked-vertex
truncation: reversing the cut of \(S\) stacks \(s\) strictly
beyond the marked facet of \(S_0\).
The factor-plane signs of \(G_1\) are now exactly those in
Theorem~\ref{thm:coplanar-extremizer}: zero at the
\(\rho\) adjacent \(D\)-planes and the five \(B\)-planes
through \(b_0\), strictly beneath the others, and beyond only
\(\sigma*B\). The complete ridge-normal argument there fixes
all seventeen facet vertex sets. Thus \(P^*\) has the face
lattice of the corresponding \(Q^{\mathrm{cop}}_{T_0,v}\).

Conversely, in the polar of each constructed \(Q^{\mathrm{cop}}\),
the smaller dual Hall face has counts \((8,6)\) and its dual
has counts \((10,7)\). The six retained \(D\)-factor facets
give six pyramid cofaces; the three replacement facets over
\(\sigma\) give the single triple coface. These are exactly
the required allocation and counts. Connectedness and matching
number sixteen follow from
Theorem~\ref{thm:coplanar-extremizer} and polarity.
The face value \(8+10=18=17+\delta(P)\) makes the indicated
dual pair maximizing.
The seven facet-size profiles in
Theorem~\ref{thm:sixteen-extremizers} are distinct, and no
original type is isomorphic to any polar. Hence exactly seven
combinatorial types satisfy the hypotheses.
\end{proof}

\section{Factor certificates and the certified seven-dimensional minimum}
\label{sec:balanced}

Corollary~\ref{cor:stacking-transfer} gives an arithmetic condition for
a surviving free-join factor to certify failure of an assignment. We establish it
first, then use three-dimensional face bounds to determine the certified
minimum \(\eta_{7,1}=33\) in dimension seven.

\begin{proposition}[Arbitrary-factor criterion]\label{prop:factor-criterion}
Let \(A\) and \(B\) have dimensions \(r\) and \(m\), with
\[
 (f_0(A),f_{r-1}(A))=(a,\alpha),\qquad
 (f_0(B),f_{m-1}(B))=(b,\beta).
\]
Let \(F\) be a facet of \(A\) with \(\rho\) ridges, and obtain \(Q\) by
strictly stacking beyond the facet \(F*B\) of \(A*B\).  Then
\begin{equation}\label{eq:general-join-stack-counts}
 f_0(Q)=a+b+1,\qquad
 f_{r+m}(Q)=\alpha+2\beta+\rho-1,
\end{equation}
and the surviving copy of \(A\) gives
\begin{equation}\label{eq:general-factor-bound}
 \nu(\N(Q))\le b+\alpha+\beta+\rho.
\end{equation}
This bound is smaller than both class sizes in
\eqref{eq:general-join-stack-counts} if and only if
\begin{equation}\label{eq:general-factor-criterion}
 a-\alpha\ge\beta+\rho,
 \qquad
 \beta-b\ge2.
\end{equation}
Whenever these conditions hold, \(r,m\ge3\), and for
\(d=r+m+1\) one necessarily has
\begin{equation}\label{eq:factor-size-floor}
 f_0(Q)\ge2d+6,\qquad f_{d-1}(Q)\ge2d+6.
\end{equation}
\end{proposition}

The floor \eqref{eq:factor-size-floor} is not sharp;
Proposition~\ref{prop:sharp-floor} improves it to \(2d+8\).

\begin{proof}
The selected facet \(F*B\) has \(\rho+\beta\) ridges.  Strict stacking
therefore gives \eqref{eq:general-join-stack-counts}.  The copy of \(A\)
has \(a\) vertices and is contained in precisely the \(\beta\) retained
facets induced by \(B\).  Corollary~\ref{cor:stacking-transfer}, applied
to \(A*B\), with \(a+b\) vertices and \(\alpha+\beta\) facets, and to
this face, gives \eqref{eq:general-factor-bound}.  Comparing its right-hand side with the two
quantities in \eqref{eq:general-join-stack-counts} gives exactly
\eqref{eq:general-factor-criterion}; integrality accounts for the weak
inequalities.

The first condition forces \(a>\alpha\), and the second forces
\(\beta>b\).  Polytopes of dimension at most two have equally many vertices
and facets, so \(r,m\ge3\).  Every \(r\)-polytope has
\(\alpha\ge r+1\), while its \((r-1)\)-dimensional facet \(F\) has
\(\rho\ge r\).  The second condition excludes an \(m\)-simplex; hence
\(b\ge m+2\) and \(\beta\ge b+2\ge m+4\).  Finally, the first condition
gives \(a\ge\alpha+\beta+\rho\), and substitution in
\eqref{eq:general-join-stack-counts} yields
\[
 f_0(Q)\ge(r+1)+(m+4)+r+(m+2)+1=2d+6
\]
and
\[
 f_{d-1}(Q)\ge(r+1)+2(m+4)+r-1=2d+6.
\]
\end{proof}

\subsection{Balanced linear-size witnesses}

For the seven-dimensional minimum and explicit examples in every
dimension \(d\ge7\), put \(m=d-4\). We use the following rational
\(m\)-polytope:

\begin{itemize}
  \item for \(m=3\), let \(B_m\) be the octahedron, with
        \((f_0,f_{m-1})=(6,8)\);
  \item for \(m\ge4\), start with the four-dimensional bipyramid over a
        tetrahedron and take \(m-4\) successive pyramids.  The resulting
        \(B_m\) has
        \((f_0,f_{m-1})=(m+2,m+4)\).
\end{itemize}

Write \(b=f_0(B_m)\) and \(\beta=f_{m-1}(B_m)\), so that always
\(\beta-b=2\).  Set \(t=\beta+4\), and let \(A_m\) be obtained by truncating
one vertex of a prism over a rational \(t\)-gon.  This is a rational simple
three-polytope with
\begin{equation}\label{eq:a-counts}
 f_0(A_m)=2t+2=2\beta+10,\qquad
 f_2(A_m)=t+3=\beta+7,
\end{equation}
and the truncation creates a triangular facet \(\Delta\).

Form \(P=A_m*B_m\), stack strictly beyond the facet
\(S=\Delta*B_m\), and denote the resulting \(d\)-polytope by \(R_d\).  All
polytopes used before the last step are rational.  A rational stacking point
exists because the points strictly beyond \(S\) and strictly beneath all
other facets form a nonempty open set defined by rational inequalities.

\begin{theorem}[Balanced obstruction]\label{thm:balanced}
For every \(d\ge7\), the rational \(d\)-polytope \(R_d\) has equally many
vertices and facets, namely
\[
 s_d:=f_0(R_d)=f_{d-1}(R_d)=
 \begin{cases}
   33,&d=7,\\
   3d+9,&d\ge8.
 \end{cases}
\]
Its nonincidence graph is connected and has matching number exactly
\(s_d-1\).  Consequently, \(R_d\) has no vertex--facet assignment and admits
no nontrivial combinatorial free-join decomposition.
\end{theorem}

\begin{proof}
The facet \(S=\Delta*B_m\) has \(3+\beta\) ridges.  Strict stacking therefore
gives
\begin{align*}
 f_0(R_d)&=(2\beta+10)+b+1=3\beta+9,\\
 f_{d-1}(R_d)&=(\beta+7)+\beta-1+(3+\beta)=3\beta+9,
\end{align*}
where \(b=\beta-2\).  For \(m=3\), this common value is \(33\); for
\(m\ge4\), it is \(3(m+4)+9=3d+9\).

The copy of \(A_m\) is a face of \(P\) with \(2\beta+10\) vertices and is
contained in exactly the \(\beta\) facets induced by \(B_m\).  It is not
contained in \(S\).  Corollary~\ref{cor:stacking-transfer} gives connectedness
and the upper bound
\[
 \nu(\N(R_d))
 \le 2(3\beta+9)-(2\beta+10)-\beta=3\beta+8=s_d-1.
\]

It remains to attain the bound.  Label the vertices of the original
\(t\)-prism by \(A_{j,\pm}\), with \(j\in\mathbb Z/t\mathbb Z\), where the
second sign records the two bases.  The lateral facet \(X_j\) lies over the
polygon edge from \(j\) to \(j+1\), and we truncate the vertex \(A_{0,+}\).
The bases are \(T_+,T_-\), and the three vertices of the new truncation facet
are denoted \(U_v,U_1,U_{t-1}\), according to the vertical edge and the two
top edges at the truncated vertex.  All subscripts \(j\) below are understood
modulo \(t\).
Match the \(t+2\) retained facets induced by \(A_m\) as follows:
\[
 U_v\longleftrightarrow T_+,\qquad
 U_1\longleftrightarrow T_-,\qquad
 U_{t-1}\longleftrightarrow X_1,\qquad
 A_{j+2,-}\longleftrightarrow X_j
 \quad(j\in\mathbb Z/t\mathbb Z,\ j\ne1).
\]
Every displayed pair is a nonincidence, and all its vertices and facets are
distinct.  Among the vertices of \(A_m\) outside \(\Delta\), exactly \(t\)
remain unused: \(A_{3,-}\) and the \(t-1\) untruncated top vertices.  There
are \(3+\beta=t-1\) new cone facets, and every one is nonincident to every
vertex outside \(\Delta\).  Match these cone facets bijectively to any
\(t-1\) of the unused vertices.  We have now matched
\((t+2)+(t-1)=2t+1=2\beta+9\) vertices of \(A_m\), leaving one unmatched.

Finally, \(\N(B_m)\) has a matching covering all \(b\) vertices while leaving
two facets unused.  For the octahedron \(\conv\{\pm e_1,\pm e_2,\pm e_3\}\),
index the facets by sign vectors, so that \(\pm e_i\) misses exactly the
facets whose \(i\)th sign is \(\mp\). The pairs
\[
 \pm e_1\longleftrightarrow(\mp,\pm,\pm),\qquad
 \pm e_2\longleftrightarrow(\pm,\mp,\pm),\qquad
 \pm e_3\longleftrightarrow(\pm,\pm,\mp)
\]
use six distinct facets and leave \((+,+,+)\) and \((-,-,-)\) unused.  For the four-dimensional bipyramid, let
\(w_1,\ldots,w_4\) be the tetrahedron vertices and \(u_+,u_-\) the two
apices, and write its facets as
\(F_{i,\sigma}\), where \(i\) is the omitted tetrahedron vertex and
\(\sigma\) is the included apex; the pairs
\[
 w_1\longleftrightarrow F_{1,+},\quad
 w_2\longleftrightarrow F_{2,+},\quad
 w_3\longleftrightarrow F_{3,-},\quad
 w_4\longleftrightarrow F_{4,-},\quad
 u_+\longleftrightarrow F_{1,-},\quad
 u_-\longleftrightarrow F_{3,+}
\]
have the required property.  Pyramiding adds an isolated apex--base edge to
the nonincidence graph, so the claim persists through all successive
pyramids.  Use this matching on the \(\beta\) retained facets induced by
\(B_m\), and match the stacking apex to one of the two unused facets.  This
adds \(b+1\) edges.  The resulting matching has
\[
 (2\beta+9)+(b+1)=3\beta+8=s_d-1
\]
edges, proving equality.  The remaining assertions follow from connectedness
and Lemma~\ref{lem:join-disconnects}.
\end{proof}

Within strict join--stack realizations of the construction suggested
by Jahn and Winter, with a three-dimensional factor as Hall witness,
the polytopes \(R_d\) are optimal. This optimality refers only to that
class; the unrestricted connected minimum in dimension seven is smaller
(Section~\ref{sec:shared-coordinate}).

\begin{corollary}[Three-dimensional witness optimum]
\label{cor:class-optimal}
Let \(A\) be a three-polytope with \(a\) vertices and \(\alpha\) facets, let
\(B\) be an \(m\)-polytope with \(b\) vertices and \(\beta\) facets, and let
\(F\) be a \(q\)-gonal facet of \(A\).  Form \(A*B\), stack strictly
beyond \(F*B\) and strictly beneath every other facet, and suppose that the surviving copy of \(A\), through the
face-neighborhood bound, certifies failure of a vertex--facet assignment.
Then the resulting \((m+4)\)-polytope has at least \(33\) vertices and at
least \(33\) facets when \(m=3\), and at least \(3m+21\) vertices and at
least \(3m+21\) facets when \(m\ge4\).  The polytopes \(R_d\) attain these
bounds.
\end{corollary}

\begin{proof}
Proposition~\ref{prop:factor-criterion}, applied with \(\rho=q\), says that
the required strict inequalities are equivalent to
\begin{equation}\label{eq:class-conditions}
 a-\alpha\ge q+\beta,\qquad \beta-b\ge2.
\end{equation}
Every three-polytope satisfies
\(\alpha\ge a/2+2\), by Euler's formula and the minimum vertex degree three.
Combining \(\alpha\le a-\beta-q\) with this bound and \(q\ge3\) gives
\[
 a\ge2\beta+2q+4\ge2\beta+10,\qquad \alpha\ge a/2+2\ge\beta+7.
\]

The second inequality excludes the \(m\)-simplex, hence \(b\ge m+2\) and
\(\beta\ge m+4\).  When \(m=3\), the sharper planar bound
\(\beta\le2b-4\) excludes \(b=5\), so \(b\ge6\) and \(\beta\ge8\).
Writing \(V\) and \(F_Q\) for the numbers of vertices and facets of the
stacked polytope, we obtain for \(m=3\)
\[
 V\ge(2\cdot8+10)+6+1=33,\qquad
 F_Q\ge(8+7)+3+2\cdot8-1=33.
\]
For \(m\ge4\),
\[
 V\ge2(m+4)+10+(m+2)+1=3m+21,
\]
and
\[
 F_Q\ge(m+4+7)+3+2(m+4)-1=3m+21.
\]
The equalities in Theorem~\ref{thm:balanced} show that all bounds are sharp.
\end{proof}

In dimension seven the class of Corollary~\ref{cor:class-optimal} is
exhaustive: by \eqref{eq:general-factor-criterion}, a witness factor must have more
vertices than facets and the complementary factor at least two more facets
than vertices, and neither is possible in dimension at most two.  Since
\(\dim A+\dim B=6\), both factors of the free join are forced to be
three-dimensional, so every member of \(\mathcal C_{7,1}\) belongs to the
class. Consequently \(\eta_{7,1}=33\).
For \(d\ge8\) the restriction to a
three-dimensional witness factor is essential.

\begin{remark}[Higher-dimensional witness factors]\label{rem:higher-witness}
Let \(d=8\), let \(A\) be a prism over a truncated tetrahedron, a simple
rational four-polytope with \(24\) vertices and \(10\) facets, and let
\(B\) be the octahedron.  One facet \(F\) of \(A\) is a triangular prism,
and the facet \(F*B\) of the free join \(A*B\) has \(5+8=13\) ridges.
Stacking strictly beyond \(F*B\) gives a rational eight-polytope with
\(31\) vertices and \(30\) facets.  The surviving copy of \(A\) is
contained in exactly the eight facets induced by \(B\) and not in
\(F*B\), so Corollary~\ref{cor:stacking-transfer} yields connectedness and
\[
 \nu\le30+18+13-24-8=29<30.
\]
This eight-polytope is therefore a counterexample beyond free joins that is
smaller than \(R_8\) in both classes.  Appendix~\ref{sec:asymptotic}
systematizes this observation, and Appendix~\ref{sec:constant} reduces the
asymptotic additive remainder to a constant. Section~\ref{sec:sharp}
determines the exact optimum for every \(d\ge12\), and
Section~\ref{sec:lowdim} settles the remaining dimensions.
\end{remark}

\section{Sharp sizes in the factor-certified class}
\label{sec:sharp}

We determine both the sharp size floor and the first dimension in which
it can be attained, for every prescribed deficiency. Throughout this section, the
certificate is the surviving factor \(A\) in
Proposition~\ref{prop:factor-criterion}.

\begin{definition}[Certified deficiency classes]\label{def:certified-class}
For positive integers \(d,h\), let \(\mathcal C_{d,h}\) consist of strict
stacks beyond \(F*B\) in \(A*B\), of dimension \(d=r+m+1\) and with the
notation of Proposition~\ref{prop:factor-criterion}, for which
the inherited bound \(L=b+\alpha+\beta+\rho\) satisfies
\[
 L\le\min\{a+b+1,\alpha+2\beta+\rho-1\}-h.
\]
Equivalently, the factor data satisfy
\begin{equation}\label{eq:h-certificate}
 a-\alpha\ge\beta+\rho+h-1,\qquad \beta-b\ge h+1.
\end{equation}
Let \(\eta_{d,h}\) be the minimum of
\(\max\{f_0(Q),f_{d-1}(Q)\}\) over \(Q\in\mathcal C_{d,h}\), with value
\(+\infty\) if the class is empty. Write \(\eta_d=\eta_{d,1}\).
A \emph{certified presentation} of \(Q\in\mathcal C_{d,h}\) is a triple
\((A,B,F)\) as above whose strict stack is combinatorially equivalent
to \(Q\) and whose data satisfy \eqref{eq:h-certificate}. An
\emph{extremizer} is a member of \(\mathcal C_{d,h}\) whose larger class
size equals \(\eta_{d,h}\).
\end{definition}

These classes specify a certificate, not just a matching number. Every
member has deficiency at least \(h\), but a polytope of deficiency \(h\)
need not have this form or this certificate.

\subsection{The sharper lower bound}

The bounds \(\alpha\ge r+1\) and \(\rho\ge r\) alone give
\(2d+2h+4\); the two additional units come from the structure of
polytopes with few facets. McMullen's classification
\cite[Section~3]{mcmullen1971} (see also
\cite[Section~6.1]{gruenbaum2003} for the dual form), restated with its facet types in
\cite[Lemma~1 and Remark~6]{pinedaTritamaWangYost2026}, says that an
\(r\)-polytope with \(r+2\) facets is combinatorially
\[
 \operatorname{Pyr}^{t}(\Delta_p\times\Delta_q),
 \qquad p,q\ge1,\quad p+q+t=r.
\]
It has \(r+1+pq\) vertices. The facets arise from a facet of either
simplex factor or from omitting one pyramid apex.

\begin{lemma}[Witness-facet floor]\label{lem:witness-facet-floor}
For factor data satisfying \eqref{eq:h-certificate},
\[
 \alpha+\rho\ge2r+3.
\]
\end{lemma}

\begin{proof}
The inequalities imply \(a>\alpha\), so \(A\) is not a simplex and
\(\alpha\ge r+2\). Since \(\rho\ge r\), the claim is immediate unless
\(\alpha=r+2\) and \(\rho=r\); assume this.
Then \(F\) is a simplex. In the classification above, the facet types are
\[
 \operatorname{Pyr}^{t}(\Delta_{p-1}\times\Delta_q),\quad
 \operatorname{Pyr}^{t}(\Delta_p\times\Delta_{q-1}),\quad
 \operatorname{Pyr}^{t-1}(\Delta_p\times\Delta_q),
\]
the last occurring only if \(t>0\). A pyramid over a nonsimplex is not a
simplex, and a product of two positive-dimensional simplices is not a
simplex. Thus a simplex facet forces \(p=1\) or \(q=1\). Consequently
\(pq\le r-1\), so \(a=r+1+pq\le2r\). On the other hand,
\eqref{eq:h-certificate} gives
\(a\ge(r+2)+r+\beta+h-1>2r\), a contradiction.
\end{proof}

\begin{proposition}[Sharp certified size floor]\label{prop:sharp-floor}
Every \(Q\in\mathcal C_{d,h}\) satisfies
\begin{equation}\label{eq:sharp-floor}
 f_0(Q)\ge2d+2h+6,\qquad f_{d-1}(Q)\ge2d+2h+6.
\end{equation}
\end{proposition}

\begin{proof}
As in Proposition~\ref{prop:factor-criterion}, both factors have dimension
at least three. The second inequality of \eqref{eq:h-certificate} excludes
an \(m\)-simplex, hence
\(b\ge m+2\) and \(\beta\ge m+h+3\).
Lemma~\ref{lem:witness-facet-floor} and the first inequality give
\begin{align*}
 a+b+1&\ge(2r+3)+(m+h+3)+(h-1)+(m+2)+1,\\
 \alpha+2\beta+\rho-1&\ge(2r+3)+2(m+h+3)-1.
\end{align*}
Both right-hand sides are \(2d+2h+6\).
\end{proof}

\begin{proposition}[Structure at balanced equality]\label{prop:sharp-equality}
Suppose \(Q\in\mathcal C_{d,h}\) has both class sizes \(2d+2h+6\).
Then
\[
 b=m+2,\quad\beta=m+h+3,\quad\alpha+\rho=2r+3,\quad
 a=\alpha+\rho+\beta+h-1.
\]
Furthermore there are positive integers \(u,v\) with \(uv=h+2\) and
\(u+v\le m\) such that
\[
 B^*\cong\operatorname{Pyr}^{m-u-v}(\Delta_u\times\Delta_v).
\]
In particular, if \(h+2\) is prime, then \(m\ge h+3\) and \(B^*\)
is an iterated pyramid over \(\Delta_1\times\Delta_{h+2}\).
\end{proposition}

\begin{proof}
Equality in the facet bound of Proposition~\ref{prop:sharp-floor} forces
\(\alpha+\rho=2r+3\) and \(\beta=m+h+3\). The bound
\(\beta\ge b+h+1\), together with \(b\ge m+2\), then gives \(b=m+2\).
The vertex equality forces the stated value of \(a\).
Thus \(B^*\) has exactly \(m+2\) facets. By the same McMullen
classification it has the displayed form and \(m+1+uv\) vertices.
Equating this with \(\beta=m+h+3\) gives \(uv=h+2\).
The prime case forces \(\{u,v\}=\{1,h+2\}\).
\end{proof}

\begin{lemma}[Few facets and a simplex facet]
\label{lem:simplex-facet-capacity}
An \(r\)-polytope with \(r+3\) facets and a simplex facet has at most
\[
 2r+\left\lfloor\frac{r^2}{4}\right\rfloor
\]
vertices.
\end{lemma}

\begin{proof}
Let \(D\) be the polar of such a polytope, and let \(x\) correspond to
the simplex facet. Then \(D\) has \(r+3\) vertices and \(x\) belongs to exactly
\(r\) facets. Delete \(x\) and put
\(D'=\conv(\V(D)\setminus\{x\})\).
This polytope is full-dimensional. Otherwise \(D\) would be a pyramid
with apex \(x\); its base would have \(r+2\) vertices but only \(r\)
facets, since \(x\) is simple. An \((r-1)\)-polytope with \(r\)
facets is a simplex, contradicting the vertex count.

Every facet of \(D\) not containing \(x\) is a facet of \(D'\).
Since \(x\notin D'\), the point \(x\) lies strictly beyond at least one
facet of \(D'\), and such a facet is not a facet of \(D\).
The polar of \(D'\) has \(r+2\) facets. McMullen's classification gives
\[
 f_{r-1}(D')\le r+1+\left\lfloor\frac{r^2}{4}\right\rfloor,
\]
since \(pq\le\lfloor r^2/4\rfloor\) whenever \(p+q\le r\).
Adding the \(r\) facets through \(x\) and omitting at least the visible
facet proves the claim.
\end{proof}

\subsection{A matching with a prescribed excess}

For \(h\ge1\), define the arithmetic parameter
\begin{equation}\label{eq:def-sigma}
 \sigma_h=\min\{u+v:u,v\ge1\text{ integers},\ uv=h+2\}.
\end{equation}
In particular \(\sigma_h\le h+3\), with equality if and only if \(h+2\)
is prime: if \(h+2=uv\) with \(2\le u\le v\), then
\(u+v\le2+(h+2)/2<h+3\).

\begin{lemma}[Complementary factors for arbitrary deficiency]
\label{lem:arbitrary-complement}
Let \(u,v\ge1\) satisfy \(uv=h+2\), and let \(m\ge u+v\).
The rational polytope
\[
 B=\left(\operatorname{Pyr}^{m-u-v}
            (\Delta_u\times\Delta_v)\right)^*
\]
has \(b=m+2\) vertices and \(\beta=m+h+3\) facets.
Its nonincidence graph has a matching covering all \(b\) vertices and
leaving \(h+1\) facets unmatched.
\end{lemma}

\begin{proof}
The product has vertices \(w_{ij}\), \(0\le i\le u\), \(0\le j\le v\).
Its \(u+1\) row facets \(H_i\) omit row \(i\), and its \(v+1\)
column facets \(V_j\) omit column \(j\). The following distinct pairs
are nonincidences:
\[
 \begin{aligned}
 H_0&\leftrightarrow w_{00},& H_1&\leftrightarrow w_{11},&
 H_i&\leftrightarrow w_{i0}\quad(2\le i\le u),\\
 V_0&\leftrightarrow w_{10},&&&
 V_j&\leftrightarrow w_{0j}\quad(1\le j\le v).
 \end{aligned}
\]
Empty index ranges contribute no pairs. They cover all \(u+v+2\) facets
and leave \(uv-1=h+1\) vertices.
Each pyramid adds one vertex, one facet and an isolated apex--base
nonincidence edge. After \(m-u-v\) pyramids the counts are
\(m+h+3\) vertices and \(m+2\) facets. Polarity exchanges the two classes
and their matchings. A rational polar is obtained after translating a
rational interior point to the origin.
\end{proof}

\subsection{A controlled cut attains the floor}

Simple-vertex truncation is a classical source of few-facet polytopes;
see \cite{mcmullen1971} and the constructions in
\cite[Section~4]{pinedaUgonYost2019} and
\cite[Section~1]{pinedaTritamaWangYost2026}; truncations of pyramids
over \(\Delta_2\times\Delta_{s-2}\) also appear in Section~6(2) of
Pineda-Villavicencio and Wang's preprint arXiv:2512.07456v1 (2025),
identified in the introduction.
Here the simplex facet created by the truncation permits a sharp Hall
certificate, and a direct matching determines the deficiency exactly.
Letting the cutting hyperplane pass through \(e\) chosen neighbors of
the cut vertex lowers the vertex count by exactly \(e\), which a cut
close to the vertex cannot achieve.

\begin{lemma}[Controlled product truncation]\label{lem:controlled-cut}
Let \(p,q\ge2\), \(t\ge0\), and \(0\le e\le p+q-2\) be integers,
and put \(r=p+q+t\). There is a rational \(r\)-polytope \(A\) with a
simplex facet \(F\) such that
\begin{equation}\label{eq:sharp-witness-data}
 a=2r+pq-e,\qquad \alpha=r+3,\qquad \rho=r.
\end{equation}
The facet \(F\) has vertices \(u_1,\ldots,u_r\) that can be matched
to \(r\) old facets of \(A\); the other two old facets can be matched
to distinct vertices outside \(F\).
\end{lemma}

\begin{proof}
Use a rational product realization of
\(E=\operatorname{Pyr}^{t}(\Delta_p\times\Delta_q)\), with product
vertices \(a_{ij}\). The original vertex \(v=a_{00}\) is simple.
Its \(r\) neighbors are the \(p+q\) product neighbors
\(a_{i0},a_{0j}\) and the \(t\) pyramid apices. The edge directions
at \(v\) are linearly independent.
Choose \(e\) product neighbors, excluding \(a_{01}\) and \(a_{10}\).
There are \(p+q-2\) eligible neighbors. Define a rational affine
functional \(\lambda\) by giving it value zero at \(v\), value one
at the chosen neighbors, and value two at every other neighbor; it is
well defined because \(v\) and its \(r\) neighbors are affinely
independent. For \(i,j>0\), the relation
\(a_{ij}+a_{00}=a_{i0}+a_{0j}\) gives
\[
 \lambda(a_{ij})=\lambda(a_{i0})+\lambda(a_{0j})\ge2.
\]
Thus \(v\) is the only old vertex with \(\lambda<1\), and the
chosen \(e\) neighbors are the only old vertices with \(\lambda=1\).
Set \(A=E\cap\{\lambda\ge1\}\) and \(F=E\cap\{\lambda=1\}\).

The \(r\) edges at \(v\) meet this hyperplane in \(r\) affinely
independent points: \(e\) are old endpoints and \(r-e\) are new
vertices. All other old vertices remain. Every old facet survives:
since \(p,q\ge2\), a facet through \(v\) contains a rectangle
\(v,a_{i0},a_{0j},a_{ij}\) with \(i,j>0\), so its vertices other than
\(v\) still affinely span it, also after pyramiding.
Facets not through \(v\) are unchanged. Consequently
\[
 f_0(A)=(r+1+pq)-1+(r-e)=2r+pq-e,
 \qquad f_{r-1}(A)=r+3,
\]
and \(F\) is a simplex with \(r\) ridges.

Index the old facets through \(v\) as \(K_1,\ldots,K_r\), so the
\(i\)-th edge misses \(K_i\). Its cut point \(u_i\), even when
it is an old endpoint, misses \(K_i\). These give \(r\) distinct
nonincidences. The remaining old facets are the extensions of the
row facet \(H_0\) and column facet \(V_0\). Pair them with
\(a_{01}\) and \(a_{10}\), respectively. Both protected vertices
have \(\lambda=2\), so lie outside \(F\).
\end{proof}

\begin{theorem}[Sharp prescribed-deficiency constructions]
\label{thm:sharp-construction}
Let \(h\ge1\), \(p,q\ge2\), \(t\ge0\), and
\(0\le e\le p+q-2\) be integers with
\(pq-e\ge2h+5+\sigma_h\). Put
\[
 r=p+q+t,\qquad m=pq-e-2h-5,\qquad d=r+m+1.
\]
There is a rational \(Q\in\mathcal C_{d,h}\) such that
\begin{equation}\label{eq:sharp-counts}
 f_0(Q)=f_{d-1}(Q)=2d+2h+6,\qquad
 \nu(\N(Q))=2d+h+6.
\end{equation}
Its nonincidence graph is connected, its deficiency is exactly \(h\),
and its face lattice has no nontrivial free-join decomposition.
\end{theorem}

The construction matches the lower bound by separating three tasks:
the simplex-facet vertices cover the old witness facets, vertices outside the
new simplex cover the cone facets, and the complementary factor supplies
the remaining matching edges.

\begin{proof}[Proof of Theorem~\ref{thm:sharp-construction}]
\emph{Step 1: the witness.}
Use \(A,F\) from Lemma~\ref{lem:controlled-cut}, with the data in
\eqref{eq:sharp-witness-data}. The cut replaces one vertex by a
simplex facet, with \(e\) vertices already present in the original
polytope.

\emph{Step 2: stacking and the upper matching bound.}
Use the polytope \(B\) of Lemma~\ref{lem:arbitrary-complement}; it exists
using a minimizing factorization in \eqref{eq:def-sigma}, since
\(m\ge\sigma_h\). Stack a rational apex \(z\) strictly beyond \(F*B\)
and strictly beneath all other facets of \(A*B\). Equations
\eqref{eq:general-join-stack-counts} and \eqref{eq:sharp-witness-data}
give the two equal counts in \eqref{eq:sharp-counts}. Moreover,
\[
 a-\alpha-(\beta+\rho+h-1)
   =pq-e-m-2h-5=0,\qquad \beta-b=h+1.
\]
Hence \(Q\in\mathcal C_{d,h}\), and its inherited bound is
\(\nu(\N(Q))\le f_0(Q)-h\).

\emph{Step 3: a matching attaining the bound.}
Lift the \(r+2\) nonincidences in Lemma~\ref{lem:controlled-cut}
to the retained join facets induced by \(A\). This uses all \(r\)
vertices of \(F\) and two vertices outside \(F\).

There are \(r+\beta\) new cone facets. Every vertex of \(A\)
outside \(F\) misses every such facet: the old part of each cone lies
in a ridge of \(F*B\). After the two preceding pairs, the number of
available vertices outside \(F\) is
\[
 a-r-2=r+pq-e-2=r+\beta+h\ge r+\beta.
\]
Use any \(r+\beta\) of them to cover the cone facets.

Finally, the matching of Lemma~\ref{lem:arbitrary-complement} pairs all
\(b\) vertices of \(B\) with retained facets of the form \(A*G\) and
leaves \(h+1\) such facets unused. The apex \(z\), being strictly
beneath every retained facet, can be paired with one of them.
The three blocks use disjoint vertices and disjoint facet types. Their
total size is
\[
 (r+2)+(r+\beta)+(b+1)=2r+\beta+b+3=f_0(Q)-h.
\]
This attains the upper bound. Strict stacking gives connectedness by
Corollary~\ref{cor:stacking-transfer}, and
Lemma~\ref{lem:join-disconnects} excludes a nontrivial free join.
Rational truncation, polarity and a rational point in the open strict
stacking region establish rationality throughout.
\end{proof}

\begin{theorem}[Exact onset of the extremal floor]\label{cor:stable-h}
For \(h\ge1\), set
\[
 K_h=2h+5+\sigma_h,\qquad
 R_h=\left\lceil2\sqrt{K_h}\right\rceil,\qquad
 D_h=\sigma_h+R_h+1.
\]
Then
\begin{equation}\label{eq:stable-h}
 \eta_{d,h}=2d+2h+6\quad\Longleftrightarrow\quad d\ge D_h.
\end{equation}
For every \(d\ge D_h\), a balanced rational polytope with connected
nonincidence graph and deficiency exactly \(h\) attains the minimum.
Moreover \(D_h=\sigma_h+O(\sqrt h)\); in particular
\(D_h\le h+O(\sqrt h)\) as \(h\to\infty\).
\end{theorem}

\begin{proof}
First suppose that the floor is attained. Both class sizes then equal
the floor by Proposition~\ref{prop:sharp-floor}. By
Proposition~\ref{prop:sharp-equality}, \(m\ge\sigma_h\) and
\[
 a=2r+m+2h+5,\qquad \alpha+\rho=2r+3.
\]
Since \(A\) is not a simplex and \(F\) is an \((r-1)\)-polytope,
\(\alpha\ge r+2\) and \(\rho\ge r\); hence there are just two
possibilities. If \(\alpha=r+3\) and \(\rho=r\),
then \(F\) is a simplex and Lemma~\ref{lem:simplex-facet-capacity}
gives \(m+2h+5\le\lfloor r^2/4\rfloor\). If
\(\alpha=r+2\) and \(\rho=r+1\), McMullen's classification gives
\(a\le r+1+\lfloor r^2/4\rfloor\), and hence the stronger inequality
\(r+m+2h+4\le\lfloor r^2/4\rfloor\). In both cases
\[
 r\ge\left\lceil2\sqrt{m+2h+5}\right\rceil\ge R_h,
 \qquad d=r+m+1\ge D_h.
\]

Conversely, take \(d\ge D_h\), and set
\[
 p=\left\lfloor R_h/2\right\rfloor,\quad
 q=\left\lceil R_h/2\right\rceil,\quad
 e=pq-K_h,\quad t=d-D_h.
\]
Then \(p+q=R_h\), \(pq=\lfloor R_h^2/4\rfloor\), and \(e\ge0\).
Since \(R_h-1<2\sqrt{K_h}\), we have
\(\lfloor(R_h-1)^2/4\rfloor\le K_h-1\), and therefore
\[
 e\le\left\lfloor R_h^2/4\right\rfloor
       -\left\lfloor(R_h-1)^2/4\right\rfloor-1
   =\left\lfloor R_h/2\right\rfloor-1\le p+q-2,
\]
using \(\lfloor n^2/4\rfloor-\lfloor(n-1)^2/4\rfloor=\lfloor n/2\rfloor\).
Every factorization \(uv=h+2\ge3\) has \(u+v\ge4\), so \(\sigma_h\ge4\),
hence \(K_h\ge11\) and \(R_h\ge7\), so
\(p\ge3\) and \(q\ge4\). Use
Theorem~\ref{thm:sharp-construction} with \(m=\sigma_h\).
Its dimensions are \(r=p+q+t\) and \(r+m+1=d\), and it attains
the floor with connected nonincidence graph and deficiency exactly
\(h\). Finally
\(K_h\le3h+8\) gives \(R_h=O(\sqrt h)\).
\end{proof}

\begin{corollary}[Quadratic deficiency in linear dimension]
\label{cor:quadratic-deficiency}
For every integer \(k\ge2\), set \(h=k^2-2\). The size floor of
\(\mathcal C_{d,h}\) is first attained in dimension
\[
 2k+1+\left\lceil2\sqrt{2k^2+2k+1}\right\rceil
\]
by polytopes from Theorem~\ref{thm:sharp-construction}, which have
matching deficiency \(k^2-2\), connected nonincidence graph and
optimal certified size.
\end{corollary}

\begin{proof}
The factorization \(h+2=k\cdot k\) minimizes the sum by the
arithmetic--geometric mean inequality, so \(\sigma_h=2k\).
Substitute into the exact threshold of Theorem~\ref{cor:stable-h}.
\end{proof}

\begin{corollary}[Deficiency-one optimum]\label{cor:stable-one}
For every \(d\ge12\),
\[
 \eta_d=2d+8.
\]
Both class sizes equal \(2d+8\) in an attaining rational construction,
whose nonincidence graph is connected and has matching number \(2d+7\).
\end{corollary}

\begin{proof}
Here \(\sigma_1=4\), \(K_1=11\), \(R_1=7\), and \(D_1=12\).
Apply Theorem~\ref{cor:stable-h}. Concretely, the attaining parameters
are \(p=3\), \(q=4\), \(e=1\), \(m=4\), and \(t=d-12\).
\end{proof}

\begin{corollary}[An eleven-dimensional benchmark]\label{cor:eleven}
There is a rational member of \(\mathcal C_{11,1}\) with \(32\) vertices,
\(32\) facets, connected nonincidence graph and matching number \(31\).
Consequently \(31\le\eta_{11}\le32\).
\end{corollary}

\begin{proof}
Use Lemma~\ref{lem:controlled-cut} with \(p=3\), \(q=4\), \(e=1\)
and \(t=0\). It gives \(r=7\), \(a=25\), \(\alpha=10\) and
\(\rho=7\). Take \(B\) to be the octahedron, with \(b=6\) and
\(\beta=8\), and stack beyond \(F*B\). The resulting dimension is
\(7+3+1=11\), the two class sizes are \(32\), and the inherited
bound is \(31\). The matching assembly of
Theorem~\ref{thm:sharp-construction} covers \(9\) retained witness
facets, \(15\) cone facets, and \(7\) vertices from \(B\) and the
apex. For the cones there are \(a-r-2=16\) available witness vertices.
The octahedral matching in Theorem~\ref{thm:balanced} leaves two
retained complementary facets, one of which receives the apex.
Thus the matching has \(31\) edges, and strict stacking gives
connectedness. The lower bound follows from the floor \(30\) and
its nonattainment in dimension eleven in Theorem~\ref{cor:stable-h}.
\end{proof}

No member of the certified class attains \(2d+8\) below dimension twelve.
Section~\ref{sec:lowdim} determines the remaining values with different
witness geometries. In particular, a three-factor product improves the
eleven-dimensional benchmark above from thirty-two to the optimum
thirty-one.

\subsection{Seven-dimensional extrema for every deficiency}
\label{subsec:seven-deficiency}

The first admissible dimension has a different extremal geometry from
the stable range. Euler's inequalities determine both factors, and an
inspection of the three-faces gives the exact deficiency of every
extremizer.

\begin{theorem}[Seven-dimensional certified extremizers]
\label{thm:seven-deficiency}
For every integer \(h\ge1\),
\begin{equation}\label{eq:seven-deficiency}
 \eta_{7,h}=7h+26.
\end{equation}
Every \(Q\in\mathcal C_{7,h}\) has at least \(7h+26\) vertices
and at least \(7h+26\) facets. For a certified presentation \((A,B,F)\)
of \(Q\), the larger class size of \(Q\) equals \(7h+26\) if and only if
\(F\) is triangular, \(A\) is a simple three-polytope with counts
\((6h+20,3h+12)\), and \(B\) is a simplicial three-polytope with
counts \((h+5,2h+6)\).
Every extremizer is balanced, has connected nonincidence graph,
and satisfies
\[
 \delta(Q)=h,\qquad \nu(\N(Q))=6h+26.
\]
The copy of \(A\) is the unique face \(\sigma\) of \(Q\) with
\(f_0(\sigma)+\kappa_Q(\sigma)=H(Q)\).
Rational extremizers exist for every \(h\).
\end{theorem}

\begin{proof}[Proof of Theorem~\ref{thm:seven-deficiency}]
\emph{Lower bound.}
Both factor dimensions are at least three by
Proposition~\ref{prop:factor-criterion}, so \(r=m=3\).
Euler's inequalities give \(\beta\le2b-4\) and
\(a\le2\alpha-4\), with equality precisely for simplicial \(B\)
and simple \(A\), respectively. Since \(F\) is a polygon,
\(\rho\ge3\). By \eqref{eq:h-certificate},
\[
 b\ge h+5,\qquad \beta\ge2h+6,\qquad
 \alpha\ge\beta+\rho+h+3\ge3h+12.
\]
The first certificate inequality further gives
\(a\ge\alpha+\beta+\rho+h-1\ge6h+20\). Consequently
\[
 a+b+1\ge7h+26,\qquad \alpha+2\beta+\rho-1\ge7h+26.
\]
If the larger count is \(7h+26\), both inequalities are equalities.
Equality in the facet count forces
\(\rho=3\), \(\alpha=3h+12\), and \(\beta=2h+6\);
equality in the vertex count gives \(a=6h+20\), \(b=h+5\).
Euler equality then gives precisely the asserted factor conditions.
Conversely, those conditions give both certificate inequalities and
balanced counts \(7h+26\).

\emph{Equality structure and exact deficiency.}
Assume these factor conditions and write \(p\) for the stacking point.
Let \(\tau\) be a proper face with \(q\) facets of its own. Each facet
of \(\tau\) is \(\tau\cap G\) for some facet \(G\) of \(Q\) not
containing \(\tau\), and distinct facets of \(\tau\) arise from distinct
such \(G\). Thus
\[
 f_6(Q)\ge q+\kappa_Q(\tau).
\]
Call \(f_0(\tau)+\kappa_Q(\tau)\) the \emph{Hall value} of \(\tau\).
In particular, if \(q\ge f_0(\tau)\), its Hall value is at most
\(f_6(Q)\). Applied to the dual face \(\tau^\Delta\) in a polar of \(Q\), which has
the same Hall value, the same inequality bounds the Hall value by
\(f_0(Q)\) whenever \(\tau^\Delta\) has at least as many facets as
vertices. This is the
pointwise inequality in Jahn and Winter's proof of
\cite[Theorem~7]{jahnWinter2021}.
Except for three-faces, every nonempty proper face of a seven-polytope
has itself or its dual of dimension at most two, hence equally many
vertices and facets (for a vertex we use the convention that its only
facet is the empty face). The empty and full faces have Hall values
\(f_6(Q)\) and \(f_0(Q)\).

It remains to examine the three-faces. By the beneath--beyond
description \cite[Section~5.2]{gruenbaum2003}, a face of \(Q\) not
containing \(p\) is a face of \(A*B\). Such a three-face is \(A\), \(B\),
a pyramid over a polygonal facet of \(A\), the join of two edges,
or a pyramid over a triangular facet of \(B\).
Other than \(A\), each has at least as many own facets as vertices:
this holds for simplicial \(B\), a polygonal pyramid has equal
counts, and the remaining types are tetrahedra.

A three-face containing \(p\) is a pyramid over a two-face of
\(F*B\): every face of \(Q\) containing \(p\) is an intersection of
cone facets \(\conv(p,R)\) from Lemma~\ref{lem:stacking}, hence equals
\(\conv(p,G)\) for a proper face \(G\) of \(F*B\)
\cite[Section~5.2]{gruenbaum2003}. Every two-face of this join is
triangular: it is \(F\), an edge of \(F\) joined to a vertex of
\(B\), a vertex of \(F\) joined to an edge of \(B\), or a facet
of \(B\). Hence every new three-face is a tetrahedron.
The pointwise inequality therefore bounds the Hall value of every
face other than \(A\) by \(7h+26\).
The surviving factor \(A\) lies in exactly \(\beta\) facets,
as in Proposition~\ref{prop:factor-criterion}, so
\[
 f_0(A)+\kappa_Q(A)=a+\beta=8h+26.
\]
It is the unique maximizing face. The exact face formula
(Theorem~\ref{thm:exact-defect}) gives deficiency \(h\) and
matching number \(6h+26\). Connectedness follows from
Corollary~\ref{cor:stacking-transfer}.

\emph{Rational attainment.}
Take a rational prism over a \((3h+9)\)-gon and strictly truncate
one vertex. A simple three-dimensional vertex cut adds two vertices
and one facet, and its new facet \(F\) is triangular
\cite{ziegler1995}. The resulting simple \(A\) has counts
\((6h+20,3h+12)\). Take \(B\) to be a rational bipyramid over
an \((h+3)\)-gon; it is simplicial with counts \((h+5,2h+6)\).
The strict stacking region beyond \(F*B\) is nonempty and open,
so it contains a rational point. This produces a rational extremizer
with all the asserted properties and proves \eqref{eq:seven-deficiency}.
\end{proof}

\subsection{Eight-dimensional extrema and four-dimensional capacity}
\label{subsec:eight-deficiency}

The next dimension replaces Euler equality for the witness by a
four-dimensional capacity bound. Deleting the vertex dual to the
stacking facet isolates that bound; the sufficiency part of the
\(g\)-theorem then supplies the required face numbers.
Bounds for faces through a distinguished vertex appear in
\cite[Theorem~2]{billeraLee1980}; the elementary deletion estimates
below suffice for both the size bound and its equality analysis.

\begin{lemma}[A distinguished vertex in dimension four]
\label{lem:four-vertex-capacity}
Let \(D\) be a four-polytope with \(n\) vertices and a vertex \(x\)
of degree \(k\). If
\(D'=\conv(\V(D)\setminus\{x\})\) is four-dimensional, then
\[
 f_3(D)\le\frac{(n-1)(n-4)}2+2k-5.
\]
If \(k>4\), then
\(f_3(D)\le\frac{(n-1)(n-4)}2+2k-6\).
\end{lemma}

\begin{proof}
For any four-polytope \(E\), every facet has at least four ridges,
and every ridge belongs to two facets. Euler's equation gives
\[
 f_2(E)\ge2f_3(E),\qquad
 f_3(E)=f_2(E)-f_1(E)+f_0(E),
\]
so
\(f_3(E)\le f_1(E)-f_0(E)\le f_0(E)(f_0(E)-3)/2\).
Equality forces every facet to have four ridges, hence to be a
tetrahedron, and every vertex pair to be an edge.
Apply this to \(D'\). The three-dimensional vertex figure at \(x\)
has \(k\) vertices and at most \(2k-4\) facets. Thus at most
\(2k-4\) facets of \(D\) contain \(x\).
Every facet not containing \(x\) survives in \(D'\), whereas at
least one facet of \(D'\) is strictly visible from \(x\), and such a
facet is not a facet of \(D\). This proves the first bound.

Suppose that equality holds in the first bound. Then \(D'\)
attains its facet bound and exactly one of its facets, say \(T\),
is strictly visible from \(x\). In particular \(T\) is tetrahedral.
All other supporting halfspaces of \(D'\) contain \(x\), so
\[
 D=D'\cup\conv(x,T).
\]
Indeed a segment from \(x\) to \(D'\) first crosses the boundary
through \(T\). Since \(x\notin D'\), near \(x\) the polytope \(D\)
coincides with the simplex \(\conv(x,T)\), even if \(x\) lies on
another supporting hyperplane of \(D'\). Hence every edge of \(D\) at
\(x\) points towards a vertex of \(T\); that vertex lies on the edge
and is a vertex of \(D\), so it is the other endpoint. Therefore
\(k\le4\), proving the strengthened bound.
\end{proof}

\begin{proposition}[Matching saturation at a tetrahedral facet]
\label{prop:tetrahedral-saturation}
Let \(A\) be a four-polytope with \(a\) vertices, \(\alpha\)
facets and a tetrahedral facet \(F\). Let \(B\) have \(b\)
vertices and \(\beta\) facets, and suppose its nonincidence graph
has a matching covering all its vertices. Assume
\[
 \beta\ge b+1,\qquad a\ge3\alpha-7,\qquad
 a\ge\alpha+\beta+3.
\]
If \(Q\) is obtained by strict stacking beyond \(F*B\) in
\(A*B\), then
\[
 \nu(\N(Q))=\alpha+\beta+b+4,\qquad
 \delta(Q)=\min\{a-\alpha-\beta-3,\,\beta-b-1\}.
\]
\end{proposition}

\begin{proof}[Proof of Proposition~\ref{prop:tetrahedral-saturation}]
Each of the four triangular ridges of \(F\) belongs to a unique
other facet of \(A\). These facets are distinct, since two facets of
a polytope share at most one ridge, and each omits
the opposite vertex of \(F\). Match them to those four vertices.
Any remaining facet \(H\) of \(A\) is a three-polytope with at
most \(\alpha-1\) facets: each ridge of \(H\) is its intersection
with a distinct other facet of \(A\). Euler's inequality gives
\[
 f_0(H)\le2f_2(H)-4\le2\alpha-6.
\]
Thus \(H\) omits at least \(a-2\alpha+2\ge\alpha-5\)
vertices outside \(F\). Greedily match the remaining
\(\alpha-5\) retained facets of type \(H*B\) to distinct such
vertices. This covers all \(\alpha-1\) retained facets of type \(H*B\),
leaving \(a-\alpha+1\ge\beta+4\) unused outside vertices.
Every new cone facet omits every vertex of \(A\) outside \(F\),
so use these vertices to cover the \(\beta+4\) cone facets.
Lift the assumed matching of \(B\) to \(b\) distinct retained
facets of type \(A*G\). Since \(\beta>b\), an unused facet of
this type remains for the new stacking point.
The three blocks are disjoint and contain
\[
 (\alpha-1)+(\beta+4)+(b+1)=\alpha+\beta+b+4
\]
edges. Conversely, the neighbors of the \(a\) vertices of \(A\)
are precisely the retained facets of type \(H*B\) and the new cone
facets,
of total number \(\alpha+\beta+3\). At most \(b+1\) further
edges can use vertices outside \(A\), proving the upper bound.
Finally \(Q\) has \(a+b+1\) vertices and \(\alpha+2\beta+3\)
facets. Subtracting the matching number from the smaller count
gives the displayed deficiency.
\end{proof}

\begin{theorem}[Eight-dimensional prescribed-deficiency optimum]
\label{thm:eight-deficiency}
For an integer \(h\ge1\), put
\[
 t_h=\left\lceil\frac{7+\sqrt{24h+81}}2\right\rceil.
\]
Then
\begin{equation}\label{eq:eight-deficiency}
 \eta_{8,h}=4h+15+t_h.
\end{equation}
Every member of \(\mathcal C_{8,h}\) has at least this many vertices
and facets. A certified presentation attains this larger class size if
and only if \(r=4\), \(m=3\), \(F\) is a tetrahedron, and the
factor counts are
\[
 (a,\alpha)=(t_h+3h+9,t_h),\qquad
 (b,\beta)=(h+5,2h+6).
\]
In particular \(B\) is simplicial, and every extremizer is balanced,
has connected nonincidence graph, has deficiency exactly \(h\), and has matching number
\(3h+15+t_h\). There are rational connected attainers with \(A\)
simple.
\end{theorem}

\begin{proof}[Proof of Theorem~\ref{thm:eight-deficiency}]
Write \(t=t_h\), \(N=f_0(Q)\), and \(M=f_7(Q)\).
Since the product of the roots of \(t^2-7t-6h-8\) is negative,
the integer \(t\) is the smallest positive integer with
\[
 t^2-7t-6h-8\ge0.
\]
Both factors have dimension at least three, so there are two
dimension splits.

\emph{Lower bound for a four-dimensional witness.}
Suppose \(r=4\), \(m=3\), and set \(u=\alpha+\rho\).
Euler's inequality for \(B\) and \eqref{eq:h-certificate} give
\[
 b\ge h+5,\quad \beta\ge2h+6,\quad
 a\ge u+\beta+h-1,\quad N,M\ge u+4h+11.
\]
Let \(D=A^*\) and let \(x\) correspond to \(F\).
Then \(D\) has \(\alpha\) vertices and \(a\) facets, and
\(x\) has degree \(\rho\). Its deleted hull \(D'\) is
four-dimensional. Otherwise \(D\) is a pyramid with apex \(x\),
\(\rho=\alpha-1\), and \(a\le2\alpha-5\) by Euler's inequality
for the base. This contradicts
\(a\ge\alpha+\beta+\rho+h-1\).
If \(\rho=4\), Lemma~\ref{lem:four-vertex-capacity} gives
\[
 a\le\frac{\alpha^2-5\alpha+10}{2},\qquad
 a\ge\alpha+3h+9,
\]
hence \(\alpha\ge t\) and \(u\ge t+4\).
If \(\rho\ge5\), the general four-dimensional bound gives
\[
 u+3h+5\le a\le\frac{\alpha(\alpha-3)}2
                   \le\frac{(u-5)(u-8)}2.
\]
Thus \(u^2-15u+30-6h\ge0\). Here \(\alpha\ge6\) because \(A\) is not
a simplex, so \(u-4\ge7\). Since
\((u-4)^2-7(u-4)-6h-8=u^2-15u+36-6h>0\), the minimality of \(t\)
gives \(u\ge t+4\). In either case \(N,M\ge4h+15+t\).

\emph{A three-dimensional witness.}
Suppose \(r=3\), \(m=4\). The nonsimplex factor \(B\) has
\(b\ge6\), since a four-polytope with five vertices is a simplex, and
\[
 \beta\le\frac{b(b-3)}2,\qquad
 h\le\frac{b^2-5b-2}2.
\]
Euler's inequality for \(A\), together with \(\rho\ge3\), gives
\(\alpha\ge\beta+\rho+h+3\),
\(a\ge2\beta+2\rho+2h+2\), and therefore
\[
 N,M\ge3b+4h+11.
\]
For \(s=3b-5\) the preceding bound on \(h\) yields
\[
 s^2-7s-6h-8\ge6(b-3)^2+4>0.
\]
Since \(s\ge13\), it follows that \(t\le3b-5\).
Both counts are consequently at least \(4h+16+t\).

\emph{Equality factor data.}
At the stated optimum the first split must hold, and the two count
bounds force
\[
 u=t+4,\quad \beta=2h+6,\quad b=h+5,\quad a=t+3h+9.
\]
If \(\rho\ge5\), the strengthened part of
Lemma~\ref{lem:four-vertex-capacity} would require
\[
 \alpha^2-7\alpha+2\rho-6h-18\ge0.
\]
At \(\rho=5\), \(\alpha=t-1\), its left-hand side is
\((t-1)^2-7(t-1)-6h-8<0\). Raising \(\rho\) by one and lowering
\(\alpha\) by one changes that expression by \(10-2\alpha<0\), since
\(\alpha\ge6\) throughout the admissible range. Hence the required
inequality fails for every \(\rho\ge5\). Thus \(\rho=4\), \(\alpha=t\), and \(F\) is
tetrahedral. Euler equality makes \(B\) simplicial.
Conversely these factor data give both certificate inequalities and
balanced size \(4h+15+t\).

\emph{Rational attainment.}
We have \(t\ge9\) and \(t\le h+8\), the latter because
\((h+8)^2-7(h+8)-6h-8=h(h+3)>0\).
Set \(q=3h+19-2t\). Then
\[
 0\le q\le\binom{t-5}{2};
\]
the upper inequality is exactly \(t^2-7t-6h-8\ge0\).
The sufficiency theorem of Billera and Lee
\cite[Theorem~1]{billeraLee1980} provides a simplicial four-polytope
\(E\) with \(t-1\) vertices and truncated \(g\)-vector
\((1,t-6,q)\). Indeed these are the degree-zero, degree-one and
degree-two dimensions of a monomial order ideal: choose any \(q\)
of the \(\binom{t-5}{2}\) quadratic monomials in \(t-6\) variables.
Its \(h\)-vector is
\((1,t-5,t-5+q,t-5,1)\), giving
\[
 f_3(E)=3(t-1)-10+q=t+3h+6.
\]
A sufficiently small rational perturbation preserves a simplicial
polytope's face lattice, so \(E\) may be chosen rational.
Strict stacking on one tetrahedral facet gives a simplicial
four-polytope with \(t\) vertices and \(t+3h+9\) facets.
Its polar \(A\) is simple and has the required counts and a
tetrahedral facet \(F\). Equivalently \(A\) is obtained by a strict
simple-vertex cut of \(E^*\).
Let \(B\) be a rational bipyramid over an \((h+3)\)-gon and
strictly stack beyond \(F*B\). All choices, including polarity
after rational interior translation and the final stacking point,
can be made rational. Corollary~\ref{cor:stacking-transfer} gives
connectedness.

\emph{Matching saturation for every extremizer.}
For any certified presentation of an extremizer, \(B\) is simplicial and
has a matching covering its vertices
\cite[Corollary~9]{jahnWinter2021}. Its counts satisfy
\[
 a-\alpha-\beta-3=\beta-b-1=h,\qquad
 a-(3\alpha-7)=3h+16-2t\ge h>0,
\]
where the last inequality uses \(t\le h+8\).
Proposition~\ref{prop:tetrahedral-saturation} therefore gives
\(\nu(\N(Q))=3h+15+t\) and \(\delta(Q)=h\) for every
extremizer, including those whose witness \(A\) is not simple.
\end{proof}

\subsection{Nine-dimensional extrema and coupled capacities}
\label{subsec:nine-deficiency}

In dimension nine both factors can be four-dimensional. Their two
quadratic capacity bounds now interact: the smallest complementary
vertex count sets the capacity required of the witness. The result
involves two nested integer thresholds, and the leading coefficient
of \(h\) drops from four to two.

\begin{theorem}[Nine-dimensional prescribed-deficiency optimum]
\label{thm:nine-deficiency}
For an integer \(h\ge1\), put
\[
 c_h=\left\lceil\frac{5+\sqrt{8h+33}}2\right\rceil,
 \qquad
 s_h=\left\lceil\frac{7+\sqrt{16h+8c_h+41}}2\right\rceil.
\]
Then
\begin{equation}\label{eq:nine-deficiency}
 \eta_{9,h}=2h+5+2c_h+s_h.
\end{equation}
Every member of \(\mathcal C_{9,h}\) has at least this many vertices
and facets. A certified presentation attains this larger class size if
and only if \(r=m=4\), \(F\) is a tetrahedron, and
\[
 (a,\alpha)=(s_h+c_h+2h+4,s_h),\qquad
 (b,\beta)=(c_h,c_h+h+1).
\]
Every extremizer is balanced, has connected nonincidence graph,
has deficiency exactly \(h\), and has matching number
\(h+5+2c_h+s_h\). There are rational attainers with \(A\) simple
and \(B\) simplicial. In particular, as \(h\to\infty\),
\[
 \eta_{9,h}=2h+(2+2\sqrt2)\sqrt h+O(1).
\]
\end{theorem}

\begin{proof}[Proof of Theorem~\ref{thm:nine-deficiency}]
Write \(c=c_h\), \(s=s_h\), and \(n_0=2h+5+2c+s\).
Since in each quadratic below the product of the roots is negative,
the thresholds are the smallest positive integers satisfying
\begin{equation}\label{eq:nine-thresholds}
 c^2-5c-2h-2\ge0,\qquad
 s^2-7s-2c-4h+2\ge0.
\end{equation}
In particular \(c\ge6\) and \(s\ge9\). We will use
\begin{equation}\label{eq:nine-margins}
 2c\le h+11,\quad s\le h+8,
 \qquad s\le h+7\quad(h\ge2).
\end{equation}
For the first inequality set \(k_0=\lfloor(h+11)/2\rfloor\).
The first polynomial in \eqref{eq:nine-thresholds} at \(k_0\) is
\((h+1)(h+3)/4\) for odd \(h\) and \((h-2)(h+4)/4\) for even
\(h\), so it is nonnegative. The second polynomial at \(h+8\)
is at least \(h^2+4h-1>0\), using \(2c\le h+11\).
For \(h\ge2\), the first polynomial at \(h+4\) is
\((h-2)(h+3)\), so \(c\le h+4\). The second polynomial at
\(h+7\) is then at least \((h-2)(h+3)\).

Let \(N=f_0(Q)\), \(M=f_8(Q)\), and \(u=\alpha+\rho\).
Both factor dimensions are at least three, so only the splits
\(4+4\), \(3+5\), and \(5+3\) are possible.

\emph{The two four-dimensional factors.}
Suppose \(r=m=4\), and write \(\beta=b+h+1+z\), where
\(z\ge0\). The general four-dimensional bound
\(\beta\le b(b-3)/2\) gives \(b\ge c\).
The certificate inequalities yield
\[
 a\ge u+b+2h+z,\qquad
 N\ge u+2b+2h+z+1,\qquad
 M\ge u+2b+2h+2z+1.
\]
The polar \(D=A^*\) has \(\alpha\) vertices and \(a\) facets;
the vertex \(x\) corresponding to \(F\) has degree \(\rho\).
Its deleted hull is four-dimensional. Otherwise \(D\) is a pyramid
with apex \(x\), \(\rho=\alpha-1\), and Euler's inequality for
its base gives \(a\le2\alpha-5\), contradicting the certificate.
If \(\rho=4\), Lemma~\ref{lem:four-vertex-capacity} gives
\[
 a\le\frac{\alpha^2-5\alpha+10}{2},\qquad
 a\ge\alpha+b+2h+4.
\]
Thus \(\alpha\ge s\) and \(u\ge s+4\).
If \(\rho\ge5\), then \(\alpha\ge6\), \(u\ge11\), and
\[
 u+b+2h\le a\le\frac{\alpha(\alpha-3)}2
                  \le\frac{(u-5)(u-8)}2.
\]
Consequently
\[
 (u-4)^2-7(u-4)-2b-4h-4\ge0.
\]
This implies the second inequality of \eqref{eq:nine-thresholds}
with \(s\) replaced by \(u-4\); hence again \(u\ge s+4\).
It follows that both \(N\) and \(M\) are at least \(n_0\).

\emph{Excluding the other dimension splits.}
If \(r=3\), \(m=5\), Euler's inequality for \(A\) gives
\[
 \alpha\ge\beta+\rho+h+3,\qquad
 a\ge2\beta+2\rho+2h+2.
\]
Since \(B\) is not a simplex, \(b\ge7\); moreover \(\beta\ge b+h+1\)
and \(\rho\ge3\). Hence
both \(N,M\ge3b+4h+11\ge4h+32\).
By \eqref{eq:nine-margins}, \(n_0\le4h+24\), so this split
is strictly worse.
If \(r=5\), \(m=3\), Euler's inequality for \(B\) gives
\(b\ge h+5\), \(\beta\ge2h+6\).
Lemma~\ref{lem:witness-facet-floor} gives \(u\ge13\), whence
\(N,M\ge u+4h+11\ge4h+24\).
For \(h\ge2\), \eqref{eq:nine-margins} gives
\(n_0\le4h+23\), excluding equality. For \(h=1\),
\(c=6\), \(s=9\), and \(n_0=28\).
An attaining presentation would require \(u=13\) and \(a\ge21\).
The only possibilities for \((\alpha,\rho)\) are \((8,5)\)
and \((7,6)\). In the first, \(\rho=5\) and a four-polytope with five
facets is a simplex, so \(F\) is a simplex facet and
Lemma~\ref{lem:simplex-facet-capacity} gives \(a\le16\).
In the second, the few-facet classification \cite{mcmullen1971} gives
\(a=6+pq\le12\) with \(p+q\le5\).
Both contradict \(a\ge21\).

\emph{Equality factor data.}
An attaining presentation therefore has \(r=m=4\).
Its count bounds force
\[
 b=c,\quad z=0,\quad u=s+4,\quad a=s+c+2h+4.
\]
If \(\rho\ge5\), the strengthened part of
Lemma~\ref{lem:four-vertex-capacity} would require
\[
 \alpha^2-7\alpha+2\rho-2c-4h-8\ge0.
\]
At \(\rho=5\), \(\alpha=s-1\), this is the second polynomial
of \eqref{eq:nine-thresholds} at \(s-1\), which is negative.
Increasing \(\rho\) by one and decreasing \(\alpha\) by one
changes it by \(10-2\alpha<0\), throughout the admissible range
\(\alpha\ge6\). Thus \(\rho=4\), \(\alpha=s\), and \(F\)
is tetrahedral. Conversely the stated factor data give both
certificate inequalities and \(N=M=n_0\).

\emph{Rational attainment.}
Set \(q_B=h+11-2c\) and \(q_A=c+2h+14-2s\).
Inequalities \eqref{eq:nine-thresholds} and \eqref{eq:nine-margins}
give
\[
 0\le q_B\le\binom{c-4}{2},\qquad
 0\le q_A\le\binom{s-5}{2}.
\]
The sufficiency theorem of Billera and Lee
\cite[Theorem~1]{billeraLee1980} gives simplicial four-polytopes
\(B,E\) with truncated \(g\)-vectors
\((1,c-5,q_B)\), \((1,s-6,q_A)\), respectively. As in the
eight-dimensional construction, these are monomial order-ideal
dimensions: choose the specified number of quadratic monomials.
The face numbers are
\[
 f_0(B)=c,\quad f_3(B)=3c-10+q_B=c+h+1,
\]
\[
 f_0(E)=s-1,\quad f_3(E)=3(s-1)-10+q_A=s+c+2h+1.
\]
Strict stacking on a tetrahedral facet of \(E\) adds one vertex
and three facets. Its polar \(A\) is simple with the required
counts, and has a tetrahedral facet \(F\) from the corresponding
simple-vertex cut. Rational perturbations preserve the simplicial
face lattices; polarity after rational interior translation and
both strict stacking choices can be made rational.
Stacking beyond \(F*B\) now attains \(n_0\).

\emph{Exact deficiency for every extremizer.}
Every four-polytope admits a vertex--facet assignment
\cite[Corollary~8]{jahnWinter2021}. Thus every optimal complementary
factor \(B\), not only the simplicial realizations just constructed,
has a matching covering its vertices, since \(\beta>b\).
The equality data give
\[
 a-\alpha-\beta-3=\beta-b-1=h,
 \qquad a-(3\alpha-7)=c+2h+11-2s\ge c-5>0.
\]
Proposition~\ref{prop:tetrahedral-saturation} therefore gives
\(\nu(\N(Q))=h+5+2c+s\) and \(\delta(Q)=h\).
Corollary~\ref{cor:stacking-transfer} gives connectedness.
Finally \(c=\sqrt{2h}+O(1)\), and the formula for \(s\) gives
\(s=2\sqrt h+O(1)\), proving the asymptotic statement.
\end{proof}

\subsection{Facet saturation and growth in every fixed dimension}
\label{subsec:fixed-deficiency}

The matching argument can be separated from the dimension of the witness.
A matching inside the stacking facet covers its adjacent facets;
vertices outside that facet cover the remaining retained facets and the
new cone facets. The following criterion states the two required counting
conditions.

\begin{proposition}[Facet matching saturation]
\label{prop:facet-saturation}
Let \(A\) be an \(r\)-polytope, \(r\ge3\), with \(a\) vertices
and \(\alpha\) facets. Let \(F\) be a facet with \(v\) vertices
and \(\rho\) ridges. Suppose \(\N(F)\) has a matching covering
all its \(\rho\) facet nodes. Let \(B\) have \(b\) vertices
and \(\beta\ge b+1\) facets, and suppose \(\N(B)\) has a
matching covering its vertices. Put
\[
 \theta=\alpha-\rho-1,\qquad
 \lambda=\max_{\substack{H\in\F(A),\ H\ne F\\
                         H\text{ not adjacent to }F}}
            |\V(H)\setminus\V(F)|,
\]
where adjacent facets share a ridge and an empty maximum is zero.
If
\begin{equation}\label{eq:facet-saturation-conditions}
 a-v-\lambda\ge\theta,\qquad a\ge v+\alpha+\beta-1,
\end{equation}
then strict stacking beyond \(F*B\) in \(A*B\) gives a polytope
\(Q\) with
\[
 \nu(\N(Q))=\alpha+\beta+b+\rho,
 \qquad
 \delta(Q)=\min\{a-\alpha-\beta-\rho+1,\ \beta-b-1\}.
\]
\end{proposition}

\begin{proof}
Each ridge \(R\) of \(F\) belongs to a unique other facet
\(H_R\) of \(A\). Different ridges give different \(H_R\),
because \(H_R\cap F\) is a proper face of \(F\) and so contains at most
one ridge of \(F\); for the same reason \(H_R\cap F=R\). A vertex of
\(F\) omitted by \(R\) is therefore also omitted by \(H_R\).
The matching in \(\N(F)\) therefore pairs \(\rho\) distinct
vertices of \(F\) with the retained facets \(H_R*B\).

There are \(\theta\) other retained facets of the form \(H*B\).
Each omits at least \(a-v-\lambda\ge\theta\) vertices of \(A\)
outside \(F\). Greedy selection pairs these \(\theta\) facets with
distinct outside vertices. At least
\[
 a-v-\theta\ge\beta+\rho
\]
outside vertices remain. All \(\beta+\rho\) cone facets omit
every such vertex, so they can be paired with distinct remaining
outside vertices. Finally, lift the \(b\) matching edges of
\(\N(B)\) to retained facets \(A*G\). Since \(\beta>b\),
one such facet remains for the new stacking vertex. These four
blocks are disjoint and contain
\(\alpha-1+\beta+\rho+b+1\) edges.

Conversely the neighbors of all vertices of \(A\) lie among the
\(\alpha-1\) retained facets of type \(H*B\) and the
\(\beta+\rho\) cone facets. At most \(b+1\) matching edges
can use vertices outside \(A\). This bounds the matching number
by the size already attained. Subtracting it from
\(a+b+1\) and \(\alpha+2\beta+\rho-1\) proves the
deficiency formula.
\end{proof}

For a simplex \(F\), the required facet matching pairs each vertex
with its opposite facet. Proposition~\ref{prop:tetrahedral-saturation}
is recovered by taking \(r=4\) and \(v=\rho=4\): a facet
of \(A\) has at most \(\alpha-1\) facets, so Euler's inequality
gives \(\lambda\le2\alpha-6\). The two inequalities in
\eqref{eq:facet-saturation-conditions} then follow from
\(a\ge3\alpha-7\) and \(a\ge\alpha+\beta+3\).

\begin{theorem}[Rational families in every fixed dimension]
\label{thm:fixed-deficiency}
Let \(h\ge1\) and \(k\ge0\) be integers, and set \(d=9+k\).
With \(c_h\) as in Theorem~\ref{thm:nine-deficiency}, put
\[
 s_{h,k}=\left\lceil
          \frac{7+\sqrt{16h+8c_h+8k+41}}2\right\rceil,
 \qquad n_{h,k}=2h+2k+5+2c_h+s_{h,k}.
\]
There is a rational balanced polytope in \(\mathcal C_{d,h}\)
with \(n_{h,k}\) vertices and facets, connected nonincidence graph
and matching number \(n_{h,k}-h\). In particular,
\begin{equation}\label{eq:fixed-deficiency-bounds}
 2d+2h+6\le\eta_{d,h}\le n_{h,d-9}\qquad(d\ge9).
\end{equation}
\end{theorem}

\begin{proof}
Write \(c=c_h\), \(s=s_{h,k}\), and \(\ell=\lceil k/2\rceil\).
The integer \(s\) is the smallest positive one for which
\[
 s^2-7s-2c-4h-2k+2\ge0.
\]
Using \(2c\le h+11\) from \eqref{eq:nine-margins}, the polynomial with
\(k=0\) at \(h+8\) equals \(h^2+5h+10-2c\ge h^2+4h-1>0\).
For general \(k\), its value at \(h+8+\ell\) exceeds this by
\(\ell(2h+9)+\ell^2-2k\), which is nonnegative because \(2\ell\ge k\).
Thus
\begin{equation}\label{eq:fixed-deficiency-margin}
 s\le h+8+\ell.
\end{equation}
Set
\[
 q_B=h+11-2c,\qquad q_A=c+2h+k+14-2s.
\]
The nine-dimensional threshold gives
\(0\le q_B\le\binom{c-4}{2}\). Also
\(q_A\ge c-2+k-2\ell\ge c-3\ge0\), and the defining
inequality for \(s\) is exactly \(q_A\le\binom{s-5}{2}\).
As before, Billera--Lee sufficiency \cite[Theorem~1]{billeraLee1980}
realizes simplicial four-polytopes \(B,E\) with truncated
\(g\)-vectors \((1,c-5,q_B)\) and \((1,s-6,q_A)\).
Their vertex--facet counts are
\[
 (c,c+h+1),\qquad (s-1,s+c+2h+k+1).
\]
Stack beyond a tetrahedral facet of \(E\) and take the polar \(A_0\).
This simple four-polytope has a tetrahedral facet \(F_0\) and counts
\[
 (a_0,\alpha_0)=(s+c+2h+k+4,s).
\]
Now take \(A=\operatorname{Pyr}^{k}(A_0)\) and
\(F=\operatorname{Pyr}^{k}(F_0)\). The facet \(F\) is a simplex;
\(r=4+k\), \(v=\rho=4+k\), \(a=a_0+k\), and
\(\alpha=s+k\). Strict stacking beyond \(F*B\) gives
\[
 f_0(Q)=f_{d-1}(Q)=n_{h,k},\qquad
 a-\alpha-\beta-\rho+1=\beta-b-1=h.
\]
Hence \(Q\in\mathcal C_{d,h}\). Rational realizations follow
from the same perturbation, polarity and strict stacking argument
as in Theorem~\ref{thm:nine-deficiency}, with rational pyramid apices.

To apply Proposition~\ref{prop:facet-saturation}, observe that every
facet of \(A\) not adjacent to \(F\) is
\(\operatorname{Pyr}^{k}(H_0)\), with \(H_0\) a facet of
\(A_0\) not adjacent to \(F_0\). Each pyramid facet omitting one
apex is adjacent to \(F\). The outside-\(F\) vertices of the
former facets are precisely the outside-\(F_0\) vertices of
\(H_0\). Since \(H_0\) has at most \(s-1\) facets, Euler's
inequality yields \(\lambda\le2s-6\). Moreover,
\[
 a_0-(3s-7)=c+2h+k+11-2s
           \ge c-5+k-2\ell\ge c-6\ge0.
\]
This implies \(a-v-\lambda\ge\alpha-\rho-1\), and
\(a-(v+\alpha+\beta-1)=h\). The simplex facet has the required
matching; the four-dimensional complement has a vertex-covering
assignment by \cite[Corollary~8]{jahnWinter2021}, since \(\beta>b\).
Saturation gives \(\delta(Q)=h\). Strict stacking gives connectedness
by Corollary~\ref{cor:stacking-transfer}. Finally,
Proposition~\ref{prop:sharp-floor} gives the lower bound in
\eqref{eq:fixed-deficiency-bounds}.
\end{proof}

\begin{corollary}[Leading deficiency growth]
\label{cor:deficiency-growth}
For each fixed integer \(d\ge7\),
\[
 \lim_{h\to\infty}\frac{\eta_{d,h}}h=
 \begin{cases}
 7,&d=7,\\
 4,&d=8,\\
 2,&d\ge9.
 \end{cases}
\]
For fixed \(d\ge9\), more precisely,
\[
 2d+2h+6\le\eta_{d,h}
 \le2h+(2+2\sqrt2)\sqrt h+O_d(1);
\]
for \(d\ge13\), Theorem~\ref{thm:second-order} gives a sharper upper
bound.
\end{corollary}

\begin{proof}
For fixed \(k\), \(c_h=\sqrt{2h}+O(1)\) and
\(s_{h,k}=2\sqrt h+O_k(1)\). Apply
Theorem~\ref{thm:fixed-deficiency} for \(d\ge9\), and
Theorems~\ref{thm:seven-deficiency} and~\ref{thm:eight-deficiency}
for the first two dimensions.
\end{proof}

\subsection{A bounded-error growth law from paired capacities}
\label{subsec:second-order}

The preceding construction uses four-dimensional capacity in every
dimension. Allowing both factor dimensions to vary gives the optimal
remainder, including every unbounded term up to a bounded additive error.
The realization input is classical:
vary the last coordinate of a simplicial \(g\)-vector, as in the
proofs of \cite[Theorems~3.2 and~3.3]{sjobergZiegler2020}.
The additional step is to realize these counts with a simplex stacking
facet and a saturated matching.

For \(s\ge3\), let \(U_s(n)\) denote the number of facets of the
cyclic \(s\)-polytope with \(n\) vertices. For \(p\ge1\), and
regarding both sides as polynomials in \(n\),
\begin{equation}\label{eq:cyclic-capacity}
 \begin{aligned}
 U_{2p}(n)&=\binom{n-p}{p}+\binom{n-p-1}{p-1},\\
 U_{2p+1}(n)&=2\binom{n-p-1}{p}.
 \end{aligned}
\end{equation}
McMullen's Upper Bound Theorem \cite{mcmullen1970} bounds the facets
of any \(s\)-polytope with \(n\) vertices by \(U_s(n)\).
Write \(p_s=\lfloor s/2\rfloor\) and
\[
 \kappa_{2p}=\frac1{p!},\qquad \kappa_{2p+1}=\frac2{p!}.
\]
Then \(U_s(n)=\kappa_s n^{p_s}+O_s(n^{p_s-1})\).

\begin{lemma}[Classical simplicial interpolation]
\label{lem:capacity-interpolation}
Fix \(s\ge4\), put \(p=\lfloor s/2\rfloor\), and let
\(n\ge s+2\), \(L=n-s-1\). Set
\[
 T_s(n)=\sum_{i=0}^{p-1}(s-2i+1)\binom{L+i-1}{i}.
\]
For every integer \(f\) with \(T_s(n)\le f\le U_s(n)\),
there is a rational simplicial \(s\)-polytope with \(n\) vertices
and \(f\) facets, provided \(f\) is even when \(s\) is odd.
Moreover \(T_s(n)=O_s(n^{p-1})\).
\end{lemma}

\begin{proof}
The sufficiency part of the \(g\)-theorem
\cite[Theorem~1]{billeraLee1980} realizes the degree counts of a
monomial order ideal as a truncated simplicial \(g\)-vector.
Take all monomials of degrees less than \(p\) in \(L\) variables,
and any \(z\) monomials of degree \(p\), where
\(0\le z\le\binom{L+p-1}{p}\). This is an order ideal and has
\[
 g_i=\binom{L+i-1}{i}\quad(0\le i<p),\qquad g_p=z.
\]
Here \(g_0=1\) and \(g_1=L\), so the vertex count is \(n\).
The standard \(h\)-vector identities give
\[
 f_{s-1}=\sum_{i=0}^{p}(s-2i+1)g_i
         =T_s(n)+(s-2p+1)z.
\]
At the largest \(z\) this is \(U_s(n)\), the cyclic value.
The step is one in even dimension and two in odd dimension; in
the latter case every coefficient is even. This proves the stated
interpolation, also used in \cite[Section~3]{sjobergZiegler2020}.
The degree bound gives the estimate for \(T_s(n)\).
Finally, a sufficiently small rational perturbation of a simplicial
realization preserves its face lattice; see also \cite[p.~237]{stanley1980}.
\end{proof}

\begin{proposition}[Paired capacity realization]
\label{prop:paired-capacity}
Fix integers \(r,m\ge4\), put \(d=r+m+1\),
\(p=\lfloor r/2\rfloor\), and \(q=\lfloor m/2\rfloor\).
For every sufficiently large integer \(h\), there is a rational
balanced member of \(\mathcal C_{d,h}\) with connected nonincidence
graph, exact deficiency \(h\), and common class size
\begin{equation}\label{eq:paired-capacity-size}
 n_h=2h+\alpha_h+2b_h+r+1,
 \qquad
 \begin{aligned}
 \alpha_h&=\left(\frac{2h}{\kappa_r}\right)^{1/p}+O_{r,m}(1),\\
 b_h&=\left(\frac{h}{\kappa_m}\right)^{1/q}+O_m(1).
 \end{aligned}
\end{equation}
Here \(r\) and \(m\) are fixed, and implied constants depend only on
the indicated parameters.
\end{proposition}

\begin{proof}
\emph{Step 1: the complementary capacity.}
Choose the smallest integer \(b=b_h\ge m+2\) such that
\(U_m(b)\ge b+h+1\). If \(m\) is odd, restrict to integers
with \(b+h+1\) even. The polynomial \(U_m(x)-x-1\) has leading
term \(\kappa_m x^q\) and error \(O_m(x^{q-1})\), since \(q\ge2\).
For \(X=(h/\kappa_m)^{1/q}\), evaluating this polynomial at
\(X+c\), with a sufficiently large fixed positive or negative
\(c\), brackets the threshold \(h\). The polynomial is eventually
increasing, and rounding on the admissible lattice costs at most two.
Thus \(b=X+O_m(1)\). Minimality also gives
\[
 U_m(b)=b+h+1+O_m(b^{q-1})=h+o(h).
\]
Moreover \(T_m(b)=O_m(b^{q-1})=o(h)\), so
\(T_m(b)\le b+h+1\le U_m(b)\) for large \(h\).
Lemma~\ref{lem:capacity-interpolation} therefore
realizes a rational simplicial \(B\) with
\((b,\beta)=(b,b+h+1)\). Since \(\beta>b\), its assignment
covers all vertices \cite[Corollary~9]{jahnWinter2021}.

\emph{Step 2: a witness with a simplex facet.}
Choose the smallest \(\alpha=\alpha_h\ge r+3\) such that
\[
 U_r(\alpha-1)\ge\alpha+b+2h+1.
\]
If \(r\) is odd, require \(\alpha+b\) odd, so the right-hand
side is even. Put \(Y=(2h/\kappa_r)^{1/p}\).
Since \(p,q\ge2\), \(b=O_m(\sqrt h)\) and \(Y=O_r(\sqrt h)\).
For a fixed additive displacement of \(\alpha\) from \(Y\), the
polynomial \(U_r(\alpha-1)\) changes by order \(h^{1-1/p}\).
This order is at least \(\sqrt h\), so a sufficiently large fixed
displacement absorbs both its lower-degree terms and the target
perturbation \(\alpha+b+1\). Positive and negative displacements
bracket the threshold, including the boundary case \(p=q=2\).
Eventual monotonicity and the parity lattice therefore give
\[
 \alpha=Y+O_{r,m}(1),\qquad
 U_r(\alpha-1)=2h+O_{r,m}(h^{1-1/p}).
\]
For sufficiently large \(h\), interpolation supplies a rational
simplicial \(r\)-polytope \(E\) with \(\alpha-1\) vertices
and \(\alpha+b+2h+1\) facets. Stack a simplex facet of \(E\),
then take a polar. The result \(A\) has
\[
 a=\alpha+b+2h+r\text{ vertices},\qquad \alpha\text{ facets},
\]
and a simplex facet \(F\), dual to the new simple vertex.
Thus \(v=\rho=r\). The two certificate margins are exactly \(h\),
and stacking beyond \(F*B\) gives the common count
\(2h+\alpha+2b+r+1\).

\emph{Step 3: saturation and rationality.}
Every facet \(H\) of \(A\) has at most \(\alpha-1\) facets.
Duality and the Upper Bound Theorem therefore bound its vertices
by \(U_{r-1}(\alpha-1)\). For the parameter \(\lambda\) in
Proposition~\ref{prop:facet-saturation}, this gives
\[
 \lambda\le U_{r-1}(\alpha-1)=
 \begin{cases}
 o(h),&r\text{ even},\\
 h+o(h),&r\text{ odd}.
 \end{cases}
\]
In odd dimension the leading coefficient is half that of \(U_r\),
and the minimality of \(\alpha\) in Step~2 gives
\(U_r(\alpha-1)=2h+o(h)\).
Since \(a=2h+o(h)\) and \(\alpha=o(h)\), eventually
\(a-\lambda\ge\alpha-1\), equivalently
\(a-r-\lambda\ge\alpha-r-1\). The other saturation condition has
the exact margin
\[
 a-(r+\alpha+\beta-1)=h.
\]
The simplex \(F\) has a facet-covering matching, and \(B\) has
the vertex-covering matching from Step~1. Saturation now gives
matching number \(n_h-h\). Strict stacking gives connectedness by
Corollary~\ref{cor:stacking-transfer}.
Rational interior translation and polarity preserve rationality,
and rational points exist in both open strict stacking regions.
This proves all assertions.
\end{proof}

\begin{theorem}[Bounded-error deficiency growth]
\label{thm:second-order}
Fix an integer \(d\ge9\), and put \(t=\lfloor(d-1)/4\rfloor\).
As \(h\to\infty\) through positive integers,
\begin{equation}\label{eq:second-order-growth}
 \eta_{d,h}=2h+K_d h^{1/t}+L_d h^{1/(t+1)}+O_d(1),
\end{equation}
where
\[
 K_d=(t!)^{1/t}
 \begin{cases}
 2^{1/t}+2,&d=4t+1,\\
 2^{1/t}+2^{1-1/t},&d=4t+2,\\
 2^{1/t},&d=4t+3,\\
 1,&d=4t+4,
 \end{cases}
\]
and
\[
 L_d=
 \begin{cases}
 0,&d=4t+1\text{ or }4t+2,\\
 2\bigl((t+1)!\bigr)^{1/(t+1)},&d=4t+3\text{ or }4t+4.
 \end{cases}
\]
Balanced rational connected examples of exact deficiency \(h\)
attain this growth profile within a bounded additive error. In particular
\(\eta_{d,h}-2h=\Theta_d(h^{1/t})\).
\end{theorem}

The remainder \(\eta_{d,h}-2h\) contains \(\alpha\) once and \(b\)
twice; this weighting selects the optimal ordered split
\((r,m)=(\dim A,\dim B)\).

\begin{proof}[Proof of Theorem~\ref{thm:second-order}]
\emph{Step 1: a lower bound for every presentation.}
For any \(Q\in\mathcal C_{d,h}\), write
\(\beta=b+h+1+z\), \(z\ge0\). The certificate gives
\(a\ge\alpha+b+2h+\rho+z\), so both class sizes are at least
\[
 2h+\alpha+2b+\rho+1.
\]
For its dimensions \(r,m\ge3\), set
\(p=\lfloor r/2\rfloor\), \(q=\lfloor m/2\rfloor\).
Applying the Upper Bound Theorem to \(A^*\) and \(B\), and
using \(a\ge2h\), \(\beta\ge h\), yields
\begin{equation}\label{eq:capacity-root-lower}
 \alpha\ge(2/\kappa_r)^{1/p}h^{1/p}-O_r(1),\qquad
 b\ge(1/\kappa_m)^{1/q}h^{1/q}-O_m(1).
\end{equation}
These inverse bounds follow directly from the leading term and
degree-\((p-1)\) or degree-\((q-1)\) error in
\eqref{eq:cyclic-capacity}; no simplicity or balance is assumed.

\emph{Step 2: the dimension comparison.}
Since \(r+m=d-1\), one has \(\min\{p,q\}\le t\).
If this minimum is smaller than \(t\),
\eqref{eq:capacity-root-lower} makes the remainder grow faster
than \(h^{1/t}\), whereas Step~3 below, which uses only
Proposition~\ref{prop:paired-capacity}, gives
\(\eta_{d,h}-2h=O_d(h^{1/t})\). For the remaining splits, the coefficient is
\[
 \mathbf1_{p=t}(2/\kappa_r)^{1/t}
 +2\mathbf1_{q=t}(1/\kappa_m)^{1/t}.
\]
Here \(\mathbf1\) denotes the indicator of the subscripted equality.
There are finitely many splits for fixed \(d\). Substitution of
the two values of \(\kappa\) gives the minima in Table~\ref{tab:capacity-splits}:
\begin{datadisplay}{Optimal dimension splits and second-order coefficients}\label{tab:capacity-splits}
\[
 \begin{array}{c|c|c}
 d & (r,m) & K_d/(t!)^{1/t}\\ \hline\noalign{\vskip 3pt}
 4t+1&(2t,2t)&2^{1/t}+2\\
 4t+2&(2t,2t+1)&2^{1/t}+2^{1-1/t}\\
 4t+3&(2t,2t+2)&2^{1/t}\\
 4t+4&(2t+1,2t+2)&1
 \end{array}
\]
\end{datadisplay}
In the first row \((2t,2t)\) is the only split with
\(\min\{p,q\}=t\). In the
second row the alternative \((2t+1,2t)\) has coefficient \(3\), and
\(2^{1/t}+2^{1-1/t}<3\) for \(t\ge2\).
In the third row the splits \((2t+1,2t+1)\) and \((2t+2,2t)\) have
coefficients \(1+2^{1-1/t}\) and \(2\), both larger than \(2^{1/t}\).
In the fourth row the splits \((2t,2t+3)\), \((2t+2,2t+1)\) and
\((2t+3,2t)\) have coefficients \(2^{1/t}\), \(2^{1-1/t}\) and \(2\),
all larger than one. Since these inequalities are strict, for each
fixed \(C\) and all sufficiently large \(h\), no other split has size
at most \(\eta_{d,h}+C\).

\emph{Step 3: attainment.}
For each displayed split, the lower bounds in Step~1 retain the full sum
\[
 2h+(2/\kappa_r)^{1/p}h^{1/p}
       +2(1/\kappa_m)^{1/q}h^{1/q}-O_d(1).
\]
Each displayed split has \(r,m\ge4\).
Proposition~\ref{prop:paired-capacity} supplies balanced rational
connected exact-\(h\) examples. Substituting its factor counts into
\eqref{eq:paired-capacity-size} gives the same sum with error
\(O_d(1)\). In the first two rows \(p=q=t\), so the two
contributions combine into \(K_d h^{1/t}\). In the last two rows
\(p=t\), \(q=t+1\) and \(m=2t+2\), so the second contribution
is \(L_d h^{1/(t+1)}\). This proves
\eqref{eq:second-order-growth} with its bounded remainder.
\end{proof}

\begin{corollary}[Bounded near-extremal budgets]
\label{cor:capacity-stability}
Fix \(d\ge9\) and \(C\ge0\). For all sufficiently large \(h\),
every member \(Q\in\mathcal C_{d,h}\) with
\(\max\{f_0(Q),f_{d-1}(Q)\}\le\eta_{d,h}+C\)
has, in every certified presentation, the ordered factor dimensions
\((r,m)=(\dim A,\dim B)\) listed in Table~\ref{tab:capacity-splits}.
For these dimensions put \(p=\lfloor r/2\rfloor\),
\(q=\lfloor m/2\rfloor\). Then
\[
 \alpha=(2h/\kappa_r)^{1/p}+O_{d,C}(1),\qquad
 b=(h/\kappa_m)^{1/q}+O_{d,C}(1).
\]
Moreover the ridge count \(\rho\), the two certificate slacks
\[
 z=\beta-b-h-1,\qquad
 w=a-\alpha-\beta-\rho-h+1,
\]
and the difference \(\lvert f_0(Q)-f_{d-1}(Q)\rvert\)
are bounded by constants depending only on \(d,C\).
\end{corollary}

\begin{proof}
The strict dimension comparisons in Theorem~\ref{thm:second-order}
exclude every nonoptimal split, even with the tolerance \(C\).
For any remaining presentation, the certificate gives \(z,w\ge0\).
Writing \(N=f_0(Q)\), \(M=f_{d-1}(Q)\), the count identities give
\[
 \max\{N,M\}
 =2h+\alpha+2b+\rho+1+\max\{z+w,2z\},\qquad N-M=w-z.
\]
Combine this equality with \eqref{eq:capacity-root-lower} and the
bounded-error profile of Theorem~\ref{thm:second-order}.
Each inverse-capacity count can exceed its lower bound by only a
bounded amount; the nonnegative quantities \(\rho,z,w\) have the
same bounded budget. The asserted estimates and imbalance bound follow.
\end{proof}

The square-root order therefore persists precisely for
\(9\le d\le12\); it becomes cube-root order in dimensions
\(13\) through \(16\), fourth-root order in dimensions \(17\)
through \(20\), and continues in the same way. Dimensions
\(4t+3\) and \(4t+4\) also have a positive correction of order
\(h^{1/(t+1)}\); both coefficients are optimal. The exact integer
formula in dimension nine agrees with
\(\eta_{9,h}=2h+(2+2\sqrt2)\sqrt h+O(1)\).

\subsection{Exact integer capacities with even optimal factors}
\label{subsec:even-capacity}

The bounded-error profile leaves an integer optimization problem.
When both optimal factors have even dimension, their facet capacities
can be interpolated without a parity restriction. A deleted-polar
bound then resolves the rounding and constrains every extremizer.

\begin{lemma}[Deleted-polar capacity]\label{lem:deleted-polar}
Let \(A\) be an \(r\)-polytope, \(r\ge4\), with \(a\)
vertices and \(\alpha\) facets. Let \(F\) be a facet with
\(\rho\) ridges, and let \(x\) be its dual vertex in \(D=A^*\).
If \(\conv(\V(D)\setminus\{x\})\) is full-dimensional, then
\begin{equation}\label{eq:deleted-polar}
 a\le U_r(\alpha-1)+U_{r-1}(\rho)-1.
\end{equation}
In particular, if \(\rho=r\), then
\(a\le U_r(\alpha-1)+r-1\).
\end{lemma}

\begin{proof}
The vertex figure at \(x\) has \(\rho\) vertices, since its vertices
correspond to the edges of \(D\) at \(x\), hence to the ridges of
\(A\) in \(F\). The Upper Bound
Theorem \cite{mcmullen1970} bounds its facets, and hence the facets
of \(D\) through \(x\), by \(U_{r-1}(\rho)\).
Every facet of \(D\) not through \(x\) remains a facet of the
deleted hull \(D'\). Since \(x\notin D'\), at least one facet
of \(D'\) is strictly visible from \(x\) and is not among these
retained facets. The remaining vertices are still extreme, so
\(D'\) has \(\alpha-1\) vertices. Another application of the
Upper Bound Theorem gives \eqref{eq:deleted-polar}.
When \(\rho=r\), the vertex figure is a simplex with \(r\) facets.
\end{proof}

\begin{theorem}[Eventual exact integer optimum]
\label{thm:even-integer}
Fix \(t\ge2\), and let
\[
 (d,r,m)=(4t+1,2t,2t)\quad\text{or}\quad(4t+3,2t,2t+2).
\]
Define integer thresholds
\begin{align}
 b_{m,h}&=\min\{b\ge m+2:U_m(b)\ge b+h+1\},\label{eq:even-b}\\
 s_{r,m,h}&=\min\{s\ge r+3:
       U_r(s-1)\ge s+b_{m,h}+2h+1\}.\label{eq:even-s}
\end{align}
For all sufficiently large \(h\), depending on \(d\),
\begin{equation}\label{eq:even-integer}
 \eta_{d,h}=2h+r+1+2b_{m,h}+s_{r,m,h}.
\end{equation}
Balanced rational connected polytopes attain this value with deficiency
exactly \(h\). Every extremizer is balanced, has connected
nonincidence graph, and has matching number \(\eta_{d,h}-h\).

In every certified presentation of an extremizer the factor dimensions
are \((r,m)\), and
\begin{equation}\label{eq:even-equality}
 \begin{gathered}
 b=b_{m,h},\qquad \beta=b_{m,h}+h+1,\\
 \alpha+\rho=s_{r,m,h}+r,\qquad
 a=s_{r,m,h}+r+b_{m,h}+2h,\qquad \rho\in\{r,r+1\}.
 \end{gathered}
\end{equation}
For an already realized certified presentation with these dimensions,
the data in \eqref{eq:even-equality} are also sufficient for
extremality. When \(r=4\), necessarily \(\rho=4\), so the
stacking facet is a simplex. For \(r\ge6\), we do not know whether the alternative
\(\rho=r+1\) occurs.
\end{theorem}

\begin{proof}
Write \(b_0=b_{m,h}\), \(s_0=s_{r,m,h}\), and
\(E_h=2h+r+1+2b_0+s_0\). Both thresholds exist because the
capacities have degree at least two.

\emph{Step 1: attainment and localization.}
Since \(r,m\) are even, the two interpolation lattices in
Proposition~\ref{prop:paired-capacity} have step one. For all
sufficiently large \(h\), take exactly \(b=b_0\) and
\(\alpha=s_0\) there. The construction gives
\[
 (a,\alpha)=(s_0+b_0+2h+r,s_0),\qquad
 (b,\beta)=(b_0,b_0+h+1),\qquad v=\rho=r.
\]
The lower interpolation endpoints \(T_m(b_0)\) and \(T_r(s_0-1)\) are
\(o(h)\), and \eqref{eq:even-b}--\eqref{eq:even-s} say exactly that the
required facet counts do not exceed the upper endpoints \(U_m(b_0)\)
and \(U_r(s_0-1)\).
Thus the classical realization and saturation argument in that
proposition applies and gives a rational connected exact-\(h\)
polytope with both counts equal to \(E_h\).

The inverse-capacity estimates in the same proposition show that
\(E_h\) has the profile in Theorem~\ref{thm:second-order} up to
\(O_d(1)\). Hence \(0\le E_h-\eta_{d,h}=O_d(1)\), without
assuming \eqref{eq:even-integer}. Corollary~\ref{cor:capacity-stability}
therefore applies to every \(Q\in\mathcal C_{d,h}\) whose
larger class size is at most \(E_h\). In every certified presentation,
the dimensions are \((r,m)\), the ridge count \(\rho\) is
uniformly bounded, and \(\alpha,b\) tend to infinity.
Since a simplex cannot satisfy \(\beta>b\), we have \(b\ge m+2\);
the certificate and the Upper Bound Theorem give
\(U_m(b)\ge\beta\ge b+h+1\), so \(b\ge b_0\).

\emph{Step 2: the integer lower bound.}
The deleted polar is full-dimensional for these localized presentations.
Otherwise \(D\) is a pyramid with apex \(x\), giving
\(\rho=\alpha-1\), contrary to bounded \(\rho\).
Put \(u=\alpha+\rho\). The certificate gives
\(a\ge u+b_0+2h\). If \(\rho=r\),
Lemma~\ref{lem:deleted-polar} implies
\(U_r(\alpha-1)\ge\alpha+b_0+2h+1\), whence
\(\alpha\ge s_0\) and \(u\ge s_0+r\).

Suppose instead that \(\rho>r\) and \(u<s_0+r\).
Then \(\alpha\le s_0-2\). The polynomial
\(P(x)=U_r(x-1)-x\) is eventually increasing. Since
\(\alpha\) tends to infinity, \eqref{eq:deleted-polar} yields
\[
 \rho+b_0+2h
 \le U_r(s_0-3)-(s_0-2)+U_{r-1}(\rho)-1.
\]
Minimality of \(s_0\) gives
\(U_r(s_0-2)\le s_0+b_0+2h-1\). Combining these inequalities,
\[
 U_r(s_0-2)-U_r(s_0-3)\le U_{r-1}(\rho)-\rho.
\]
The left side tends to infinity because \(r/2\ge2\), whereas
the right side is bounded. Thus \(u\ge s_0+r\) in all cases.
With the nonnegative slacks \(z,w\) of
Corollary~\ref{cor:capacity-stability}, the exact identity is
\[
 \max\{f_0(Q),f_{d-1}(Q)\}
 =2h+u+2b+1+\max\{z+w,2z\}\ge E_h.
\]
Polytopes of larger size already satisfy the inequality. Together with
Step~1 this proves \eqref{eq:even-integer}.

\emph{Step 3: equality.}
Equality forces \(b=b_0\), \(u=s_0+r\), and \(z=w=0\),
giving the displayed factor counts and balance. If \(\rho\ge r+2\),
then \(\alpha\le s_0-2\), and the divergent-gap argument in
Step~2 again gives a contradiction. Hence \(\rho=r\) or \(r+1\).
Conversely, the displayed data give both counts equal to \(E_h\)
in any realized certified presentation.
If \(r=4\) and \(\rho=5\), then \(\alpha=s_0-1\).
The strengthened bound in Lemma~\ref{lem:four-vertex-capacity} gives
\[
 a\le U_4(s_0-2)+4\le s_0+b_0+2h+3,
\]
contrary to \(a=s_0+b_0+2h+4\). This excludes the second
alternative when \(r=4\).

\emph{Step 4: matching saturation for every extremizer.}
Write \(r=2p\), \(m=2q\). A vertex of \(B\) belongs to
at most \(U_{m-1}(b_0-1)=O_d(b_0^{q-1})=o(h)\) facets,
by applying the Upper Bound Theorem to its vertex figure.
Consequently every vertex omits at least \(\beta-o(h)\ge b_0\)
facets for large \(h\). Greedy selection gives a matching covering
all vertices of \(B\), without a simplicity assumption.

The facet \(F\) has dimension \(r-1\) and \(r\) or
\(r+1\) facets. In the first case it is a simplex. In the second,
the classical few-facet classification \cite{mcmullen1971} writes
it as an iterated pyramid over \(\Delta_i\times\Delta_j\),
\(i,j\ge1\). This product has a facet-covering nonincidence
matching: with row facets \(H_k\), column facets \(V_\ell\),
and vertices \(w_{k\ell}\), use
\[
 \begin{gathered}
 H_0\leftrightarrow w_{00},\quad H_1\leftrightarrow w_{11},\quad
 H_k\leftrightarrow w_{k0}\ (2\le k\le i),\\
 V_0\leftrightarrow w_{10},\quad
 V_\ell\leftrightarrow w_{0\ell}\ (1\le\ell\le j).
 \end{gathered}
\]
Empty ranges are omitted; the pairs remain valid when \(ij=1,2\).
Each pyramid adds an isolated apex--base edge. Thus \(\N(F)\)
has a matching covering its facets in either case.

Finally \(v=f_0(F)\le U_{r-1}(\rho)=O_d(1)\).
Each facet of \(A\) has at most \(\alpha-1\) ridges, so
polarity and the Upper Bound Theorem bound its vertices by
\(U_{r-1}(\alpha-1)=O_d(\alpha^{p-1})=o(h)\).
The outside-facet parameter of Proposition~\ref{prop:facet-saturation}
therefore satisfies \(\lambda=o(h)\). Since
\(a=2h+o(h)\), \(\alpha,b=o(h)\), and \(v,\rho=O_d(1)\),
\[
 a-v-\lambda\ge\alpha-\rho-1,\qquad
 a-(v+\alpha+\beta-1)=\rho+h-v\ge0
\]
for large \(h\). Proposition~\ref{prop:facet-saturation} applies
to every extremizer, giving \(\delta(Q)=h\) and \(\nu(\N(Q))=E_h-h\).
Connectedness follows from Corollary~\ref{cor:stacking-transfer}.
\end{proof}

The theorem determines factor counts, not face lattices. Its
integer formula is eventual; the onset is not evaluated here.
An odd optimal complement adds a face-number feasibility constraint,
which we isolate next.

\subsection{Odd complementary factors and extremal plateaus}
\label{subsec:odd-capacity}

An odd-dimensional complement changes the integer problem. Its simplicial
facet counts are even, and some odd counts immediately below the cyclic
capacity are impossible. We isolate the single face-number pair that
affects the optimum, then determine the equality slacks and their effect
on the actual matching deficiency.

Fix \(t\ge2\) and put \((d,r,m)=(4t+2,2t,2t+1)\). Define
\begin{equation}\label{eq:odd-baseline}
 b_h^0=\min\{b\ge m+2:U_m(b)\ge b+h+1\},\qquad
 \beta_h^0=b_h^0+h+1.
\end{equation}
Let \(\chi_h\) be zero if an \(m\)-polytope with \(b_h^0\)
vertices and \(\beta_h^0\) facets exists, and one otherwise.
Existence here is over the real numbers. Set
\begin{align}
 \beta_h^*&=\beta_h^0+\chi_h,\label{eq:odd-beta}\\
 s_h^*&=\min\{s\ge r+3:U_r(s-1)\ge s+\beta_h^*+h\}.
 \label{eq:odd-s}
\end{align}
For sufficiently large \(h\), classical interpolation gives
\(\chi_h=0\) when \(\beta_h^0\) is even. If \(\beta_h^0\)
is odd and \(U_m(b_h^0)-\beta_h^0<t\), then \(\chi_h=1\).
Indeed, a simplicial odd-dimensional polytope has an even number of
facets. A nonsimplicial \(m\)-polytope with \(b\) vertices has at
most \(U_m(b)-\lfloor m/2\rfloor\) facets
\cite[proof of Theorem~3.3]{sjobergZiegler2020}. This bound follows
from the almost-simplicial upper bound
\cite[Theorem~1.2 and Proposition~4.2]{nevoPinedaUgonYost2018}
by pulling while preserving one nonsimplex facet, as explained in
\cite{sjobergZiegler2020}.

\begin{theorem}[An odd-complement integer reduction]
\label{thm:odd-integer}
With the preceding definitions, for all sufficiently large \(h\),
\begin{equation}\label{eq:odd-integer}
 \eta_{d,h}=E_h:=s_h^*+r+2\beta_h^*-1.
\end{equation}
There are connected balanced attainers of deficiency exactly \(h\),
and connected exact-\(h\) attainers whose complement has \(b_h^0\)
vertices. These attainers can be chosen rational if \(\chi_h=1\)
or \(\beta_h^0\) is even. In the remaining case the constructions
use a realization of the specified real face-number pair.

In every certified presentation of an extremizer, the factor dimensions
are \((r,m)\), and, writing \(z=\beta-b-h-1\) and
\(w=a-\alpha-\beta-\rho-h+1\), the factor data satisfy
\begin{equation}\label{eq:odd-equality}
 \begin{gathered}
 \beta=\beta_h^*,\quad \alpha+\rho=s_h^*+r,\quad
 \rho\in\{r,r+1\},\quad b_h^0\le b\le\beta_h^*-h-1,\\
 z\in\{0,1\},\quad 0\le w\le z,\quad
 a=s_h^*+r+\beta_h^*+h-1+w.
 \end{gathered}
\end{equation}
When \(r=4\), necessarily \(\rho=4\). For an already realized
certified presentation of these dimensions, the data in
\eqref{eq:odd-equality} are sufficient for extremality. Every extremizer
has connected nonincidence graph and
\begin{equation}\label{eq:odd-matching}
 \delta(Q)=h+w,\qquad \nu(\N(Q))=\alpha+\rho+\beta_h^*+b,
 \qquad f_0(Q)-f_{d-1}(Q)=w-z.
\end{equation}
\end{theorem}

The reduction leaves only the feasibility of the one pair in
\eqref{eq:odd-baseline}; it does not require a classification of all
odd-dimensional vertex--facet pairs. Unlike
Theorem~\ref{thm:even-integer}, equality can retain one unit of
certificate slack. The next corollary realizes this alternative.

\begin{proof}[Proof of Theorem~\ref{thm:odd-integer}]
Write \(b_0=b_h^0\), \(\beta_0=\beta_h^0\),
\(\beta_*=\beta_h^*\), and \(s_*=s_h^*\).

\emph{Step 1: complementary realization.}
The argument of Step~1 of Proposition~\ref{prop:paired-capacity},
without its parity restriction, gives
\[
 b_0=(h/\kappa_m)^{1/t}+O_d(1),\qquad U_m(b_0)=h+o(h),
 \qquad T_m(b_0)=o(h).
\]
If \(\beta_0\) is even, Lemma~\ref{lem:capacity-interpolation}
realizes the pair \((b_0,\beta_0)\) rationally. If \(\beta_0\)
is odd, the even cyclic capacity satisfies
\(\beta_0+1\le U_m(b_0)\), so the same lemma realizes
\((b_0,\beta_0+1)\). When \(\chi_h=0\), use a realization of
\((b_0,\beta_0)\) instead. Thus \(\beta_*\) is attained with
\(b=b_0\) in every case.
The value \(\beta_*\) is the least facet count of a complement of
dimension \(m\) in any certified presentation: the certificate and Upper Bound Theorem give
\(b\ge b_0\) and \(\beta\ge b+h+1\ge\beta_0\).
If \(\beta=\beta_0\), then \(b=b_0\), so that very complement
realizes the pair defining \(\chi_h\).

For balanced attainment we also need
\(b_*=\beta_*-h-1\). This is \(b_0\) if \(\chi_h=0\).
Otherwise it is \(b_0+1\), and the even count \(\beta_*\) is
between \(T_m(b_0+1)=o(h)\) and \(U_m(b_0+1)\).
Interpolation supplies a rational complement with this pair.

\emph{Step 2: an exact upper bound.}
Again by inverse capacity,
\(s_*=(2h/\kappa_r)^{1/t}+O_d(1)\).
Since \(r\) is even, interpolation realizes a rational simplicial
\(r\)-polytope with \(s_*-1\) vertices and \(s_*+\beta_*+h\)
facets. Stack a simplex facet and polarize. The witness then has
\[
 \alpha=s_*,\quad v=\rho=r,\quad
 a=s_*+\beta_*+h+r-1.
\]
For either complement just constructed the certificate margins are
\(h\) and \(h+z\), where \(z=\beta_*-b-h-1\). The strict stack
therefore has
\[
 M:=f_{d-1}(Q)=E_h,\qquad N:=f_0(Q)=E_h-z\le E_h.
\]

To obtain the matching, bound each complementary vertex figure by
the Upper Bound Theorem. For \(b=b_0\) or \(b_0+1\),
the actual cyclic formulas give
\[
 U_{m-1}(b-1)=\tfrac12 U_m(b)+O_d(b^{t-1})
             =\tfrac12 h+o(h).
\]
Every vertex of \(B\) consequently omits at least
\(\beta_*-h/2-o(h)\ge b\) facets for large \(h\). Greedy
selection gives a matching covering its vertices, even for a
nonsimplicial realization. The simplex \(F\) has a facet-covering
matching. The witness-facet bound is
\(\lambda\le U_{r-1}(\alpha-1)=o(h)\), while \(a=2h+o(h)\).
Hence
\[
 a-v-\lambda\ge\alpha-\rho-1,\qquad
 a-(v+\alpha+\beta_*-1)=h.
\]
Proposition~\ref{prop:facet-saturation} gives deficiency exactly
\(h\). Connectedness follows from Corollary~\ref{cor:stacking-transfer}.
Taking \(b=b_*\) makes \(z=0\) and gives balance. Rational
stacking preserves the rationality in the stated branches.

\emph{Step 3: the integer lower bound.}
The value \(E_h\) has the independently established profile in
Theorem~\ref{thm:second-order} up to \(O_d(1)\). Thus
\(0\le E_h-\eta_{d,h}=O_d(1)\), without assuming the desired
equality. Corollary~\ref{cor:capacity-stability} localizes every
presentation with larger class size at most \(E_h\): its dimensions
are \((r,m)\), its ridge count is bounded, and \(\alpha,b\) tend
to infinity. Its actual \(\beta\) is \(h+o(h)\).

For this \(\beta\), put
\[
 s(\beta)=\min\{s\ge r+3:U_r(s-1)\ge s+\beta+h\},
 \qquad u=\alpha+\rho.
\]
The deleted polar is full-dimensional; otherwise \(D=A^*\) is a
pyramid with apex \(x\), the vertex dual to \(F\), so \(\rho=\alpha-1\), contrary to localization.
Lemma~\ref{lem:deleted-polar} and the certificate imply
\[
 u+\beta+h-1\le a\le U_r(\alpha-1)+U_{r-1}(\rho)-1.
\]
For \(\rho=r\) this gives \(\alpha\ge s(\beta)\).
For \(\rho>r\), a putative \(u<s(\beta)+r\) would give
\(\alpha\le s(\beta)-2\). As in Step~2 of
Theorem~\ref{thm:even-integer}, eventual monotonicity of
\(U_r(x-1)-x\) and minimality of \(s(\beta)\) would imply
\[
 U_r(s(\beta)-2)-U_r(s(\beta)-3)
       \le U_{r-1}(\rho)-\rho.
\]
The left side tends to infinity and the right side is bounded,
uniformly for these presentations. Thus \(u\ge s(\beta)+r\).
Step~1 gives \(\beta\ge\beta_*\), and \(s(\beta)\) is
nondecreasing. Consequently
\[
 M=u+2\beta-1\ge s(\beta)+r+2\beta-1\ge E_h.
\]
Presentations of larger size already satisfy the bound. This proves
\eqref{eq:odd-integer}.

\emph{Step 4: equality and matching.}
Equality forces \(\beta=\beta_*\), \(u=s_*+r\), and \(N\le M\).
Since \(N-M=w-z\), it follows that \(0\le w\le z\).
The complement satisfies \(b_0\le b\le\beta_*-h-1\le b_0+1\),
so \(z\in\{0,1\}\). If \(\rho\ge r+2\), then
\(\alpha\le s_*-2\), and the same divergent-gap argument excludes
it. For \(r=4,\rho=5\), Lemma~\ref{lem:four-vertex-capacity}
would give
\[
 a\le U_4(s_*-2)+4\le s_*+\beta_*+h+2,
\]
whereas the certificate gives \(a\ge s_*+\beta_*+h+3\).
Thus \(\rho=4\) when \(r=4\).

For general \(r\), the few-facet classification and explicit
product matching in Step~4 of Theorem~\ref{thm:even-integer}
give a facet-covering matching for \(F\), with
\(v\le U_{r-1}(\rho)=O_d(1)\). The complementary vertex-figure
argument of Step~2 still applies, since \(b\) is \(b_0\) or
\(b_0+1\). Finally \(\lambda=o(h)\), \(a=2h+o(h)\), and
\[
 a-(v+\alpha+\beta_*-1)=\rho+h+w-v\ge0.
\]
Both saturation inequalities hold, and
Proposition~\ref{prop:facet-saturation} yields
\(\delta(Q)=\min\{h+w,h+z\}=h+w\) and the stated matching
number. Connectedness again follows from strict stacking.
Conversely, any realized certified presentation with
\eqref{eq:odd-equality} has \(M=E_h\) and \(N\le E_h\),
and is therefore extremal.
\end{proof}

\begin{corollary}[Ten-dimensional parity plateaus]
\label{cor:ten-plateau}
For every sufficiently large integer \(k\), set
\[
 h_k=k^2-8k+10,\qquad n_k=2k^2-12k+23.
\]
Then
\begin{equation}\label{eq:ten-plateau}
 \begin{aligned}
 \eta_{10,h_k-1}&=n_k-4,\\
 \eta_{10,h_k}=\eta_{10,h_k+1}&=n_k,\\
 \eta_{10,h_k+2}&=n_k+5.
 \end{aligned}
\end{equation}
Each value has a rational balanced connected attainer of exactly
the requested deficiency. At \(h=h_k\), put
\(a_k=2k^2-13k+21\). In every certified presentation of an extremizer
at \(h=h_k\),
\(\alpha=2k-4\), \(\rho=v=4\), and \(\beta=(k-3)(k-4)\).
The possible count and deficiency data are exactly the three rows below;
each row in Table~\ref{tab:ten-plateau-data} is rationally realized.
\begin{datadisplay}{The three ten-dimensional plateau configurations}\label{tab:ten-plateau-data}
\[
\begin{array}{c|c|c|c|c}
 b&a&(f_0(Q),f_9(Q))&\delta(Q)&\nu(\N(Q))\\ \hline
 k+1&a_k&(n_k,n_k)&h_k&n_k-h_k\\
 k&a_k&(n_k-1,n_k)&h_k&n_k-h_k-1\\
 k&a_k+1&(n_k,n_k)&h_k+1&n_k-h_k-1
\end{array}
\]
\end{datadisplay}
\end{corollary}

\begin{proof}[Proof of Corollary~\ref{cor:ten-plateau}]
Write \(h=h_k\) and \(U=U_5(k)=(k-3)(k-4)\).
At \(h\), the baseline in \eqref{eq:odd-baseline} is \(b_0=k\)
and \(\beta_0=U-1\): at \(k-1\), the capacity
\(U_5(k-1)=U-2k+8\) is smaller than the required \(U-2\)
for large \(k\). The count \(U-1\) is odd and violates the
nonsimplicial bound \(U_5(k)-2\), so \(\chi_h=1\) and
\(\beta_*=U\).
After multiplication by \(2\), the condition in \eqref{eq:odd-s} becomes
\[
 s^2-7s-4k^2+30k-40\ge0.
\]
At \(s=2k-4\) the left side is \(4\); at \(s-1=2k-5\)
it is \(20-4k<0\). Hence \(s_*=2k-4\), the capacity slack
\(U_4(s_*-1)-(s_*+\beta_*+h)\) is \(2\), and Theorem~\ref{thm:odd-integer} gives \(\eta_{10,h}=n_k\).

At \(h-1\) and \(h+1\), the baseline remains \(k\), with
even facet counts \(U-2\) and \(U\). The right side
\(s+\beta_*+h\) of \eqref{eq:odd-s} changes by \(-3\) and \(+1\), so the same \(s_*=2k-4\) remains minimal
for large \(k\). The respective slacks are \(5\) and \(1\),
giving \(n_k-4\) and \(n_k\).
At \(h+2\), the baseline is \(k+1\), with even count \(U+2\).
The right side of \eqref{eq:odd-s} increases by \(4\) from its value at \(h\).
It fails at \(2k-4\) and succeeds at \(2k-3\), where its slack
is \(2k-9\). The optimum is therefore \(n_k+5\).
All four values of \(\beta_*\) are even, so the theorem gives rational balanced
connected exact-deficiency attainers.

It remains to realize the three equality rows explicitly. For \(w=0,1\),
take a rational simplicial four-polytope \(E_w\) with \(2k-5\)
vertices and \(2k^2-13k+18+w\) facets. Its truncated \(g\)-vector is
\[
 (1,2k-10,2k^2-19k+43+w).
\]
The last coordinate is \(2-w\) below the Macaulay maximum
\(\binom{2k-9}{2}\); all coordinates are nonnegative for large
\(k\). Billera--Lee sufficiency \cite[Theorem~1]{billeraLee1980}
and rational perturbation realize this vector. Stack a new vertex beyond one facet of \(E_w\)
and polarize to obtain \(A_w\) with \(\alpha=2k-4\),
\(a=a_k+w\), and tetrahedral \(F\).
For \(b=k\) use the cyclic five-polytope. For \(b=k+1\),
the rational simplicial five-polytope with \(\beta=U\) has vector
\[
 \left(1,k-5,\frac{k^2-11k+26}{2}\right).
\]
Its last coordinate is a nonnegative integer, \(k-3\) below
\(\binom{k-4}{2}\). The pairs \((A_0,B_{k+1})\),
\((A_0,B_k)\), and \((A_1,B_k)\) have certificate margins
\((h,h)\), \((h,h+1)\), and \((h+1,h+1)\), respectively.
Saturation and strict stacking give the displayed counts and matchings.
The equality data of Theorem~\ref{thm:odd-integer}, with \(r=4\),
exhaust all rows.
\end{proof}

\subsection{Parity repair for an odd-dimensional witness}
\label{subsec:odd-witness}

An odd witness and an odd complement have different roles in the
integer problem. The complement must realize a prescribed face-number
pair. The witness instead needs a distinguished simplex facet after
polarity. Classical coplanar stacking removes its parity restriction
while preserving that facet.

\begin{lemma}[Simple-vertex parity repair]\label{lem:coplanar-stacking}
Let \(E\) be a rational simplicial \(r\)-polytope, \(r\ge3\),
with \(n\) vertices and \(f\) facets. For each
\(f'\in\{f+r-1,f+r-2\}\) there is a rational \(r\)-polytope with
\(n+1\) vertices and \(f'\) facets, of the form \(\conv(E\cup\{x\})\)
with the new vertex \(x\) simple.
\end{lemma}

\begin{proof}
The first count comes from strict stacking on a simplex facet.
For the second, use the coplanar stacking construction of
Sj\"oberg and Ziegler \cite[proof of Theorem~3.3]{sjobergZiegler2020}.
Choose adjacent simplex facets \(T,G\), with common ridge \(R\),
and a rational point \(p\in\operatorname{relint}R\).
Their supporting normals are independent, so there is a rational
direction \(q\) strictly outward from \(T\) and parallel to \(G\).
For sufficiently small positive rational \(\eps\), the point
\(x=p+\eps q\) is beyond \(T\), in \(\aff G\), and strictly
beneath every other facet of \(E\).

Every old vertex remains extreme. For each vertex \(v\), choose a
strictly exposing linear functional \(\ell_v\). Since
\(p\in E\setminus\{v\}\), we have \(\ell_v(p)<\ell_v(v)\); these finitely many inequalities persist
for small \(\eps\). The point \(x\) is extreme because it lies
beyond \(T\).

In \(\conv(E,x)\), the facet \(T\) disappears, \(G\) becomes
\(\conv(G,x)\), and the other old facets persist, by the
beneath--beyond theorem \cite[Theorem~5.2.1]{gruenbaum2003}; here
the restriction of the defining form of \(T\) to \(\aff G\) defines the
facet \(R\) of \(G\), so \(x\) lies beyond \(R\) within \(\aff G\). The cone over
\(R\) is contained in this extended facet. The other \(r-1\)
ridges of \(T\) give \(r-1\) new cone facets. The net increase
is \(r-2\). The new vertex belongs precisely to the extended
facet and these \(r-1\) cones, so it is simple.
\end{proof}

\begin{theorem}[Exact integer optimum with an odd witness]
\label{thm:odd-witness-integer}
Fix \(t\ge2\), and let \((d,r,m)=(4t+4,2t+1,2t+2)\).
Define \(b_{m,h}\) and \(s_{r,m,h}\) by the same inverse capacities
as in \eqref{eq:even-b}--\eqref{eq:even-s}, now with this odd \(r\).
For all sufficiently large \(h\), depending on \(d\),
\begin{equation}\label{eq:odd-witness-integer}
 \eta_{d,h}=2h+r+1+2b_{m,h}+s_{r,m,h}.
\end{equation}
Balanced rational connected polytopes attain this value with deficiency
exactly \(h\). Every extremizer is balanced, has connected nonincidence
graph, and has matching number \(\eta_{d,h}-h\).

In every certified presentation of an extremizer the dimensions are
\((r,m)\), and
\begin{equation}\label{eq:odd-witness-equality}
 \begin{gathered}
 b=b_{m,h},\qquad \beta=b_{m,h}+h+1,\\
 \alpha+\rho=s_{r,m,h}+r,\qquad
 a=s_{r,m,h}+r+b_{m,h}+2h,\qquad \rho\in\{r,r+1\}.
 \end{gathered}
\end{equation}
For an already realized certified presentation of these dimensions,
these data are also sufficient for extremality. We do not know whether
the case \(\rho=r+1\) is realized.
\end{theorem}

\begin{proof}
Put \(b_0=b_{m,h}\), \(\beta_0=b_0+h+1\),
\(s_0=s_{r,m,h}\), \(L_0=s_0+\beta_0+h\), and
\(E_h=s_0+r+2\beta_0-1\). The inverse-capacity estimates in
Proposition~\ref{prop:paired-capacity} give
\[
 b_0=(h/\kappa_m)^{1/(t+1)}+O_d(1),\qquad
 s_0=(2h/\kappa_r)^{1/t}+O_d(1).
\]
In particular, \(U_m(b_0)=h+o(h)\) and
\(U_r(s_0-1)=2h+o(h)\).

\emph{Step 1: parity-free attainment.}
Since \(m\) is even, Lemma~\ref{lem:capacity-interpolation}
gives a rational simplicial complement with
\((b,\beta)=(b_0,\beta_0)\). If \(L_0\) is even, realize a
rational simplicial \(r\)-polytope \(E\) with \(s_0-1\)
vertices and \(L_0\) facets, then stack strictly.
If \(L_0\) is odd, realize \(E\) with \(L_0+1\) facets,
then use the coplanar branch of Lemma~\ref{lem:coplanar-stacking}.
The capacity \(U_r(s_0-1)\) is even, so the latter target
is still within capacity. Both targets exceed the \(o(h)\)
lower interpolation endpoint for large \(h\).

Polarizing the result gives a witness \(A\) with
\[
 \alpha=s_0,\qquad a=L_0+r-1,\qquad v=\rho=r,
\]
where \(F\) is the simplex facet dual to the new simple vertex.
Both certificate margins equal \(h\). The final operation is
strict stacking beyond \(F*B\) in \(A*B\); the coplanarity used
to construct \(A\) does not affect this final strict operation.
Its two class sizes equal \(E_h\).

The simplex \(F\) has a facet-covering matching, and \(B\), being
simplicial with \(\beta>b\), has a vertex-covering matching
\cite[Corollary~9]{jahnWinter2021}. A facet of \(A\) has at most
\(\alpha-1\) ridges, giving
\[
 \lambda\le U_{r-1}(\alpha-1)=h+o(h).
\]
Here \(r-1=2t\); the leading term is \(h\), not \(o(h)\).
Since \(a=2h+o(h)\) and \(\alpha=o(h)\),
\(a-r-\lambda\ge\alpha-r-1\) eventually. The other saturation
margin is \(a-(r+\alpha+\beta_0-1)=h\).
Proposition~\ref{prop:facet-saturation} gives deficiency \(h\)
and matching number \(E_h-h\). Corollary~\ref{cor:stacking-transfer}
gives connectedness. All these constructions are rational.

\emph{Step 2: a lower bound for every presentation.}
The independent capacity profile in Theorem~\ref{thm:second-order}
and the inverse estimates above give
\(0\le E_h-\eta_{d,h}=O_d(1)\). Thus
Corollary~\ref{cor:capacity-stability} localizes every certified
presentation of size at most \(E_h\): its dimensions are \((r,m)\),
\(\rho\) is bounded, \(\alpha,b\) tend to infinity, and
\(\beta=h+o(h)\). The complementary certificate and Upper Bound
Theorem imply \(b\ge b_0\) and \(\beta\ge\beta_0\).
For the actual \(\beta\), set
\[
 s(\beta)=\min\{s\ge r+3:U_r(s-1)\ge s+\beta+h\},
 \qquad u=\alpha+\rho.
\]
The deleted polar is full-dimensional: otherwise \(D=A^*\) is a
pyramid with apex \(x\), the vertex dual to \(F\), so \(x\) has degree \(\alpha-1\); but its vertex
figure is polar to \(F\), so its degree is \(\rho\), which is bounded. The certificate and Lemma~\ref{lem:deleted-polar} give
\[
 u+\beta+h-1\le a
 \le U_r(\alpha-1)+U_{r-1}(\rho)-1.
\]
For \(\rho=r\), this implies \(\alpha\ge s(\beta)\).
If \(\rho>r\) and \(u<s(\beta)+r\), then
\(\alpha\le s(\beta)-2\). Eventual monotonicity of
\(U_r(x-1)-x\) and minimality of \(s(\beta)\), as in
Step~2 of Theorem~\ref{thm:even-integer}, yield
\[
 U_r(s(\beta)-2)-U_r(s(\beta)-3)
 \le U_{r-1}(\rho)-\rho.
\]
The left side diverges because \(U_r\) has degree \(t\ge2\);
the right side is bounded. Therefore \(u\ge s(\beta)+r\), and
since \(s(\beta)\) is nondecreasing in \(\beta\), \(s(\beta_0)=s_0\),
and \(\beta\ge\beta_0\), the facet count satisfies
\[
 M=u+2\beta-1\ge s(\beta)+r+2\beta-1
 \ge s_0+r+2\beta_0-1=E_h.
\]
Larger presentations already satisfy the bound, proving the formula.

\emph{Step 3: equality and its matching number.}
Equality forces \(\beta=\beta_0\), \(u=s_0+r\), and \(N\le M\).
The certificate slacks \(z=\beta-b-h-1\) and
\(w=a-\alpha-\beta-\rho-h+1\) are nonnegative.
Since \(b\ge b_0\), we have \(b=b_0\) and \(z=0\).
The identity \(N-M=w-z\) forces \(w=0\), balance, and the
displayed value of \(a\). If \(\rho\ge r+2\), then
\(\alpha\le s_0-2\), and the same divergent-gap argument excludes
this possibility. Conversely, the displayed count data give
\(N=M=E_h\) in any realized presentation of the stated dimensions.

For \(\rho=r,r+1\), the few-facet classification
\cite{mcmullen1971} makes \(F\) a simplex or an iterated pyramid
over \(\Delta_i\times\Delta_j\). The explicit product matching
in Step~4 of Theorem~\ref{thm:even-integer}, together with the
apex--base edges, covers its facets. Also
\(v\le U_{r-1}(\rho)=O_d(1)\). Every vertex of an arbitrary
complement \(B\) belongs to at most
\(U_{m-1}(b_0-1)=o(h)\) facets. It therefore omits at least
\(\beta_0-o(h)\ge b_0\) facets, and greedy selection gives a
vertex-covering matching. Finally \(\alpha=s_0+O_d(1)\) implies
\(\lambda\le U_{r-1}(\alpha-1)=h+o(h)\), whereas
\(a=2h+o(h)\). Hence
\[
 a-v-\lambda\ge\alpha-\rho-1,\qquad
 a-(v+\alpha+\beta_0-1)=h+\rho-v\ge0
\]
eventually. Proposition~\ref{prop:facet-saturation} applies to
every extremizer, giving deficiency exactly \(h\) and matching
number \(E_h-h\). Corollary~\ref{cor:stacking-transfer} gives
connectedness.
\end{proof}

Thus witness parity causes no loss in the eventual certified optimum.
This differs from the odd-complement feasibility constraint in
Theorem~\ref{thm:odd-integer}. The formula remains eventual, and
the count data do not constitute a face-lattice classification.

\subsection{Parity staircases at almost every capacity threshold}
\label{subsec:plateau-density}

The ten-dimensional plateau is not isolated. In every dimension
\(d=4t+2\) with \(t\ge3\), where the optimal complement is odd, the
same facet gap produces \(\lfloor t/2\rfloor\) consecutive plateaus
near almost every complementary-capacity threshold. The rounding
mechanism is different: the ratio \(4^{1/t}\) of the linear root scales
of the two capacity polynomials is irrational.

Fix \(t\ge3\) and \((d,r,m)=(4t+2,2t,2t+1)\). Put
\begin{align}
 U(k)&=U_m(k)=2\binom{k-t-1}{t},& P(s)&=U_r(s-1)-s,
 \label{eq:density-capacities}\\
 h_k&=U(k)-k-2,& C_k&=2U(k)-k-2,\notag\\
 s_k&=\min\{s\ge r+3:P(s)\ge C_k\},&
 n_k&=s_k+r+2U(k)-1.\label{eq:density-thresholds}
\end{align}
Here \(P(s)\) denotes a polynomial, rather than a polytope.
Let \(J=\lfloor t/2\rfloor\), \(A_t=\max\{3,4(J-1)\}\), and define
\begin{equation}\label{eq:density-good}
 \mathcal K_t=\{k\ge m+2:
       P(s_k-1)<C_k-A_t,\quad P(s_k)\ge C_k+4\}.
\end{equation}
Only sufficiently large indices will be used, so the finitely many
initial capacity values have no effect.

\begin{theorem}[Density-one parity staircases]
\label{thm:density-plateaus}
Let \(t\ge3\) and \((d,r,m)=(4t+2,2t,2t+1)\), with the notation
\eqref{eq:density-capacities}--\eqref{eq:density-good}. The set \(\mathcal K_t\) has natural density one:
\[
 \lim_{K\to\infty}\frac{|\mathcal K_t\cap\{1,\ldots,K\}|}{K}=1.
\]
For every sufficiently large \(k\in\mathcal K_t\),
\begin{equation}\label{eq:density-four}
 \begin{aligned}
 \eta_{d,h_k-1}&=n_k-4,\\
 \eta_{d,h_k}=\eta_{d,h_k+1}&=n_k,\\
 \eta_{d,h_k+2}&=n_k+4.
 \end{aligned}
\end{equation}
Moreover, for the same \(k\) and \(0\le j<J\), set
\(h_{k,j}=h_k-2j\) and \(n_{k,j}=n_k-4j\). Then
\begin{equation}\label{eq:density-pairs}
 \eta_{d,h_{k,j}}=\eta_{d,h_{k,j}+1}=n_{k,j}.
\end{equation}
Every minimum in these formulas has a balanced rational attainer
with connected nonincidence graph and exactly the requested deficiency.
At \(h=h_{k,j}\), writing \(n=n_{k,j}\), the possible extremal
count data are exactly those in Table~\ref{tab:density-plateau-data}:
\begin{datadisplay}{Face counts and matchings at a capacity plateau}\label{tab:density-plateau-data}
\[
\begin{array}{c|c|c|c}
 f_0(Q)&f_{d-1}(Q)&\delta(Q)&\nu(\N(Q))\\ \hline
 n&n&h&n-h\\
 n-1&n&h&n-h-1\\
 n&n&h+1&n-h-1
\end{array}
\]
\end{datadisplay}
Each row is realized by a rational polytope with connected
nonincidence graph.
\end{theorem}

The density statement concerns the capacity indices \(k\), not all
deficiencies: \(h_k\) grows as a degree-\(t\) polynomial. The table
classifies count data, rather than face lattices. The proof first
evaluates \(\chi_h\), and then shows that the witness threshold
\(s_h^*\) equals \(s_k\) for every deficiency in the statement. Classical irrational-rotation
equidistribution supplies the density statement.

\begin{proof}[Proof of Theorem~\ref{thm:density-plateaus}]
\emph{Step 1: the complementary facet minima.}
The consecutive differences of \(U(k)-k-1\) tend to infinity.
Thus, at each fixed shift \(h=h_k+\ell\) with \(\ell\le1\)
occurring in the statement, the baseline in
\eqref{eq:odd-baseline} is \(b_h^0=k\) for sufficiently large
\(k\). Indeed, the inequality at \(k\) is equivalent to
\(\ell\le1\), whereas the preceding capacity is eventually too
small. At \(h=h_k+2\), the baseline is \(k+1\).

For \(h=h_k-2j\), the baseline facet count is
\(U(k)-(2j+1)\). It is odd, and \(2j+1<t\) for \(0\le j<J\).
Simplicial odd-dimensional polytopes have an even number of facets;
the nonsimplicial bound \(U_m(k)-t\) rules out this pair as well
\cite[proof of Theorem~3.3]{sjobergZiegler2020}.
Hence \(\chi_h=1\) and \(\beta_h^*=U(k)-2j\).
At \(h+1\), the baseline is already this even count, so
\(\chi_{h+1}=0\) and the same facet minimum results.
Similarly, the minima at \(h_k-1\) and \(h_k+2\) are
\(U(k)-2\) and \(U(k)+2\), respectively. All the even pairs
are rationally realized by Lemma~\ref{lem:capacity-interpolation}:
their facet counts have order \(k^t\), while the lower interpolation
endpoint has order \(k^{t-1}\).

\emph{Step 2: the common witness threshold.}
The inequality defining \(s_h^*\) is
\(P(s)\ge\beta_h^*+h\). For the four consecutive deficiencies
in \eqref{eq:density-four}, its targets are
\[
 C_k-3,\quad C_k,\quad C_k+1,\quad C_k+4.
\]
For the \(j\)-th pair in \eqref{eq:density-pairs}, they are
\(C_k-4j\) and \(C_k-4j+1\). Every target lies between
\(C_k-A_t\) and \(C_k+4\). Condition~\eqref{eq:density-good}
and eventual monotonicity of \(P\) show that the minimal witness
threshold is \(s_k\) throughout. Theorem~\ref{thm:odd-integer}
therefore gives all the displayed minima and balanced rational
connected attainers of exact deficiency.

\emph{Step 3: density of the good indices.}
Both \(P(x)\) and \(C(x)=2U_m(x)-x-2\) are degree-\(t\)
polynomials; their leading coefficients are \(1/t!\) and \(4/t!\).
Let \(R_k\) be the unique large positive solution of
\(P(R_k)=C_k\). Comparing the first two coefficients gives
\begin{equation}\label{eq:density-root}
 R_k=\theta k+c_t+O_t(k^{-1}),\qquad \theta=4^{1/t},
\end{equation}
where \(c_t\) is a fixed real number. To justify the error, write
\[
 P(x)=\mu_0x^t+\mu_1x^{t-1}+O_t(x^{t-2}),\qquad
 C_k=\gamma_0k^t+\gamma_1k^{t-1}+O_t(k^{t-2}).
\]
Take \(\theta=(\gamma_0/\mu_0)^{1/t}\) and
\(c_t=(\gamma_1-\mu_1\theta^{t-1})/(t\mu_0\theta^{t-1})\).
Substitution leaves an error \(O_t(k^{t-2})\); the derivative
of \(P\) is comparable to \(k^{t-1}\) near \(\theta k\).
The mean value theorem proves \eqref{eq:density-root}.
On the increasing tail, \(s_k=\lceil R_k\rceil\).

If the upper slack in \eqref{eq:density-good} fails, then
\(s_k-R_k=O_t(k^{1-t})\). If its lower gap fails, then
\(R_k-(s_k-1)=O_t(k^{1-t})\). These estimates follow from
the same derivative bound on the intervening unit interval and include
integral roots. Consequently every sufficiently large bad index satisfies
\[
 \operatorname{dist}(\theta k+c_t,\mathbb Z)=O_t(k^{-1}).
\]
For \(t\ge3\), \(\theta\) is irrational: a rational root of
\(x^t-4\) would be an integer, and no integer has \(t\)-th
power four. The sequence \(\theta k+c_t\) is uniformly distributed
modulo one by Weyl's criterion \cite{weyl1916}: for every nonzero
integer \(\ell\), the number \(\ell\theta\) is not an integer, so the
finite geometric-series formula bounds
\(\bigl|\sum_{k=1}^K e^{2\pi i\ell(\theta k+c_t)}\bigr|\) by
\(1/\lvert\sin(\pi\ell\theta)\rvert\), and the averages tend to zero. For any
\(0<\eps<1/2\), the bad indices eventually lie in the fixed
\(\eps\)-neighborhood of the integers, of limiting proportion
\(2\eps\). Letting \(\eps\) tend to zero proves that their
natural density is zero.

\emph{Step 4: equality counts and their rational realization.}
At \(h=h_{k,j}\), put \(\beta_*=U(k)-2j\).
Then \(\beta_*-h-1=k+1\).
The equality data in Theorem~\ref{thm:odd-integer} give
\[
 b\in\{k,k+1\},\quad z=k+1-b,\quad 0\le w\le z,
 \qquad \alpha+\rho=s_k+r.
\]
They imply \(f_{d-1}(Q)=n\), \(f_0(Q)=n+w-z\),
\(\delta(Q)=h+w\), and \(\nu(\N(Q))=n-h-z\).
The three choices \((z,w)=(0,0),(1,0),(1,1)\) give exactly
the table.

To realize them, take \(\alpha=s_k\) and \(\rho=v=r\).
For \(w=0,1\), interpolate a rational simplicial \(r\)-polytope
\(E_w\) with \(s_k-1\) vertices and
\[
 s_k+\beta_*+h+w=s_k+C_k-4j+w
\]
facets. The upper slack in \eqref{eq:density-good} keeps both counts
at most \(U_r(s_k-1)\); their order \(k^t\) dominates the lower
interpolation endpoint. Since \(r\) is even, every intervening
integer facet count is available. Strictly stack a simplex facet of
\(E_w\) and polarize to obtain \(A_w\), with a simplex stacking
facet and
\[
 a=s_k+\beta_*+h+r-1+w,\qquad \alpha=s_k.
\]
Choose rational complements \(B_k,B_{k+1}\) with their indicated
vertex counts and \(\beta_*\) even facets. Their interpolation
conditions hold by Step~1 and monotonicity of the upper capacity.
Use \((A_0,B_{k+1})\), \((A_0,B_k)\), and
\((A_1,B_k)\), respectively. The certificate margins are
\(h+w\) and \(h+z\), and final strict stacking gives the
three count rows. Theorem~\ref{thm:odd-integer} supplies their
exact matching and deficiency; strict stacking supplies connectedness.
\end{proof}

For \(t=2\), the scale ratio \(4^{1/t}=2\) is rational, so
the density argument is replaced by the explicit calculation in
Corollary~\ref{cor:ten-plateau}. Its final jump is five, whereas
\eqref{eq:density-four} has a final jump of four. In higher dimensions
the parity obstruction persists while almost every capacity index
avoids the narrow witness-rounding boundary.

\section{The complete deficiency-one optimum}
\label{sec:lowdim}

The exact onset theorem leaves only four dimensions between the
seven-dimensional construction and the stable range. Two additional
witness geometries settle them: the polar of a two-facet extension of
the cyclic polytope \(C(7,4)\), truncated at a vertex, and a product of
three simplices with a controlled vertex cut. Throughout this section,
\(\eta_d\) refers to the certified class of
Definition~\ref{def:certified-class}.

\subsection{Four-dimensional capacity and witnesses}

\begin{lemma}[Four-dimensional capacity]\label{lem:four-capacity}
A four-polytope with at most seven facets has at most fourteen vertices.
A four-polytope with eight facets and a tetrahedral facet has at most
seventeen vertices.
\end{lemma}

\begin{proof}
McMullen's Upper Bound Theorem \cite{mcmullen1970} gives at most
\(n(n-3)/2\) facets for a four-polytope with \(n\) vertices; the
cyclic polytope attains this number. Polarity gives the first assertion.

For the second, take the polar \(D\), with eight vertices, and the
simple vertex \(x\) dual to the tetrahedral facet. Deleting \(x\)
leaves a full-dimensional four-polytope \(D'\): otherwise \(D\)
would be a pyramid with simple apex, whose base must be a tetrahedron,
contrary to its vertex count. The facets of \(D\) not through \(x\)
persist in \(D'\). At least one facet of \(D'\) is strictly visible
from \(x\), and is not such a persistent facet. Since \(x\) lies on
four facets,
\[
 f_3(D)\le4+f_3(D')-1\le4+14-1=17.
\]
\end{proof}

\begin{lemma}[A cyclic witness]\label{lem:cyclic-witness}
There is a rational four-polytope \(A\) with twenty-one vertices,
nine facets and a tetrahedral facet \(F\). Its eight other facets
can be matched to distinct nonincident vertices, using all four
vertices of \(F\) and four vertices outside \(F\).
\end{lemma}

\begin{proof}
Start with the rational cyclic polytope \(C(7,4)\), with fourteen
tetrahedral facets \cite{mcmullen1970,ziegler1995}. Choose two adjacent
facets meeting in a triangular ridge, and a rational point in the
relative interior of that ridge. A sufficiently small rational
displacement of this point lies strictly beyond precisely those two
facets and strictly beneath all others. Indeed, the two supporting normals are
independent, and every other facet inequality has a strict margin.
Choose the displacement small enough to preserve every old vertex:
a linear functional strictly exposing any old vertex has a strict
gap at the selected ridge point, and the finitely many gaps persist.

The two visible tetrahedra have six horizon triangles. Their replacement
by cones gives a simplicial four-polytope \(D\) with eight vertices
and \(14-2+6=18\) facets. This is the strict visible/retained extension
description of Section~\ref{sec:preliminaries}. Its polar \(E\) is
simple with eighteen vertices and eight facets. A sufficiently close
rational truncation of one vertex gives \(A\) with twenty-one vertices
and nine facets; the new facet \(F\) is a tetrahedron.

The four vertices of \(F\), one on each old edge, match to the four
old facets through the truncated vertex: each edge point misses the
unique one of these facets not containing its edge. There are four
old facets not through that vertex. Each is a simple three-polytope
with at most seven facets, hence at most ten vertices by Euler's
formula. Among the seventeen retained old vertices, each therefore
misses at least seven. Greedily choosing four distinct missed vertices
matches these remaining facets outside \(F\).
\end{proof}

\begin{proposition}[Dimensions eight and nine]\label{prop:eight-nine}
\(\eta_8=\eta_9=28\). Each value is attained by a balanced rational
polytope with connected nonincidence graph and deficiency one.
\end{proposition}

\begin{proof}
For attainment, use the witness of Lemma~\ref{lem:cyclic-witness}:
\(r=4\), \(a=21\), \(\alpha=9\), \(\rho=4\).
Take \(B\) to be the octahedron for \(d=8\), and the four-dimensional
tetrahedral bipyramid for \(d=9\). Both have \(b=6\), \(\beta=8\)
and an assignment covering their six vertices. In either case,
\[
 f_0(Q)=f_{d-1}(Q)=28,\qquad L=b+\alpha+\beta+\rho=27.
\]
The inherited certificate gives \(\nu(\N(Q))\le27\).
Use the eight witness pairs of Lemma~\ref{lem:cyclic-witness}, then
match the twelve cone facets to twelve of the thirteen unused witness
vertices outside \(F\). These vertices miss every cone facet.
Match the six vertices of \(B\) to retained facets of type \(A*G\),
and the new apex to one of the two unused such facets. The disjoint
blocks have sizes \(8,12,7\), yielding twenty-seven edges. Strict
stacking gives connectedness by Corollary~\ref{cor:stacking-transfer}.

For the lower bound in dimension eight, write \(r+m=7\).
Both dimensions are at least three. If \(r=3,m=4\),
Corollary~\ref{cor:class-optimal} gives size at least thirty-three.
Otherwise \(r=4,m=3\). Euler's formula and \(\beta-b\ge2\)
give \(b\ge6\), \(\beta\ge8\). Put \(u=\alpha+\rho\); then \(u\ge11\)
by Lemma~\ref{lem:witness-facet-floor}.
If the larger class size were at most twenty-seven, then
\(u+2\beta-1\le27\), forcing \(\beta=8\), \(u\le12\),
and \(b=6\). The certificate requires \(a\ge u+8\ge19\).
If \(\alpha\le7\), Lemma~\ref{lem:four-capacity} gives \(a\le14\).
If \(\alpha=8\), then \(\rho=4\) and \(F\) is a tetrahedron,
giving \(a\le17\). Larger \(\alpha\) is incompatible with
\(\rho\ge4\), \(u\le12\). All cases contradict \(a\ge19\).

For dimension nine, suppose the larger class size is at most twenty-seven.
Write
\[
 b=m+2+s,\quad \beta=m+4+s+z,\quad
 \alpha+\rho=2r+3+w,
\]
with \(s,z,w\ge0\) integers, by Lemma~\ref{lem:witness-facet-floor}
and \(\beta-b\ge2\); the latter also shows that \(B\) is not a simplex. The facet count is
\(26+w+2s+2z\le27\), so \(s=z=0\), \(w\le1\).
Now \(B^*\) has \(m+2\) facets and \(m+4\) vertices.
McMullen's few-facet classification \cite{mcmullen1971} gives
\(B^*\cong\operatorname{Pyr}^{m-u-v}(\Delta_u\times\Delta_v)\) with
\(m+1+uv=m+4\), so \(\{u,v\}=\{1,3\}\) and \(m\ge4\). Since \(r+m=8\), either \(r=3,m=5\),
where Corollary~\ref{cor:class-optimal} gives size at least thirty-six,
or \(r=m=4\). In the latter case \(\alpha+\rho\le12\),
\(\beta=8\), and \(a\ge\alpha+\rho+8\ge19\), which contradicts
Lemma~\ref{lem:four-capacity} by the case analysis used in dimension
eight. Neither lower bound assumes balance.
\end{proof}

\subsection{Products of three simplices}

\begin{lemma}[A three-factor controlled cut]\label{lem:triple-cut}
Let \(p,q,s\ge1\), \(r=p+q+s\), and
\(T=pq+ps+qs+pqs\). For every integer \(0\le e\le r\), there
is a rational \(r\)-polytope \(A\) with
\[
 a=2r+T-e,\qquad \alpha=r+4,
\]
and a simplex facet \(F\) with \(\rho=r\) ridges. The \(r+3\)
other facets match to distinct nonincident vertices, using all \(r\)
vertices of \(F\) and three vertices outside it.
\end{lemma}

\begin{proof}
Use \(E=\Delta_p\times\Delta_q\times\Delta_s\), with vertices
\(a_{ijk}\), based at \(a_{000}\). This simple polytope has
\(r+3\) facets and \((p+1)(q+1)(s+1)=r+1+T\) vertices.
Choose \(e\) of the \(r\) edge-neighbors of \(a_{000}\).
Assign an affine rational functional \(\lambda\) value zero at
\(a_{000}\), one at the chosen neighbors and two at the others.
Additivity across the product defines it uniquely. Every other old
vertex has value at least one, and one varying in at least two factors
has value at least two. Cut by \(\lambda\ge1\).

Only \(a_{000}\) is removed. The cut meets every incident edge,
at \(e\) old endpoints and \(r-e\) new points; the new facet is a
simplex with \(r\) vertices. Every old facet through \(a_{000}\)
contains a product rectangle through it, since the two other factors
have positive dimension. The other three vertices of that rectangle
span \(a_{000}\) affinely. Removing it thus preserves the affine
span of every old facet, and all \(r+3\) old facets survive. This
gives the stated counts.

The \(r\) cut-facet vertices match to the \(r\) old facets through
\(a_{000}\), using the unique missed facet on each edge.
The remaining three facets are the opposite facets in the three
factors. They match respectively to \(a_{011},a_{101},a_{110}\).
These are distinct, nonincident to their assigned facets, and lie
outside \(F\) because \(\lambda\ge2\).
\end{proof}

\begin{proposition}[Dimensions ten and eleven]\label{prop:ten-eleven}
\(\eta_{10}=29\) and \(\eta_{11}=31\). Both values have balanced
rational attainers with connected nonincidence graph and deficiency one.
\end{proposition}

\begin{proof}
Take \(B\) to be the four-dimensional tetrahedral bipyramid, with
\(b=6\), \(\beta=8\). Use Lemma~\ref{lem:triple-cut} with
\((p,q,s,e)=(1,2,2,0)\) and \((1,2,3,5)\), respectively. The
resulting data are listed in Table~\ref{tab:lowdim-witnesses}:
\begin{datadisplay}{Witness data in dimensions ten and eleven}\label{tab:lowdim-witnesses}
\[
\begin{array}{c|rrrr|rr}
 d&r&a&\alpha&\rho&f_0(Q)=f_{d-1}(Q)&L\\ \hline
 10&5&22&9&5&29&28\\
 11&6&24&10&6&31&30.
\end{array}
\]
\end{datadisplay}
The factor certificate gives \(\nu(\N(Q))\le L\).
Match the retained witness facets using Lemma~\ref{lem:triple-cut}.
This gives eight and nine pairs, respectively, and leaves
\(a-r-3=14\) and \(15\) outside witness vertices. Match the
\(\beta+\rho=13\) and \(14\) cone facets to distinct such vertices.
Finally match the six complementary vertices and the apex as in
Proposition~\ref{prop:eight-nine}. The three blocks attain \(L\).
Strict stacking gives connectedness, so both upper bounds hold.
The certified floors are twenty-eight and thirty, and
Theorem~\ref{cor:stable-h} excludes their attainment below dimension
twelve. By integrality, \(\eta_{10}\ge29\) and \(\eta_{11}\ge31\).
\end{proof}

\begin{theorem}[Complete deficiency-one optimum]\label{thm:complete-one}
The class \(\mathcal C_{d,1}\) is empty for \(d<7\). For every
\(d\ge7\),
\begin{equation}\label{eq:complete-one}
 \eta_d=
 \begin{cases}
 33,&d=7,\\
 28,&d=8,9,\\
 29,&d=10,\\
 31,&d=11,\\
 2d+8,&d\ge12.
 \end{cases}
\end{equation}
Every finite value is attained by a balanced rational polytope with
connected nonincidence graph and matching number \(\eta_d-1\).
\end{theorem}

\begin{proof}
Both factors have dimension at least three
(Proposition~\ref{prop:factor-criterion}); this excludes \(d<7\) and
forces \(r=m=3\) when \(d=7\). Combine Corollary~\ref{cor:class-optimal},
Propositions~\ref{prop:eight-nine} and~\ref{prop:ten-eleven}, and
Corollary~\ref{cor:stable-one}.
\end{proof}

\section{The infinite rational family}
\label{sec:construction}

Fix integers \(d\ge7\) and \(n\ge2^{d-4}+6\), and set \(m=d-4\ge3\), so that
\(n\ge2^m+6\ge14\). Throughout Sections~\ref{sec:construction}
and~\ref{sec:matching}, \(n\) and \(m\) denote these parameters, not
the class sizes of Section~\ref{sec:defect}, and the prism \(A_n\) and
polygon \(C_n\) are specific to this coordinate family. Let
\[
  C_n=\conv\{(k,k^2):0\le k\le n-1\}\subset\R^2
\]
and let \(A_n=C_n\times[-1,1]\subset\R^3\).  We denote its vertices by
\[
  A_{k,s}=(k,k^2,s),
  \qquad 0\le k\le n-1,\quad s\in\{-1,1\}.
\]
The polygon vertices lie on a strictly convex parabola, so their displayed
order is their cyclic order.  The facets of \(A_n\) are the top and bottom
facets \(z\le1\), \(-z\le1\), together with
\begin{align}
  (2k+1)x-y &\le k(k+1), &&0\le k\le n-2,\label{eq:lateral}\\
  -(n-1)x+y &\le0. &&\label{eq:closing}
\end{align}
We index the lateral facets cyclically as \(X_0,\ldots,X_{n-1}\), with
\(X_{n-1}\) denoting \cref{eq:closing}, and write \(T_+\) and \(T_-\) for
the top and bottom facets.

Let \(B=\conv\{\pm e_1,\ldots,\pm e_m\}\subset\R^m\) be the
\(m\)-dimensional cross-polytope; for \(d=7\) it is the octahedron.
In coordinates \((x,y,z,u_1,\ldots,u_m,t)\) on \(\R^d\), embed the vertex
sets of \(A_n\) and \(B\) as
\begin{align*}
 \widetilde A_n
   &=\{(k,k^2,s,0,\ldots,0,0):
       0\le k\le n-1,\ s\in\{-1,1\}\},\\
 \widetilde B
   &=\{(0,0,0,\sigma e_i,1):
       1\le i\le m,\ \sigma\in\{-1,1\}\},
\end{align*}
and set
\begin{equation}\label{eq:free-join}
  P_{d,n}=\conv(\widetilde A_n\cup\widetilde B),
\end{equation}
so that \(P_{d,n}\) is the free join \(A_n*B\) in the realization of
Section~\ref{sec:preliminaries}.
We use the same symbols \(A_{k,s}\) for the first vertex class and write
\(B_{i,\sigma}\) for the second.

Every facet inequality \(\alpha\cdot(x,y,z)\le c\) of \(A_n\) induces
\begin{equation}\label{eq:x-lift}
  \alpha\cdot(x,y,z)+ct\le c
\end{equation}
on \(P_{d,n}\).  For \(\eps=(\eps_1,\ldots,\eps_m)\in\{-1,1\}^m\), the
\(2^m\) facets induced by \(B\) are
\begin{equation}\label{eq:y-facets}
  Y_\eps:\qquad \eps_1u_1+\cdots+\eps_mu_m-t\le0.
\end{equation}
We abbreviate sign tuples in subscripts; for example,
\(Y_{--+}=Y_{(-1,-1,1)}\) when \(m=3\).
Thus \(P_{d,n}\) has \(2n+2m\) vertices, dimension \(d\), and \(n+2+2^m\)
facets.

Let \(R_0=\conv\{A_{0,-},A_{0,+},A_{1,-},A_{1,+}\}\).  We stack beyond the
facet \(S=X_0=R_0*B\), supported by
\[
  g_0=x-y\le0.
\]
The point
\[
  q=\left(\frac15,\frac15,0,\ldots,0,\frac35\right)
   =\frac25\left(\frac12,\frac12,0,0,\ldots,0,0\right)
   +\frac35(0,0,0,0,\ldots,0,1)
\]
is a convex combination, with positive weights, of the centers of
\(R_0\) and of the copy of \(B\); hence it lies in the relative interior
of \(S=R_0*B\).  Set
\begin{equation}\label{eq:apex}
  p=q+\Bigl(\frac1{100},-\frac1{100},0,\ldots,0,0\Bigr).
\end{equation}

\begin{proposition}[Uniform visibility]\label{prop:visibility}
The point \(p\) lies strictly beyond \(S\) and strictly beneath every other
facet of \(P_{d,n}\).
\end{proposition}

\begin{proof}
The value of \(g_0\) at \(p\) is \(1/50>0\).  For
\(1\le k\le n-2\), the value at \(p\) of
\((2k+1)x-y+k(k+1)t-k(k+1)\), the lifted form of \cref{eq:lateral}, is
\begin{equation}\label{eq:lateral-slack}
  -\frac{2k^2}{5}+\frac{2k+2}{100}<0;
\end{equation}
after multiplication by \(100\), the left-hand side is
\(-40k^2+2k+2\le-36\).  The value of the closing form is
\begin{equation}\label{eq:closing-slack}
  -\frac{n-2}{5}-\frac{n}{100}<0.
\end{equation}
The top and bottom forms both take the value \(-2/5\), and every
\(Y_\eps\) form takes the value \(-3/5\).  By the facet description of
free joins, these forms exhaust the facets of \(P_{d,n}\), and none of
the values depends on \(m\).
\end{proof}

Define
\begin{equation}\label{eq:qn}
  Q_{d,n}=\conv(P_{d,n}\cup\{p\}).
\end{equation}

\begin{proposition}[Complete facet count]\label{prop:facet-count}
The polytope \(Q_{d,n}\) is rational and \(d\)-dimensional, with
\(2n+2d-7\) vertices and \(n+5+2^{d-3}\) facets.
\end{proposition}

\begin{proof}
The selected facet \(S\) is the free join of a quadrilateral and the
cross-polytope \(B\).  Its ridges correspond to the four edges of the
quadrilateral and the \(2^m\) facets of \(B\), so it has \(4+2^m\) ridges.
Proposition~\ref{prop:visibility} and Lemma~\ref{lem:stacking} show that
stacking adds one vertex, removes \(S\), and adds \(4+2^m\) cone facets.
Therefore
\[
  f_0(Q_{d,n})=2n+2m+1=2n+2d-7,
  \qquad
  f_{d-1}(Q_{d,n})=n+2+2^m-1+4+2^m=n+5+2^{d-3}.
\]
All input coordinates are rational.  The old facets are supported by
\cref{eq:x-lift,eq:y-facets}, and the new facets by \cref{eq:cone-form}, so
the description is complete and exact.
\end{proof}

For later use, the four old facets meeting \(S\) in the ridges
\(e*B\), with \(e\) an edge of \(R_0\), are \(T_+,T_-,X_1,X_{n-1}\).
Together with the \(2^m\) facets \(Y_\eps\), which meet \(S\) in the
ridges \(R_0*G\), with \(G\) a facet of \(B\), these are the \(4+2^m\) retained facets adjacent
to \(S\).  We denote the
corresponding new cone facets by
\[
 \widehat T_+,\ \widehat T_-,\ \widehat X_1,\
 \widehat X_{n-1},\ \widehat Y_\eps.
\]

\section{Exact matching deficiency}
\label{sec:matching}

The copy of \(A_n\) in \(P_{d,n}\) survives the stacking operation as a face
of \(Q_{d,n}\): it is the intersection of \(Q_{d,n}\) with the \(2^m\)
supporting hyperplanes of the retained facets \(Y_\eps\), and the apex
satisfies all \(2^m\) inequalities strictly.  No retained \(X_j\) or
\(T_\pm\) contains all vertices of \(A_n\).  Every new cone facet meets the
old \(A_n\)-vertex class inside \(R_0\), so no new facet contains \(A_n\).
Hence precisely the \(2^m\) facets \(Y_\eps\) contain this face.  This is
also the situation covered by Corollary~\ref{cor:stacking-transfer}.

\begin{proposition}\label{prop:matching-upper}
Every matching in \(\N(Q_{d,n})\) has at most \(n+2^{d-4}+2d-2\) edges.
Since
\[
  f_0(Q_{d,n})-\bigl(n+2^{d-4}+2d-2\bigr)=n-2^m-5\ge1
\]
and
\[
  f_{d-1}(Q_{d,n})-\bigl(n+2^{d-4}+2d-2\bigr)=2^m-2m-1\ge1,
\]
the polytope \(Q_{d,n}\) has no vertex--facet assignment.
\end{proposition}

\begin{proof}
In \(P_{d,n}\), the face \(A_n\) has \(r=2n\) vertices and is contained in
\(c=2^m\) facets.  The polytope \(P_{d,n}\) has \(v=2n+2m\) vertices and
\(f=n+2+2^m\) facets, the selected facet does not contain \(A_n\), and it
has \(\rho=4+2^m\) ridges.  The transfer bound in
  Corollary~\ref{cor:stacking-transfer} therefore gives
\[
 \nu\bigl(\N(Q_{d,n})\bigr)
 \le (2n+2m)+(n+2+2^m)+(4+2^m)-2n-2^m
 =n+2m+2^m+6,
\]
which equals \(n+2^{d-4}+2d-2\).  The two displayed differences are
\(n-2^m-5\) and \(2^m-2m-1\); the first is at least one because
\(n\ge2^m+6\), and the second is at least one because \(m\ge3\).  An
assignment would be a matching of size
\(\min\{f_0,f_{d-1}\}>n+2^{d-4}+2d-2\), which the bound rules out.
\end{proof}

The bound is exact: we now exhibit a matching of size \(n+2m+2^m+6\) by a
uniform rule.  The incidence facts used are the following; each one reads
off directly from
\cref{eq:x-lift,eq:y-facets,eq:cone-form}.
\begin{enumerate}
  \item \(X_j\) contains all \(2m\) vertices of the \(B\)-class and the two
        vertical pairs \(A_{j,\pm}\), \(A_{j+1,\pm}\) over the endpoints of
        its polygon edge (indices modulo \(n\)).
  \item \(T_s\) contains all \(B\)-vertices and the \(n\) vertices
        \(A_{k,s}\).
  \item \(Y_\eps\) contains all \(2n\) vertices of the \(A\)-class and the
        \(m\) vertices \(B_{i,\eps_i}\); the apex \(p\) lies on no
        \(Y_\eps\).
  \item \(\widehat T_s\) contains \(p\), all \(B\)-vertices, and
        \(A_{0,s},A_{1,s}\); the facet \(\widehat X_1\) contains \(p\), all
        \(B\)-vertices, and \(A_{1,\pm}\); the facet \(\widehat X_{n-1}\)
        contains \(p\), all \(B\)-vertices, and \(A_{0,\pm}\); and
        \(\widehat Y_\eps\) contains \(p\), the four vertices of \(R_0\),
        and the \(m\) vertices \(B_{i,\eps_i}\).  In particular, no cone
        facet contains a vertex \(A_{k,s}\) with \(k\ge2\).
\end{enumerate}

\begin{proposition}\label{prop:matching-lower}
The graph \(\N(Q_{d,n})\) contains the following matching of size
\(n+2m+2^m+6\):
\begin{enumerate}[label=\textup{(M\arabic*)}]
  \item\label{m:lateral}
        \(A_{j-1,-}\longleftrightarrow X_j\) for \(1\le j\le n-1\);
  \item\label{m:top}
        \(A_{n-1,-}\longleftrightarrow T_+\) and
        \(A_{0,+}\longleftrightarrow T_-\);
  \item\label{m:cone}
        enumerate the \(2^m+4\) cone facets in any fixed order as
        \(G_1,\ldots,G_{2^m+4}\) and take
        \(A_{i+1,+}\longleftrightarrow G_i\) for \(1\le i\le 2^m+4\);
  \item\label{m:cross}
        \(B_{i,+}\longleftrightarrow Y_\eps\) with \(\eps_i=-1\) and all
        other entries \(+1\), and
        \(B_{i,-}\longleftrightarrow Y_\eps\) with \(\eps_i=+1\) and all
        other entries \(-1\), for \(1\le i\le m\); and
        \(p\longleftrightarrow Y_{(1,\ldots,1)}\).
\end{enumerate}
Consequently \(\nu\bigl(\N(Q_{d,n})\bigr)=n+2^{d-4}+2d-2\).
\end{proposition}

\begin{proof}
Each pair is a nonincidence.  In \ref{m:lateral}, the vertex \(A_{j-1,-}\)
lies on the lateral facets \(X_{j-2}\) and \(X_{j-1}\) only (indices modulo
\(n\)), hence not on \(X_j\).  In \ref{m:top} the signs disagree with
rule~(2).  In \ref{m:cone}, the used vertices are \(A_{k,+}\) with
\(2\le k\le2^m+5\le n-1\), and by rule~(4) no cone facet contains such a
vertex.  In \ref{m:cross}, the vertex \(B_{i,\sigma}\) lies on \(Y_\eps\)
exactly when \(\eps_i=\sigma\), and each listed pattern was chosen with
\(\eps_i=-\sigma\); the apex lies on no \(Y_\eps\).

All vertices used are distinct: the minus class contributes
\(A_{0,-},\ldots,A_{n-2,-}\) in \ref{m:lateral} and \(A_{n-1,-}\) in
\ref{m:top}; the plus class contributes \(A_{0,+}\) in \ref{m:top} and
\(A_{2,+},\ldots,A_{2^m+5,+}\) in \ref{m:cone}; the vertices in
\ref{m:cross} are the \(B\)-class and the apex.  All facets used are
distinct: the retained facets \(X_1,\ldots,X_{n-1},T_+,T_-\), all
\(2^m+4\) cone facets, and \(2m+1\) pairwise distinct facets \(Y_\eps\) in
\ref{m:cross}, since for \(m\ge3\) the sign patterns with exactly one
\(-1\), those with exactly one \(+1\), and the all-plus pattern are pairwise
distinct.  The total count is
\[
  (n-1)+2+(2^m+4)+(2m+1)=n+2m+2^m+6.
\]
Together with Proposition~\ref{prop:matching-upper}, this proves
\begin{equation}\label{eq:matching-number}
  \nu\bigl(\N(Q_{d,n})\bigr)=n+2m+2^m+6=n+2^{d-4}+2d-2.
\end{equation}
\end{proof}

The matching misses exactly the \(2^m-2m-1\) facets \(Y_\eps\) whose sign
pattern is neither all-plus, nor has exactly one \(-1\), nor has exactly one
\(+1\); for \(d=7\) the unique missed facet is \(Y_{(-1,-1,-1)}\).
Corollary~\ref{cor:stacking-transfer} gives connectedness of
\(\N(Q_{d,n})\) and excludes every nontrivial free-join decomposition.
Together with Proposition~\ref{prop:facet-count} and
\cref{eq:matching-number}, this proves Theorem~\ref{thm:main}; the case
\(d=7\), where \(m=3\), is Corollary~\ref{cor:seven}. Consequently
\[
 \delta(Q_{d,n})=\min\{n-2^{d-4}-5,\;2^{d-4}-2d+7\},
\]
which equals \(1\) when \(d=7\).

\section{Extremal geometry and further questions}
\label{sec:discussion}

Connectedness changes the geometry of an assignment obstruction.
The obstruction can no longer be explained solely by opposite
imbalances in free-join components; the ambient incidences of a Hall
face must support the deficiency within one connected graph. In
dimension seven, the sharp value \(\omega_{7,1}=17\) gives the exact
cost in both class sizes. Shared-coordinate hulls and coplanar
extensions attain it by different mechanisms, and the resulting at
least sixteen equality types show that minimal size does not force
a unique combinatorial structure.

Rigidity instead appears after fixing the Hall-face data.
A zero-slack maximizing pair determines the whole face lattice up to
polarity (Theorem~\ref{thm:zero-seventeen-rigidity}). In the
triple-coface branch, extension rows determine a marked six-vertex seed,
and the seven marked seeds give the complete converse
(Theorem~\ref{thm:triple-global-classification}). The incidence-only
support injection explains one obstruction common to the lower-bound
arguments: a nonzero extension vector at a simple base vertex determines
an outside neighbor, and distinct base vertices require distinct
neighbors (Theorem~\ref{thm:incidence-support-injection}).
Its three-dimensional consequence is \(\alpha-3\le s\).
The ambient nonvanishing arguments in
Propositions~\ref{prop:one-simple-double-slack}
and~\ref{prop:almost-simple-unit-slack} are therefore as important
as the local row-space calculation.

\subsection{Prescribed deficiency and construction constraints}

The unrestricted construction realizes every positive deficiency in
dimension seven, but existence alone does not determine its minimum
size. The certified problem fixes more geometry: an inherited free-join
witness and a strict stacking facet. Admissible factor counts must then
be realized with that facet, and the surviving-face bound must be
attained by a matching. Solving these three requirements together gives
exact optima for the certified class. The distinction is already visible
at deficiency one: \(\omega_{7,1}=17\), whereas \(\eta_{7,1}=33\).

The exact first-attainment dimension makes the arithmetic of the
certificate visible. The factorization \(uv=h+2\), followed by the
witness capacity, determines \(D_h\). For example, if \(h+2\) is prime,
\[
 D_h=h+4+\left\lceil2\sqrt{3h+8}\right\rceil;
\]
if \(h=k^2-2\), then
\[
 D_h=2k+1+\left\lceil2\sqrt{2k^2+2k+1}\right\rceil
       =(2+2\sqrt2)k+O(1).
\]
Thus the square sequence attains the certified size floor with
quadratic deficiency in linear dimension. This is an assertion about
the sharp certified floor, rather than the dimension needed for
unrestricted existence.

For fixed dimension, classical face-vector realization
\cite{billeraLee1980,stanley1980,sjobergZiegler2020} supplies the
capacity inputs. The distinguished stacking facet and matching
saturation turn these counts into assignment obstructions.
Theorem~\ref{thm:second-order} determines all divergent terms of the
resulting size law, and Corollary~\ref{cor:capacity-stability}
bounds the remaining near-extremal budgets. At the integer level,
coplanar stacking removes witness parity, giving exact eventual
formulas for even-dimensional complements
(Theorems~\ref{thm:even-integer} and~\ref{thm:odd-witness-integer}).
Odd-dimensional complements retain a face-pair feasibility condition
(Theorem~\ref{thm:odd-integer}). Their classical facet gaps cause
plateaus: an example of deficiency \(h+1\) can already minimize size
among examples of deficiency at least \(h\). The explicit
ten-dimensional plateaus and the density-one families in
Corollary~\ref{cor:ten-plateau} and
Theorem~\ref{thm:density-plateaus} distinguish this integer phenomenon
from a bounded asymptotic error.

\subsection{The next numerical and structural problems}

For \(h\ge2\), the unrestricted seven-dimensional problem remains
between
\[
 3h+14\le\omega_{7,h}\le5h+12.
\]
Determining its leading growth constant is a natural first step.
In dimensions eight and nine, the shared-coordinate examples give
additional unrestricted upper bounds; in each fixed dimension at
least nine, certified constructions are smaller for sufficiently
large deficiency. A better unrestricted construction must use
Hall-face geometry beyond an optimization of either existing family.

The lower boundary also poses a classification problem.
For \(h\ge2\), Corollary~\ref{cor:final-seven-boundary} forces simple
dual Hall faces with count slacks \((2,1)\), up to polarity.
Theorem~\ref{thm:triple-coface-exclusion} excludes assigning both
surplus vertices to one coface. The remaining allocations are two
hidden vertices, a hidden vertex and a double coface, or two double
cofaces. Deleting the columns of exceptional inducing facets may lower
the rank of the restricted slack matrix, so an argument for these
allocations must retain ambient information; compare
\cite[Section~4.2]{gouveiaEtAl2019}. The same three allocations remain
at the minimum size seventeen; the two converse theorems classify
only their specified branches, not the whole equality set.

Within the certified class, the remaining integer questions concern
the exceptional odd face-pair bit, the finite range before the eventual
formulas, an effective onset, and realization of the remaining
\(\rho=r+1\) equality data. The known odd-dimensional gaps
\cite[Theorem~3.3]{sjobergZiegler2020} do not decide every intervening
pair, and a realization must retain both the stacking facet and the
matching-saturation inequalities.

More generally, the horizon transfer theorem separates graph
connectedness from preservation of a deficient face for extensions
with several visible facets. In the family \(Q_{d,n}\), the two class
sizes exceed the matching number by
\[
 n-2^{d-4}-5\qquad\text{and}\qquad 2^{d-4}-2d+7,
\]
respectively, so the condition \(n\ge2^{d-4}+6\) is exactly where the
displayed certificate begins.
Extending the construction to smaller unrestricted examples therefore
requires more than improving a face-count budget. The sharp
seven-dimensional argument identifies the additional constraint:
the proposed Hall faces and their cofaces must fit into an ambient
incidence structure that preserves both deficiency and connectedness.

\appendix
\section{Asymptotically extremal factor witnesses}
\label{sec:asymptotic}

This appendix records a prismatic alternative to the sharp construction
of Section~\ref{sec:sharp}. Its small-dimensional examples are
superseded by the exact values of Section~\ref{sec:lowdim}. Its asymptotic
estimate also explains what is gained by the simplex-facet witness in
the sharp construction.

For \(k\ge4\), let \(D_k\) be obtained by truncating one vertex of a prism
over a rational \((k-1)\)-gon.  It is a rational simple three-polytope with
\[
 f_0(D_k)=2k,\qquad f_2(D_k)=k+2,
\]
and the truncation produces a triangular facet \(\Delta\).  For
\(\ell\ge0\), put
\[
 A_{\ell,k}=D_k\times[0,1]^\ell,
 \qquad F_{\ell,k}=\Delta\times[0,1]^\ell.
\]
Then \(A_{\ell,k}\) has dimension \(\ell+3\), and
\begin{equation}\label{eq:prismatic-witness-data}
 a=2^{\ell+1}k,\qquad
 \alpha=k+2+2\ell,
 \qquad
 \rho=3+2\ell,
\end{equation}
where \(\rho\) is the number of ridges of \(F_{\ell,k}\).

For \(m\ge3\), use the rational \(m\)-polytope \(B_m\) from
Section~\ref{sec:balanced}, and write
\begin{equation}\label{eq:beta-m}
 \beta_m=f_{m-1}(B_m)=
 \begin{cases}
 8,&m=3,\\
 m+4,&m\ge4,
 \end{cases}
 \qquad f_0(B_m)=\beta_m-2.
\end{equation}
Stack strictly beyond \(F_{\ell,k}*B_m\) in
\(A_{\ell,k}*B_m\), and denote the result by
\(S_{\ell,m,k}\).

\begin{theorem}[Prismatic witnesses]\label{thm:prismatic-witness}
Let \(\ell\ge0\), \(m\ge3\), and \(k\ge4\).  If
\begin{equation}\label{eq:prismatic-condition}
 (2^{\ell+1}-1)k\ge\beta_m+4\ell+5,
\end{equation}
then \(S_{\ell,m,k}\) is a rational polytope of dimension
\(d=\ell+m+4\) with connected nonincidence graph and no vertex--facet
assignment.  More precisely,
\begin{align}
 f_0(S_{\ell,m,k})&=2^{\ell+1}k+\beta_m-1,
 \label{eq:prismatic-vertices}\\
 f_{d-1}(S_{\ell,m,k})&=k+2\beta_m+4\ell+4,
 \label{eq:prismatic-facets}
\end{align}
the facet class is no larger than the vertex class, and
\begin{equation}\label{eq:prismatic-defect}
 \nu(\N(S_{\ell,m,k}))
 \le f_{d-1}(S_{\ell,m,k})-1.
\end{equation}
Consequently its face lattice has no nontrivial free-join decomposition.
\end{theorem}

\begin{proof}
The data in \eqref{eq:prismatic-witness-data} give
\[
 a-\alpha-\rho
 =(2^{\ell+1}-1)k-(4\ell+5).
\]
Thus \eqref{eq:prismatic-condition} is exactly the first inequality in
\eqref{eq:general-factor-criterion}; the second follows from
\(\beta_m-f_0(B_m)=2\).  Proposition~\ref{prop:factor-criterion} proves the
obstruction and, after substitution, the two displayed face numbers.  It also
gives
\[
 \nu(\N(S_{\ell,m,k}))
 \le f_0(B_m)+\alpha+\beta_m+\rho
 =k+2\beta_m+4\ell+3,
\]
which is \eqref{eq:prismatic-defect}.  The difference between the two class
sizes is
\[
 f_0(S_{\ell,m,k})-f_{d-1}(S_{\ell,m,k})
 =(2^{\ell+1}-1)k-(\beta_m+4\ell+5),
\]
so it is nonnegative.  Rationality follows from the rational input polytopes
and a rational choice in the nonempty open stacking region.  Connectedness
and free-join indecomposability follow from
Corollary~\ref{cor:stacking-transfer}.
\end{proof}

The parameter inequality is elementary enough to optimize directly.
Table~\ref{tab:prismatic-records} records choices minimizing the larger of the two displayed
class sizes within this prismatic family for the first few dimensions.  These
are construction records, not claims of global minimality.

\begin{table}[htbp]
\centering
\caption{Prismatic construction records: dimension \(d\), parameters
\(\ell,k\), and vertex--facet counts, with \(m=d-\ell-4\)}
\label{tab:prismatic-records}
\begin{tabular}{c|c|c|c|c}
\(d\)&\(\ell\)&\(k\)&\(f_0\)&\(f_{d-1}\)\\ \hline
8&1&6&31&30\\
9&1&6&31&30\\
10&1&6&32&32\\
11&1&7&37&35\\
12&1&7&38&37\\
13&1&7&39&39\\
14&2&4&43&40\\
15&2&4&44&42\\
16&2&4&45&44\\
17&2&4&46&46
\end{tabular}
\end{table}

In particular, dimensions eight through ten improve the corresponding
balanced counts \(33,36,39\) to \((31,30),(31,30),(32,32)\),
respectively; Section~\ref{sec:lowdim} improves these to the exact
values \(28,28,29\).

Recall \(\eta_d=\eta_{d,1}\) from Definition~\ref{def:certified-class}.

\begin{corollary}[Prismatic benchmark]\label{cor:asymptotic-optimal}
As \(d\to\infty\),
\[
 \eta_d=2d+O(\sqrt d).
\]
In particular, the prismatic family already attains the best possible
leading constant two within the certified classes \(\mathcal C_{d,1}\).
Theorem~\ref{thm:constant-witness} sharpens the remainder to \(O(1)\),
and Corollary~\ref{cor:stable-one} determines the stable exact value.
\end{corollary}

\begin{proof}
The lower bound \(\eta_d\ge2d+8\) is
Proposition~\ref{prop:sharp-floor} with \(h=1\).  For the upper bound, take \(d\) sufficiently
large and choose \(\ell\) so that \(L=2^{\ell+1}\) is the largest power of
two not exceeding \(\sqrt d\).  Put \(m=d-\ell-4\), so \(m\ge4\) and
\(\beta_m=d-\ell\), and choose
\[
 k=\max\left\{4,
 \left\lceil\frac{\beta_m+4\ell+5}{L-1}\right\rceil\right\}.
\]
Then \eqref{eq:prismatic-condition} holds.  Moreover
\(L=\Theta(\sqrt d)\), \(\ell=O(\log d)\), and
\(k=O(\sqrt d)\).  Equation~\eqref{eq:prismatic-facets} therefore gives
\(f_{d-1}=2d+O(\sqrt d)\).  If
\[
 e=(L-1)k-(\beta_m+4\ell+5),
\]
then \(0\le e<L-1\) once the ceiling term is at least four, and the difference
formula in the proof of Theorem~\ref{thm:prismatic-witness} gives
\(f_0=f_{d-1}+e=2d+O(\sqrt d)\).  Finitely many remaining dimensions are
absorbed by the constant implicit in the notation.
\end{proof}

\section{Constant-remainder factor witnesses}
\label{sec:constant}

The prismatic family of Appendix~\ref{sec:asymptotic} balances an exponential
vertex multiplier against the complementary dimension and leaves an
\(O(\sqrt d)\) remainder.  We now replace that multiplier by the vertex
count of a product of two simplices.  The resulting family has a bounded
remainder while retaining a connected nonincidence graph and deficiency
exactly one.  The product witness \(A_s=\Delta_2\times\Delta_s\) below is
independent of the truncated prisms \(A_m\) and \(A_{\ell,k}\) of the
preceding constructions; the cut cube \(C_4\) is likewise independent of
the polygonal base \(C_n\) of Section~\ref{sec:construction}.

\subsection{A dual-excess factor}

Let \(C_3=[-1,1]^3\) be the cube.  Let
\begin{equation}\label{eq:c4-definition}
 C_4=[-1,1]^3\cap\{(x,y,z):2x+2y+z\le3\}.
\end{equation}
Both are rational three-polytopes, and
\begin{equation}\label{eq:cc-data}
 \bigl(f_0(C_c),f_2(C_c)\bigr)=(c+5,c+3),
 \qquad c\in\{3,4\}.
\end{equation}

\enlargethispage{\baselineskip}
\begin{lemma}[The cut cube]\label{lem:cut-cube}
The polytope \(C_4\) has the nine vertices
\[
 \{-1,1\}^3\setminus\{(1,1,1)\}
 \quad\cup\quad\{(1,0,1),(0,1,1)\}
\]
and seven facets.  The six coordinate facets of the cube survive, and the
seventh is the triangle
\[
 \conv\{(1,1,-1),(1,0,1),(0,1,1)\}.
\]
Moreover, for each \(c\in\{3,4\}\), the nonincidence graph of \(C_c\)
has a matching that covers all \(c+3\) facets and leaves exactly two
vertices unmatched.
\end{lemma}

\begin{proof}
The cutting plane in \eqref{eq:c4-definition} separates only the cube vertex
\((1,1,1)\).  The plane meets the three edges incident with that vertex at
\((1,1,-1)\), \((1,0,1)\), and \((0,1,1)\), respectively.  The first is
already a cube vertex and the other two are new.  Thus the displayed nine
points are precisely the vertices.  Each coordinate supporting plane still
contains a two-dimensional face, while the cut plane supports the displayed
triangle; hence these are precisely the seven facets.

For \(C_4\), name the coordinate facets \(X^\pm,Y^\pm,Z^\pm\) and the cut
facet \(T\).  The following are seven distinct nonincidences:
\[
\begin{array}{c|c@{\qquad}c|c}
 X^+&(-1,-1,-1)&X^-&(0,1,1)\\
 Y^+&(-1,-1,1)&Y^-&(-1,1,-1)\\
 Z^+&(1,-1,-1)&Z^-&(-1,1,1)\\
 T&(1,-1,1)&&
\end{array}
\]
and they leave \((1,0,1)\) and \((1,1,-1)\) unmatched.  For \(C_3\), use
\[
\begin{array}{c|c@{\qquad}c|c}
 X^+&(-1,-1,-1)&X^-&(1,1,1)\\
 Y^+&(-1,-1,1)&Y^-&(-1,1,-1)\\
 Z^+&(1,-1,-1)&Z^-&(-1,1,1).
\end{array}
\]
\end{proof}

For \(q\ge3\), put
\begin{equation}\label{eq:bqc-definition}
 D_{q,c}=\operatorname{Pyr}^{q-3}(C_c),
 \qquad B_{q,c}=D_{q,c}^{*},
\end{equation}
where the superscript star denotes a combinatorial dual (not the free
join \(*\)).  A rational realization of the
dual is obtained by translating a rational interior point to the origin and
taking the polar.

\begin{lemma}[The complementary factor]\label{lem:bqc}
The \(q\)-polytope \(B_{q,c}\) satisfies
\begin{equation}\label{eq:bqc-data}
 b:=f_0(B_{q,c})=q+c,\qquad
 \beta:=f_{q-1}(B_{q,c})=q+c+2.
\end{equation}
Its nonincidence graph has a matching that covers all \(b\) vertices and
leaves exactly two facets unmatched.
\end{lemma}

\begin{proof}
Each pyramid step adds one vertex and one facet.  Starting from
\eqref{eq:cc-data}, \(D_{q,c}\) therefore has \(q+c+2\) vertices and
\(q+c\) facets; duality gives \eqref{eq:bqc-data}.

If \(D'=\operatorname{Pyr}(D)\), its nonincidence graph is the disjoint
union of \(\N(D)\) and the single edge joining the new apex to the base
facet.  Indeed, the apex belongs to every side facet and the base contains
every old vertex.  The matching in Lemma~\ref{lem:cut-cube}, followed by
this isolated edge at each pyramid step, covers all facets of \(D_{q,c}\)
and leaves two vertices.  Swapping the two sides under duality proves the
claim for \(B_{q,c}\).
\end{proof}

\subsection{The product-of-simplices witness}

For \(s\ge1\), let
\[
 A_s=\Delta_2\times\Delta_s.
\]
Index its vertices as \(a_{ij}\), where \(0\le i\le2\) and
\(0\le j\le s\).  Its three row facets and \(s+1\) column facets are
\[
 H_i=\conv\{a_{rj}:r\ne i\},\qquad
 V_j=\conv\{a_{ik}:k\ne j\}.
\]
Thus
\begin{equation}\label{eq:as-data}
 \dim A_s=s+2,\qquad
 a:=f_0(A_s)=3s+3,\qquad
 \alpha:=f_{s+1}(A_s)=s+4.
\end{equation}
Select the facet
\[
 F_s=H_2\cong\Delta_1\times\Delta_s.
\]
Its facets, hence its ridges as a facet of \(A_s\), are the two row faces
\[
 R_i^{\rm row}=\conv\{a_{ij}:0\le j\le s\},
 \qquad i=0,1,
\]
and the \(s+1\) column faces
\[
 R_j^{\rm col}=
 \conv\{a_{ik}:i\in\{0,1\},\ k\ne j\},
 \qquad 0\le j\le s.
\]
In particular,
\begin{equation}\label{eq:fs-ridges}
 \rho=f_s(F_s)=s+3.
\end{equation}

Form \(A_s*B_{q,c}\), stack strictly beyond the facet
\(F_s*B_{q,c}\), and denote the resulting polytope by \(T_{s,q,c}\).
The input polytopes are rational, and a rational point exists in the
nonempty open region strictly beyond the selected facet and strictly beneath
all other facets.

\begin{theorem}[Constant-remainder witnesses]
\label{thm:constant-witness}
Let
\begin{equation}\label{eq:constant-parameters}
 c=\begin{cases}3,&d\text{ even},\\4,&d\text{ odd},\end{cases}
 \qquad
 s=\frac{d+c+3}{2},
 \qquad
 q=\frac{d-c-9}{2}.
\end{equation}
For even \(d\ge18\) and odd \(d\ge19\), these are integers with \(q\ge3\).
The rational \(d\)-polytope
\[
 T_d:=T_{s,q,c}
\]
has connected nonincidence graph and balanced class sizes
\begin{equation}\label{eq:constant-size}
 f_0(T_d)=f_{d-1}(T_d)=N_d=2d+2c+4
 =\begin{cases}
 2d+10,&d\ge18\text{ even},\\
 2d+12,&d\ge19\text{ odd}.
 \end{cases}
\end{equation}
Its matching number and deficiency are
\begin{equation}\label{eq:constant-matching}
 \nu(\N(T_d))=N_d-1,\qquad \delta(T_d)=1.
\end{equation}
Consequently, \(T_d\) has no vertex--facet assignment and admits no
nontrivial combinatorial free-join decomposition.
\end{theorem}

\begin{proof}
The parameter identities are equivalent to
\[
 q=s-c-6,\qquad d=s+q+3=2s-c-3.
\]
The selected facet \(F_s*B_{q,c}\) has \(\rho+\beta\) ridges.  Strict
stacking and \eqref{eq:bqc-data}, \eqref{eq:as-data}, and
\eqref{eq:fs-ridges} give
\begin{align}
 f_0(T_d)
 &=a+b+1=3s+q+c+4,\label{eq:constant-vcount}\\
 f_{d-1}(T_d)
 &=\alpha+2\beta+\rho-1
   =2s+2q+2c+10.\label{eq:constant-fcount}
\end{align}
Substituting \(q=s-c-6\) makes both quantities \(4s-2\), which is
\(2d+2c+4\).

By Corollary~\ref{cor:stacking-transfer}, since
\(A_s\not\subseteq F_s*B_{q,c}\), the copy of \(A_s\) survives as a face
and is contained in precisely the
\(\beta=q+c+2=s-4\) retained facets induced by \(B_{q,c}\).  Therefore
\[
v(A_s)+\kappa_{T_d}(A_s)=(3s+3)+(s-4)=4s-1=N_d+1.
\]
Theorem~\ref{thm:exact-defect} gives
\(\nu(\N(T_d))\le N_d-1\).

We attain this bound by an explicit matching.  Write \(\widehat R\) for the
new cone facet over the ridge \(R*B_{q,c}\) of \(F_s*B_{q,c}\), where \(R\)
is a facet of \(F_s\), and write \(\widehat G\) for the
new cone facet over \(F_s*G\) when \(G\) is a facet of \(B_{q,c}\).  First
match all retained facets induced by \(A_s\):
\begin{equation}\label{eq:constant-matching-aold}
 H_0\longleftrightarrow a_{00},\quad
 H_1\longleftrightarrow a_{10},\quad
 V_0\longleftrightarrow a_{20},\quad
 V_j\longleftrightarrow a_{0j}\quad(1\le j\le s).
\end{equation}
For the new cone facets over the ridges of \(F_s\), take
\begin{align}
 \widehat R_0^{\rm row}&\longleftrightarrow a_{11},&
 \widehat R_1^{\rm row}&\longleftrightarrow a_{21},
 \label{eq:constant-matching-rows}\\
 \widehat R_j^{\rm col}&\longleftrightarrow a_{1j}
   &&(2\le j\le s),\label{eq:constant-matching-cols}\\
 \widehat R_0^{\rm col}&\longleftrightarrow a_{2,s-2},&
 \widehat R_1^{\rm col}&\longleftrightarrow a_{2,s-1}.
 \label{eq:constant-matching-cols01}
\end{align}
Finally, enumerate the \(\beta=s-4\) facets \(G\) of \(B_{q,c}\), and
match the corresponding cone facets over \(F_s*G\) bijectively to
\begin{equation}\label{eq:constant-matching-bg}
 a_{2,2},a_{2,3},\ldots,a_{2,s-3}.
\end{equation}
Every pair in
\eqref{eq:constant-matching-aold}--\eqref{eq:constant-matching-bg} is a
nonincidence: \(a_{ij}\) lies on \(H_r*B_{q,c}\) if and only if \(i\ne r\),
on \(V_k*B_{q,c}\) if and only if \(j\ne k\), and on \(\widehat R\) if and
only if \(a_{ij}\in R\); no vertex in row \(2\) lies on a cone facet.  The vertices and facets are distinct, and the construction
uses \(3s+2=a-1\) vertices of \(A_s\), leaving only \(a_{2,s}\).

By Lemma~\ref{lem:bqc}, match all \(b\) vertices of \(B_{q,c}\) to
retained facets induced by \(B_{q,c}\), leaving two such facets.  Match the
stacking apex to one of the two.  Altogether this gives
\[
 (a-1)+(b+1)=a+b=N_d-1
\]
edges, proving \eqref{eq:constant-matching}.

Connectedness follows from
Corollary~\ref{cor:stacking-transfer}.  Lemma~\ref{lem:join-disconnects}
then excludes every nontrivial free-join decomposition.
\end{proof}

Recall that \(\eta_d\) is the minimum of the larger class size over the
certified class \(\mathcal C_{d,1}\) of
Definition~\ref{def:certified-class}. The family \(T_d\) belongs to that
class; it is additionally balanced and has exact deficiency one.

\begin{corollary}[Constant extremal remainder]
\label{cor:constant-optimal}
For every even \(d\ge18\) and odd \(d\ge19\),
\begin{equation}\label{eq:constant-eta-bounds}
 2d+8\le\eta_d\le
 \begin{cases}
 2d+10,&d\ge18\text{ even},\\
 2d+12,&d\ge19\text{ odd}.
 \end{cases}
\end{equation}
In particular,
\begin{equation}\label{eq:constant-eta}
 \eta_d=2d+O(1).
\end{equation}
\end{corollary}

\begin{proof}
The lower bound is Proposition~\ref{prop:sharp-floor} with \(h=1\);
the upper bounds are \eqref{eq:constant-size}.
\end{proof}

By Corollary~\ref{cor:stable-one}, the lower bound is the exact value
for \(d\ge12\). The family \(T_d\) exceeds this optimum by two or four,
but both of its factors are given explicitly.

\section{Exact verification}
\label{sec:verification}

The constructions and proofs above are self-contained and do not depend
on computation. Six Python programs, using only the standard library,
provide independent checks (Online Resource 1). They are included as
ancillary files in the source package and cover the balanced,
prismatic and constant-remainder families, the coordinate family
\(Q_{d,n}\), and the cut-cube seed. They are not used as premises in
any proof. The shared-coordinate and coplanar constructions, the sharp
prescribed-deficiency family of Section~\ref{sec:sharp}, and the
low-dimensional witnesses of Section~\ref{sec:lowdim} are established
analytically and are not covered by these programs. The balanced,
prismatic and constant-remainder programs check the incidence structure
given by the join and stacking descriptions; coordinate realizations
are checked by \texttt{anc/verify\_extremal.py},
\texttt{anc/verify\_family.py} and \texttt{anc/verify\_c4.py}.

The program \texttt{anc/verify\_balanced\_family.py} reconstructs the entire
vertex--facet incidence relation of \(R_d\) from \(d\).  For every requested
dimension it checks the two class sizes, the inherited Hall witness, a
maximum matching of size \(s_d-1\), and connectedness of all \(2s_d\) nodes.
A typical call is
\texttt{python3 anc/verify\_balanced\_family.py 7 8 20 50}.
The companion \texttt{anc/verify\_extremal.py} specializes to \(R_7\).  In
addition to the combinatorial checks, it constructs rational coordinates for
the truncated prism and the stacked free join, verifies every supporting
inequality and facet equality set in exact \texttt{fractions.Fraction}
arithmetic, checks affine dimensions and incident-normal ranks, and prints an
explicit 32-edge matching certificate.

The program \texttt{anc/verify\_prismatic\_family.py} audits the prismatic
family of Section~\ref{sec:asymptotic}.  For each requested triple
\((\ell,m,k)\), passed as \texttt{ell:m:k}, it rebuilds every vertex--facet
incidence of \(S_{\ell,m,k}\), checks the two face-count formulas, the
containment of the witness in exactly \(\beta_m\) facets, the Hall bound
\(f_{d-1}-1\), a maximum matching, and connectedness.  With no arguments it
first re-derives the record table of Section~\ref{sec:asymptotic} by brute
force over \(\ell\) and \(k\) and then verifies all table rows, the
consistency case \(0{:}3{:}13\), which reproduces the combinatorics of
\(R_7\), and two larger parameter choices. The program prints the
computed maximum matching number and asserts that it does not exceed
\(f_{d-1}-1\).

The program \texttt{anc/verify\_constant\_family.py} audits the
constant-remainder construction.  It builds the two base factors \(C_3,C_4\),
performs the pyramid and dual incidence transformation, reconstructs
\(A_s=\Delta_2\times\Delta_s\) and all strict-stack facets, and checks the
face intersection \(\bigcap\{G:A_s\subseteq G\}=A_s\).  It then checks the
class counts, the two equalities in the factor criterion, the displayed
explicit matching of size \(N_d-1\), a maximum matching, that the
witness value \(|\V(A_s)|+\beta\) equals
\(|\V(T_d)|+|\F(T_d)|-\nu(\N(T_d))\), and graph connectedness.
The program \texttt{anc/verify\_c4.py} is an independent geometric check of
the cut cube: it solves every triple of defining equations in exact
rational arithmetic, checks the nine vertices, seven genuine
two-dimensional facets and full dimension. It then finds a seven-edge
nonincidence matching that covers all seven facets of the cut cube and
leaves two vertices unmatched (in the polar, it covers all seven
vertices and leaves two facets unmatched).

Finally, \texttt{anc/verify\_family.py} audits the family \(Q_{d,n}\) of
Section~\ref{sec:construction}. For each
requested pair \((d,n)\) with \(d\ge7\) and \(n\ge2^{d-4}+6\), it checks:
\begin{enumerate}
  \item the apex lies strictly beyond \(X_0\) and strictly beneath every
        other facet;
  \item every listed inequality is valid and has a
        \((d-1)\)-dimensional equality set;
  \item the vertices span dimension \(d\), and at every listed point the
        incident facet normals have rank \(d\), so every listed point is a
        vertex;
  \item the face \(A_n\) is contained in exactly \(2^{d-4}\) facets and
        gives the upper bound \(n+2^{d-4}+2d-2\), which is strictly smaller
        than both the number of vertices and the number of facets;
  \item the uniform matching in Section~\ref{sec:matching} has
        \(n+2^{d-4}+2d-2\) distinct nonincidence edges; and
  \item the full nonincidence graph is connected.
\end{enumerate}
Completeness of the facet list is not tested; it follows from
Lemma~\ref{lem:stacking}, and items (4)--(6) refer to the listed facets.
Parameters are passed as \texttt{d:n} pairs; a bare integer is read as a
value of \(n\) in dimension seven. From the source directory, for example,
the command is
\texttt{python3 anc/verify\_family.py 7:14 7:100 8:22 10:70}.
Without arguments, the four parameterized programs run deterministic
suites; \texttt{verify\_extremal.py} and \texttt{verify\_c4.py} take no
arguments. For \texttt{verify\_family.py} the suite is
\(7{:}14\), \(7{:}15\), \(7{:}16\), \(7{:}20\), \(7{:}50\), \(7{:}100\),
\(8{:}22\), \(8{:}30\), \(9{:}38\), \(10{:}70\), which includes the
threshold parameter \(n=2^{d-4}+6\) of every listed dimension.  The balanced
program checks \(d=7,8,9,10,20,50\); the constant-remainder program checks
\(d=18,20,22,19,21,23,30,31\) with the parity-selected \(c\). No program uses floating-point arithmetic. The balanced, prismatic and
constant-remainder programs generate incidences from the construction
rules (the vertex set of the cut cube used there is derived separately by
\texttt{verify\_c4.py}); the other three programs compute incidences from
exact coordinates, and \texttt{verify\_extremal.py} compares the two.

\paragraph{Online Resource 1.}
Exact rational and combinatorial verification programs for the
balanced, prismatic and constant-remainder families, the coordinate
family \(Q_{d,n}\), and the cut-cube seed.

\section*{Statements and Declarations}

\subsection*{Funding}

No specific funding was received for this work.

\subsection*{Competing interests}

The author has no relevant financial or non-financial interests to disclose.

\subsection*{Data and code availability}

All construction data are given in the paper. Six verification programs
are provided in Online Resource 1 and as ancillary files in the source
package. They check the rational and combinatorial certificates for
the families described in Appendix~\ref{sec:verification}. The
shared-coordinate and coplanar examples, the sharp prescribed-deficiency
family and the low-dimensional witnesses are proved analytically in the
text and are not covered by the programs.

\section*{Disclosure of automated assistance}

OpenAI GPT-5.6 Sol, Anthropic Claude Fable 5, xAI Grok 4.6, OpenAI
GPT-6.1 Sol, OpenAI GPT-6 Astra, and Anthropic Claude Opus 5.5
contributed to the development and preparation of this paper.
Their combined work included generating mathematical ideas and proof
strategies, constructing examples, developing extremal and classification
arguments; checking proofs; checking exact incidence, matching and
connectivity certificates; comparing prior work; and drafting and
revising the manuscript.

\bibliographystyle{amsplain-fullnames}
\enlargethispage{-3\baselineskip}
\bibliography{references}

\end{document}